\documentclass[11pt,a4paper,twoside,openright]{book}

\pdfoutput=1 

\usepackage[english]{babel}
\usepackage{ae,aecompl}
\usepackage{fullpage}
\usepackage{graphicx}
\usepackage{amssymb,amsmath,amsthm}
\usepackage{pictex,dcpic124}
\usepackage{longtable}

\usepackage{hyperref}
\hypersetup{%
 pdftitle={Noncommutative Geometry and Quantum Groups},
 pdfauthor={Francesco D'Andrea},
 pdfsubject={PhD Thesis},
 plainpages=false,
 colorlinks=true,
 pdfstartpage=1,
 pdfstartview=FitH,
 bookmarksopen=true,
 bookmarksnumbered=true,
 bookmarksopenlevel=0,
 linkcolor=black,
 anchorcolor=black,
 citecolor=black,
 urlcolor=black}

\newlength{\tabella}
\newtheorem{thm}{Theorem}[chapter]
\newtheorem{lemma}[thm]{Lemma}
\newtheorem{prop}[thm]{Proposition}
\newtheorem{cor}[thm]{Corollary}

\newtheorem{df}[thm]{Definition}

\allowdisplaybreaks[4]
\numberwithin{equation}{chapter}

\newcommand{\afo}[3]{
\hfill\begin{minipage}{#1}
 \noindent{\textit{``#2''}}, 
                     #3  
\end{minipage} \vspace{1.5cm} }

\newenvironment{prova}{\begin{proof}\parindent=0in}{\end{proof}}

\newcommand{\arxiv}[1]{[arxiv:\htmladdnormallink{#1}{http://arxiv.org/abs/#1}]}
\newcommand{\doi}[1]{[doi:\htmladdnormallink{#1}{http://dx.doi.org/#1}]}

\newcommand{\mc}{\mathcal}
\newcommand{\mf}{\mathfrak}
\newcommand{\mr}{\mathrm}

\newcommand{\A}{\mathcal{A}}
\newcommand{\B}{\mathcal{B}}
\newcommand{\HH}{\mathcal{H}}
\newcommand{\N}{\mathbb{N}}
\newcommand{\Z}{\mathbb{Z}}
\newcommand{\R}{\mathbb{R}}
\newcommand{\C}{\mathbb{C}}
\newcommand{\D}{D\mkern-11.5mu/\,}
\newcommand{\tr}{\mathrm{Tr}}
\newcommand{\ma}[2]{\left(\rule{0pt}{#1}\!\smash[b]{\smash[t]{\begin{array}{cccc}#2 \end{array}}}\!\right)}
\newcommand{\sma}[1]{\textrm{\scriptsize$\left(\!\begin{array}{cccc}#1\end{array}\!\right)$}}
\newcommand{\de}{\mathrm{d}}
\newcommand{\inner}[1]{\left<#1\right>}
\newcommand{\nint}{\int\mkern-19mu-\;}
\newcommand{\op}{\mathrm{OP}^{-\infty}}
\newcommand{\opz}{\mathrm{OP}^0}
\newcommand{\az}{\triangleright}
\newcommand{\za}{\triangleleft}
\newcommand{\aaz}{\,\textrm{\footnotesize$\blacktriangleright$}\,}
\newcommand{\ket}[1]{\left|#1\right>}
\newcommand{\kkett}[1]{\,\left|\mkern-1mu\left|\smash[b]{\smash[t]{#1}}\right>\mkern-3mu\right>}

\newcommand{\sqbn}[2]{\Big[\textrm{\footnotesize$\!\!\begin{array}{c}#1 \\ #2\end{array}\!\!$}\Big]}
\newcommand{\CP}{\mathbb{C}\mathrm{P}}
\newcommand{\oh}{\smash[t]{\smash[b]{\tfrac{1}{2}}}}
\newcommand{\mate}{\textsc{Mathematica}$^\textrm{\tiny{\textregistered}}\,$5}
\newcommand{\GT}[1]{\mathbf{#1}}
\newcommand{\qan}[1]{1-q^{2(#1)}}
\newcommand{\dN}{\!\cdot\!}

\begin{document}

\frontmatter

\pagestyle{headings}


\typeout{Frontmatter}

\begin{titlepage}
\begin{center}

\smallskip

\begin{center}
\includegraphics[width=5cm]{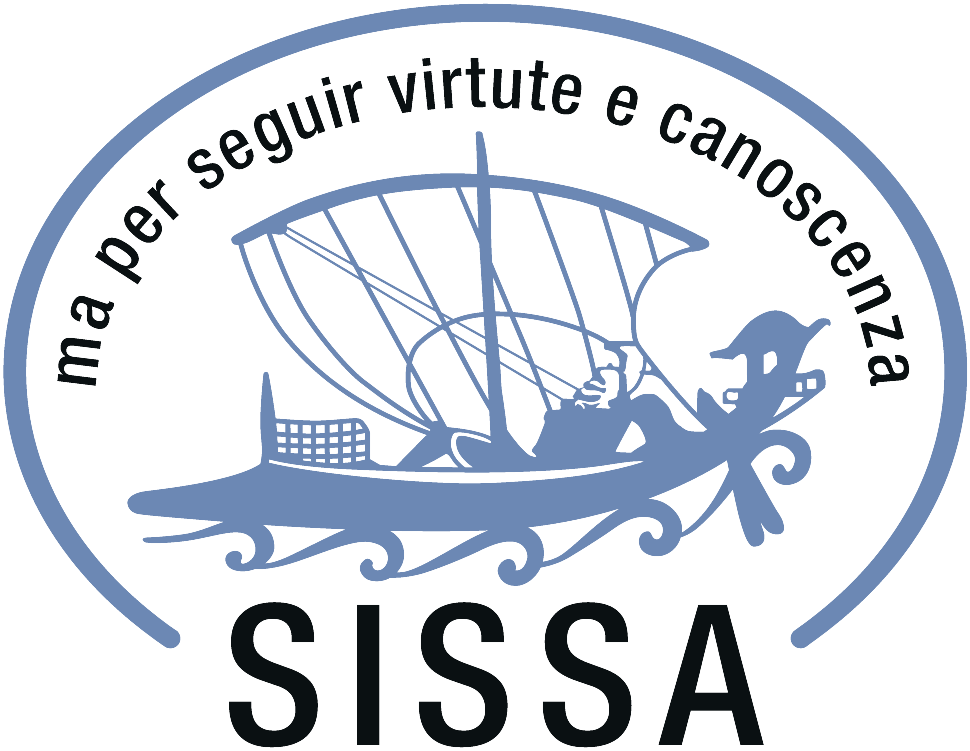}
\end{center}

\textsf{\large International School for Advanced Studies}

\medskip
\textsf{\large Mathematical Physics Sector}

\medskip
\textsf{\large Academic Year 2006-2007}

\vspace{2.5cm}
\textsc{\Large Ph.D.~thesis}

\medskip
\medskip

\begin{center}
{\huge\textsc{
\rule{0pt}{30pt} Noncommutative Geometry and \\
\rule{0pt}{30pt} Quantum Group Symmetries
 }}
\end{center}

\medskip

\vspace{1.5cm}

\textsc{\large candidate}

\medskip

\href{mailto:francesco.dandrea@uclouvain.be}{\textsc{\Large Francesco D'Andrea}}

\vspace{4cm}

\begin{tabular}{c p{1.5cm} c}
\textsc{\large supervisors} & & \textsc{\large external referees}\\
\rule{0pt}{22pt}
\href{mailto:dabrow@sissa.it}{\Large\textsc{Prof.~Ludwik D\k{a}browski}}
& &
\href{mailto:rieffel@math.berkeley.edu}{\Large\textsc{Prof.~Marc A.~Rieffel}} \\
\rule{0pt}{20pt}
\href{mailto:landi@univ.trieste.it}{\Large\textsc{Prof.~Giovanni Landi}}
& &
\href{mailto:varilly@cariari.ucr.ac.cr}{\Large\textsc{Prof.~Joseph C.~V{\'a}rilly}}
\end{tabular}

\end{center}
\end{titlepage}

\newpage
\thispagestyle{empty}
\cleardoublepage


\typeout{Preface}

\phantomsection
\addcontentsline{toc}{chapter}{Preface}

\section*{\LARGE Preface}
\thispagestyle{plain}

This dissertation is based on research done at the Mathematical Physics
sector of the International School for Advanced Studies of Trieste, during
the period from October 2003 to June 2007. It is divided into two parts:
the first part (Chapters \ref{chap:intro}-\ref{chap:NCG}) is an account of
the general theory and a collection of some general notions and results;
the second part (Chapters \ref{chap:S2qs}-\ref{chap:CPlq}) contains the
original work, carried out under the supervision of Prof.~Ludwik D\k{a}browski
and Prof.~Giovanni Landi. Part of the original material presented here has been
published or submitted as a preprint in the following papers:
\begin{itemize}
\item \textit{Local Index Formula on the Equatorial Podle\'s Sphere}, \\
      with L.~D\k{a}browski, \\
      Lett. Math. Phys. \textbf{75} (2006) 235--254. \cite{DD06}
\item \textit{Dirac operators on all Podle\'s spheres}, \\
      with L.~D\k{a}browski, G.~Landi and E.~Wagner, \\
      J.~Noncomm. Geom. \textbf{1} (2007) 213--239. \cite{DDLW07}
\item \textit{The Isospectral Dirac Operator on the $4$-dimensional Orthogonal Quantum Sphere}, \\
      with L.~D\k{a}browski and G.~Landi, \\
      Commun. Math. Phys. \textbf{279} (2008) 77--116 \cite{DDL06}
\end{itemize}


\vfill

\begin{flushright}
\textbf{\Large \ldots acknowledgements}
\end{flushright}

\noindent
I would like to express my gratitude to my supervisors, L.~D{\k a}browski and G.~Landi,
for their help and guidance. Special thanks go to E.~Wagner for our successful
collaboration, and to S.~Brain and to the referees for their valuable comments
on a preliminary version of this dissertation. I would like to thank also
H.~Grosse, P.M.~Hajac, U.~Krahmer, A.~Michelangeli, A.~Sitarz, J.~Varilly,
and many others for interesting and useful conversations.
Finally, my thanks go to my parents, my brother, and my sister
for their support and care.


\pagebreak
\section*{\LARGE Notations}

\begin{longtable}{rp{13cm}}
$\N$ & natural numbers including $0$. \\
$\Z$ & integer numbers. \\
$\mc{S}(\Z)$ & rapid decay sequences: $\{a_n\}_{n\in\Z}\in\mc{S}(\Z)\,\,$
if{}f $\,\,\lim_{n\to\pm\infty}|n|^ka_n=0\,\;\forall\;k\in\N$. \\
$\mc{S}$ & rapid decay matrices: $((a_{j,k}))_{j,k\in\N}\in\mc{S}\,\,$
if{}f $\,\,\lim_{j,k\to\infty}j^mk^n a_{j,k}=0\,\;\forall\;n,m\in\N$. \\
$[x]$ & $q$-analogue of $x\in\C$, given by $[x]:=(q^x-q^{-x})/(q-q^{-1})$.\\
$[n]!$ & $q$-factorial of $n\in\N$, $[n]!:=[1][2][3]\cdots [n]$ if $n\geq 1$, $[0]!:=1$.\\
$[n]!!$ & $q$-double factorial of $n\in\N-1$, given by $[n]!!:=[1][3][5]\cdots [n]$ if $n>0$
    is odd, $[n]!!=[2][4][6]\cdots [n]$ if $n>0$ is even, $[0]!!:=1$ and $[-1]!!:=1$.\\
$\sqbn{n}{m}$ & $q$-binomial of $n,m\in\N$,
    $\sqbn{n}{m}:=\,${\footnotesize${\displaystyle\frac{[n]!}{[m]![n-m]!}}\,$}.\\
$\dot{\otimes}$ & the dotted tensor product of two matrices $A=((A_{ij}))$ and
    $B=((B_{ij}))$ is the matrix $A\,\,\dot{\otimes}\,B$ with entries
    $(A\,\,\dot{\otimes}\,B)_{ij}=\sum_kA_{ik}\otimes B_{kj}$.\\
$*$-\emph{algebra} & $=$ associative unital involutive algebra over $\C$.\\
$\sum_{n\in\emptyset}$ & an empty sum will always mean $0$.\\
$\prod_{n\in\emptyset}$ & an empty product will always mean $1$.\\

~\\

\multicolumn{2}{l}{\bf Landau's `big-O' and `little-o' notations} \\
$a_n=o(b_n)$ & whenever $\{a_n\}$ and $\{b_n\}$ are two sequences and
   $\lim_{n\to\infty}a_n/b_n=0$. \\
$a_n=O(b_n)$ & whenever $\{a_n\}$ and $\{b_n\}$ are two sequences and
   there is a constant $C>0$ such that $|a_n|\leq C|b_n|$. \\
$a_n\lesssim b_n$ & another notation for $a_n=O(b_n)$. \\
\end{longtable}


\cleardoublepage

\pagestyle{headings}

\phantomsection
\addcontentsline{toc}{chapter}{Contents}
\tableofcontents


\mainmatter

\typeout{Capitoli}


\typeout{Capitolo 1}

\chapter{Prologue}\label{chap:intro}

\afo{8cm}{\ldots it seems that the empirical notions on which the metrical
          determinations of space are founded [\ldots] cease to be valid
          for the infinitely small.}
         {Bernhard Riemann (1826--1866)~\cite{Rie73}.}

\noindent
It is a widespread belief that mathematics originates
from the desire to understand (and eventually to formalize) some aspects
of the real world. Quoting \cite{Man07},
$\ll\;$we are doing mathematics in order to understand, create,
and handle things, and perhaps \emph{this understanding} is mathematics$\;\gg\,$.
Let me thus begin with a brief discussion of the
physical ideas that motivated the development of Noncommutative
Geometry and Quantum Group Theory --- the areas of mathematics to
which this dissertation belongs.
Some physicists believe, and Einstein himself expressed this view
in~\cite{Ein98}, 
that physics progresses in stages: there is no `final' theory of Nature,
but simply a sequence of theories which provide more and more accurate
descriptions of the real world.
This idea was made more precise~\cite{BFF77,BFF78a,BFF78b}, by saying that
one passes from one stage to the next by a deformation, the meaning
of which is explained for example in~\cite{Fla82} as follows:
$\ll\,$%
Mathematically one starts with an algebraic structure which is e.g.~a Lie
algebra or an associative algebra and asks the question: does there exist
a $1$ (or $n$) parameter family of similar structures such that for an
initial value (say zero) of the parameter we get the structure we started
with? If such a field of structures exists, we call it a deformation of
the original structure.%
$\,\gg$ \
In this way, one can think of special relativity as a deformation of classical
mechanics with deformation parameter the speed of light $c$: in a
range of velocities much smaller than $c$, which corresponds \emph{formally}
to an initial value $c^{-1}=0$ (although we recall that $c$ is a universal
constant, not a parameter) the Lorentz transformations reduce to the Galilean
ones and special relativity reduces to classical mechanics.
Similarly, quantum mechanics can be viewed as a deformation of classical
mechanics in the framework of associative algebras, with deformation
parameter the reduced Planck constant $\hbar$~\cite{BFF77}.
A typical example is the quantization of a one-dimensional system:
in the classical picture, observables are smooth real functions of
the coordinate $x$ and momentum $p$ and the dynamic is dictated
by the Poisson structure defined by the rule
$$
\{x,p\}=1
$$
and by a preferred function called the Hamiltonian, depending on
the system under consideration and controlling the time evolution.
In the quantum picture, functions are
replaced by elements of an algebra \cite{BJ25,BHJ26} generated
by two real elements $\hat{x}$ and $\hat{p}$ subject to the relation
\begin{equation}\label{eq:QM}
[\hat{x},\hat{p}]:=\hat{x}\hat{p}-\hat{p}\hat{x}=i\hbar\;,
\end{equation}
and the crucial property of this algebra is that its elements
do not commute, unless `$\hbar=0$' which formally corresponds to the
classical limit. The phase space of quantum mechanics is the archetype
of a noncommutative space: a `space' which is no longer described as
a set of points with some additional structure, but by a noncommutative
algebra replacing (or better, deforming) the coordinate algebra of a
manifold.

Of course, the deformations which give rise to special relativity and
quantum mechanics do not ``fall from heaven''. The development of
special relativity originated with Einstein's thought (or ``gedanken'')
experiment, which derived the main properties of the theory (time
dilation/length contraction) from the assumption that the speed
of light is finite and has the same value in every inertial frame
(cf.~\cite{Ein05} or Chapter~IX and X of \cite{Ein85}). On the other hand,
one of the motivations behind the formalism of quantum mechanics (say,
in the Heisenberg formulation) is Heisenberg's gedanken experiment (see
e.g.~Chapter~2 of \cite{Hei30}), which imposes limits on the localizability
of a particle in phase space. This gives rise to the celebrated Heisenberg
uncertainty relations, first presented in \cite{Hei27}, which are
implemented mathematically by the assumption that position and momentum
belong to a noncommutative algebra like (\ref{eq:QM}).
In the same way, a refined version of Heisenberg's gedanken experiment,
suitably generalized to include the effects of gravity, imposes limitations
upon the localizability of a particle~\cite{SW58,Mea64,Pad85,Pad87} and
suggests that at very small length scales the notion of space as a
set of points is inappropriate and that space--time coordinates should
be replaced by noncommuting operators~\cite{DFR94,DFR95,AC97}.

While the physical ideas which led to the study of noncommutative
space--times seem to be still in an embryonic stage (due to a lack
of an experimental evidence which could support or falsify such ideas),
the mathematics behind these ideas is now a well--established independent
field of research in its own right. 
The geometry of noncommutative spaces has been intensively studied
during recent decades by Connes and many others, and most of the tools
of differential geometry have been generalized to the noncommutative framework.
Until now, I have given only a vague idea of what a noncommutative space is,
and one might ask for an exact definition. In fact, there is no general
consensus on the meaning of these words, and often it depends on the context.
For example, the equivalence between the category of locally compact
Hausdorff topological spaces and the category of commutative $C^*$-algebras
(via the Gelfand--Na{\u\i}mark theorem) suggests the idea 
that the theory of $C^*$-algebras should be regarded as noncommutative generalization
of topology. Similarly, from the work of Connes on von Neumann algebras, it comes the
habit to consider them as the noncommutative analogue of measure spaces.
Spectral triples (discussed in Chapter~\ref{chap:NCG}) replace Riemannian
manifolds, and real spectral triples~\cite{Con95} are the algebraic counterpart
of spin manifolds. Now as with any theory, one has to test the postulates
of noncommutative geometry 
by considering a certain number of examples.
This is the purpose of the research described in this dissertation, and
as a source of examples, we focus on Quantum Group Theory.

What then is a `quantum group'? In the middle of the `80s, Drinfeld proposed
a program of quantization of Poisson--Lie groups~\cite{Dri86}. Quantized
Poisson--Lie groups, or simply \emph{quantum groups}, are deformations
of Poisson--Lie groups in the intuitive sense explained above (which
we will make more precise in Chapter~\ref{chap:QG}). They provide a large class
of examples of algebras which for many reasons we would like to interpret
as `coordinate algebras' over noncommutative spaces (with some additional
structure). One interesting feature of these examples is that it seems very
difficult (if not impossible) to fulfill all the axioms of a real spectral triple,
at least in their original formulation.
Thus, a natural question is: are quantum groups `noncommutative spin manifolds'?
In other words: do they satisfy the conditions needed for a real spectral triple,
either in the original or in some modified form?
The aim of this work is to shed some light on these issues.

The plan of the dissertation is the following.
Chapter~\ref{chap:NCG} is a survey of the basic notions of Noncommutative Geometry
which are needed in order to discuss the examples presented in the remaining part
of the dissertation.
In Chapter~\ref{chap:S2qs} we study the family of noncommutative $2$-spheres
introduced by Podle\'s in \cite{Pod87} as quantum homogeneous spaces for
Woronowicz' quantum $SU(2)$ group \cite{Wor87}.
In Chapter~\ref{chap:S4q} we study the $4$-dimensional quantum orthogonal sphere
of~\cite{RTF90}.
Finally in Chapter~\ref{chap:CPlq} we generalize some of these results
to the whole family of odd-dimensional quantum spheres~\cite{VS91}.

Throughout this thesis, unless otherwise stated, the term \emph{$*$-algebras}
is used to refer to complex involutive unital algebras, and the symbol $q$
indicates a fixed parameter in the open interval $]0,1[$.


\typeout{Capitolo 2}

\chapter{Noncommutative Geometry and Quantum Groups}\label{chap:NCG}

\afo{8cm}{The correspondence between geometric spaces and commutative algebras
          is a familiar and basic idea of algebraic geometry. The purpose of this
          book is to extend the correspondence to the noncommutative case in the
          framework of real analysis.}
         {Alain Connes~\cite{Con94}.}

\noindent
As in many theoretical physics' stories, all begins with
the old dream of unifying relativity with quantum mechanics. In an attempt to
find a relativistic version of Schr{\"o}dinger equation, which for a free particle
in $3$ dimensions is
\begin{equation}\label{eq:Schro}
\big(i\partial_t+{\textstyle\sum_{i=1}^3}\partial_{x_i}^2\big)\psi=0 \;,
\end{equation}
physicists noticed that the differential operator $i\partial_t+\sum_{i=1}^3
\partial_{x_i}^2$ does not commute with the natural action of the Poincar\'e Lie algebra
on $C^\infty(\R^{1,3})$ via vector fields, that is: the operator is not \emph{Poincar\'e
invariant}. As covariance is one of the guiding principles of relativity, it became clear that
a relativistic version of (\ref{eq:Schro}) should have involved a Poincar{\'e}--invariant
differential operator. Klein-Gordon proposal\footnote{See~\cite{Kra84} for a historical
review; however, it seems to me that Klein-Gordon equation appeared for the first time
in 1923 in~\cite{Gor23} and not in 1926 as claimed in~\cite{Kra84}.} was to use
the \emph{Laplace--Beltrami operator}, which on $\R^{1,3}$ is
$$
\Delta=-\partial_t^2+{\textstyle\sum_{i=1}^3}\partial_{x_i}^2 \;,
$$
and write the relativistic version of (\ref{eq:Schro}) as $\Delta\psi=m\psi$, for
a free particle of mass $m$. Since $\Delta$ is not of first order in the time
derivative, the knowledge of a solution $\psi$ of Klein-Gordon equation at $t=0$
does not completely determine its time evolution. Motivated by this observation,
Dirac's idea was to search for a first order differential operator $\D$ acting
on vector valued functions and with square $\D^2=\Delta$; the solution for $\R^{1,3}$
first appeared in 1928 in~\cite{Dir28a,Dir28b}. The operator $\D$, which is
called \emph{Dirac operator}, is given on $\R^{1,3}$ by
$$
\D=\gamma_0\tfrac{\partial}{\partial t}+{\textstyle\sum_{i=1}^3\gamma_i
\frac{\partial}{\partial x_i}} \;,
$$
with $\gamma_\mu$ the well known Dirac's gamma matrices.
As everybody knows Dirac equation $\D\psi=m\psi$ provides a description
of elementary particles with spin $1/2$ (and mass $m$), such as electrons,
which is consistent with both the principles of quantum mechanics and special
relativity.

We stress that the demand of symmetries played a crucial role in Dirac's construction;
we will discuss the role of symmetries in noncommutative geometry in Sections
\ref{chap:HA}--\ref{chap:QG}. Another important point is the quest for
a first order differential operator: a square root $\Delta^{1/2}$ of the
Laplacian is well defined (at least in the Riemannian case), but in general
is not a differential operator. We'll return on this when discussing the
`first order condition'. A difference between $\Delta^{1/2}$ and $\D$
is that the former acts on scalar functions, the latter on vector valued
functions (in general, on sections of the spinor bundle, or simply `spinors')
whose components are interpreted as probability amplitudes for a couple of
spin $1/2$ particle--antiparticle; the prediction of the existence of antiparticles,
confirmed in 1933 by the experimental discovery of the positron (the antiparticle
of the electron), was one of the greatest triumphs of modern
theoretical physics and decreed the success of the Dirac operator in physics.

\bigskip

Although the importance of the Dirac operator was immediately recognized
in physics, its importance in mathematics became clear much later, with
Atiyah-Singer local index formula. The generalization of the Dirac operator
to `noncommutative spaces' is in part motivated by its applications to
index theory, which we briefly recall in the following.

Consider an even--dimensional compact Riemannian spin manifold $M$ with fixed
spin structure, let $\HH=\HH_+\oplus\HH_-$ the Hilbert space of
$L^2$-spinors and $\D$ the Dirac operator on $\HH$. Let $\gamma$ be the
natural grading operator on $\HH$ and $\D^+:=\frac{1}{4}(1-\gamma)\D(1+\gamma)$.
If $V\to M$ is any smooth vector bundle and $\D^+_V$ a lift of $\D^+$ to $V$,
the map sending the isomorphism class of $V$ to the index of $\D^+_V$ gives
a group homomorphism from $K^0(M)\simeq K_0(C^\infty(M))$ to the integers $\Z$.
In the odd--dimensional case a similar construction gives a homomorphism from
$K^1(M)\simeq K_1(C^\infty(M))$ to $\Z$.
The natural generalization of this `index map' $K_\bullet(C^\infty(M))\to\Z$
to the $K$-groups of a generic complex involutive unital algebra $\A$
(`$*$-algebra', for short) is possible if there exists a pre--Fredholm module
over it.

\begin{df}
A \emph{pre--Fredholm module} over
a $*$-algebra $\A$ is the data $(\HH,\pi,F)$ of a Hilbert space
$\HH$, a $*$-representation $\pi:\A\to\B(\HH)$ and a bounded
operator $F$ on $\HH$ such that $F-F^*$, $F^2-1$ and $[F,\pi(a)]$,
for all $a\in\A$, are compact operators. The pre--Fredholm module is \emph{even},
or also `graded', if there exists a grading operator $\gamma$ on
$\HH$ (i.e.~a bounded operator satisfying $\gamma=\gamma^*$, $\gamma^2=1$)
commuting with $\pi(a)$ for any $a\in\A$ and anticommuting with $F$.
Otherwise it is odd.
A \emph{Fredholm module} is a pre--Fredholm module satisfying $F=F^*$
and $F^2=1$.
\end{df}

Usually pre--Fredholm modules are called Fredholm module `tout court',
but since we need to distinguish between the two cases in which
$F-F^*$ and $F^2-1$ are zero or just compact, we follow the terminology
of \cite{GVF01} and use the prefix \emph{pre} in the latter case.

Equivalence classes of pre--Fredholm modules form a homology theory called
$K$-homology and dual to $K$-theory~\cite{Kas75}.
Given an even pre--Fredholm module over $\A$, the index map $K_0(\A)\to\Z$ is
constructed as follows. A class $[e]\in K_0(\A)$ is represented by an idempotent
$e=e^2\in\mr{Mat}_N(\A)$ for some $N\in\N$; we can lift in the obvious way the
$*$-representation $\pi$ to a $*$-representation $\tilde{\pi}$ of
$\mr{Mat}_N(\A)$ over $\HH\otimes\C^N$ and form the Fredholm operator
$F^+_e:\HH_+\otimes\C^N\to\HH_-\otimes\C^N$,
with $\HH_\pm:=(1\pm\gamma)\HH$, defined by
\begin{equation}\label{eq:Fpe}
F^+_e:=\tfrac{1}{4}(1-\gamma)\tilde{\pi}(e)F\tilde{\pi}(e)(1+\gamma) \;.
\end{equation}
The index map is just $[e]\mapsto\mr{Index}(F^+_e)$, and gives a pairing
between $K$-theory and $K$-homology.

To perform computations it is crucial that $F-F^*$ and $F^2-1$ are zero (not
just compact) and that the Fredholm module is finite-summable: we shall assume that
this is always the case.
A Fredholm module is \emph{$p$-summable}
if $p\geq 1$ and $[F,\pi(a)]$ is in the $p$-th Schatten-von Neumann ideal
$\mc{L}^p(\HH)$ for all $a\in\A$ (i.e.~$[F,\pi(a)]^p$ is of trace class),
and in this case the index of $F^+_e$ can be evaluated through the formula~\cite{Con85}
\begin{equation}\label{eq:ind}
\mr{Index}(F^+_e)=\tfrac{1}{2}(-1)^m\,\tr_{\HH\otimes\C^N}(\gamma F[F,\tilde{\pi}(e)]^{2m+1})
\end{equation}
for any $m\in\N$ such that $2m+1\geq p$.

The canonical example is the Fredholm module associated with a spin
structure by taking $F:=\D|\D|^{-1}$, the sign of the Dirac operator $\D$
(which we assume for simplicity to be invertible). The fact that in this
case $F$ comes from an unbounded operator $\D$ is what allows to express
the index (\ref{eq:ind}) through local formul{\ae}, which considerably simplify
the computations. This motivates, in the noncommutative case, the passage
from Fredholm modules to unbounded Fredholm modules
(to the best of my knowledge, this name was first used in~\cite{Con89}),
later called \emph{spectral triples}~\cite{Con95}.

\begin{df}\label{def:ST}
A \emph{spectral triple} $(\A,\HH,D)$ over a $*$-algebra $\A$
consists of
\begin{itemize}
\item[(i)] a Hilbert space $\HH$ with a $*$-representation $\pi:\A\to\B(\HH)$;
\item[(ii)] an unbounded selfadjoint operator $D$ with dense domain in $\HH$
 such that
\begin{itemize}
\item[(1)] $\,D+i$ has a compact inverse;
\item[(2)] $\,[D,\pi(a)]$ extends to a bounded operator on $\HH$ for all $a\in\A$.
\end{itemize}
\end{itemize}
The spectral triple is \emph{even} if there exists a
grading operator $\gamma$ on $\HH$ such that $D\gamma=-\gamma D$
and $\pi(a)\gamma=\gamma\pi(a)$ for all $a\in\A$.
\end{df}

\noindent
The reason for the name `unbounded Fredholm module' is that the data $(\HH,\pi,F)$,
with $F:=D(1+D^2)^{-1/2}$, give a pre--Fredholm module over $\A$.
In fact, it was shown in \cite{BJ83} (although they referred to them
by a different name), that any pre--Fredholm module class can be
constructed in this way.

The canonical example of spectral triple comes from the spin
structure over a spin manifold, where $D=\D$ is just the Dirac operator;
Definition \ref{def:ST} is the starting point for a noncommutative
generalization of the theory of spin structures, which is usually
called `noncommutative differential geometry'~\cite{Con94}.
With the additional conditions of finite metric dimension and `regularity'
(which will be recalled later)
it is possible to develop an abstract theory of pseudo--differential operators,
to give meaning to the word `local' in the noncommutative setting, and to
write a local formula for the index (\ref{eq:ind}) which generalizes
Atiyah-Singer index formula~\cite{CM95}.

\bigskip

In the remaining part of this chapter we recall some general definitions of Connes'
noncommutative differential geometry, following \cite{Con94}, 
and collect some theorems which we shall need later.
In Section~\ref{sec:2.1} we describe how to compute the index
(\ref{eq:ind}) using an unbounded Fredholm operator.
In Section~\ref{sec:2.2} we recall some analytic
properties of spectral triples which are needed in Section~\ref{sec:2.3}
to discuss local index formul{\ae} in noncommutative geometry.
In Section~\ref{sec:2.5} we analyze in detail some ideals which will be used
in the sequel. Section~\ref{sec:2.6} is devoted to the noncommutative analogue
of smooth functions. In Section~\ref{sec:2.7} we collect some considerations
about the `reality' and `first order' axioms. Section~\ref{sec:2.8}
is a digression on the topological dimension of quantum spaces.
Section~\ref{sec:2.8r}--\ref{sec:2.4} contain some results which are
used in the examples. Finally, in Sections~\ref{chap:HA} and
\ref{chap:QG} we discuss noncommutative spaces whose symmetries
are described by Hopf algebras, and introduce the examples which
will be studied in the subsequent chapters.

\section{Connection theory and index computations}\label{sec:2.1}
As discussed before, in the commutative case one constructs the
index map by lifting $D$ to vector bundles. We now explain the exact
meaning of this statement by considering a more general (possibly noncommutative)
case, described by a spectral triple $(\A,\HH,D)$.
We first need some preliminary notions about differential calculi.

\bigskip

To any spectral triple is associated a representation of the
\emph{universal differential calculus} (UDC) of $\A$ (cf.~Section III.1 of~\cite{Con94}).
The UDC $(\Omega^\bullet,\de)$ over $\A$ is the data of the graded
algebra $\Omega^\bullet$ generated by degree $0$ elements $a\in\A$
and degree $1$ elements $\de a$, $a\in\A$, with $\de:\A\to\Omega^1$
linear and satisfying the relation
$$
\de(ab)=(\de a)b+a(\de b)\;\;\forall\;a,b\in\A\;.
$$
The involution of $\A$ extends to an involution on $\Omega^\bullet$
via the rule $(\de a)^*=-\de(a^*)$, while the differential $\de$ on
$\Omega^\bullet$ can be defined unambiguously as
$$
\de(a_0\de a_1\ldots\de a_n)=\de a_0\de a_1\ldots\de a_n
$$
and satisfies $\de^2=0$ and the graded Leibniz rule
$$
\de(\omega_1\omega_2)=(\de\omega_1)\omega_2+(-1)^p\omega_1(\de\omega_2)
$$
for $\omega_i\in\Omega$ and with $p$ the degree of $\omega_1$.

Any other differential calculus can be constructed as a
quotient of the universal differential calculus by a differential
ideal. In particular, if $(\A,\HH,D)$ is a spectral triple
the representation $\pi:\A\to\B(\HH)$ can be extended to
an algebra morphism $\pi_D:\Omega\to\B(\HH)$ as follows
$$
\pi_D(a_0\de a_1\ldots\de a_n)=\pi(a_0)[D,\pi(a_1)]\ldots[D,\pi(a_n)] \;.
$$
We denote $\de_Da:=\pi_D(\de a)=[D,\pi(a)]$. \\
This is what is usually called \emph{Connes' differential calculus}.
Quotienting $\pi_D\Omega^\bullet$ by the differential ideal
$J^\bullet:=\pi_D(\de\ker\pi_D)$,
called `Junk ideal' (cf.~\cite{Con94}), one gets a differential
calculus $(\Omega_D^\bullet,\de_D)$ which in the commutative case
coincides with de Rham calculus.
In particular there are no junk $1$-forms, since the Connes' differential
of any $a\in\ker\pi$ is zero; hence $\Omega^1_D=\pi_D\Omega^1$.

\bigskip

The next notion we need is the analogue of Hermitian vector bundles.
If $M$ is a compact $C^\infty$-manifold, then by Serre-Swan theorem there
is a functorial correspondence between smooth vector bundles over $M$ and
finitely generated projective $C^\infty(M)$-modules. The latter replace
vector bundles in the algebraic setting. We recall that any finitely
generated (right) projective $\A$-module $\mc{E}$ is of the form
$$
\mc{E}=e\A^N:=\Big\{\xi=(\xi_1,\ldots,\xi_N)\in\A^N\,\Big|\,e\xi:=\xi\;\mr{i.e.}\;
\sum\nolimits_{j=1}^Ne_{ij}\xi_j=\xi_i\;\forall\;i=1,\ldots,N\Big\}
$$
for $e=e^2\in\mr{Mat}_N(\A)$ an idempotent and for a suitable $N\in\N$.
The group $K_0(\A)$ is used to attach an invariant to any finitely generated
projective $\A$-module, i.e.~to any idempotent $e\in\bigcup_{N\in\N}\mr{Mat}_N(\A)$.

An \emph{Hermitian structure} on $\mc{E}$
is a sesquilinear map $(\,,\,):\mc{E}\times\mc{E}\to\A$ satisfying
$\,(\eta a,\xi b)=a^*(\eta,\xi)b\,$ and $\,(\eta,\eta)\geq 0\,$,
for all $\eta,\xi\in\mc{E}$ and $a,b\in\A$, and such that $\mc{E}$
is selfdual with respect to $(\,,\,)$.\footnote{This is sometimes
called a \emph{non-degenerate} Hermitian structure, see e.g.~Section 4.3
of \cite{Lan97}.} Hermitian structures exist on
any finitely generated projective module $\mc{E}$, and when $e=e^*$
is a projection (this can be always assumed, see e.g.~pages 88--89
of \cite{GVF01}) they are all isomorphic to the one obtained by
restricting to $\mc{E}$ the Hermitian structure on $\A^N$ given by
$$
(\eta,\xi)=\sum\nolimits_{i=1}^N\eta_i^*\xi_i
$$
for $\eta=(\eta_1,\ldots,\eta_N)$ and $\xi=(\xi_1,\ldots,\xi_N)\in\A^N$.
A \emph{connection} on $\mc{E}$ is a linear map
$\nabla:\mc{E}\to\mc{E}\otimes_{\A}\Omega^1_D$ which satisfies the
Leibniz rule
$$
\nabla(\eta a)=(\nabla\eta)a+\eta\otimes\de_Da
$$
and the condition
$$
(\nabla\eta,\xi)-(\eta,\nabla\xi)=\de_D(\eta,\xi)
$$
for all $\eta,\xi\in\mc{E}$ and $a\in\A$. Here if $\eta_{1,2}\in\mc{E}$ and
$\omega_{1,2}\in\Omega^1_D$ we define
$$
(\eta_1\otimes\omega_1,\eta_2\otimes\omega_2):=\omega_1^*(\eta_1,\eta_2)\omega_2\,\in\Omega^1_D \;.
$$

On any $\mc{E}$ there always exists a connection $\nabla_e$ called \emph{Grassmannian
connection}, given by
$$
\nabla_e\eta=e\de_D\eta\quad\mr{for}\;\eta=(\eta_1,\ldots,\eta_N)\in e\A^N\;
\mr{and}\;\de_D\eta=(\de_D\eta_1,\ldots,\de_D\eta_N)\;.
$$
Any other connection differs from $\nabla_e$ by an element in
$\mr{Hom}_{\A}(\mc{E},\mc{E}\otimes_{\A}\Omega^1_D)$.

By completing $\mc{E}\otimes_{\A}\HH$ with respect to the inner product
$$
\inner{\eta_1\otimes v_1,\eta_2\otimes v_2}:=
\inner{v_1,(\eta_1,\eta_2)v_2}
$$
with $\eta_i\in\mc{E}$, $v_i\in\HH$ and $\inner{\,,\,}$ the inner product
of $\HH$, we get (eventually quotienting out the norm zero tensors, if the
Hermitian structure is degenerate) a Hilbert space $\HH_{\mc{E}}$ canonically
isomorphic to $e(\C^N\otimes\HH)$.
The lift $D_{\mc{E},\nabla}$ of $D$ to $\HH_{\mc{E}}$ is defined
by
$$
D_{\mc{E},\nabla}:=1\otimes D+\nabla\otimes 1 \;,
$$
and gives a spectral triple $(\A,\HH_{\mc{E}},D_{\mc{E},\nabla})$.
The operator $D_{\mc{E},\nabla}$ is an unbounded Fredholm operator
and, in the even--dimensional case, a group homomorphism $\mr{Ind}_D:K_0(\A)\to\Z$,
called simply \emph{index map}, is defined by
\begin{equation}\label{eq:imap}
\mr{Ind}_D([e])=\mr{Index}(D_{\mc{E},\nabla}^+) \;,
\end{equation}
with $D_{\mc{E},\nabla}^+:=\frac{1}{2}(1+\gamma)D_{\mc{E},\nabla}(1-\gamma)$
and $\gamma$ the grading operator.
The right hand side of (\ref{eq:imap}) is independent of the particular connection
chosen, and depends on $\mc{E}$ only through its equivalence class in $K_0(\A)$.
Using the Grassmannian connection $\nabla_e$ one proves that (for a
Lipschitz-regular spectral triple, cf.~\cite{Mos97})
$$
\mr{Ind}_D([e])=\mr{Index}(F^+_e) \;,
$$
where $F=D|D|^{-1}$ and $F^+_e$ is given by (\ref{eq:Fpe}).

\section{Analytic properties of spectral triples}\label{sec:2.2}
The next task is to express (\ref{eq:ind}) as a pairing between a
homology theory and the dual cohomology. Before this, we need to
recall some properties of spectral triples. Let
$(\A,\HH,D)$ be a spectral triple, with $D$ invertible.
We also assume that $\pi$ is faithful (otherwise we can replace
$\A$ by $\A/\ker\pi$), we shall identify $\A$ with $\pi(\A)$ and
omit the symbol $\pi$ when there is no risk of confusion.

Such a spectral triple is called \emph{Lipschitz-regular}
if $[|D|,a]$ is bounded for all $a\in\A$. If $(\A,\HH,D)$ is a
Lipschitz-regular spectral triple and $F:=D|D|^{-1}$, by construction
$$
[F,a]=([D,a]-F[|D|,a])|D|^{-1}
$$
is compact for any $a\in\A$, being the product of a bounded operator
with the compact operator $|D|^{-1}$. Hence $(\A,\HH,F)$ is a Fredholm
module. On the other hand the condition that $[F,a]$ is compact for
all $a\in\A$ does not guarantee that a spectral triple is
Lipschitz-regular; by
$$
[|D|,a]=[F,a]D+F[D,a]
$$
the spectral triple is Lipschitz-regular, i.e.~$[|D|,a]\in\B(\HH)$ for
all $a\in\A$, if and only if the operators $[F,a]|D|$ are bounded.
This happens when $[F,a]$ are operators of `order $-1$', a stronger
condition than compactness, the meaning of which we explain below.

For each $s\in\R$ we let $\HH^s$ be the domain of $|D|^s$
(notice that $\HH^s\subset\HH^{s'}$ for $s\geq s'$) and
$\HH^\infty=\bigcap_{s\geq 0}\HH^s$. The space $\HH^s$ is
closed with respect to the norm (on the right hand side
there is the norm of $\HH$)
$$
||v||_s^2=||v||^2+||D^sv||^2 \;.
$$
In the canonical commutative example, $\HH^s$ are Sobolev spaces
and $\HH^\infty$ smooth sections of the spinor bundle.
For each $r\in\R$ a class of operators $\mr{op}^r$
is defined as the subset of (unbounded) operators with domain containing
$\HH^\infty$ which extend to continuous operators $\HH^s\to\HH^{s-r}$
for all $s\in\R$.
We have $\mr{op}^r\subset\mr{op}^{r'}$ for each $r\leq r'$.
As the class $\mr{op}^r$ is difficult to deal with, a smaller class
denoted $\mr{OP}^r$ was introduced in~\cite{CM95} as follows.
Let $\delta$ be the unbounded derivation on $\B(\HH)$ given by
$$
\delta(T):=[|D|,T] \;,
$$
for all $T\in\B(\HH)$. Then, \emph{order zero} operators
are defined as the elements of the class
$$
\opz:=\bigcap\nolimits_{j\in\N}\mr{dom}\,\delta^j \;,
$$
which is a $*$-subalgebra of $\B(\HH)$. 
We remark that by definition, $T\in\mr{dom}\,\delta$ if (i)
$T\in\B(\HH)$, (ii) $T\cdot\mr{dom}\,|D|\subset\mr{dom}\,|D|$
(a necessary condition for $\delta(T)=|D|T-T|D|$ to be defined on $\mr{dom}\,|D|=\HH^1$),
(iii) $\delta(T)$, which is \emph{a priori} defined on $\mr{dom}\,|D|$,
extends to a bounded operator on $\HH$.
By iterating, $T\in\mr{dom}\,\delta^j$ if{}f $T\in\mr{dom}\,\delta^{j-1}$
and $\delta^{j-1}(T)\in\mr{dom}\,\delta$.
\\
Operators of \emph{order $r$}, denoted by $\mr{OP}^r$, are defined by
$$
\mr{OP}^r:=\big\{T\in\B(\HH)\,\big|\,|D|^{-r}T\in\opz\big\} \;.
$$
One has $\mr{OP}^r\subset\mr{op}^r$ for any $r$, and
$\mr{OP}^r\subset\mr{OP}^{r'}$ for all $r'\geq r$. Moreover
if $r\leq 0$, then $\mr{OP}^r$ is a $*$-algebra (of bounded operators)
and a two-sided $*$-ideal in $\mr{OP}^{r'}$, for all $r\leq r'\leq 0$
(cf.~Appendix B of~\cite{CM95}). These classes of operators are the
starting point of an abstract theory of pseudo--differential operators,
which is beyond the scope of this introductory chapter.

\bigskip

For a spectral triple, a stronger requirement than Lipschitz-regularity is the
`regularity' condition~\cite{CM95}.
Classically, while Lipschitz-regular elements of $C(M)$ are Lipschitz-class
functions, regular elements are $C^\infty$ functions.
A spectral triple is called \emph{regular},
or $QC^\infty$, if the following set inclusion holds,
$$
\A\cup[D,\A]\subset\opz\;,
$$
where $[D,\A]=\{[D,a]\,|\,a\in\A\}$.
For a regular spectral triple, the algebra $\Psi^0$ of `pseudo--differential
operators' of order (less than or equal to) zero is defined as the algebra generated
by $\bigcup_{k\in\N}\delta^k(\A\cup[D,\A])$.

Recall that, with $n\in\R^+$, $D$ is called \emph{$n^+$-summable}
if the operator $|D|^{-n}$ is in the Dixmier ideal $\mc{L}^{1,\infty}(\HH)$,
and $n$ is called the \emph{metric dimension} of the spectral triple
if $D$ is $n^+$-summable and $\tr_\omega(|D|^{-n})>0$ (where $\tr_\omega$
is the Dixmier trace, which we recall later)\footnote{This definition is
reminiscent of Hausdorff dimension in fractal geometry~\cite{Fal52}.
In that case one has a one-parameter family of measures $\mu^s$, $s\in\R^+$,
and the dimension of a fractal set $\mc{F}$ is the real number $\,\dim\mc{F}\,$
separating the open interval $\{s\in\R^+\,|\,\mu^s(\mc{F})=0\}$ from the open interval
$\{s\in\R^+\,|\,\mu^s(\mc{F})=+\infty\}$: that is $\mu^s(\mc{F})=0$ for $s<\dim \mc{F}$
and $\mu^s(\mc{F})=+\infty$ for $s>\dim \mc{F}$. The metric dimension is defined
analogously, with $\tr_\omega(|D|^{-s})$ replacing $\mu^s(\mc{F})$.}.
If a spectral triple is regular and $n$-dimensional, the `zeta-type'
function
$$
\zeta_T(s):=\tr_{\HH}(T|D|^{-s})
$$
associated to $T\in\Psi^0$ is defined (and holomorphic) for $s\in\C$ with
$\mr{Re}\,s>n$ and the following definition makes sense.
If it happens that
\begin{itemize}
\item[(i)] for all $T\in\Psi^0$ the function $\zeta_T(s)$ has
a meromorphic extension to $\C$ whose only singularities are poles,
\item[(ii)] the union of such singularities is a countable set $\Sigma\subset\C$,
\end{itemize}
then we call $\Sigma$ the \emph{dimension spectrum} of the spectral triple.
The meaning of $\Sigma$ is that its points are associated with states
on the algebra. In particular if $\Sigma$ is made only of simple poles,
one can prove that the Wodzicki-type residue functionals
\begin{equation}\label{eq:nint}
\tau_m(T)=\nint T|D|^{-m}:=\mr{Res}_{s=m}\tr(T|D|^{-s})
\end{equation}
are defined on $\Psi^0$ for all $m\in\Sigma$, and tracial on suitable
subsets. For example, the top residue $\tau_n$ is tracial on all $\Psi^0$
and there coincides with the Dixmier trace (cf.~Proposition 4.2.4 of~\cite{Con94}).

As in ordinary pseudodifferential calculus, computations can be simplified
by neglecting \emph{smoothing operators}, which
in this abstract setting are defined as the elements of the following class
$$
\op:=\bigcap\noindent_{r\geq 0}\mr{OP}^{-r}=\big\{T\in\B(\HH)\,\big|
\,|D|^kT\in\opz\;\forall\;k\in\N\big\}\;.
$$
The class $\op$ is a two-sided $*$-ideal in the $*$-algebra $\opz$
(actually, in any $\mr{OP}^{-r}$ with $r\geq 0$).
Moreover, if $T$ is a smoothing operator, $\zeta_T(s)$ is holomorphic on $\C$
and the `integral' (\ref{eq:nint}) vanishes.
In this spirit, we say that a multilinear functional on $\Psi^0$ is
\emph{local} if it is insensitive to smoothing perturbations
of its arguments; in particular, residues of zeta-functions are local expressions.
The general philosophy is that locality makes complicated expressions computable,
by ignoring irrelevant details.

\section{Spectral triples and homology theories}\label{sec:2.3}
Let $(\HH,F)$ be a $p$-summable Fredholm module over a $*$-algebra $\A$,
and $n=2m\geq p-1$ an even integer. We can get (\ref{eq:ind}) in two
steps; first we call
\begin{align}
& \mr{ch}^F_n(a_0,\ldots,a_n)=\tfrac{(n/2)!}{2\cdot n!}\,
      \tr(\gamma F[F,a_0]\ldots[F,a_n])\;,\quad\qquad a_i\in\A\,, \label{eq:chF} \\
& \mr{ch}^n(e)=\sum_{k_0,k_1,\ldots,k_n}(e_{k_0k_1}-\tfrac{1}{2}c\,\delta_{k_0k_1})
      \otimes e_{k_1k_2}\otimes\ldots\otimes e_{k_nk_0} \;, \label{eq:che}
\end{align}
(here $c\in\C$ is arbitrary, as $[F,c]=0$) and then write
\begin{equation}\label{eq:indch}
\mr{Index}(F^+_e)=\inner{\mr{ch}^F_n,\mr{ch}^n(e)}
\end{equation}
where
\begin{equation}\label{eq:pair}
\inner{\tau,a_0\otimes\ldots\otimes a_n}:=(-1)^{n/2}\tfrac{n!}{(n/2)!}\,\tau(a_0,\ldots,a_n)
\end{equation}
for any multilinear functional $\tau:\A^{n+1}\to\C$ and any $a_0\otimes\ldots\otimes a_n\in\A^{n+1}$.
The advantage of writing (\ref{eq:ind}) in this way is that (\ref{eq:chF})
is a cyclic cocycle, (\ref{eq:che}) is a cyclic cycle, and (\ref{eq:pair})
is a pairing between the associated cohomology/homology groups.
The element (\ref{eq:che}) is usually called Chern-Connes character.

We recall~\cite{Lod97} that Hochschild $n$-chains $C_n(\A)$
are elements of $\A^{n+1}$, Hochschild $n$-cochains $C^n(\A)$ are multilinear functionals
$\A^{n+1}\to\C$, the Hochschild boundary $b:C_{n+1}(\A)\to C_n(\A)$ is given by
\begin{align*}
b(a_0\otimes a_1\otimes\ldots\otimes a_{n+1}):= \;&
b'(a_0\otimes a_1\otimes\ldots\otimes a_{n+1}) \\ &+
(-1)^{n+1}a_{n+1}a_0\otimes a_1\otimes\ldots\otimes a_n \;,
\end{align*}
with
$$
b'(a_0\otimes a_1\otimes\ldots\otimes a_{n+1}):=
\sum_{j=0}^n(-1)^ja_0\otimes\ldots\otimes a_ja_{j+1}\otimes\ldots\otimes a_{n+1} \;,
$$
and the Hochschild coboundary is its pullback $b^*:C^n(\A)\to C^{n+1}(\A)$, given by
\begin{align*}
(b^*\tau)(a_0,\ldots,a_{n+1})= \;&
\sum_{j=0}^n(-1)^j\tau(a_0,\ldots,a_ja_{j+1},\ldots,a_{n+1}) \\ &+
(-1)^{n+1}\tau(a_{n+1}a_0,a_1,\ldots,a_n) \;.
\end{align*}
The \emph{Hochschild homology} $H\!H_{\bullet}(\A)$ is defined as the homology
of the complex $\{C_n(\A),b\}_{n\geq 0}$, while the dual \emph{Hochschild cohomology}
$H\!H^{\bullet}(\A)$ is defined as the cohomology of the complex $\{C^n(\A),b^*\}_{n\geq 0}$.
If $\A=\A(X)$, with $X$ a smooth affine variety, by Hochschild-Kostant-Rosenberg
theorem (generalized by Connes to smooth manifolds in~\cite{Con95}) there is an
isomorphism of $\A(X)$-modules $\,H\!H_n(\A(X))\;\tilde{\to}\;
\Omega^n_{\mr{de Rham}}(X)\,$, defined on representatives by
\begin{equation}\label{eq:HKR}
a_0\otimes a_1\otimes\ldots\otimes a_n\mapsto \tfrac{1}{n!}\,a_0\de a_1\ldots\de a_n\;.
\end{equation}
This motivates the general philosophy for which in the noncommutative framework
differential forms are replaced by Hochschild chains and de Rham currents by
Hochschild cochains.
The image of (\ref{eq:che}) under (\ref{eq:HKR}), for $c=1$, is the $n$th Chern
character of the Grassmannian connection of $e$; this explains the name `Chern-Connes
character' for $\mr{ch}^n(e)$.

Now, while the map $\mr{ch}^F_n$ in (\ref{eq:chF}) is a Hochschild $n$-cocycle,
that is $b\,\mr{ch}^F_n=0$, the element $\mr{ch}^n(e)$ in (\ref{eq:che}) is
not always a Hochschild cycle; so, in general (\ref{eq:indch}) is not a
pairing between Hochschild homology and the dual cohomology.
To understand what is the appropriate (co)homology theory, one introduces
the \emph{permutation operator} $\lambda_n:\A^{n+1}\to\A^{n+1}$,
$$
\lambda_n(a_0\otimes a_1\otimes\ldots\otimes a_n):=(-1)^na_n\otimes a_0\otimes\ldots\otimes a_{n-1}\;,
$$
its pull-back $\lambda^*_n:C^{n+1}(\A)\to C^{n+1}(\A)$,
$$
(\lambda^*_n\tau)(a_0,a_1,\ldots,a_n)=(-1)^n\tau(a_n,a_0,\ldots,a_{n-1})\;,
$$
and notice that \fbox{$\lambda^*_n\mr{ch}^F_n=\mr{ch}^F_n$} . We call \emph{cyclic
$n$-cochains} $C^n_\lambda(\A)$ those Hochschild $n$-cochains $\tau$ which satisfy
$\lambda^*_n\tau=\tau$. Since
\begin{equation}\label{eq:bbp}
b(1-\lambda_{n+1})=(1-\lambda_n)b'
\end{equation}
we have $\,(1-\lambda_{n+1}^*)b^*=b'^*(1-\lambda_n^*)\,$ and $b^*$ maps the kernel
of $\,1-\lambda_n^*\,$ into the kernel of $\,1-\lambda_{n+1}^*\,$, that is it sends
cyclic cochains into cyclic cochains. The \emph{cyclic cohomology}
$H\!C^\bullet(\A)$ is the cohomology of the complex $\{C_\lambda^n(\A),b^*\}$. 
If $\tau\in C_\lambda^n(\A)$ and $\omega\in C_n(\A)$, the value of $\tau(\omega)$
depends only on the class of $\omega$ in $C_n(\A)/\mr{Im}(1-\lambda_n)$. This suggests
to define \emph{cyclic homology} $H\!C_\bullet(\A)$ as follows. Chains are elements
of $C_n^\lambda(\A):=C_n(\A)/\mr{Im}(1-\lambda_n)$. By (\ref{eq:bbp})
the map $b$ sends the image of $\,1-\lambda_{n+1}\,$ into the image of $\,1-\lambda_n\,$,
hence it descends to a map $b:C^\lambda_{n+1}(\A)\to C^\lambda_n(\A)$.
The \emph{cyclic homology} $H\!H_\bullet(\A)$ is the homology of the complex
$\{C^\lambda_n(\A),b\}$. Remark: the one discussed here is the version
of cyclic homology presented in~\cite{Con85} (cf.~also~\cite{Lod97}).

If we fix $c=0$ in (\ref{eq:che}), then $\mr{ch}^n(e)$ is the
representative of a cyclic cocycle, that is
$$
b\,\mr{ch}^n(e)=
\sum\nolimits_{k_1,\ldots,k_n}e_{k_1k_2}\otimes\ldots\otimes e_{k_nk_1}
=\tfrac{1}{2}(1-\lambda_{n-1})
\sum\nolimits_{k_1,\ldots,k_n}e_{k_1k_2}\otimes\ldots\otimes e_{k_nk_1}
$$
is in the image of $1-\lambda_{n-1}$.
Thus (\ref{eq:indch}) is a pairing between a cyclic cycle and a cyclic cocycle, and this
means that we can replace $\mr{ch}^F_n$ in (\ref{eq:indch}) with any cyclic cycle in the
same cyclic cohomology class, if this helps to simplify computations. If we are in
the fortunate case in which $\mr{ch}^n(e)$ is a Hochschild cycle, i.e.~$b\,\mr{ch}^n(e)=0$,
then we can replace $\mr{ch}^F_n$ in (\ref{eq:indch}) with any Hochschild cocycle in
the same Hochschild cohomology class. Of course, the goal is to replace $\mr{ch}^F_n$
with a local functional.

If $(\A,\HH,D)$ is an even $k$-dimensional spectral triple ($k\in 2\N$),
a Hochschild cocycle is given by
\begin{equation}\label{eq:loc-Hc}
\Psi_D(a_0,\ldots,a_k):=\frac{(k/2)!}{k\cdot k!}\,\nint
\gamma a_0\,\de_D\hspace{1pt}a_1\ldots \de_D\hspace{1pt}a_k\,|D|^{-k} \;,
\end{equation}
where we recall that $\,\de_D\hspace{1pt}a:=[D,a]\,$.
For such a spectral triple the associated Fredholm module is at least $k+1$-summable,
and the cocycles (\ref{eq:loc-Hc}) and (\ref{eq:chF}) have the same class
in Hochschild homology
(Hochschild character theorem, cf.~\cite{Con94}).
So, when $\mr{ch}_n(e)$ is a Hochschild cycle, we can use $\Psi_D$
to compute the index instead of $\mr{ch}^F_n$.
In the general case, more complicated tools are needed.
The idea is to look at the class of $\mr{ch}_n^F$ in a different
cohomology theory, called \emph{periodic
cyclic cohomology}, $H\!P^\bullet$.
Similarly to $K$-theory, $H\!P^\bullet$ is made of two groups,
$H\!P^{\mr{ev}}$ and $H\!P^{\mr{odd}}$. The relevant group for
us is $H\!P^{\mr{ev}}$. An even periodic cyclic cocycle $\phi$ is
a sequence
$$
\phi=(\phi_0,\phi_2,\phi_4,\ldots)
$$
of multilinear functionals $\phi_{2k}:\A^{2k+1}\to\C$ which are all
zero but for finitely many $k$ (and such that $\phi$ is in the kernel
of a suitable coboundary operator, cf.~\cite{Lod97}). An integer--valued
pairing between the class $[\phi]\in H\!P^{\mr{ev}}(\A)$ of $\phi$
and the class $[e]\in K_0(\A)$ of $e$ is defined by
\begin{align*}
\inner{[\phi],[e]} &=\phi_0(e)+\sum\nolimits_{k\geq 1}(-1)^k\,\tfrac{(2k)!}{k!}\,
                   \phi_{2k}(e-\tfrac{1}{2},e,\ldots,e) \;.
\end{align*}
The trick now is to notice that any cyclic $n$-cocycle $\tau$ determines
a periodic cyclic $n$-cocycle $\phi=(0,\ldots,0,\tau,0,\ldots)$,
which has $\tau$ in the $n$th position and zero everywhere else.
In particular, $\mr{ch}^F_n$ determines a periodic cyclic cocycle,
denoted $\mr{ch}^F$ without subscript $n$, which is independent
of $n$ (for even $n$ no less than the metric dimension).
We can rewrite (\ref{eq:indch}) as
$$
\mr{Index}(F^+_e)=\inner{[\mr{ch}^F],[e]}
$$
and although it is always the same expression, in this form it is
clearly invariant if we replace $\mr{ch}^F$ by an element with
the same class in $H\!P^{\mr{ev}}$.

A general theorem (cf.~Theorem II.3 in~\cite{CM95}) gives a
local representative of $[\mr{ch}^F]$ which can be used to
compute the index of $F_e^+$. We quote this theorem in the
particular case we are interested in.

\begin{thm}[Connes-Moscovici]\label{thm:CM}
Let $(A,\HH,D)$ be a regular, even spectral
triple (with finite metric dimension $n$), with dimension spectrum $\Sigma$
made of simple poles. Then, the following formul{\ae} define a periodic cyclic
cocycle in the same class of $\mr{ch}^F$:
\begin{align*}
\varphi_0(a_0) &=\mr{Res}_{s=0}s^{-1}\tr(\gamma a_0|D|^{-2s}) \;\;,\\
\varphi_j(a_0,\ldots,a_j) &=\frac{1}{2}\sum_{k\in\N^j}\frac{(-1)^k}{k_1!\ldots k_j!}\,\alpha_k\,
\nint\gamma a_0[D,a_1]^{(k_1)}\ldots[D,a_j]^{(k_j)}|D|^{-(2|k|+j)}\;\;.
\end{align*}
Here $a_i\in\A\,$, $\,|k|:=k_1+k_2+\ldots+k_j\,$,
$\,\alpha_k^{-1}:=(k_1+1)(k_1+k_2+2)\ldots(k_1+\ldots+k_j+j)\,$,
$\,T^{(0)}:=T\,$ and $\,T^{(i+1)}:=[D^2,T^{(i)}]\;\forall\;i\in\N\,$.
\end{thm}

\noindent
The finiteness of the metric dimension guarantees that all terms in the sum with
$\,|k|+j>n\,$ are zero, hence all sums are finite and only finitely many $\varphi_j$
are different from zero. \\
All the terms $\varphi_j$ with $j>0$ are local (i.e.~they are independent of smoothing
contributions). The unique non-local term is $\varphi_0$.

We close this section with a remark. In the commutative case, when
$\A$ is some algebra of smooth functions on a \emph{connected} compact smooth manifold
$M$, there is always a Hochschild cocycle associated with any projection $p=p^*=p^2$
with elements $p_{ij}\in\A$. Indeed, for $p$ a projection-valued function on $M$,
the map sending $x\in M$ to the matrix trace $\tr\,p(x)$ is a continuous map $M\to\Z$,
then it is constant provided $M$ is connected (its geometrical meaning is the rank of
the bundle underlying the projective module $p\A^k$). Thus, $\tr(p)$ belongs to
the subalgebra $\C\subset\A$, and the tensor
$$
\omega_p=\tr(p-\tfrac{1}{2})\,\dot{\otimes}\,p\,\dot{\otimes}\,(p-1)
$$
has vanishing Hochschild boundary
$$
b\omega_p=\tfrac{1}{2}\tr(p\otimes 1-1\otimes p)=
\tfrac{1}{2}\big\{\tr(p)\otimes 1-1\otimes\tr(p)\big\}=0 \;,
$$
the tensor product being over $\C$ and $\tr(p)\in\C$.
Here the dotted tensor product of two matrices $A=((A_{ij}))$ and
$B=((B_{ij}))$ is the matrix $A\,\,\dot{\otimes}\,B$ with entries
$(A\,\,\dot{\otimes}\,B)_{ij}=\sum_kA_{ik}\otimes B_{kj}$.
Under the map (\ref{eq:HKR}) the cycle $\omega_p$ becomes
the 1st Chern character (modulo a normalization constant)
of the Grassmannian connection on the module $p\A^k$.
For a generic $*$-algebra, a sufficient condition for $\omega_p$
to be a Hochschild cycle is that $\tr(p)\in\C$; in the quantum
group case, it is often possible to deform the matrix trace
in such a way that `$\tr_q$'$(p)$ is a scalar:
the resulting tensor $\omega_2$ will be a cycle in some twisted
Hochschild homology (more on this in Section~\ref{sec:2.8}).

\section{Nested ideals in Hilbert spaces}\label{sec:2.5}
In this section we explain how to concretely work with order $r$ operators
and other relevant ideals. In all this section, two-sided $*$-ideals will
be called simply \emph{ideals}.

Let $\HH$ be a \emph{separable} Hilbert space, and denote by $\mu_j(T)$ the
\emph{$j$th singular value} of the bounded operator $T$, which we recall
is defined by
$$
\mu_j(T):=\inf_{\substack{V\subset\HH, \\ \dim V=j-1}}\;\sup_{v\perp V,\;||v||=1}||Tv|| \;.
$$
We have $\mu_j(T)\geq\mu_{j+1}(T)$, and $\mu_1$ is the operator norm.
Recall that $T\in\B(\HH)$ is \emph{compact} if and only if $\mu_j(T)\to 0$
for $j\to\infty$, and let $\mc{K}$ be the class of compact operators on $\HH$.
This is a closed ideal in $\B(\HH)$, and most ideals we are going to discuss
are dense subspaces of $\mc{K}$. Firstly, with the inequalities (cf.~\cite{Sim05})
\begin{align}
\mu_{j+k}(T_1+T_2)\leq\mu_j(T_1)+\mu_k(T_2)\quad &\textrm{for all}\;T_1,T_2\in\mc{K}\;, \notag\\
\rule{0pt}{16pt}
\mu_j(ST)\leq ||S||\mu_j(T)\quad\mr{and}\quad\mu_j(TS)\leq ||S||\mu_j(T)\quad
&\textrm{for all}\;T\in\mc{K}\;\mr{and}\;S\in\B(\HH)\;,\label{eq:muineq}
\end{align}
one proves that the Schatten-von Neumann classes $\mc{L}^p(\HH)$,
$$
\mc{L}^p(\HH):=\left\{T\in\B(\HH)\,\big|\,
\tr(|T|^p):=\sum\nolimits_{n=1}^\infty\mu_n(T)^p<\infty\right\}\;,
$$
are ideals in $\B(\HH)$. Clearly $\mc{L}^p(\HH)\subsetneq\mc{L}^{p'}(\HH)$
for all $1\leq p<p'<\infty$. In particular, $\mc{L}^1(\HH)$ are the operators
of \emph{trace class}.

For the \emph{spectral} or \emph{quantized calculus}~\cite{Con94} the following
ideals~\cite{Voi79,Voi81} are crucial ($p\geq 1$)
$$
\mc{L}^{(p,\infty)}(\HH) :=
\big\{T\in\B(\HH)\,\big|\,\mu_n(T)=O(n^{-1/p}),\;\mr{as}\;n\to\infty\big\} \;.
$$
We have $\mc{L}^{p}\subsetneq\mc{L}^{(p,\infty)}\subsetneq\mc{L}^{p'}$ for all
$1\leq p<p'<\infty$. In particular, the \emph{Dixmier ideal}
$\mc{L}^{(1,\infty)}$ is the domain of the Dixmier trace, which vanishes on $\mc{L}^1$.
For a positive $T\in\mc{L}^{(1,\infty)}$ the Dixmier trace $\tr_\omega$ is given
by a suitable limit of the bounded sequence~\cite{Dix66}
\begin{equation}\label{eq:limDix}
\tfrac{1}{\log N}\sum\nolimits_{n=1}^N\mu_N(T)
\end{equation}
and extended to all $\mc{L}^{(1,\infty)}$ by linearity. Most important,
if the ordinary limit of (\ref{eq:limDix}) exists, then $T$ is called
\emph{measurable} and its Dixmier trace (if $T\geq 0$) is simply given by
$$
\tr_\omega(T)=\lim_{N\to\infty}\tfrac{1}{\log N}\sum\nolimits_{n=1}^N\mu_N(T) \;.
$$
Elements of $\mc{L}^{(1,\infty)}$ are called for obvious reasons
\emph{infinitesimals of order $1$}, being in the domain of the Dixmier
trace which plays the role of a `noncommutative integral'.
More generally, infinitesimals of order $\alpha\geq 0$ are defined by
$$
\mc{I}_\alpha=\big\{T\in\B(\HH)\,\big|\,\mu_n(T)=O(n^{-\alpha})\big\} \;,
$$
and $\mc{I}_\alpha=\mc{L}^{(p,\infty)}$ if $\alpha=1/p\leq 1$. Infinitesimals
behaves as one expects~\cite{Lan97}: $\mc{I}_{\alpha_1}\supset\mc{I}_{\alpha_2}$ for all
$\alpha_1\leq\alpha_2$, $\mc{I}_{\alpha_1}\cdot\mc{I}_{\alpha_2}\subset
\mc{I}_{\alpha_1+\alpha_2}$ (which in particular means that $\mc{I}_\alpha$
is an ideal in $\B(\HH)$), and the Dixmier trace vanishes on infinitesimals
of order greater than $1$. We let
$$
\mc{I}_\infty:=\bigcap\nolimits_{\alpha\geq 0}\mc{I}_\alpha
=\big\{T\in\B(\HH)\,\big|\,\lim_{n\to\infty}n^\alpha\mu_n(T)=0
\;\forall\;\alpha\geq 0\big\}
$$
be the ideal of infinitesimals of arbitrary high order.

\subsection{Rapid decay matrices}
Up to now, we have discussed subspaces of $\mc{K}$ which are ideals in $\B(\HH)$.
We now introduce subspaces of $\mc{K}$ which are ideals in algebras smaller than $\B(\HH)$.

Let $\{\ket{n}\}_{n\in\N}$ be an orthonormal basis of $\HH$ and
$$
\mc{S}:=\Big\{T\in\B(\HH)\,\Big|\,\lim_{n,m\to\infty}n^\alpha m^\beta\inner{n|T|m}=0\;
\forall\;\alpha,\beta\geq 0\Big\}
$$
the non-unital involutive algebra of \emph{rapid decay matrices}
relative to this basis.

\begin{lemma}
We have $\mc{S}\subset\mc{I}_\infty$.
\end{lemma}
\begin{prova}
We define a selfadjoint positive operator $N$ 
by $N\ket{n}:=(n+1)\ket{n}$, and note that $N^{-\alpha}$ is compact
if $\alpha>0$. For any $T\in\mc{S}$ and any $\alpha>0$, the product
$T N^\alpha$ is bounded. So by (\ref{eq:muineq})
$$
\mu_k(T)\leq ||T N^\alpha||\mu_k(N^{-\alpha})=||T N^\alpha|| k^{-\alpha}
$$
and $k^\beta\mu_k(T)$ goes to zero for any $\beta>0$ (just chose $\alpha>\beta$
in previous inequality). This proves that $T\in\mc{I}_\infty$.
\end{prova}

\begin{lemma}
The algebra $\mc{S}$ is not an ideal in $\B(\HH)$. In particular,
the inclusion $\mc{S}\subset\mc{I}_\infty$ is proper.
\end{lemma}

\begin{prova}
Let $L_q\in\mc{S}$, $0<q<1$, be the operator given by
$$
L_q\ket{n}:=q^n\ket{n} \;.
$$
We now construct a $T\in\B(\HH)$ such that $L_qT\notin\mc{S}$, thus proving
that $\mc{S}$ is not an ideal in $\B(\HH)$. Since $\mc{I}_\infty$ \emph{is}
an ideal in $\B(\HH)$, we conclude that $\mc{S}$ does not coincide with
$\mc{I}_\infty$.

Let $\ket{v}:=\frac{\sqrt{6}}{\pi}\sum_n(n+1)^{-1}\ket{n}$ (it is in $\HH$
since by Euler's formula $||v||^2=\frac{6}{\pi^2}\sum_{k\geq 1}k^{-2}=1$),
and call $T$ the partial isometry $T:=\ket{0}\!\left<v\right|$. Although
$T$ is of rank $1$, it is not a rapid decay matrix. Indeed
$$
\lim_{n\to\infty}n^\alpha\inner{0|T|n}=\tfrac{\sqrt{6}}{\pi}\lim_{n\to\infty}n^{\alpha-1}
$$
is zero if and only if $\alpha<1$, while it should
be zero for all $\alpha$ if $T$ were in $\mc{S}$.
The observation that $L_qT=T$ concludes the proof.
\end{prova}

\subsection{Properties of finite-dimensional spectral triples}
We want now to discuss the relation between `infinitesimals' and the classes
of order $r$ operators defined in Section \ref{sec:2.2}, in the particular
case in which $(\A,\HH,D)$ is a finite-dimensional spectral triple.
We assume for simplicity that $D$ is invertible.

\begin{lemma}
If $D$ is $d^+$-summable, then $\mr{OP}^{-rd}\subset\mc{I}_r$ for
all $r>0$.
\end{lemma}
\begin{prova}
For $T\in\mr{OP}^{-rd}$, $|D|^{rd}T$ is bounded and
$\mu_k(T)\leq ||\,|D|^{rd}T||\mu_k(|D|^{-rd})=O(k^{-r})$.
Thus $T\in\mc{I}_r$ and this concludes the proof.
\end{prova}

Since $|D|^{-1}$ is positive and compact we can find
an orthonormal basis $\{\ket{n}\}_{n\in\N}$ of $\HH$ such that
$$
|D|^{-1}\ket{n}=\lambda_{n+1}^{-1}\ket{n}
$$
with $\lambda_k^{-1}$ ($k\geq 1$) the eigenvalues of $|D|^{-1}$ in decreasing order,
counting multiplicities. As for all positive compact operators, we have
$\mu_k(|D|^{-1})=\lambda_k^{-1}$ and so by $d^+$-summability
$\lambda_k^{-1}=O(k^{-1/d})$.
Let $\mc{S}$ denote rapid decay matrices relative to the eigenbasis
of $D$.

\begin{lemma}\label{lemma:SOP}
For a finite-dimensional spectral triple, $\mc{S}\subset\op$.
\end{lemma}
\begin{prova}
Let $d\in\R^+$ be the metric dimension.
The sequence $\lambda^{-1}_n/n^{-1/d}$ must have a non-zero $n\to\infty$ limit,
otherwise $\lambda_n^{-d}$ would be $o(n^{-1})$,
and $\tr_\omega(|D|^{-d})$ would be zero.
Therefore the inverse sequence $\lambda_n n^{-1/d}$ must be bounded
by a constant $C$. If $T\in\mc{S}$, the sequence
$$
\big|\!\big<n\big||D|^jT|D|^k\big|m\big>\!\big|=
\lambda_n^j\lambda_m^k\,|\!\inner{n|T|m}\!|\leq
C^{j+k}n^{j/d}m^{k/d}\,|\!\inner{n|T|m}\!|\;\;,
$$
is of rapid decay in $n,m$. Thus for all $j,k$, $|D|^jT|D|^k$ belongs
$\mc{S}$ which in particular means that it is bounded.
To conclude the proof now we need just the following Lemma.

\begin{lemma}\label{lemma:uno}
$\;\{T\in\op\}\iff\big\{|D|^jT|D|^k\in\B(\HH)\;\forall\;j,k\in\N\big\}\;$.
\end{lemma}
\noindent\textit{Proof:}\vspace{5pt} \\
\hspace*{-4mm}
\begin{tabular}[t]{rp{144mm}}
(``$\Leftarrow$'') &
By hypothesis, $\delta^j(|D|^rT)=\sum_{l=0}^j(-1)^l\binom{j}{l}|D|^{j+r}T|D|^{l-j}\in\B(\HH)$
for all $r\geq 0$, hence $|D|^rT\in\opz$ for all $r\geq 0$ and $T\in\op$. \\
\rule{0pt}{20pt}
(``$\Rightarrow$'') &
$|D|^jT|D|^k=\sum_{l=0}^k\binom{k}{l}(-1)^{k-l}|D|^{j+k-l}\delta^l(T)\;$.
But $\delta^l(T)\in\op$ since $\op$ is $\delta$-invariant. So
all factors in the right hand side are bounded operators, and then
$\,|D|^jT|D|^k\in\B(\HH)$.
\hfill$\square$
\end{tabular}

\bigskip

This concludes the proof of Lemma \ref{lemma:SOP}.
\end{prova}

\subsection{The dimension spectrum}
Let $(\A,\HH,D)$ be a $d$-dimensional spectral triple, and assume for simplicity
that $D$ is invertible.
If $T\in\mr{OP}^{-r}$, the zeta-function $\zeta_T(s)=\tr(T|D|^{-s})$ is defined
and holomorphic in the half plane $\mr{Re}\,s>d-r$. In particular, if $T\in\op$
the associated zeta-function is holomorphic on the full complex plane and
gives no contribution to the dimension spectrum. As the proof is simple and
instructive, we report it in the following.

\begin{lemma}\label{lemma:seq}
Let $\{a_n\}$ and $\{b_n\}$ two sequences such that
\begin{enumerate}
\item $0<a_n<Cn^{-1/\ell}$, $C\in\R^+$;
\item $|b_n|\leq C'\in\R^+$ is a bounded sequence.
\end{enumerate}
Then $\sum_{n\in\N}a^z_nb_n$ converges (pointwise) to a holomorphic function
in the half plane $\mr{Re}\,z>\ell$.
\end{lemma}
\begin{prova}
The partial sums $f_N(z)=\sum_{n\leq N}a_n^zb_n$ are holomorphic on
$\C$ and dominated by
$$
|f_N(z)|\leq C^{\mr{Re}\,z}C'\sum\nolimits_{n\leq N}n^{-\mr{Re}\,z/\ell} \;.
$$
Thus $f_N(z)$ converges pointwise in the half plane $\mr{Re}\,z>\ell$.
Let $E$ be a bounded closed domain in the half plane $\mr{Re}\,z>\ell$,
and call
$$
k_E=\sup_{z\in E}C^{\mr{Re}\,z}C'\;,\qquad
z_0=\inf_{z\in E}\mr{Re}\,z \;.
$$
Notice that $z_0>\ell$. We have the inequality
$$
|f(z)-f_N(z)| \leq C^{\mr{Re}\,z}C'\sum_{n>N}n^{-\mr{Re}\,z/\ell}\leq k_E
                \sum_{n>N}n^{-z_0/\ell} \;,\qquad\forall\;z\in E\;.
$$
The limit $N\to\infty$ of the right hand side is zero, since the series is
convergent. So
$$
||f-f_N||_\infty \leq k_E\sum_{n>N}n^{-z_0/\ell}\stackrel{N\to\infty}
 {\longrightarrow} 0 \;,\qquad \forall\;z\in E\;,
$$
and $f_N$ converges uniformly to $f$ in $E$. By Weierstrass theorem~\cite[pag.~640]{PF55}
the uniform limit of a sequence of holomorphic functions is holomorphic.
Thus, $f$ is holomorphic in any bounded closed domain $E$ of the half
plane $\mr{Re}\,z>\ell$. We conclude that $f$ is holomorphic on the
entire half plane $\mr{Re}\,z>\ell$.
\end{prova}

For $T\in\mr{OP}^{-r}$, $|D|^rT$ is bounded and so $\inner{n||D|^rT|n}=:b_n$
is a bounded sequence. Called $a_n=\lambda_n^{-1}$ the $n$th eigenvalue
of $|D|^{-1}$, applying Lemma \ref{lemma:seq} we prove that
$\tr(T|D|^{-z+r})=\sum_{n\in\N}a_n^zb_n$ converges to a function
$\zeta_T(z-r)$ which is holomorphic in $\mr{Re}\,z>d$. With the
substitution $s=z-r$, we have proved that $\zeta_T(s)$ is holomorphic
on $\mr{Re}\,s>d-r$ for any $T\in\mr{OP}^{-r}$ (and for all $r\in\R$,
not necessarily positive). For $T\in\op$ previous claim holds for any
$r>0$, and then $\zeta_T(s)$ is holomorphic on all $\C$.

We remark that the dimension spectrum $\Sigma$ contains always
the metric dimension $d$, and any point $s\in\Sigma$ has $\mr{Re}\,s\leq d$
since $\zeta_T(s)$ is holomorphic on $\mr{Re}\,s>d$ for any $T$ bounded.

\subsection{Residues of zeta-type functions}
Motivated by the discussion in previous section, we introduce the following
classes of operators
\begin{equation}\label{eq:Jr}
\mc{J}_r:=\big\{ T\in\opz\,\big|\, T|D|^{-p}\in\mc{L}^1(\HH)\;\forall\;p>r \big\}
\end{equation}
which in the finite-dimensional case contains $\mr{OP}^{-(n-r)}$, with $n$
the metric dimension, and can be neglected when computing residues of
zeta-functions in poles with real part bigger than $r$.

\begin{lemma}
The collection $\mc{J}_r$ is a two-sided ideal in $\opz$.
\end{lemma}

\begin{prova}
Clearly $\mc{J}_r$ is a vector space: if $T_1,T_2\in\mc{J}_r$, that is
$T_1|D|^{-p}\in\mc{L}^1(\HH)$, $T_2|D|^{-p}\in\mc{L}^1(\HH)$
for all $p>r$, then $T_1|D|^{-p}+T_2|D|^{-p}\in\mc{L}^1(\HH)$ for all $p>r$,
which means $T_1+T_2\in\mc{J}_r$.

That $\mc{J}_r$ is a left ideal is straightforward. 
Since $\mc{L}^1(\HH)$ is a two-sided ideal in $\B(\HH)$,
if $T_1\in\opz$ and $T_2\in\mc{J}_r$, for all $p>r$ we have that 
$T_1\cdot T_2|D|^{-p}$ is the product of a  bounded operator,
$T_1$, with a trace class one, $T_2|D|^{-p}$, thus it is of
trace class, and $T_1T_2\in\mc{J}_r$. 

From Appendix B of~\cite{CM95} for any $p>0$, we know that the
bounded operator $|D|^{-p}$ maps $\HH$ to $\HH^p:=\mr{dom}\,|D|^p$,
that $T\in\opz\subset\mr{op}^0$ is a bounded operator $\HH^p\to\HH^p$,
and finally that $|D|^p$ is bounded from $\HH^p$ to $\HH$.
Thus, for $T\in\opz$, the product $|D|^pT|D|^{-p}$ is a bounded
operator on $\HH$. 
Now, if $T_1\in\opz$ and $T_2\in\mc{J}_r$, for all $p>r$ we can
write $T_2T_1|D|^{-p}=T_2 |D|^{-p} \cdot |D|^p T_1 |D|^{-p}$
as the product of a bounded operator, $|D|^{p}T_1|D|^{-p}$,
with a trace class one, $T_2|D|^{-p}$; thus $T_2T_1|D|^{-p}$
is of trace class; therefore $T_2T_1\in\mc{J}_r$ and $\mc{J}_r$
is also a right ideal.
\end{prova}

Of course $\opz=\mc{J}_n$.
One might ask whether in the finite-dimensional case $\mr{OP}^{-r}$ and $\mc{J}_{n-r}$
are the same ideal. In general they are not, and the following is a counterexample.

Let $\HH$ be the Hilbert space with orthonormal basis $\ket{l,m}$, with $l\in\N$
and $l-|m|\in\N$ (e.g.~think of $L^2(S^2)$ in the basis of Laplace spherical
harmonics), and define a positive $D$ and a bounded $\tilde{L}_q$ by
$$
D\ket{l,m}=(l+1)\ket{l,m} \;,\qquad
\tilde{L}_q\ket{l,m}=q^{l+m}\ket{l,m} \;.
$$
We have $\tilde{L}_q\in\mc{J}_1$, since
$$
0\leq\tr(\tilde{L}_qD^{-s})=\sum\nolimits_{l\in\N}(l+1)^{-s}\tfrac{1-q^{2l+1}}{1-q}
\leq\tfrac{1}{1-q}\zeta(s)\;,\qquad\forall\;s>1\;,
$$
where $\zeta(s)$ is the Riemann zeta-function.
On the other hand
$$
||D\tilde{L}_q||^2\geq\sup\nolimits_{l,m}\big<l,m\big|D\tilde{L}_q\big|l,m\big>
\geq\sup\nolimits_{l}\big<l,-l\big|D\tilde{L}_q\big|l,-l\big>
=\sup\nolimits_{l}(l+1)=+\infty\;,
$$
and then $D\tilde{L}_q$ is not bounded and $\tilde{L}_q\notin\mr{OP}^{-1}$.

\section{Smooth `functions' on noncommutative spaces}\label{sec:2.6}
As continuous functions carry all topological information about a space,
$C^\infty$-functions carry all relevant information about its smooth structure.
In the canonical commutative example of a spectral triple, we can replace
$C^\infty(M)$ with any subalgebra and still have a regular spectral triple.
Thus, regularity condition is not sufficient to select $C^\infty(M)$ against
proper subalgebras. The characteristic property which enables to select
$C^\infty(M)$ is its stability under holomorphic functional calculus.

Let $\mc{C}$ be a $C^*$-algebra. If a function $\psi$
is holomorphic in a neighborhood of the spectrum of $a\in\mc{C}$,
we define an element $\psi(a)\in\mc{C}$ as the integral
\begin{equation*}
\psi(a)=\oint_\Gamma\psi(z)(z\cdot 1-a)^{-1}\de z
\end{equation*}
on a rectifiable contour $\Gamma$ winding once around $\mr{Sp}(a)$. The map
$a\mapsto\psi(a)$ is called \emph{holomorphic functional calculus}, and a
dense involutive subalgebra $\mc{A}\subset\mc{C}$ is called \emph{pre-$C^*$-algebra}
if $\psi(a)\in\A$ whenever $a\in\A$ and $\psi(a)$ is well-defined in $\mc{C}$.

Smooth functions on a compact manifold are characterized also by another
property: they form a Frech\'et space, that is a locally convex topological
space closed with respect to a countable family of semi-norms. Actually, they
form a Frech\'et algebra, that is an algebra and a Frech\'et space such that
the semi-norms are submultiplicative.
Given a regular spectral triple $(\A,\HH,D)$ and with the notation $\delta(a):=[|D|,a]$,
the algebra $\A^\infty$, completion of $\A$ in the topology defined by the family
of seminorms
$$
||a||_{1,n}:=||\delta^n(a)||\;,\qquad
||a||_{2,n}:=||\delta^n([D,a])||\;,
$$
is a Frech\'et pre-$C^*$-algebra~\cite{RV06} and $(\A^\infty,\HH,D)$ is
still a regular spectral triple.


\section{Reality and first order condition}\label{sec:2.7}

%
%
%
A real structure on a spectral triple is a generalization of the charge
conjugation operator acting on spinors~\cite{Con95}, has application in
physics to the construction of the spectral action \cite{CC97} (see also
\cite{CCM06} and \cite{Kas00}), and its existence is one of a list of conditions
which should characterize a noncommutative spin manifold \cite{Con96}
(cf.~\cite{RV06} for an up-to-date account of this line of research).
When testing these axioms on quantum groups, or quantum homogeneous
spaces, things start to become complicated exactly when one arrives to the
reality axiom.

Let $(\A,\HH,D)$ be a spectral triple.
A \emph{real structure} is an antilinear isometry $J$ on
$\HH$ such that for all $a,b\in\A$:
$$
J^2=\pm 1\;,\qquad
JD=\pm\,DJ\;,\qquad
[a,Jb^*J^{-1}]=0\;,\qquad
\boxed{[[D,a],Jb^*J^{-1}]=0}\;.
$$
If the spectral triple is even, we require also that $J\gamma=\pm\,\gamma J$. The
signs `$\pm$' determine the \mbox{$KO$-dimension} of the spectral triple (cf.~Table 9.45
of~\cite{GVF01}), which classically corresponds to the metric dimension modulo $8$.
The framed condition is usually called \emph{first order axiom}.
In the canonical commutative example $(C^\infty(M),\HH,\D)$, with $M$
a compact Riemannian spin manifold, $\HH$ the square integrable spinors
and $\D$ the Dirac operator, $J$ is the charge conjugation and since
$Jb^*J^{-1}=b$ the 3rd condition trivially follows from commutativity
of the algebra. For $\D$ a 
pseudo--differential operator, the
condition $[\D,a]\in\B(\HH)$ implies that $\D$ is of order \mbox{$\leq 1$},
while the further condition $[[\D,a],b]=0$ is satisfied if $\D$ is
a first-order \emph{differential} operator. A clarifying example
is $(C^\infty(M),L^2(M),\Delta^{1/2})$,
with $M$ compact Riemannian manifold and $\Delta$ the Laplace-Beltrami
operator; it is a regular spectral triple, but not real since in general
$[[\Delta^{1/2},a],b]$ is not zero.
In the noncommutative case, the 3rd condition is inspired by von Neumann
algebras, for which Tomita-Takesaki theorem ensures the existence of a $J$
such that $b\mapsto Jb^*J^{-1}$ is an isomorphism from $\A$ into the
opposite algebra, concretely represented inside the commutant~$\A'$.

The first order condition is what sometimes creates problems.
Typical examples are quantum groups (and quantum homogeneous spaces),
where the commutators $[[D,a],JbJ^{-1}]$ vanish only modulo some ideal of
`infinitesimals'. Motivated by examples coming from quantum group theory,
we try with the following definition.

\begin{df}
A \emph{weak real structure} is an antilinear
isometry $J$ on $\HH$ such that $\forall\;a,b\in\A$,
$$
J^2=\pm 1\;,\qquad
JD=\pm DJ\;,\qquad
[a,JbJ^{-1}]\in\mc{I}\;,\qquad
[[D,a],JbJ^{-1}]\in\mc{I}\;,
$$
where the set $\mc{I}$ is a suitable two-sided ideal made of infinitesimals.
If the spectral triple is even with grading $\gamma$, we impose the further
relation $J\gamma=\pm\gamma J$.
\end{df}

\noindent
A real structure corresponds to the choice $\mc{I}=0$;
while in examples coming from quantum groups~\cite{DLPS05,DLS05,DDLW07,DDL06}
one usually takes $\mc{I}=\op$. We stress that the term `weak real structure'
is not standard and is provisionally used in this dissertation to simplify
the discussion.

A real spectral triple of dimension $n$ gives an element
of $K\!R^n(\A\otimes\A^{\mr{op}},\C)$, and this relation with real
$K$-homology is not clear if $J$ is a weak real structure.

The second remark is that this weaker condition might be too weak.
The ideal $\mc{I}$ should be chosen with some care.
For example, let $M$ be a compact Riemannian manifold, $\nabla$ the Levi-Civita connection,
$\Delta=\nabla^*\nabla$ the Laplacian, and denote $\delta(f)=[\Delta^{1/2},f]$
(bounded if $f$ is Lipschitz). The commutator
$$
[[\Delta,f],g]=2(\de f,\de g)
$$
(on the right hand side we have the Hermitian structure on $1$-forms) is bounded
for all $f,g\in C^\infty(M)$, and then from $\Delta=\Delta^{1/2}\Delta^{1/2}$
we deduce that the operator
$$
2\Delta^{1/2}[\delta(f),g]=[[\Delta,f],g]-2\delta(f)\delta(g)+[\delta^2(f),g]
$$
is bounded too. We have both $\Delta^{-1/2}\in\mc{L}^{(n,\infty)}(\HH)=\mc{I}_{1/n}$
(by Weyl estimates where $n:=\dim M$, we assume for simplicity that $\Delta$
is invertible) and $\Delta^{-1/2}\in\mr{OP}^{-1}$. Thus, the commutator
$$
[[\Delta^{1/2},f],g]=[\delta(f),g]=\Delta^{-1/2}\cdot(\textrm{bounded operator})
$$
is compact (it is even an infinitesimal of order $1/n$ and an element of $\mr{OP}^{-1}$).
This proves that the complex conjugation $J$ is a weak real structure for
$(C^\infty(M),L^2(M),\Delta^{1/2})$, with $\mc{I}$ any
ideal of compact operators containing either $\mr{OP}^{-1}$ or infinitesimals
of order $1/n$, even if $M$ may not be a spin manifold. This is a situation
that one probably would like to avoid, by requesting for example that elements
of $\mc{I}$ are infinitesimals of order greater than $1$.


\section{The homological dimension of quantum spaces}\label{sec:2.8}
One can generalize Hochschild homology by considering Hochschild homology with
coefficients in a bimodule~\cite{Lod97}. A particular case is twisted Hochschild
homology, which we recall briefly here.
Let $\A$ be a (complex, associative) unital algebra and $\kappa\in\mr{Aut}(\A)$
a scaling automorphism (i.e.~$\A$ admits a linear basis in which
$\kappa$ is diagonal). On $\A^{\otimes(n+1)}$ define a linear map $\theta_\kappa$ by
$$
\theta_\kappa(a_0\otimes\ldots\otimes a_n)=\kappa(a_0)\otimes\ldots\otimes\kappa(a_n)\;.
$$
The set of twisted $n$-chains $C^\kappa_n$ is the subset of
elements of $\A^{\otimes(n+1)}$ in the kernel of $\theta_\kappa-1$.
Twisted Hochschild homology $H\!H^\kappa_{\bullet}(\A)$ is defined
as the homology of the complex $\{C^\kappa_n,b_\kappa\}_{n\geq 0}$, with 
boundary operator $b_\kappa:C^\kappa_{n+1}\to C^\kappa_n$
given by
\begin{align*}
b_\kappa(a_0\otimes a_1\otimes\ldots\otimes a_{n+1})= \;&
\sum_{j=0}^n(-1)^ja_0\otimes\ldots\otimes a_ja_{j+1}\otimes a_{n+1} \\ &+
(-1)^{n+1}\kappa(a_{n+1})a_0\otimes a_1\otimes\ldots\otimes a_n \;.
\end{align*}
This is not the ordinary definition, but is equivalent to it
if $\kappa$ is a scaling automorphism (see e.g.~\cite{Sit05}).
If $\kappa=id_{\A}$ one gets the `untwisted' groups $H\!H_\bullet(\A)$.

From the dual point of view, a twisted $n$-cochain is a multilinear map
$\phi:\A^{\otimes(n+1)}\to\C$ in the kernel of $\theta_\kappa^*-1$,
with $\theta_\kappa^*$ the pullback of $\theta_\kappa$,
$$
\theta_\kappa^*\phi(a_0,a_1,\ldots,a_n)=
\phi\bigl(\kappa(a_0),\kappa(a_1),\ldots,\kappa(a_n)\bigr)\;.
$$
We denote twisted $n$-cochains with $C_\kappa^n$. Let $b_\kappa^*$
be the pullback of $b_\kappa$. The twisted Hochschild cohomology
$H\!H_\kappa^{\bullet}(\A)$ is defined as the cohomology of the complex
$\{C_\kappa^n,b^*_\kappa\}_{n\geq 0}$.

The expression $\inner{[\phi],[\omega]}:=\phi(\omega)$ is a pairing
between $H\!H_\kappa^{\bullet}(\A)$ and $H\!H^\kappa_{\bullet}(\A)$, that is
it depends only on the class of $\phi$ in $H\!H_\kappa^{\bullet}(\A)$ and
by the class of $\omega$ in $H\!H^\kappa_{\bullet}(\A)$.

If $\mc{U}$ is a Hopf algebra and $\A$ a left $\mc{U}$-module algebra and $\kappa$
commutes with the action of $\mc{U}$, the action on $\A^{n+1}$ obtained by
Hopf tensor product defines a structure of left $\mc{U}$-module on $C_n^\kappa(\A)$.
Similarly, the pullback of the action on $\A^{n+1}$ on multilinear maps
$\phi:\A^{n+1}\to\C$, given by
$$
(h\az\phi)(\omega):=\phi(S^{-1}(h)\az\omega)\;,\qquad
\forall\;\omega\in\A^{n+1},\;h\in\mc{U}\,,
$$
turns $C^n_\kappa(\A)$ into a left $\mc{U}$-module.

If $\A$ is the algebra of smooth functions on a $C^\infty$ manifold, the 
metric dimension of the space equals the Hochschild dimension
$\,\sup\{n\,|\,H\!H_n(\A)\neq 0\}\,$. In the general case one defines
the twisted Hochschild dimension as
$$
\mr{Hdim}^\kappa(\A):=\sup\big\{n\,\big|\,H\!H^\kappa_n(\A)\neq 0\big\} \;.
$$
A common situation on $q$-deformations $M_q$ of homogeneous spaces $M$
is that the ordinary Hochschild dimension of $M_q$ is less than the dimension
of the commutative limit $M$, a phenomenon known as `dimension drop'.
For example the Hochschild dimension of $SL_q(N)$ equals its rank,
while it is $1$ for Podle\'s spheres, as proved in \cite{MNW91}.
On the other hand, $\mr{Hdim}^\kappa$ is believed to coincide with
the classical dimension if we let $\kappa$ be the modular automorphism
(this was proved in the case of $SL_q(N)$ in~\cite{HK05,HK06}, for Podle\'s
spheres in~\cite{Had07} and quantum hyperplanes in~\cite{Sit05}).

Suppose now we have an even Fredholm module $(\A,\HH,F,\gamma)$,
and that the automorphism $\kappa$ is implemented on $\HH$ by the adjoint
action of a (unbounded, invertible, positive, with dense domain) operator
$K$ commuting with $F$ and $\gamma$, that is
$$
\kappa(a)v=KaK^{-1}v\;,\qquad
[F,K]v=0\;,\qquad
[\gamma,K]v=0\;,
$$
for any $a\in\A$ and any $v$ in the domain of $K$. Instances of this situation
are provided by subalgebras of the `coordinate algebra' of a compact quantum
group, with $\kappa$ the modular automorphism, and $(\A,\HH,F,\gamma)$ an
equivariant even Fredholm module (cf.~Section \ref{chap:HA}).
Finally, suppose there exists $p\in[1,\infty)$ such that the operator
$$
K^{-1}[F,a_0][F,a_1]\ldots [F,a_{p-1}]
$$
is of trace class on $\HH$ for all $a_0,\ldots,a_{p-1}\in\A$.

\begin{lemma}\label{lemma:ch}
A twisted Hochschild $n$-cocycle is defined by the formula
$$
\mr{ch}^{F,\kappa}_n(a_0,\ldots,a_n):=
\tr_{\HH}(K^{-1}\gamma F[F,a_0]\ldots [F,a_n])\;,
$$
for any even $n\geq p-1$.
\end{lemma}
\begin{prova}
Using Leibniz rule
\begin{align*}
\{b^*_\kappa\mr{ch}^{F,\kappa}_n\}(a_0,\ldots,a_{n+1})&=
\tr_{\HH}(K^{-1}\gamma Fa_0[F,a_1]\ldots [F,a_{n+1}]) \\
&+\tr_{\HH}(K^{-1}\gamma F[F,a_0]\ldots [F,a_n]a_{n+1}) \\
&-\tr_{\HH}(K^{-1}\gamma F[F,\kappa(a_{n+1})a_0]\ldots [F,a_n]) \;.
\end{align*}
By hypothesis
$K^{-1}\gamma F[F,\kappa(a_{n+1})a_0]=\gamma F[F,a_{n+1}K^{-1}a_0]$.
Then using again Leibniz rule and cyclicity of the trace
we deduce that the r.h.s.~vanishes, and $\mr{ch}^{F,\kappa}_n$
is closed.

Similarly, from cyclicity of the trace we get
\begin{align*}
\mr{ch}^{F,\kappa}_n(a_0,\ldots,a_{n+1}) &=
\tr_{\HH}([F,a_n]K^{-1}\gamma F[F,a_0]\ldots [F,a_{n-1}]) \\ &=
\tr_{\HH}(K^{-1}[F,\kappa(a_n)]\gamma F[F,a_0]\ldots [F,a_{n-1}])\;.
\end{align*}
But $F^2=1$ implies $[F,\kappa(a_n)]F=-F[F,\kappa(a_n)]$, and
$$
\mr{ch}^{F,\kappa}_n(a_0,\ldots,a_{n+1})=
\tr_{\HH}(K^{-1}\gamma F[F,\kappa(a_n)][F,a_0]\ldots [F,a_{n-1}])\;.
$$
Thus, $\mr{ch}^{F,\kappa}_n$ is twisted cyclic
\begin{equation}\label{eq:chcyc}
\mr{ch}^{F,\kappa}_n(a_0,a_1,\ldots,a_n)=
\mr{ch}^{F,\kappa}_n(\kappa(a_n),a_0,\ldots,a_{n-1})
\end{equation}
and by iterating last identity $n$-times we prove that $\mr{ch}^{F,\kappa}_n$
is in the kernel of $\theta_\kappa^*-1$.
\end{prova}

We'll show in the case of Podle\'s spheres $S^2_{qs}$ and of the quantum
orthogonal $4$-sphere $S^4_q$ how to construct a twisted Hochschild cycle
which is a deformation of the volume form on $S^2$ and $S^4$ respectively.
We'll also observe that on $S^2_{qs}$ and $S^4_q$, the operators
$K^{-1}[F,a_0][F,a_1]\ldots [F,a_{p-1}]$ are of trace class for $p>2$ and $p>4$
respectively, as in the $q\to 1$ limit. This is in contrast with the fact that
Fredholm modules are usually at most $1$-summable on $q$-deformed spaces
(for $q$ real).

\section{Some general results}\label{sec:2.8r}
Since in all our examples there is a standard strategy to prove that a
collection $(\A,\HH,D,\gamma)$ is a (even, regular) spectral triple,
we explain it here.

Assume $\A$ is a $*$-algebra of bounded operators on a separable
graded Hilbert space $\HH=\HH_+\oplus\HH_-$, and call
$\pi_\pm(a):=\frac{1}{2}(1\pm\gamma)a$ with $\gamma$ the natural
grading. Let $\ket{n,k}_\pm$ be an orthonormal basis of $\HH_\pm$,
with $n\in\N$ and integer $k$ in some finite index set $I_n$ depending
on $n$ (or eventually a multi-index, it makes no difference), and
assume the cardinality of $I_n$ grows at most polynomially in $n$.
Let $D=F|D|$ with absolute value defined by
$$
|D|\ket{n,k}_\pm=\lambda_n\ket{n,k}_\pm \;,
$$
where $\lambda_n>0$ and $\lambda_n^{-1}\to 0$ for $n\to\infty$,
and with sign
$$
F\ket{n,k}_\pm=\ket{n,k}_\mp \;.
$$
It is a typical situation to have $D$ defined by its eigenvalues
as above. Such a $D$ is selfadjoint with domain the set of $L^2$-vectors
$\psi=\sum_{n,k}c_{n,k}^\pm\ket{n,k}_\pm$ such that $|\lambda_nc_{n,k}^\pm|^2$
is summable.
We'll assume also that $\A$ has a finite number of generators $\{x_i\}$,
and that each generator is represented by a \emph{finite} sum
$\pi_\pm(x_i)=\sum_{j,l}x_{i,\pm}^{j,l}$ of bounded shift operators
$x_{i,\pm}^{j,l}$,
$$
x_{i,\pm}^{j,l}\ket{n,k}_\epsilon=\delta_{\epsilon,\pm}
C_{i,\epsilon}^{j,l}(n,k)\ket{n+j,k+l}_\epsilon
$$
with $\{C_{i,\pm}^{j,l}(n,k)\}_{n,k}$ a bounded sequence (for all
fixed $i,j,l,\pm$). Finally, we assume that
\begin{equation}\label{eq:ineqqn}
|C_{i,+}^{j,l}(n,k)-C_{i,-}^{j,l}(n,k)|\leq K_i^{j,l}q^n
\end{equation}
for some $K_i^{j,l}>0$ and a fixed $q$, $0<q<1$.

\begin{lemma}\label{lemma:5.1}
The data $(\A,\HH,F,\gamma)$ is a $1$-summable even Fredholm module.
\end{lemma}
\begin{prova}
By construction
\begin{equation}\label{eq:gFFx}
\gamma F[F,x_{i,+}^{j,l}+x_{i,-}^{j,l}]\ket{n,k}_\pm=
\big\{C_{i,+}^{j,l}(n,k)-C_{i,-}^{j,l}(n,k)\big\}\ket{n+j,k+l}_{\pm}
\end{equation}
and then the commutators
$$
[F,x_i]=(\gamma F)^{-1}\sum\nolimits_{j,l}\gamma F[F,x_{i,+}^{j,l}+x_{i,-}^{j,l}]
$$
because of (\ref{eq:ineqqn}) are of trace class for all $i$. Since
trace--class operators form a (two-sided $*$-)ideal in $\B(\HH)$, from Leibniz rule
$$
[F,ab]=[F,a]b+a[F,b] \;,\qquad\forall\;a,b\in\A\,,
$$
we deduce that $[F,a]$ is of trace class for any $a\in\A$.
\end{prova}

The next quest is for a simple sufficient condition on the data
$(\A,\HH,D,\gamma)$ to have a regular even spectral triple.

\begin{lemma}\label{lemma:5.2}
If \mbox{$\,\lambda_{n+1}-\lambda_n$} is a bounded sequence, the data
$(\A,\HH,D,\gamma)$ is a even regular spectral triple.
\end{lemma}
\begin{prova}
Recall that $\delta:=[|D|,\,.\,]$ and note the identity
\begin{equation}\label{eq:iddpDF}
\delta^p([D,a])=\delta^{p+1}(a)F+|D|[F,\delta^p(a)] \;, \qquad\forall\;p\in\N\,.
\end{equation}
From
\begin{equation}\label{eq:deltap}
\delta^p(x^{j,l}_{i,\pm})\ket{n,k}_\pm=(\lambda_{n+j}-\lambda_n)^p
x^{j,l}_{i,\pm}\ket{n,k}_\pm
\end{equation}
we deduce that $\delta^p(x^{j,l}_{i,\pm})$ are bounded operators for any
$p\in\N$, since $\lambda_{n+j}-\lambda_n=O(1)$. This implies that
$\delta^p(x_i)$ are bounded operators for any $p\in\N$.

From
\begin{equation}\label{eq:DFx}
|D|[F,\delta^p(x^{j,l}_{i,\pm})]
\ket{n,k}_\pm=(\lambda_{n+j}-\lambda_n)^p
\lambda_{n+j}[F,x^{j,l}_{i,\pm}]\ket{n,k}_\pm \;,
\end{equation}
from (\ref{eq:gFFx}) and (\ref{eq:ineqqn}) we deduce that $|D|[F,\delta^p(x_i)]$
is bounded for all $p\in\N$. Since for $a=x_i$ both summands in (\ref{eq:iddpDF})
are bounded, $\delta^p([D,x_i])$ is also bounded for all $p\in\N$.
Thus, by Leibniz rule we have the set inclusion
$$
\A\cup [D,\A]\subset\opz
$$
and $(\A,\HH,D)$ is a regular spectral triple. That it is even
is obvious.
\end{prova}

As a particular case we have the following, which is a simple consequence
of the observation that the series $\sum_{n\in\N}(n+1)^{-s}$ converges
to the Riemann zeta function $\zeta(s)$ for any $s>1$, and that
$\mr{Res}_{s=1}\zeta(s)=1$.

\begin{lemma}\label{lemma:5.3}
If $\lambda_n=n+1$ and its multiplicity is a polynomial in $n$ of order
$\alpha$, then $(\A,\HH,D)$ is a regular spectral triple with metric
dimension $\alpha+1$.
\end{lemma}
\begin{prova}
Simply
$$
\tr(|D|^{-s})\propto\sum\nolimits_{n\in\N}\big\{(n+1)^{-s+\alpha}+\;
\textrm{lower order monomials in}\;n\big\} \;.
$$
Thus $|D|^{-s}$ is of trace class for $s>\alpha+1$. The Dixmier trace
of $|D|^{-\alpha-1}$ is the residue in $s=\alpha+1$ of the meromorphic
extension of $\tr(|D|^{-s})$, which is $\zeta(s-\alpha)$ plus terms which
are holomorphic for $\mr{Re}\,s>\alpha$ (coming from monomials
of order $\leq\alpha-1$ in the multiplicity of $\lambda_n$).
Then $\tr_\omega(|D|^{-\alpha-1})\propto\mr{Res}_{s=1}\zeta(s)\neq 0$.
\end{prova}

\begin{lemma}\label{lemma:5.4}
If $\lambda_n=O(q^{-n})$ and if, for all $j\neq 0$,
\begin{equation}\label{eq:Kprime}
|C_{i,\pm}^{j,l}(n,k)|\leq K'^{j,l}_{i,\pm}\,q^n \;,
\end{equation}
for some constants $K'^{j,l}_{i,\pm}>0$,
the data $(\A,\HH,D,\gamma)$ is a even
Lipschitz-regular spectral triple.  The triple
is regular if $\lambda_n$ diverges at most
polynomially.
\end{lemma}
\begin{prova}
By (\ref{eq:iddpDF}) to prove the first statement it is sufficient to
prove that $\delta(x_i)$ and $|D|[F,x_i]$ are bounded. As above using
(\ref{eq:DFx}) with $p=0$, (\ref{eq:gFFx}) and (\ref{eq:ineqqn}) one
proves that $|D|[F,x_{i,+}^{j,l}+x_{i,-}^{j,l}]$ are bounded operators,
and then $|D|[F,x_i]$ are bounded.

Now using (\ref{eq:deltap}) we prove that
$\delta^p(x^{j,l}_{i,\pm})$ is zero for $j=0$, while
for $j\neq 0$
$$
\Big|\rule{0pt}{11pt}_\pm\big<n+j,k+l\big|\delta^p(x^{j,l}_{i,\pm})
\big|n,k\big>\rule{0pt}{11pt}_\pm\Big| \leq 
K'^{j,l}_{i,\pm}(\lambda_{n+j}-\lambda_n)^p q^n\;.
$$
For $p=1$, the sequence $(\lambda_{n+j}-\lambda_n)^p q^n$ is
bounded. Thus $\delta(x_{i,\pm}^{j,l})$ are bounded operators,
and $\delta(x_i)$ are bounded too. This proves the first statement.

If $\lambda_n$ grows at most polynomially in $n$,
then $(\lambda_{n+j}-\lambda_n)^p q^n$ is a bounded sequence
for all $p$, and by previous inequality $\delta^p(x^{j,l}_{i,\pm})$
are bounded operators for all $p$.
Using (\ref{eq:DFx}), (\ref{eq:gFFx}) and (\ref{eq:ineqqn}) one
proves that also $|D|[F,\delta^p(x_{i,+}^{j,l}+x_{i,-}^{j,l})]$ are
bounded operators for all $p$, and then $|D|[F,\delta^p(x_i)]$ are bounded.
By (\ref{eq:iddpDF}), this proves that $x_i$ and $[D,x_i]$ are in
$\opz$, and then the spectral triple is regular.
\end{prova}

\section{Local index formula for quantum homogeneous spaces}\label{sec:2.4}
Typically, Fredholm modules associated to isospectral Dirac operators
on even-dimensional quantum homogeneous spaces satisfy $[F,a]\in\op$ for all $a\in\A$
(cf.~\cite{DDLW07,DDL06}), and this results in a drastic simplification of
the local index formula. The following proposition clarifies this last statement.
We use the notations of Theorem \ref{thm:CM}.

\begin{prop}
Let $(\A,\HH,D,\gamma)$ be a regular even spectral triple with finite metric dimension,
with dimension spectrum $\Sigma$ made of simple poles only, $0\notin\Sigma$ and
$[F,a]\in\op$ for all $a\in\A$ ($F=|D|D^{-1}$). Then, the Connes-Moscovici
local cocycle $\varphi$ has a unique non-vanishing component $\varphi_0=\mr{ch}^F_0$.
\end{prop}
\begin{prova}
Let $n\in 2\N$ be the metric dimension. Then
$\tr(\gamma a_0|D|^{-2s})$ is holomorphic for $\mr{Re}\,2s>n$,
and in this half plane coincides with the function $\psi(s):=\frac{1}{2}\tr(\gamma F[F,a_0]|D|^{-2s})$
(just use $F^2=1$, $[F,|D|]=0$ and the cyclicity of the trace).
But $\psi(s)$ is holomorphic on $\C$, since $[F,a_0]\in\op$, and then $\psi(s)$ is
just the holomorphic extension on $\C$ of the zeta-function $\tr(\gamma a_0|D|^{-2s})$.
This implies that
\begin{equation*}
\varphi_0(a_0):=\mr{Res}_{s=0}s^{-1}\psi(s)=\psi(0)=\mr{ch}^F_0(a_0)\;.
\end{equation*}
It remains to prove that all other components $\varphi_j$ of the cocycle $\varphi$
vanish. The generic component, for $j$ even greater than $0$, is given by:
$$
\varphi_j(a_0,\ldots,a_j)=\frac{1}{2}\sum_{k\in\N^j}\frac{(-1)^k}{k_1!\ldots k_j!}\,\alpha_k\,
\nint\gamma a_0[D,a_1]^{(k_1)}\ldots[D,a_j]^{(k_j)}|D|^{-(2|k|+j)}\;\;.
$$
Let $\chi(a)=\frac{1}{2}(a+FaF)$. Then $\chi(a)$ commutes with $F$ and $\gamma$,
and $a-\chi(a)=\frac{1}{2}F[F,a]\in\op$.

Smoothing operators do not contribute to $\varphi_j$, so we can replace
$a_i$ with $\chi(a_i)$ in the formula for $\varphi_j$ and obtain:
$$
\varphi_j(a_0,\ldots,a_j)=\frac{1}{2}\sum_{k\in\N^j}\frac{(-1)^k}{k_1!\ldots k_j!}\,\alpha_k\,
\nint\gamma F^j\chi(a_0)\delta\chi(a_1)^{(k_1)}\ldots\delta\chi(a_j)^{(k_j)}|D|^{-(2|k|+j)}\;\;.
$$
But $F^j=1$ since $j$ is even. Since $\chi(a)$ commutes with $F$ and $\gamma$,
it commutes also with $\sigma_\pm:=\frac{1}{2}(1\pm\gamma)F$.
The operators $\sigma_\pm$ commute with $|D|$ as well, so we
conclude that $\xi:=\chi(a_0)\delta\chi(a_1)^{(k_1)}\ldots\delta\chi(a_j)^{(k_j)}|D|^{-(2|k|+j)}$
commutes with $\sigma_\pm$.
Using the identity $\gamma=[\sigma_+,\sigma_-]$ we can write
$$
\varphi_j(a_0,\ldots,a_j)=\frac{1}{2}\sum_{k\in\N^j}\frac{(-1)^k}{k_1!\ldots k_j!}\,\alpha_k\,
\mr{Res}_{s=0}\tr\bigl([\sigma_+,\sigma_-\xi |D|^{-s}]\bigr)\;\;.
$$
Now the expression under the trace is of trace class for $\mr{Re}\,s>n$ and in
such a half plane is also traceless, being a commutator.
This proves that $\varphi_j(a_0,\ldots,a_n)=0$ for all even $j>0$.
\end{prova}

For some even-dimensional quantum spaces, for example quantum complex projective
spaces, there is a strong belief that the situation can drastically change
by considering twisted homologies (see for example \cite{NT05} for the case
of the standard Podle\'s sphere $\CP^1_q$).

\section{Hopf algebras and equivariance}\label{chap:HA}

An important notion, which enters in the construction of many noncommutative
spectral triples, is the one of \emph{equivariance}. It is inspired by the
example of coset spaces. Given a (compact) Riemannian $G$-symmetric space
$M$, with $G$ a Lie group and $K$ the isotropy group of a fixed point $p\in M$,
it is natural to look for the existence of a homogeneous spin structure~\cite{Fri00,CG88,BILP03},
which is a commutative diagram
\begin{center}\begin{tabular}{c}
\begindc{\commdiag}[5]
 \obj(1,11)[A]{$\mr{Spin}(T_pM)$}
 \obj(1,1)[B]{\raisebox{14pt}{$\mr{SO}(T_pM)$}}
 \obj(14,11)[C]{$K$}
 \mor{A}{B}[10,16]{\scriptsize $\lambda$}
 \mor{C}{A}{\scriptsize $\widetilde{\mr{Ad}}$}[\atright,\dasharrow]
 \mor{C}{B}{\scriptsize Ad}
\enddc
\end{tabular}\end{center}

\vspace{-10pt}

\noindent
where the arrows are group homomorphisms: $\lambda$ is the double covering of
the $SO$ group, $\mr{Ad}$ is the restriction of the adjoint action of $G$
on its Lie algebra, and the dashed arrow is the lift one looks for.
Homogeneous spin structures correspond to commutative `equivariant'
spectral triples. One identifies $C(M)$ with $C_K(G)$, and spinors
with $L^2$ sections of the bundle $G\times_K V\to M$ associated to
the representation $V$ of $K$ obtained by composing the map $\widetilde{\mr{Ad}}$
with the defining representation of the spin group. The representation
$\pi$ of $C(M)$ on spinors is given by pointwise multiplication, while
a representation $U$ of $g\in G$ is given by the pullback of left
multiplication by $g^{-1}$, $\{U(g)\psi\}(m):=\psi(g^{-1}m)$; the
representation $U$ is unitary since $G$ acts by isometries.
The map $\pi$ intertwines the representation $U$ of $G$ with
the natural action of $G$ on $C(M)$:
\begin{equation}\label{eq:Gcov}
U(g)\pi(f)U(g)^{-1}=\pi(L_gf)\;,
\end{equation}
where $(L_gf)(m):=f(g^{-1}m)$ for all $g\in G$, $f\in C(M)$ and $m\in M$.
When a representation $\pi$ satisfies (\ref{eq:Gcov}) we say that it is
\emph{$G$-covariant}. A fundamental characteristic of homogeneous spin structures
is that they have associated an invariant Dirac operator (i.e.~commuting
with $U(g)$ for all $g\in G$), and eventually an invariant grading if we
are in the even--dimensional case; the charge conjugation $J$, instead, is not
invariant but intertwines the representation $U$ of $G$ with the representation
\begin{equation}\label{eq:realequiv}
g\mapsto U(g^{-1})^\dag
\end{equation}
(in this case $U(g^{-1})^\dag=U(g)$ by unitarity). All these properties can
be formalized to get an abstract definition of equivariant spectral triple
(see e.g.~\cite{Sit03}).

Any $G$-covariant representation of $C(M)$ is equivalent to a representation
of the `transformation group $C^*$-algebra' $G\ltimes C(M)$ (replaced by
\emph{crossed product algebras} in the noncommutative case), and both are equivalent
to Mackey's system of imprimitivity (cf.~\cite{Lan06} and reference [60] therein),
thus establishing a link with quantization of mechanical systems.

By Serre-Swan theorem, $L^2(G\times_KV)$ is the completion of a finitely
generated projective $C(M)$-module. Thus, another way to construct
this equivariant spin representation $\pi$ is starting from $C(M)^N$
and looking for a suitable invariant projection, i.e.~a bounded linear map
$$
p:C(M)^N\to C(M)^N
$$
such that $p=p^2=p^*$ and $\rho(g)p=p\rho(g)$ for all $g\in G$ and with
$\rho$ a suitable representation. The Hilbert space completion of $C(M)^Np$
will carry automatically a $G$-equivariant representation of $C(M)$, and
one has just to guess the correct $p$ to get spinors.
This is the strategy we'll adopt for the quantum orthogonal $4$-sphere
$S^4_q$ in Chapter \ref{chap:S4q}.

We are going to discuss the noncommutative analogue of $G$-covariant
representations (and to define equivariant spectral triples in full
generality). Before this, we shall briefly discuss the objects which
replace groups in the noncommutative case, i.e.~Hopf algebras.
Reference textbooks for this topic are~\cite{Abe80,Swe69,Maj95,Var01,FGB05}.

\subsection{Hopf algebras and crossed products}
We recall that by a ($*$-)algebra we mean an associative algebra over $\C$
with unity (and involution). Let $\A$ be a finite-dimensional algebra.
On the dual space $\mc{C}=\mr{Hom}(\A,\C)$ we can define two operations as follows.
Writing $f(a)=:\inner{a,f}$ we call $\Delta:\mc{C}\to\mc{C}\otimes\mc{C}$
and $\epsilon:\mc{C}\to\C$ the linear maps uniquely defined by
\begin{subequations}\label{eq:De}
\begin{align}
\inner{a\otimes b,\Delta f}&=\inner{ab,f} \\
\epsilon(f)&=\inner{1,f}
\end{align}
\end{subequations}
for all $a,b\in\A$ and $f\in\mc{C}$, and known as `coproduct'
and `counit' respectively. The associativity and the unitality of $\A$
imply, respectively, the identities
\begin{subequations}\label{eq:coa}
\begin{align}
(\Delta\otimes id)\circ\Delta &=(id\otimes\Delta)\circ\Delta \;,\label{eq:coaA} \\
(\epsilon\otimes id)\circ\Delta &=(id\otimes\epsilon)\circ\Delta=id \;. \label{eq:coaB}
\end{align}
\end{subequations}
A vector space with two linear maps $(\Delta,\epsilon)$ satisfying (\ref{eq:coa})
is called a `coalgebra'. In the infinite-dimensional case, the dual space
$\mr{Hom}(\A,\C)$ of $\A$ is not always a coalgebra with the operations defined
by (\ref{eq:De}), since in general the map $\Delta$ sends $\mr{Hom}(\A,\C)$ into
$\mr{Hom}(\A\otimes\A,\C)$, which is a suitable completion of the algebraic tensor
product $\mr{Hom}(\A,\C)\otimes\mr{Hom}(\A,\C)$. This is the reason why one prefers
to talk of \emph{dual pair} $(\A,\mc{C})$ of algebra-coalgebra, meaning that $\A$
is an algebra, $\mc{C}$ a coalgebra, and there is a bilinear map $\inner{\,,\,}:
\A\otimes\mc{C}\to\C$ such that (\ref{eq:De}) is satisfied. The dual pair is
\emph{non-degenerate} whenever the pairing is non-degenerate, i.e.~if the unique
element of $\A$ (resp.~$\mc{C}$) which has vanishing pairing with all the elements
of $\mc{C}$ (resp.~$\A$) is the zero. In the finite-dimensional case there is
a unique coalgebra in non-degenerate dual pairing with any $\A$, the coalgebra
$\mr{Hom}(\A,\C)$ with operations (\ref{eq:De}); in the infinite-dimensional
case this is no more true: there can be many non-isomorphic coalgebras in
non-degenerate dual pairing with a given algebra. An example where this
situation occurs is the coordinate algebra of the quantum $SU(2)$ group
\cite{KS97}. In a dual pair $(\A,\mc{C})$, $\A$ is commutative if and only
if the coproduct of $\mc{C}$ is symmetric (we say that $\mc{C}$ is `cocommutative').

If $\A$ is an algebra and $(V,W)$ a non-degenerate dual pair of vector spaces,
any left $\A$-module structure on $W$ corresponds to a right $\A$-module structure
on $V$ (and vice versa) through the formula
\begin{equation}\label{eq:dualaz}
\inner{va,w}=\inner{v,aw}\;,\qquad\forall\;a\in\A,\;v\in V,\;w\in W\,.
\end{equation}
The vector space $\A$ carries a natural left/right $\A$-module structure
given by left/right multiplication, sometimes called `left/right regular
action'. If $\mc{C}$ is a coalgebra in non-degenerate dual pairing with $\A$,
a right/left $\A$-module structure on $\mc{C}$ dual to the left/right regular
action is defined using (\ref{eq:dualaz}). The corresponding action
is called `right/left canonical', denoted $\za$ and $\az$ respectively,
and given by ($a\in\A$, $c\in\mc{C}$)
\begin{equation}\label{eq:canaz}
c\za a=\inner{a,c_{(1)}}c_{(2)}\;,\qquad
a\az c=c_{(1)}\inner{a,c_{(2)}} \;.
\end{equation}
As usual we adopt Sweedler notation for the coproduct, $\Delta c=c_{(1)}\otimes c_{(2)}$
with summation implied.

A natural question is: if we have two representations $\rho_i:\A\to\mr{Aut}(V_i)$
of an algebra, how to construct a representation $\psi:\A\to\mr{Aut}(V_1\otimes V_2)$
on the tensor product of the vector spaces $V_i$ starting from $\rho_{1,2}$.
The tensor product $\rho_1\otimes\rho_2$ is a representation of $\A\otimes\A$
on $V_1\otimes V_2$, so all what we need is an algebra morphism $\Delta:\A\to\A\otimes\A$;
the map
\begin{equation}\label{eq:Htp}
\psi:=(\rho_1\otimes\rho_2)\circ\Delta
\end{equation}
gives the representation we were looking for, and is called \emph{Hopf-tensor
product} of $\rho_1$ and $\rho_2$. With $n$ representations $\rho_i$ we can
recursively define an $n$-fold tensor product, and the result does not depend
on the order we multiply the $\rho_i$s exactly when $\Delta$ is coassociative,
i.e.~satisfies (\ref{eq:coaA}). Similarly, two representations $\rho_{1,2}$
of $\A$ can be constructed from a representation $\psi$ of $\A\otimes\A$ if
we have an algebra morphism $\epsilon:\A\to\C$, simply putting
$\rho_1=\psi\circ(id\otimes\epsilon)$ and $\rho_2=\psi\circ(\epsilon\otimes id)$;
these operations are inverses of (\ref{eq:Htp}) if $\epsilon$ is a counit.
Motivated by these considerations, one calls \emph{bialgebra} an algebra $\A$
which is also a coalgebra and such that coproduct and counit are algebra morphisms.
Note that the dual object of a bialgebra is itself a bialgebra.

We remark that since in general $\Delta$ is not symmetric, the Hopf tensor
products $V_1\otimes V_2$ and $V_2\otimes V_1$ are different. The order
can be crucial: the difference between a first 
attempt of
constructing a spectral triple for $SU_q(2)$ in \cite{Gos01} and a later
successful construction in \cite{DLS05} was exactly the order in an
Hopf tensor product of representations.

A source of commutative examples of bialgebras is from group theory.
If $\A$ is the algebra of functions on a finite group $G$ (or also
representative functions of a compact matrix group), a coproduct and
counit on $\A$ can be defined by
$$
(\Delta f)(g,g'):=f(gg')\;,\qquad
\epsilon(f):=f(e) \;,
$$
with $f\in\A$, $g\in G$ and $e$ the neutral element.
A further operation $S:\A\to\A$ is obtained by dualizing the
operation of taking the inverse,
$$
\{S(f)\}(g)=f(g^{-1})\;,
$$
and satisfies
\begin{equation}\label{eq:HA}
S(f_{(1)})f_{(2)}=f_{(1)}S(f_{(2)})=\epsilon(f)\;,\qquad\forall\;f\in\A\,.
\end{equation}
A bialgebra $\A$ with a map $S:\A\to\A$ satisfying (\ref{eq:HA}) is called
\emph{Hopf algebra}, and $S$ is called \emph{antipode}. A remarkable fact
is that when $S$ exists it is unique, being the inverse of the identity map
under the convolution product in $\mr{End}(\A)$, cf.~equation (\ref{eq:convprod}).
Given a non-degenerate dual pair of bialgebras $(\A,\mc{U})$ and an antipode $S$
on $\A$, one defines an antipode on $\mc{U}$ as
$$
\inner{S(a),h}=\inner{a,S(h)}\;,\qquad a\in\A,\;h\in\mc{U}\,,
$$
and talks of (non-degenerate) dual pair of Hopf-algebras.

If $\A$ is a Hopf algebra, other natural actions of $\A$ on itself
besides the regular ones are the left/right adjoint actions, defined
by
$$
\mr{Ad}_L(a)b:=a_{(1)}bS(a_{(2)})\;,\qquad
\mr{Ad}_R(a)b:=S(a_{(1)})ba_{(2)}\;,
$$
for all $a,b\in\A$.

The last point to discuss is the involution. Suppose $\A$ is a Hopf algebra,
and as an algebra it possesses an involution. Let $\rho_1$ and $\rho_2$ be two
$*$-representations. Their Hopf tensor product (\ref{eq:Htp}) is
a $*$-representation (for any pair $\rho_1,\rho_2$) if and only if
\begin{subequations}\label{eq:star}
\begin{equation}
\Delta(a^*)=(a_{(1)})^*\otimes(a_{(2)})^*
\end{equation}
for all $a\in\A$.
Similarly given a $*$-representation $\psi$ of $\A\otimes\A$, using
$\epsilon$ as explained before we get two $*$-representations of $\A$
(for any $\psi$) if and only if
\begin{equation}
\epsilon(a^*)=\overline{\epsilon(a)}
\end{equation}
\end{subequations}
for all $a\in\A$. If (\ref{eq:star}) are satisfied, we call $\A$
a \emph{Hopf $*$-algebra}. From the axioms it automatically follows
that the antipode $S$ of a Hopf $*$-algebra is invertible, with
inverse
$$
S^{-1}(a)=S(a^*)^*\;,\qquad\forall\;a\in\A\;.
$$
Given a non-degenerate dual pair of Hopf algebras $(\A,\mc{U})$, if the first
is a Hopf $*$-algebra, also the second becomes a Hopf $*$-algebra with
the involution
$$
\inner{a,h^*}=\overline{\inner{S(a)^*,h}} \;,\qquad\forall\;a\in\A,\;h\in\mc{U}\,.
$$

Inspired by the (left) canonical and (left) adjoint actions of a Hopf $*$-algebra,
if $(\mc{U},\Delta,\epsilon,S)$ is a Hopf $*$-algebra and $\A$ a $*$-algebra,
we say that $\A$ is a (left) \emph{$\mc{U}$-module $*$-algebra} if there is a
(left) action `$\az$' of $\mc{U}$ on $\A$ satisfying
$$
h\az ab=(h_{(1)}\az a)(h_{(2)}\az b)\;,\qquad
h\az 1=\varepsilon(h)1\;,\qquad
h\az a^*=\{S(h)^*\az a\}^*\;,
$$
for all $h\in\mc{U}$ and $a,b\in\A$.
The (left) \emph{crossed product} $\A\rtimes\mc{U}$ is the $*$-algebra generated
by $\A$ and $\mc{U}$ with crossed commutation relations
$$
ha=(h_{(1)}\az a)h_{(2)}\;,\quad\forall\;h\in\mc{U},\;a\in\A\,.
$$
Right module $*$-algebras and right crossed products are defined in
a similar way.

\subsection{Equivariant projective modules}\label{sec:3.2}
Let $\mc{U}$ be a Hopf $*$-algebra, $\A$ be an $\mc{U}$-module $*$-algebra and
$\varphi:\A\to\C$ be an invariant faithful state (i.e.~$\varphi$ is linear,
$\varphi(a^*a)>0$ for all non-zero $a\in\A$, and $\varphi(h\az a)=\epsilon(h)
\varphi(a)$ for all $a\in\A$ and $h\in\mc{U}$). Suppose also that there exists
$\kappa\in\mr{Aut}(\A)$ such that the `twisted' cyclicity
$$
\varphi(ab)=\varphi\bigl(b\,\kappa(a)\bigr)
$$
holds for all $a,b\in\A$. Instances of this situation are provided by
subalgebras of compact quantum group algebras with $\varphi$ the Haar
state and $\kappa$ the modular automorphism. KMS states in Thermal Quantum
Field Theory provide additional examples.

For $N\in\N$, let $\A^N:=\A\otimes\C^N$ be the linear space with
elements $v=(v_1,\ldots,v_N)$, $v_i\in\A$, and 
scalar product given by
\begin{equation}\label{eq:in}
\inner{v,w}:=\sum\nolimits_{i=1}^N\varphi(v_i^*w_i)\;.
\end{equation}

\begin{lemma}\label{lemma:3.2}
Let $\sigma:\mc{U}\to\mr{Mat}_N(\C)$ be a $*$-representation.
The formul{\ae}:
\begin{equation}\label{eq:starrep}
(a\dN v)_i:=av_i\;,\qquad
(h\dN v)_i:=\sum\nolimits_{j=1}^N(h_{(1)}\az v_j)\sigma_{ij}(h_{(2)})\;,
\end{equation}
for all $a,v_i\in\A$ and $h\in\mc{U}$ (and $i=1,\dots,N$), define a $*$-representation
of the crossed product algebra $\A\rtimes\mc{U}$ on the linear space $\A^N$.
\end{lemma}
\begin{prova}
The inner product allows us to define the adjoint of an element of $\A\rtimes\mc{U}$
in the representation on $\A^N$. For $x\in\mr{End}(\A^N)$ its adjoint, denoted with
$x^\dag$, is defined by
$$
\inner{\smash[t]{x^\dag\dN v},w}:=\inner{v,x\dN w}\;, \qquad \forall \;, v,w\in\A^N\;.
$$
Recall that being a $*$-representation means that $x^\dag\dN\,v=x^*\dN\,v$ for any
operator $x$ and any $v\in\A^N$.

The non-trivial part of the proof consists in showing that $h^\dag\dN\,v=h^*\dN\,v$ for all
$h\in\mc{U}$ and $v\in\A$. It is enough to take $N=1$. For $N>1$ we
are considering the Hopf tensor product of the $N=1$ representation with
a matrix representation that is a $*$-representation by hypothesis.

The $\mc{U}$-invariance of $\varphi$ implies:
$$
\epsilon(h)\inner{v,w}=\varphi\bigl(h\az(v^*w)\bigr)=
\varphi\bigl((h_{(1)}\az v^*)(h_{(2)}\az w)\bigr)\;.
$$
But $h_{(1)}\az v^*=\{S(h_{(1)})^*\az v\}^*$ by definition of module
$*$-algebra. Then,
$$
\epsilon(h)\inner{v,w}=\inner{S(h_{(1)})^*\dN v,h_{(2)}\dN w}=
\big<v,S(h_{(1)})^{*\dag}h_{(2)}\dN w\big>\;.
$$
We deduce that for all $h\in\mc{U}$ one has that
\begin{equation}\label{eq:newS}
S(h_{(1)})^{*\dag}h_{(2)}=\epsilon(h)\;.
\end{equation}
Recall that the convolution product `$\star$' for any $F,G\in\mr{End}(\mc{U})$ 
is defined by
\begin{equation}\label{eq:convprod}
(F\star G)(h):=F(h_{(1)})G(h_{(2)})\qquad\forall\;h\in\mc{U} \;;
\end{equation}
and $(\mr{End}(\mc{U}),\star)$ is an associative algebra with unity given
by the endomorphism $h\mapsto\epsilon(h)1_{\mc{U}}$,  with $S$  a
left and right inverse for $id_{\mc{U}}$ in $(\mr{End}(\mc{U}),\star)$, that is
$$
S\star id_{\mc{U}}=1_{\mc{U}}\epsilon=id_{\mc{U}}\star S\;.
$$
Let $S'\in\mr{End}(\mc{U})$ be the composition $S':=\dag\circ *\circ S$.
Equation \eqref{eq:newS} implies that $S'$ is a left inverse for $id_{\mc{U}}$:
$$
S'\star id_{\mc{U}}=1_{\mc{U}}\epsilon\;.
$$
Applying $\star S$ to the right of both members of this equation and using
$id_{\mc{U}}\star S=1_{\mc{U}}\epsilon$ we get $S'=S$ as endomorphisms of $\mc{U}$,
i.e.~$S(h)^{*\dag}=S(h)$ for all $h\in\mc{U}$.

Now, the antipode of a Hopf $*$-algebra is invertible,
with $S^{-1}=*\circ S\circ *$, thus we arrive at $h^{*\dag}=h$ for all $h\in\mc{U}$.
Replacing $h$ with $h^*$ we prove that $h^\dag=h^*$ for all $h\in\mc{U}$, and
this concludes the proof.
\end{prova}

Now, let $e=(e_{ij})\in\mr{Mat}_N(\A)$ be an $N\times N$ matrix
with entries $e_{ij}\in\A$. Let $\pi:\A^N\to\A^N$ be the (linear)
endomorphism defined by:
\begin{equation}\label{eq:pi}
\pi(v)_j:=\sum\nolimits_{i=1}^Nv_ie_{ij}\;,
\end{equation}
for all $v\in\A^N$ and $j=1,\dots,N$.
Since $\A$ is associative, left and right multiplication commute
and $\pi(av)=a\pi(v)$ for all $a\in\A$ and $v\in\A^N$.
Thus we have the following Lemma.

\begin{lemma}
The map $\pi$ defined by \eqref{eq:pi} is an $\A$-module map.
\end{lemma}

Recall that an endomorphism $p$ of an inner
product space $V$ is a projection (not necessarily orthogonal)
if $p\circ p=p$. A projection $p$ is \emph{orthogonal}
if the image of $p$ and $id_V-p$ are orthogonal with respect
to the inner product of $V$, and this happens exactly when $p^\dag=p$.

A simple computation shows that the map $\pi$ in 
\eqref{eq:pi} is a projection if{}f $e^2=e$, that is the matrix
$e\in\mr{Mat}_N(\A)$ is an idempotent. Now we use
the twisted-cyclicity of $\varphi$ to deduce:
$$
\big<v,\pi^\dag(w)\big>=\inner{\pi(v),w}=
\sum\nolimits_{ij}\varphi(e_{ij}^*v_i^*w_j)=
\sum\nolimits_{ij}\varphi\bigl(v_i^*w_j\kappa(e_{ij}^*)\bigr) \;, 
$$
for all $v,w\in\A^N$. Hence the adjoint $\pi^\dag$ of the endomorphism $\pi$
is given by
$$
\pi^\dag(w)_i=\sum\nolimits_{j=1}^Nw_j\kappa(e_{ij}^*)\;.
$$
Let $e^*$ be the matrix with entries $(e^*)_{jk}:=e_{kj}^*$.
We have the following trivial Lemma.

\begin{lemma}\label{lemma:next}
The endomorphism $\pi$ in \eqref{eq:pi} is an orthogonal
projection if{}f $\,e^2=e=\kappa(e^*)$.
\end{lemma}

Next, we  determine a sufficient condition for the endomorphism $\pi$ to be
not only an $\A$-module map, but also an $\mc{U}$-module map.

\begin{lemma}
With `$\phantom{|}^t$' denoting transposition, if
\begin{equation}\label{eq:cov}
h\az e=\sigma(h_{(1)})^t\,e\,\sigma(S^{-1}(h_{(2)}))^t
\end{equation}
for all $h\in\mc{U}$, then the endomorphism $\pi$ in
\eqref{eq:pi} is an $\mc{U}$-module map.
\end{lemma}
\begin{prova}
Equation (\ref{eq:cov}) can be rewritten as,
$$
h\az e_{ij}=\sum\nolimits_{kl}\sigma_{ki}(h_{(1)})\,e_{kl}
\,\sigma_{jl}(S^{-1}(h_{(2)}))\;;
$$
by using it into the definition \eqref{eq:starrep}
one checks that $\pi(h\dN v)=h\dN\pi(v)$ for all $h\in\mc{U}$ and $v\in\A^N$.
Indeed
\begin{align*}
\{h\dN\pi(v)\}_i &=\sum\nolimits_{jk}(h_{(1)}\az v_je_{jk})\sigma_{ik}(h_{(2)}) \\
 &=\sum\nolimits_{jk}(h_{(1)}\az v_j)(h_{(2)}\az e_{jk})\sigma_{ik}(h_{(3)}) \\
 &=\sum\nolimits_{jklm}(h_{(1)}\az v_j)\sigma_{lj}(h_{(2)})\,e_{lm}\,
   \sigma_{km}\bigl(S^{-1}(h_{(3)})\bigr)\sigma_{ik}(h_{(4)}) \\
 &=\sum\nolimits_{jlm}(h_{(1)}\az v_j)\sigma_{lj}(h_{(2)})\,e_{lm}\,
   \sigma_{im}\bigl(h_{(4)}S^{-1}(h_{(3)})\bigr) \;, \\
 &=\sum\nolimits_{lm}(h_{(1)}\dN v)_l\,e_{lm}\,
   \sigma_{im}\bigl(h_{(3)}S^{-1}(h_{(2)})\bigr) \;, \\
\intertext{and since $\sigma_{im}\bigl(h_{(3)}S^{-1}(h_{(2)})\bigr)=\epsilon(h_{(3)})\delta_{im}$,
last equation is equal to}
 &=\sum\nolimits_{l}(h_{(1)}\dN v)_l\,e_{li}\,\epsilon(h_{(2)})=\pi(h\dN v)_i \;,
\end{align*}
which concludes the proof.
\end{prova}

When Lemma \ref{lemma:next} and equation (\ref{eq:cov}) are satisfied,
the orthogonal projections $\pi$ and \linebreak $\pi^\perp=1-\pi\,$
split $\A^N$ into the orthogonal sum of two $*$-subrepresentations
$\pi(\A^N)$ and $\pi^\perp(\A^N)$ of $\A\rtimes\mc{U}$.
The next lemma gives a (quite obvious) sufficient condition for $\pi(\A^N)$
and $\pi^\perp(\A^N)$ to be inequivalent as representations of $\A$.
Recall that an isomorphism of $\A$-modules is an invertible $\A$-linear
map, so isomorphic left modules correspond to equivalent representations.

\begin{lemma}\label{lemma:inv}
Let $(\A,\HH,F,\gamma)$ be an even Fredholm module over $\A$,
and $F^+_e$ the operator in (\ref{eq:Fpe}).
If $\;\mr{Index}(F^+_e)\neq 0\,$, the left $\A$-modules $\pi(\A^N)$ and
$\pi^\perp(\A^N)$ are not equivalent.
\end{lemma}
\begin{prova}
Suppose $\pi(\A^N)$ and $\pi^\perp(\A^N)$ are isomorphic $\A$-modules,
then $\,[e]=[1-e]\,$ and $\,\mr{Index}(F^+_{1-e})=\mr{Index}(F^+_e)\,$.
But from equation (\ref{eq:ind}), $\,\mr{Index}(F^+_{1-e})=-\,\mr{Index}(F^+_e)\,$.
Hence $\,\mr{Index}(F^+_e)=0\,$, and this concludes the proof by contradiction.
\end{prova}

\noindent
Notice that if (\ref{eq:cov}) is satisfied, the \emph{right} module
$e\A^N$ carries a left action `$\aaz$' of $\mc{U}$ given by
\begin{equation}\label{eq:RmodAz}
h\aaz v:=\sigma\bigl(S^{-1}(h_{(1)})\bigr)^t h_{(2)}\az v\;,
\end{equation}
for all $h\in\mc{U}$ and $v=(v_1,\ldots,v_N)^t\in e\A^N$ (with row by column
multiplication implied). One can easily
check that if $ev=v$, i.e.~$v\in e\A^N$, the vector
$h\aaz v$ is still in $e\A^N$, i.e.~$e(h\aaz v)=(h\aaz v)$.

\noindent
Indeed, by (\ref{eq:cov}) and (\ref{eq:RmodAz})
\begin{align*}
h\aaz (ev)&=\sigma\bigl(S^{-1}(h_{(1)})\bigr)^t h_{(2)}\az (ev) \\
&=\sigma\bigl(S^{-1}(h_{(1)})\bigr)^t (h_{(2)}\az e)(h_{(3)}\az v) \\
&=\sigma\bigl(S^{-1}(h_{(1)})\bigr)^t\sigma(h_{(2)})^t\,e\,
  \sigma\bigl(S^{-1}(h_{(3)})\bigr)^t(h_{(4)}\az v) \\
&=\sigma\bigl(h_{(2)}S^{-1}(h_{(1)})\bigr)^t\,e\,
  \sigma\bigl(S^{-1}(h_{(3)})\bigr)^t(h_{(4)}\az v) \\
&=\epsilon(h_{(1)})\,e\,\sigma\bigl(S^{-1}(h_{(2)})\bigr)^t(h_{(3)}\az v) \\
&=e\,\sigma\bigl(S^{-1}(h_{(1)})\bigr)^t(h_{(2)}\az v)=e(h\aaz v) \;,
\end{align*}
which when $ev=v$ implies $e(h\aaz v)=(h\aaz v)$.

\subsection{Equivariant spectral triples}
Let $(\mc{U},\epsilon,\Delta,S)$ be a Hopf $*$-algebra and suppose $\A$ is a
(left) $\mc{U}$-module $*$-algebra. The data $(\A,\HH,D)$ is called a
\emph{$\mc{U}$-equivariant spectral triple} if
\begin{itemize}
\item[(i)] there is a dense subspace $\mc{M}$ of $\HH$ carrying a $*$-representation $\pi$ of $\A\rtimes\mc{U}$,
\item[(ii)] $D$ is a 
selfadjoint operator with compact resolvent and with domain containing $\mc{M}$,
\item[(iii)] $\pi(a)$ and $[D,\pi(a)]$ extend to bounded operators on $\HH$ for all $a\in\A$,
\item[(iv)] $[D,\pi(h)]=0$ on $\mc{M}$ for any $h\in\mc{U}$.
\end{itemize}
An equivariant spectral triple is called \emph{even} if there exists a grading
$\gamma$ on $\HH$ (i.e.~a bounded operator satisfying $\gamma=\gamma^*$ and $\gamma^2=1$)
such that the Dirac operator is odd and the crossed product algebra is even:
$$
\gamma D+D\gamma=0\;,\quad\qquad \pi(t)\gamma=\gamma\pi(t) \quad \forall\;t\in\A\rtimes\mc{U}\;.
$$
With antilinear operators the definition of equivariance is more subtle.
Recall the representation (\ref{eq:realequiv}): in a generic Hopf algebras
the antipode plays the role of the inverse and this motivates the following
definition. An antilinear operator $T$ defined on $\mc{M}$ is called
\emph{equivariant} if it intertwines the representation $h\mapsto\pi(h)$
with the representation $h\mapsto\pi(S(h)^*)$, that is
$$
T\pi(h)v=\pi(S(h)^*)Tv
$$
for all $h\in\mc{U}$ and $v\in\mc{M}$. Notice that if $T$ is an equivariant
antilinear operator, its square $T^2$ is invariant, i.e.~commutes with all
$h\in\mc{U}$. On the other hand, since
$$
T^*T\pi(h)=T^*\pi(S(h)^*)T=\{\pi(S(h))T\}^*T=\{T\pi(S^2(h))^*\}^*T=\pi(S^2(h))T^*T \;,
$$
the linear operator $T^*T$ is not invariant, unless $S^2=id_{\mc{U}}$ (which
holds if{}f the Hopf algebra is cocommutative, i.e.~has a symmetric coproduct).

A $\mc{U}$-equivariant \emph{real} spectral triple is given by a $\mc{U}$-equivariant
spectral triple $(\A,\HH,D)$ together with a (weak) real structure $J$, with the
additional condition that $J$ is the antiunitary part of an equivariant antilinear
operator.

Equivariant spectral triples have associated (left) covariant differential
calculi, equivariant homologies, and so on. On the other hand, a basic object
like the Dixmier trace is not invariant as one would expect. Let $n$
be the metric dimension and consider the functional $f:\A\to\C$,
$$
f(a):=\tr_\omega(\pi(a)|D|^{-n}) \;.
$$
Using the equivariance 
and the cyclicity of the trace we get
$$
f(h\az a)=\tr_\omega\big(\pi(h_{(1)})\pi(a)\pi(S(h_{(2)}))|D|^{-n}\big)
=\tr_\omega\big(\pi(S(h_{(2)})h_{(1)})\pi(a)|D|^{-n}\big) \;,
$$
and in general $S(h_{(2)})h_{(1)}$ is not equal to $\epsilon(h)$ due to
the wrong order of the factors (unless $S^2=id_{\mc{U}}$). Thus,
we cannot conclude that $f(h\az a)$ and $\epsilon(h)f(a)$ are equal:
in general, the Dixmier trace is not invariant.

\section{Quantum Homogeneous Spaces}\label{chap:QG}

A source of examples of noncommutative spaces is from quantum group theory.
Quantum groups are particular examples of Hopf algebras which appear in the
theory of quantum integrable systems, in the same way as their classical
counterparts, Poisson-Lie groups, arise in the theory of classical integrable
systems. Compact quantum homogeneous spaces are described by unital algebras
$\A$ which are comodule $*$-algebras for some quantum group and such that the
trivial corepresentation appears in the decomposition of $\A$ with multiplicity
$1$ (that is, there are no other coinvariant elements in $\A$ besides the constants).
Reference textbooks for this topic are~\cite{Cha94,Kas95,KS97}.

Recall that on $C^\infty(G)$, $G$ a semisimple Lie group, one can canonically define
a Hopf algebra structure (modulo technicalities, since $\Delta$
does not map into the algebraic tensor product $C^\infty(G)\otimes C^\infty(G)$
but into some completion).
If $\mc{U}(\mf{g})$ is the universal enveloping algebra of the Lie
algebra $\mf{g}$ of $G$, a dual pairing between $\mc{U}(\mf{g})$
and $C^\infty(G)$ is defined on generators $x\in\mf{g}$ by
$$
\inner{x,f}=\tfrac{\de}{\de t}f(\exp tx)\big|_{t=0}\;,\qquad \forall\;f\in C^\infty(G)\;.
$$
The algebra $\mc{U}(\mf{g})$ is a Hopf algebra with operations dual to those
of $C^\infty(G)$, and on an element $x$ of the Lie algebra, the coproduct, counit
and antipode are explicitly given by
$$
\Delta x=x\otimes 1+1\otimes x\;,\qquad\epsilon(x)=0\;,
\qquad S(x)=-x\;.
$$
With these, one can compute explicitly the left and the right canonical actions
(\ref{eq:canaz}) of $x\in\mf{g}$, and finds
$$
(x\az f)(g)=\tfrac{\de}{\de t}f(ge^{tx})\big|_{t=0} \;,\qquad
(f\za x)(g)=\tfrac{\de}{\de t}f(e^{tx}g)\big|_{t=0} \;.
$$
Thus, if $\{x_i\}$ is a basis of $\mf{g}$ orthonormal with respect
to the Killing form, the derivations $x_i\az$ and $\za x_i$ are
a basis for left invariant and right invariant vector fields on $G$.

A Poisson-Lie group $G$ is a Lie group with Poisson structure satisfying
a suitable compatibility condition with the coproduct of $C^\infty(G)$.
If $G$ is simple, the most general Poisson bracket giving a Poisson-Lie
group structure has the form~\cite{Tak94}
\begin{equation}\label{eq:Pois}
\{f_1,f_2\}=\sum r^{ij}\big(f_1\za x_i\cdot f_2\za x_j-x_i\az f_1\cdot x_j\az f_2\big)
\end{equation}
for all $f_{1,2}\in C^\infty(G)$,
where $r^{ij}\in C^\infty(G)$ is a skew-symmetric tensor.
Calling $c=\sum x_i^2$ the quadratic Casimir of $\mc{U}(\mf{g})$ and
$$
t=\Delta c-1\otimes c-c\otimes 1\;,
$$
one has the following important characterization of simple Poisson-Lie groups,
which can be taken as a definition. Equation (\ref{eq:Pois}) gives a
Poisson-Lie group if and only if the element $r=\sum r^{ij}x_i\otimes x_j
\in\bigwedge^2\mf{g}$, named \emph{classical $r$-matrix},
satisfies the modified \emph{classical Yang-Baxter equation}
\begin{equation}\label{eq:mCYBE}
[r_{12},r_{13}]+[r_{12},r_{23}]+[r_{13},r_{23}]=s[t_{13},t_{23}]
\end{equation}
for some constant $s$, and with usual notation
$A_{12}:=\sum A_i\otimes B_i\otimes 1$,
$A_{13}:=\sum A_i\otimes 1\otimes B_i$ and
$A_{23}:=\sum 1\otimes A_i\otimes B_i$
for any $A:=\sum A_i\otimes B_i$.

If $G$ is a matrix group, i.e.~it is realized as a subgroup of $GL(n)$
for some $n$, then the algebra $\A(G)$ of \emph{representative functions}
on $G$ is generated by the functions $u^i_j$ associating to $g\in G$ its
$(i,j)$th matrix entry (in $GL(n)$), and it is dense in $C^\infty(G)$.\footnote{Coproduct
is given on generators by $\Delta u^i_j =\sum_ku^i_k\otimes u^k_j$, and
maps $\A(G)$ into $\A(G)\otimes\A(G)$, which then is a Hopf algebra.
On the other hand $C^\infty(G)$ is not, since the image of $\Delta$ is
$C^\infty(G\times G)$ and in general this is bigger than $C^\infty(G)\otimes C^\infty(G)$.}
In this case the Poisson bracket (\ref{eq:Pois}) is explicitly given on
generators by
\begin{equation}\label{eq:Lax}
\{u^i_j,u^l_m\}=[r,U\otimes U]^{i\phantom{j}l}_{\phantom{i}j\phantom{l}m}
\end{equation}
where $U=((u^i_j))$ and we identify $r\in\mf{g}\otimes\mf{g}$ with its realization
in $gl(n)\otimes gl(n)$. Now (\ref{eq:Lax}) is exactly the equation satisfied when
$U$ is a local Lax $L$-operator (or also the monodromy matrix) of an integrable
model on a lattice~\cite{Bax82}. We can interpret such $L$-operators
as coordinate functions on Poisson-Lie matrix groups. Similarly,
the quantization of these models is associated to \emph{non-commuting}
$L$-operators which can be interpreted as `coordinates' on some
virtual `quantum group'.
Quantum groups are `morally' deformation quantizations of Poisson-Lie
groups, with the additional condition that they are Hopf algebras with
undeformed coproduct and counit (the unique antipode needs to
be deformed since the product changes).
In the next section as an illustrative example we show how to get
$SU_q(2)$ by canonical quantization of the (unique) Poisson-Lie
group structure on $SU(2)$.
The problem with deformation quantization is compatibility
of the product with the undeformed coalgebra structure,
which has to be verified case by case, as well as the
existence of the antipode.

A different approach is to start directly with a noncommutative
algebra with a finite number of generators $u^i_j$. The idea is
to write commutation relations among the generators in the form
$
\sum_{j_1j_2} R^{i_1i_2}_{j_1j_2}u^{j_1}_{k_1}u^{j_2}_{k_2}=
\sum_{j_1j_2} u^{i_1}_{j_1}u^{i_2}_{j_2}R^{j_1j_2}_{k_1k_2}\,
$, or in matrix form
\begin{equation}\label{eq:commrel}
R\,U_1U_2=U_1U_2R
\end{equation}
with $U=((u^i_j))$, $U_1=U\otimes 1$ and $U_2=1\otimes U$, and to look for
a suitable $2$-tensor $R$ with matrix elements depending on a deformation parameter
$\hbar$. The condition $R=1-i\hbar r+o(\hbar)$ guarantees that the semi-classical
limit is (\ref{eq:Lax}), that is ``$\,\lim_{\hbar\to 0}[\,,\,]/i\hbar=\{\,,\,\}\,$''.
A sufficient condition (but not necessary) to get an associative algebra
is that $R$ solves the \emph{quantum Yang-Baxter equation}
\begin{equation}\label{eq:QYBE}
R_{12}R_{13}R_{23}=R_{23}R_{13}R_{12} \;.
\end{equation}
Although (\ref{eq:QYBE}) is not linear and highly non-trivial to treat,
a number of solutions were discovered in the `80s (see for
example \cite{RTF90}).

Thus, for any solution of the quantum Yang-Baxter equation 
we have an associative unital algebra $\A(R)$ generated by $u^i_j$ with
relations (\ref{eq:commrel}). This is always a bialgebra, with undeformed
coproduct and counit
$$
\Delta u^i_j =\sum_ku^i_k\otimes u^k_j\;,\qquad
\epsilon(u^i_j)=\delta^i_j \;.
$$
So, to get a noncommutative deformation of a Poisson-Lie group,
one just needs to find an antipode (and eventually an involution).
The algebra $\A(R)$ is a noncommutative deformation of
$\A(\mr{Mat}_n)$, and to get an antipode one needs to restrict
the attention at least to invertible matrices: that is, one
must quotient $\A(R)$ by a suitable two-sided ideal.
A universal way to construct quotients of $\A(R)$ which are
quantum groups is indicated in \cite[Chapter 7]{Man88}; this
approach is nevertheless quite difficult to handle.
For Drinfeld-Jimbo solutions of the quantum Yang-Baxter equation,
the construction of some noncommutative deformation of classical
Lie groups was performed explicitly in \cite{RTF90}.
It is worth mentioning a different approach due to Woronowicz \cite{Wor80},
which studies compact quantum groups from the point of view of
$C^*$-algebras.

\bigskip

The objects dual to quantum groups, which are deformations of universal
enveloping algebras and are called for this reason \emph{quantum universal
enveloping algebras}, arise for example in conformal quantum field
theory in the study of current algebras (see e.g.~\cite{AGGS89}
and~\cite{Bab87}). A canonical deformation of $\mc{U}(\mf{g})$
was constructed by Drinfeld and Jimbo for any finite-dimensional
complex semisimple Lie algebra $\mf{g}$ (cf.~\cite{Jim85,Jim86a,%
Jim86b,Dri86}).

In the following we recall the definition of those quantum groups
and quantum homogeneous spaces which we'll needed in the original
part of the work, starting with the celebrated $SU_q(2)$ and
Podle\'s quantum spheres.

\subsection{Podle\'s quantum spheres and quantum $SU(2)$}
As mentioned,
there is no general algorithm for finding solutions of the modified
classical Yang-Baxter equation. Despite this, for any simple Lie group
$G$ a canonical $r$-matrix can be constructed as follows. Let $\inner{\,,\,}$
be the Killing form on $\mf{g}$, and $P_\pm:\mf{g}\to\mf{n}_\pm$ the
orthogonal projections onto its maximal nilpotent subalgebras $\mf{n}_\pm$.
Then, an $r$-matrix solving (\ref{eq:mCYBE}) is defined by the formula
(cf.~\cite{STS94})
\begin{equation}\label{eq:STS}
r=i\sum\nolimits_j r'_j\otimes r''_j\;,\qquad
\mr{with}\quad\sum\nolimits_j r'_j\inner{r''_j,\,\circ\,}=P_+-P_-
\end{equation}
(we rescaled the $r$-matrix of \cite{STS94} by a constant, which is always
possible and simply means that the parameter $s$ in (\ref{eq:mCYBE}) has
to be rescaled accordingly). For $G=SL(2,\C)$, with $\inner{A,B}=\tr_{\C^2}(AB)$
and Pauli matrices
$$
\sigma_+=\ma{14pt}{0 & 1 \\ 0 & 0} \;,\qquad
\sigma_-=\ma{14pt}{0 & 0 \\ 1 & 0} \;,
$$
one has $P_\pm=\sigma_\pm\inner{\sigma_\mp,\,\circ\,}$ and then by (\ref{eq:STS})
$$
r=i\sigma_+\wedge\sigma_-\;.
$$
Now let
$$
U=\ma{14pt}{\alpha & \beta \\ -\beta^* & \alpha^*}\;,\qquad
\mr{with}\;\alpha\alpha^*+\beta\beta^*=1\;,
$$
where $\alpha$ and $\beta$ are the generators of the coordinate $*$-algebra of $SU(2)$;
this is obtained from $\A(SL(2,\C))$ by imposing a $*$-structure, as discussed e.g.~in
Chapter 4 of \cite{KS97}.
Thanks to the choice of the normalization, the Poisson structure on $SL(2,\C)$
associated to $r$ gives a Poisson structure on $\A(SU(2))$ which is compatible
with the real structure, i.e.
$$
\{a,b\}^*=\{a^*,b^*\}\;,\qquad\forall\;a,b\in\A(SU(2))\;.
$$
This is explicitly given on generators by
\begin{equation}\label{eq:Psu2}
\{\alpha,\beta\}=i\alpha\beta\;,\qquad
\{\alpha,\beta^*\}=i\alpha\beta^*\;,\qquad
\{\beta,\beta^*\}=0\;,\qquad
\{\alpha,\alpha^*\}=-2i\beta\beta^*\;.
\end{equation}
By \cite[Section~1.5D]{Cha94} this is the unique Poisson-Lie group structure
on $SU(2)$ up to a multiplicative constant. Symplectic leaves are open disks $D_t$,
defined for $t\in U(1)$ by $\beta=t|\beta|$ with $|\beta|\neq 0$,
and the single points at the boundary of the disks~\cite[Section~1.5C]{Cha94}.

The $2$-sphere $S^2$ can be identified with the sub-manifold of matrices
in $SU(2)$ with real off-diagonal elements, i.e.~$\A(S^2)$ is the quotient
of $\A(SU(2))$ by the ideal generated by the element $\beta-\beta^*$ (that is,
the class of $\beta$ is real in the quotient algebra). The class $A$ of $\beta$ and
the class $B$ of $\alpha$ generate $\A(S^2)$. The Poisson structure (\ref{eq:Psu2})
induces a Poisson structure on the quotient, and $S^2$ splits into two
non-trivial symplectic leaves: these are the two open hemispheres $D_\pm$, defined
by $\pm A>0$, and glued together along the equator $S^1$, defined by $A=0$.
While on the equator the Poisson structure is trivial, on each hemisphere
$D_\pm$ the generators of $\A(S^2)$ can be rewritten in terms of canonical
coordinates $x,p$, satisfying
$$
\{x,p\}=1\;,
$$
via the transform ($x\geq 0$)
\begin{equation}\label{eq:leaves}
A=\pm e^{-x} \;,\qquad
B=\sqrt{1-e^{-2x}}\,e^{ip} \;,\qquad
B^*=e^{-ip}\sqrt{1-e^{-2x}} \;.
\end{equation}
It is interesting to notice that the Lax operators for this Poisson structure
are the ones associated to the Liouville model on a lattice \cite{Tak94}.
By applying the recipe of canonical quantization, that means replacing $x$ and $p$
with operators $\hat{x}$ and $\hat{p}$ satisfying the commutation relation (\ref{eq:QM}),
$$
[\hat{x},\hat{p}]=i\hbar\;,
$$
we get a noncommutative algebra whose underlying virtual space is known as the
\emph{equatorial Podle\'s sphere}, denoted $S^2_q$ and discovered in \cite{Pod87}.
Thus, coordinates on the equatorial Podle\'s sphere (actually, on its two open
hemispheres) can be interpreted as matrix entries of Lax operators for the
quantization of the Liouville model on a lattice.
In fact, that $S^2_q$ is even a \emph{strict} deformation quantization~\cite{Rie89}
of $S^2$ was proved in~\cite{She91}.\footnote{It is worth mentioning the appendix of \cite{She91},
where it is reviewed the work in \cite{LW90} about the classification of $SU(2)$-covariant
Poisson structures on $S^2$.}

Commutation relations among generators of $\A(S^2_q)$ can be computed by replacing
$x,p$ with $\hat{x},\hat{p}$ in (\ref{eq:leaves}) (keeping the order in which the products
are written), and using the Weyl-type relation $e^{-\hat{x}}e^{i\hat{p}}=e^\hbar e^{i\hat{p}}e^{-\hat{x}}$.
They are written in the definition that follows. We set $q:=e^\hbar$.

\begin{df}\label{def:S2q}
We call \emph{equatorial Podle\'s sphere} the virtual space underlying
the $*$-algebra $\A(S^2_q)$ generated by $A=A^*$, $B$ and $B^*$ with
relations:
$$
AB=qBA \;,\qquad
BB^*+A^2=1 \;,\qquad
B^*\! B+q^2A^2=1 \;.
$$
\end{df}

As usual due to symmetries of the algebra one can focus on the case
\ \fbox{$0<q<1$}\ . In all this dissertation we'll always assume that $q$
is in this range (unless stated otherwise).

The equatorial Podle\'s sphere turns out to be the actual building block
for $q$-deformations. We'll see that working with local formul{\ae} essentially
means shifting the problem from a complicated quantum space to a number of
copies of the equatorial Podle\'s sphere.

Similar relations yield a deformation of $\A(SU(2))$ which is a Hopf
$*$-algebra with undeformed coproduct $\Delta$ and counit $\epsilon$.
This is the well known quantum group $SU_q(2)$ of Woronowicz \cite{Wor80,Wor87}.
The sphere $S^2_q$ is `embedded' in $SU_q(2)$ as an equator,
though this is not the reason for the name `equatorial' (it
is called `equatorial' since one can think of the space of
characters of the algebra, which is a circle, as the equator
of $S^2_q$).

\begin{df}\label{def:SUq2}
The algebra $\A(SU_q(2))$ of polynomial functions on $SU_q(2)$ is
the $*$-algebra generated by $\alpha,\beta$ and their
adjoints, with relations:
$$
\beta\alpha=q\alpha\beta \;,\quad
\beta^*\alpha=q\alpha\beta^* \;,\quad
[\beta,\beta^*]=0 \;,\quad
\alpha\alpha^*+\beta\beta^*=1 \;,\quad
\alpha^*\alpha+q^2\beta^*\beta=1 \;.
$$
It is a Hopf $*$-algebra with coproduct, counit and antipode given by
\begin{gather*}
\Delta\ma{14pt}{\alpha & \beta \\ -q\beta^* & \alpha^*}
=\ma{14pt}{\alpha & \beta \\ -q\beta^* & \alpha^*}\,\dot{\otimes}\,
\ma{14pt}{\alpha & \beta \\ -q\beta^* & \alpha^*}\;, \\
\rule{0pt}{25pt}
\epsilon\ma{14pt}{\alpha & \beta \\ -q\beta^* & \alpha^*}=
\ma{14pt}{1 & 0 \\ 0 & 1}\;,\qquad
S\ma{14pt}{\alpha & \beta \\ -q\beta^* & \alpha^*}=
\ma{14pt}{\alpha^* & -q\beta \\ \beta^* & \alpha}\;.
\end{gather*}
\end{df}

\noindent
Here the dotted tensor product is defined as
$(A\,\,\dot{\otimes}\,B)_{ij}=\sum_kA_{ik}\otimes B_{kj}$.
We use the same notation of~\cite{DLS05}, but for Greek letters instead
of Latin ones.
In~\cite[sec.~4.1.4]{KS97} they call $d=\alpha$, $c=-\beta$ and their
fundamental corepresentation is related to ours by a unitary equivalence:
\begin{equation*}
\ma{14pt}{d^* & -qc^* \\ c & d}=\ma{14pt}{0 & 1 \\ -1 & 0}
\ma{14pt}{\alpha & \beta \\ -q\beta^* & \alpha^*}
\ma{14pt}{0 & -1 \\ 1 & 0}\;.
\end{equation*}

Another quantum $2$-sphere can be obtained by recalling that $S^2$
is diffeomorphic to the projective line $\CP^1=SU(2)/U(1)$.
For any $q$, the algebra $\A(SU_q(2))$ carries a covariant
action of $u\in U(1)$ given by the map $\alpha\mapsto u\alpha$,
$\beta\mapsto u\beta$. The algebra $\A(\CP^1_q)$ is defined as
the fixed point algebra of $\A(SU_q(2))$ for this action of $U(1)$,
and coincides with the subalgebra generated by $A=\beta\beta^*$,
$B=\beta^*\alpha$ (we use again the same symbols $A$, $B$).
This is called \emph{standard Podle\'s sphere}, and commutation
relations among its generators can be easily computed from the ones
defining $\A(SU_q(2))$.

\begin{df}
We call coordinate algebra on the \emph{standard Podle\'s sphere}
$\CP^1_q$ the $*$-algebra generated by $A=A^*$, $B$ and $B^*$ with
relations
$$
AB=q^2BA\;, \qquad
BB^*+A(A-1)=0\;,\qquad
B^*\! B+q^2A(q^2A-1)=0\;.
$$
\end{df}

The standard and equatorial sphere are limiting cases of a family
of noncommutative spheres $S^2_{qs}$, depending on a parameter $s\in[0,1]$
($q$ is supposed to be fixed), whose elements are all quantum homogeneous
spaces for $SU_q(2)$. The spheres $\CP^1_q$ and $S^2_{q^2}$ correspond to
the `extremal' values $s=0$ and $s=1$, respectively.

\begin{df}\label{def:S2qs}
For $0<q<1$ and $0\leq s\leq 1$, the coordinate algebra on the \emph{Podle{\'s} sphere}
$S^2_{qs}$ is the $*$-algebra generated by $A=A^*$, $B$ and $B^*$ with relations
$$
AB=q^2BA\;, \qquad
BB^*+(A+s^2)(A-1)=0\;,\qquad
B^*\! B+(q^2A+s^2)(q^2A-1)=0\;.
$$
For $s=0$ we obtain the \emph{standard} Podle{\'s} sphere $S^2_{q0}=\CP^1_q$, while
for $s=1$ we obtain the \emph{equatorial} Podle\'s sphere $S^2_{q1}=S^2_{q^2}$.
\end{df}

\noindent
We are using, with minor changes, the notations of~\cite{Dab06}.
The original Podle\'s parameter is $c=(s^{-1}-s)^{-2}$, while the
original generators are $A'=(1-s^2)^{-1}A,\;B'=(1-s^2)^{-1}B^*$ if $s\neq 1$,
and $A'=A,\;B'=B^*$ if $s=1$ \cite{Pod87}.
The generators used in \cite{DDLW07} are related to ours by
$x_1:=\sqrt{1+q^2}\,B$, $x_0:=-(1+q^2)A+(1-s^2)$,
$x_{-1}:=-q^{-1}x_1^*$.

The geometrical meaning of the parameter $s$ is clear when one
takes the $q\to 1$ limit. For $q=1$, one gets the coordinate algebra
of a $2$-sphere in $\R^3$ with radius $\frac{1+s^2}{2}$ and center
$(0,0,\frac{1-s^2}{2})$, both depending on $s$.
Notice that, for any $q$ and $s$, we have an injective $*$-algebra morphism
$\A(S^2_{qs})\to\A(SU_q(2))$ defined on generators by
\begin{equation}\label{eq:SUtoS2}
A\mapsto s(\alpha\beta+\beta^*\alpha^*)+(1-s^2)\beta\beta^*\;,\qquad
B\mapsto s\bigl(\alpha^2-q(\beta^*)^2\bigr)+(1-s^2)\beta^*\alpha\;,
\end{equation}
and which is the noncommutative analogue of the principal bundle
$SU(2)\to S^2$.

The noncommutative geometry of Podle\'s spheres will be studied in Chapter
\ref{chap:S2qs}, following the papers \cite{DD06} and \cite{DDLW07}.

\subsection{The `zoo' of quantum spaces}

\afo{8cm}{We think in generalities, but we live in details}
         {Alfred North Whitehead (1861 -- 1947).}

\noindent
Generalizations of Podle\'s spheres and $SU_q(2)$ are possible in different
directions (most of the spaces we discuss here were introduced in a more or
less explicit form in \cite{RTF90}).
The equatorial Podle\'s sphere is the first of a family $S^n_q$ of quantum
orthogonal spheres, which are homogeneous spaces for the quantum orthogonal
groups but not quotient spaces (they are not fixed point subalgebras of
$\A(SO_q(n+1))$ for a coaction of $SO_q(n)$, unless $q=1$). The case $n=4$
will be studied in Chapter \ref{chap:S4q}; here we give its algebra.

\begin{df}\label{def:S4q}
We call \emph{quantum orthogonal $4$-sphere} the virtual space underlying the algebra
$\A(S^4_q)$ generated by $x_0=x_0^*$, $x_1$, $x_1^*$, $x_2$, $x_2^*$, with commutation
relations:
\begin{align*}
&x_ix_j     = q^2x_jx_i   \;, \qquad\quad\forall\;\;0\leq i<j\leq 2 \;, \\
&x_i^*x_j   =q^2x_jx_i^*  \;, \qquad\quad\forall\;\;i\neq j \;, \\
&[x_1^*,x_1]=(1-q^4)x_0^2 \;, \\
&[x_2^*,x_2]=x_1^*x_1-q^4x_1x_1^* \;, \\
&x_0^2+x_1x_1^*+x_2x_2^*=1 \;.
\end{align*}
\end{df}

\noindent
The original notations of \cite[eq.~(1.14)]{RTF90}
can be obtained by defining $x'_1:=x_2^*$, $x'_2:=x_1^*$, $x'_3:=\sqrt{q(1+q^2)}\,x_0$,
$x'_4:=x_1$, $x'_5:=x_2$ and $q':=q^2$.
The notations in~\cite[eq.~(2.1)]{HL04} can be obtained by the replacement
$x_i\mapsto x_i^*$ and $q^2\mapsto q^{-1}$.

A family of odd-dimensional quantum spheres were introduced also
in~\cite{VS91} as quantum homogeneous spaces for quantum unitary groups,
and it was later recognized in~\cite{HL04,Dab03} that these `unitary'
quantum spheres and the odd-dimensional orthogonal spheres
were the same (up to an isomorphism).

\begin{df}
For $\ell\geq 1$, we denote $\A(SU_q(\ell+1))$ the $*$-algebra generated
by the $(\ell+1)^2$ elements $u^i_j$ ($i,j=1,\ldots,\ell+1$)
with commutation relations
\begin{align*}
u^i_ku^j_k &=qu^j_ku^i_k &
u^k_iu^k_j &=qu^k_ju^k_i &&
\forall\;i<j\;, \\
[u^i_l,u^j_k]&=0 &
[u^i_k,u^j_l]&=(q-q^{-1})u^i_lu^j_k &&
\forall\;i<j,\;k<l\;,
\end{align*}
and determinant relation
$$
\sum\nolimits_{p\in S_{\ell+1}}(-q)^{||p||}
u^1_{p(1)}u^2_{p(2)}\ldots u^{\ell+1}_{p(\ell+1)}=1
$$
where the sum is over all permutations $p$ of the $\ell+1$ elements
$\{1,2,\ldots,\ell+1\}$ and $||p||$ is the (minimal) number of inversions
in $p$. The $*$-structure is given by
$$
(u^i_j)^*=(-q)^{j-i}\sum\nolimits_{p\in S_\ell}(-q)^{||p||}u^{k_1}_{p(m_1)}
u^{k_2}_{p(m_2)}\ldots u^{k_\ell}_{p(n_\ell)}
$$
with $\{k_1,\ldots,k_\ell\}=\{1,\ldots,\ell+1\}\smallsetminus\{i\}$,
$\{m_1,\ldots,m_\ell\}=\{1,\ldots,\ell+1\}\smallsetminus\{j\}$
(as ordered sets) and the sum is over all permutations $p$ of the $\ell$
elements $\{m_1,\ldots,m_\ell\}$.

This is a Hopf $*$-algebra with standard coproduct, counit and antipode:
$$
\Delta(u^i_j)=\sum\nolimits_ku^i_k\otimes u^k_j\;,\qquad
\epsilon(u^i_j)=\delta^i_j\;,\qquad
S(u^i_j)=(u^j_i)^*\;.
$$
\end{df}

\noindent
We are using the same notations of~\cite[Section 9.2]{KS97}.

The coordinate algebra of the odd-dimensional quantum spheres
$S^{2\ell+1}_q$ can be identified both with the subalgebra of
$\A(SU_q(\ell+1))$ generated by the elements $u^{\ell+1}_i$,
or with the one generated by the elements $u^1_i$. In both
cases one gets isomorphic $*$-algebras, and both fixed point
subalgebras of $\A(SU_q(\ell+1))$ for coactions of $SU_q(\ell)$
obtained with two different choices of the inclusion
$\A(SU_q(\ell))\hookrightarrow\A(SU_q(\ell+1))$.

\begin{df}
We call $\A(S^{2\ell+1}_q)$ the $*$-algebra generated by
$\{z_i,z_i^*\}_{i=1,\ldots,\ell+1}$ with relations
\begin{align*}
z_iz_j &=qz_jz_i &&\forall\;i<j \;,\\
z_i^*z_j &=qz_jz_i^* &&\forall\;i\neq j \;,\\
[z_1^*,z_1] &=0 \;,\\
[z_{i+1}^*,z_{i+1}] &=(1-q^2)\sum\nolimits_{j=1}^i z_jz_j^*
    &&\forall\;i=1,\ldots,\ell \;,\\
z_1z_1^*+z_2z_2^* &+\ldots+z_{\ell+1}z_{\ell+1}^*=1 \;.
\end{align*}
\end{df}

\noindent
The original generators of \cite{VS91} are obtained by the
replacement $z_i\to z_{\ell+1-i}$ and setting $q=e^{h/2}$.
The generators of \cite{CP07} are the same of \cite{VS91} but numbered
from $0$ to $\ell$ instead of from $1$ to $\ell+1$;
there should be a $q$ instead of $q^{-1}$ in their second relation.
The generators $x_i$ used in~\cite{HL04} are related to ours by
$x_i=z_i^*$ and by the replacement $q\to q^{-1}$.

One can check that both the maps $\,z_i\mapsto u^{\ell+1}_i\,$ and
$\,z_{\ell+2-i}^*\mapsto q^{-i+1}u^1_i\,$
define $*$-algebra morphisms. That these are also injective is
proved in \cite{VS91}. We'll spend a few words about odd-dimensional
quantum spheres in Chapter \ref{chap:CPlq}.


\typeout{Capitolo 3}

\chapter[Noncommutative geometry of Podle\'s spheres]
{Noncommutative geometry of \\ Podle\'s quantum spheres}\label{chap:S2qs}

\afo{8cm}{Mathematics is the part of physics where experiments are cheap.}
         {Vladimir I. Arnold~\cite{Arn98}.}

\noindent
It is not yet completely clear how to generalize the axioms for commutative
real spectral triples to get a theory of `noncommutative manifolds'.
Too narrow axioms would exclude a lot of interesting examples,
while too mild axioms would give a too general theory, where it is
difficult to prove interesting theorems.
To shed light on this problem, a fundamental step is to study
the properties of a certain number of examples, which is the
`experimental' part in doing mathematics.
The class of examples on which we focus is that of quantum
homogeneous spaces.

In this chapter we study regular even spectral triples on $S^2_{qs}$
following \cite{DD06,DDLW07} (recall that $0<q<1$ and $0\leq s\leq 1$).
We then focus on a particular spectral triple, which is a deformation of
the canonical one associated to the round sphere $S^2$, and whose Dirac
operator has undeformed spectrum. For this, we compute dimension spectrum
and the noncommutative integral, which will be simply the integral over
the space of characters of the algebra.
The problem of constructing a real structure is also discussed.

We recall Definition \ref{def:S2q} of the coordinate
algebra $\A(S^2_{qs})$ of the quantum Podle\'s spheres, and Definition
\ref{def:SUq2} of the $*$-Hopf algebra $\A(SU_q(2))$. The former is
generated by $A=A^*$, $B$ and $B^*$ with relations
$$
AB=q^2BA\;, \qquad
BB^*+(A+s^2)(A-1)=0\;,\qquad
B^*\! B+(q^2A+s^2)(q^2A-1)=0\;.
$$
The latter is generated as a $*$-algebra by $\alpha,\beta,\alpha^*$ and
$\beta^*$, with relations
$$
\beta\alpha=q\alpha\beta \;,\quad
\beta^*\alpha=q\alpha\beta^* \;,\quad
[\beta,\beta^*]=0 \;,\quad
\alpha\alpha^*+\beta\beta^*=1 \;,\quad
\alpha^*\alpha+q^2\beta^*\beta=1 \;.
$$

\section{The symmetry Hopf algebra $U_q(su(2))$}
We call $U_q(su(2))$ the compact real form of the Drinfeld-Jimbo deformation of
$sl(2)$; it is a real form of the Hopf algebra called $\breve{U}_q(su_2)$
in~\cite[Section~3.1.2]{KS97}. As a $*$-algebra, $U_q(su(2))$ is generated
by $K=K^*,K^{-1},E$ and $F=E^*$, with relations
$$
KEK^{-1}=qE\;,\qquad
[E,F]=\tfrac{K^2-K^{-2}}{q-q^{-1}}\;,
$$
all the other relations being determined by compatibility with the involution.
The Hopf algebra structure $(\Delta,\epsilon,S)$ of $U_q(su(2))$ is given by
$$
\begin{array}{c}
\Delta K=K\otimes K\;,\quad \Delta E=E\otimes K+K^{-1}\otimes E\;, 
\\ \rule{0pt}{3.5ex}\epsilon(K)=1\;,\quad \epsilon(E)=0\;,\\ \rule{0pt}{3.5ex}
S(K)=K^{-1}\;,\quad S(E)=-qE\;.
\end{array}
$$
Using Lemma 3.2.2 of~\cite{Maj95} we write the coproduct of $E^n$ as
\begin{equation}\label{eq:coJn}
\Delta(E^n)=\sum_{k=0}^n\sqbn{n}{k}\,E^{n-k}K^{-k}\otimes E^kK^{n-k}\;,
\end{equation}
where $q$-factorials and $q$-binomials are defined as:
\begin{equation*}
[0]!:=1\;,\qquad
[n]!:=[n][n-1]\ldots[1]\;\forall\;n\geq 1\;,\qquad
\sqbn{n}{k}:=\frac{[n]!}{[k]!\,[n-k]!}\;.
\end{equation*}

For $l\in\frac{1}{2}\N$, $\omega=\pm 1$, let $V_l^\omega$ be a complex vector space with
orthonormal basis $\ket{l,m}$, $m=-l,-l+1,\ldots,l$. Then $V_l^\omega$ carries an
irreducible $*$-representation of $U_q(su(2))$ defined by
\begin{align}
K\ket{l,m}&=\omega\,q^m\ket{l,m}\;, \notag\\
E\ket{l,m}&=\omega\sqrt{[l-m][l+m+1]}\ket{l,m+1}\;, \label{eq:rep}\\
F\ket{l,m}&=\omega\sqrt{[l+m][l-m+1]}\ket{l,m-1}\;, \notag
\end{align}
where $[x]=(q^x-q^{-x})/(q-q^{-1})$ is the $q$-analogue of $x\in\C$.
Any finite-dimensional irreducible $*$-representation of $U_q(su(2))$ is equivalent
to one of the $V_l^\omega$~\cite[Section~3.2.3]{KS97}. In the following, we'll be interested
only in $V_l:=V_l{}^+$.

For $l\in\N$, let $C$ be the \emph{antilinear} operator on $V_l$ defined by
\begin{equation}\label{eq:stru}
C\ket{l,m}=(-q)^m\ket{l,-m}\;,
\end{equation}
and which is of square equal to $1$. 
Moreover
\begin{equation}\label{eq:simi}
ChC=S(h)^*\;,\qquad\forall\;h\in U_q(su(2))\;,
\end{equation}
as one checks first taking $h=K,E,F$:
\begin{gather*}
CKC\ket{l,m}=q^{-m}\ket{l,m}=K^{-1}\ket{l,m}=S(K)^*\ket{l,m} \;, \\
CEC\ket{l,m}=-q\sqrt{[l+m][l-m+1]}\ket{l,m-1}=-qF\ket{l,m}=S(E)^*\ket{l,m} \;,\\
CFC\ket{l,m}=-q^{-1}\sqrt{[l-m][l+m+1]}\ket{l,m+1}=-q^{-1}E\ket{l,m}=S(F)^*\ket{l,m} \;,
\end{gather*}
and then noting that both members of the equation are algebra morphisms.
%
Similarly, considering first $h=K,E,F$ and then using the fact
that both member of the equation are algebra morphism,
one proves that, for all $h\in U_q(su(2))$,
\begin{equation}\label{eq:Squad}
S^2(h)=K^2hK^{-2} \;.
\end{equation}
Later we need the Clebsch-Gordan coefficients for the decomposition
$V_l\otimes V_{1/2}\simeq V_{l+1/2}\oplus V_{l-1/2}$ and
$V_{1/2}\otimes V_l\simeq V_{l+1/2}\oplus V_{l-1/2}$.
If we call $\ket{\pm}:=\ket{\oh,\pm\oh}$, then~\cite{KS97}:
\begin{subequations}\label{eq:Clebsch}
\begin{align}
\ket{l+\tfrac{1}{2},m-\tfrac{1}{2}} &=
\sqrt{\rule{0pt}{11pt}\smash[t]{\tfrac{q^{-(l+m)}[l-m+1]}{[2l+1]}}} \ket{l,m}\otimes\ket{-}
+\sqrt{\rule{0pt}{11pt}\smash[t]{\tfrac{q^{l-m+1}[l+m]}{[2l+1]}}} \ket{l,m-1}\otimes\ket{+} \;,\\
\ket{l-\tfrac{1}{2},m+\tfrac{1}{2}} &=
\sqrt{\rule{0pt}{11pt}\smash[t]{\tfrac{q^{l-m}[l+m+1]}{[2l+1]}}} \ket{l,m+1}\otimes\ket{-}
-\sqrt{\rule{0pt}{11pt}\smash[t]{\tfrac{q^{-(l+m+1)}[l-m]}{[2l+1]}}} \ket{l,m}\otimes\ket{+} \;,
\end{align}
\end{subequations}
and
\begin{subequations}\label{eq:copClebsch}
\begin{align}
\ket{l+\tfrac{1}{2},m-\tfrac{1}{2}} &=
\sqrt{\rule{0pt}{11pt}\smash[t]{\tfrac{q^{l+m}[l-m+1]}{[2l+1]}}} \ket{-}\otimes\ket{l,m}
+\sqrt{\rule{0pt}{11pt}\smash[t]{\tfrac{q^{-(l-m+1)}[l+m]}{[2l+1]}}} \ket{+}\otimes\ket{l,m-1} \;,\\
\ket{l-\tfrac{1}{2},m+\tfrac{1}{2}} &=
\sqrt{\rule{0pt}{11pt}\smash[t]{\tfrac{q^{l+m+1}[l-m]}{[2l+1]}}} \ket{+}\otimes\ket{l,m}
-\sqrt{\rule{0pt}{11pt}\smash[t]{\tfrac{q^{-(l-m)}[l+m+1]}{[2l+1]}}} \ket{-}\otimes\ket{l,m+1} \;.
\end{align}
\end{subequations}

\section{Spherical functions}
A non-degenerate dual pairing $\inner{\,,\,}:U_q(su(2))\times\A(SU_q(2))\to\C$
is defined by (cf.~Section~4.4.1 of \cite{KS97}, the generators there being
$d=\alpha$ and $c=-\beta$)
$$
\inner{K,\alpha}=q^{1/2}\;,\qquad\inner{K,\alpha^*}=q^{-1/2}\;,\qquad
\inner{E,-\beta}=\inner{F,q\beta^*}=1\;,
$$
and $\inner{h,a}=0$ for all other pairs of generators $h,a$.
A complete proof that this pairing is non-degenerate
is for example in Section VII.4 of \cite{Kas95}.
Identifying the $*$-algebra $\A(S^2_{qs})$ with its image in $\A(SU_q(2))$ under
the $*$-algebra morphism (\ref{eq:SUtoS2}) one checks that $\A(S^2_{qs})$
are invariant subspaces for the left canonical action (\ref{eq:canaz}),
which on generators is explicitly given by
\begin{align*}
K\az A&=A \;, &
K\az B\phantom{^*}&=qB \;, &
K\az B^*&=q^{-1}B^* \;,\\
E\az A&=-q^{-1/2}B \;, &
E\az B^*&=q^{1/2}(1+q^2)A-q^{1/2}(1-s^2) \;, &
E\az B\phantom{^*}&=0 \;,\\
F\az A&=q^{-3/2}B^* \;, &
F\az B\phantom{^*}&=-q^{-1/2}(1+q^2)A+q^{-1/2}(1-s^2) \;, &
F\az B^*&=0 \;.
\end{align*}
The algebra $\A(S^2_{qs})$ is then a $U_q(su(2))$-module $*$-algebra for this
action.

We can decompose $\A(S^2_{qs})$ into irreducible representations of $U_q(su(2))$,
i.e.~find a basis $\ket{l,m;\lambda}^\omega$ carrying --- for $l$, $\lambda$
and $\omega$ fixed --- a representation of $U_q(su(2))$ equivalent to $V_l^\omega$.
The index $\lambda$ counts possible multiplicities.
From (\ref{eq:rep}) we read $K\ket{l,m;\lambda}^\omega=\omega q^m\ket{l,m;\lambda}^\omega$.
Every element of $\A(S^2_{qs})$ can be written as a finite sum:
\begin{equation*}
f=\sum_{n\in\N}\big\{B^nf_n(A)+(B^*)^{n+1}f_{-n-1}(A)\big\}\;,
\end{equation*}
where $f_n$ are polynomials in $A$. Then $K\az f=-q^mf$ has no solution,
and $K\az f=q^mf$ if and only if $f=B^mf_m(A)$ if $m\in\N$ or
$f=(B^*)^{-m}f_m(A)$ if $-m-1\in\N$. Hence, only the representations $V_l^\omega$
with $\omega=1$ and $l$ integer appear in the decomposition of $\A(S^2_{qs})$;
the index $\omega$ can be omitted.

The \emph{lowest weight vectors} $\ket{l,-l;\lambda}=(B^*)^lf_{-l;\lambda}(A)$
are determined by the condition that
\begin{equation*}
F\ket{l,-l;\lambda}=q^l(B^*)^l\,F\az f_{-l;\lambda}(A)=0\;,
\end{equation*}
implying that $f_{-l;\lambda}(A)=c_l$ is a normalization constant. Thus,
for each $l\in\N$ the lowest weight vector is unique, and the representations
$V_l$ have multiplicity $1$: $\A(S^2_{qs})\simeq\bigoplus_{l\in\N}V_l$. We can
omit the index $\lambda$. The equivalence $\A(S^2_{qs})\to\bigoplus_{l\in\N}V_l$
of $U_q(su(2))$-modules is unitary if the involution $*$ is sent to the operator
$C$ on $V_l$, and this in particular means
\begin{equation*}
\ket{l,m}=(-q)^m\,\ket{l,-m}^*\;.
\end{equation*}
We need just to determine $\ket{l,m}$ for $m\geq 0$.
For all $m=-l,-l+1,\ldots,l$ we have
\begin{equation*}
\ket{l,m}=\bigl([l+m][l-m+1]\bigr)^{-1/2}\,E\ket{l,m-1}\;,
\end{equation*}
which iterated $l+m$ times becomes
\begin{equation}\label{eq:ketlm}
\ket{l,m}=\sqrt{\frac{[l-m]!}{[l+m]!\,[2l]!}}\,E^{l+m}\ket{l,-l}
         =c_l\sqrt{\frac{[l-m]!}{[l+m]!\,[2l]!}}\,E^{l+m}\az (B^*)^l \;.
\end{equation}
Last equation reduces to Rodriguez formula for Laplace spherical harmonics,
in the limit $q\to 1$.

A way to determine the coefficients $c_l$ is by imposing the normalization of $\ket{l,m}$
to $1$ with respect to the inner product associated to the Haar measure on $SU_q(2)$
and restricted to $S^2_{qs}$. We'll use a simpler (even if equivalent) method: we'll
compute $c_l$ at the end of this section by imposing that $A$ is represented by a
symmetric operator. In both cases, changing $c_l$ of a phase factor does not affect
the normalization neither the self-adjointness of operators (it correspond to a unitary
transformation), and we can choose $c_l\in\R$.

The next step is to compute the action of $E^{l+m}$.
The explicit expression for spherical functions on the generic Podle\'s
sphere already appears in~\cite{NM90}, but for completeness we prefer to
compute them from scratch.

For $l\geq 1$, we call $P_l^{\pm}$ the following polynomials:
\begin{equation*}
P_l^+=\prod\nolimits_{k=1}^l(q^{2k}A+s^2)\;,\qquad
P_l^-=\prod\nolimits_{k=1}^l(q^{2k}A-1)\;,
\end{equation*}
with the notation $\prod_{k=1}^0:=1$ and $P^{\pm}_0:=1$.
Then, a little induction yelds
\begin{equation}\label{eq:6.6}
(B^*)^lB^l=(-1)^lP_l^-P_l^+\qquad\forall\;l\geq 0\;.
\end{equation}
The inductive step comes by observing that
\begin{equation*}
(B^*)^lB^l
=-(B^*)^{l-1}(q^2A+s^2)(q^2A-1)B^{l-1}
=-(B^*)^{l-1}B^{l-1}(q^{2l}A+s^2)(q^{2l}A-1)\;.
\end{equation*}
Multiplying (\ref{eq:ketlm}) from the right by $B^l$ and using $E\az B=0$
and (\ref{eq:6.6}), we get
\begin{equation}\label{eq:ins}
q^{l(l+m)}\ket{l,m}B^l=(-1)^lc_l\sqrt{\frac{[l-m]!}{[l+m]!\,[2l]!}}\,E^{l+m}\az (P_l^-P_l^+)\;.
\end{equation}
Induction on $l$ gives
\begin{equation*}
E\az P^\pm_l =-q^{l+1/2}[l]\,P^{\pm}_{l-1}B\qquad\forall\;l\geq 0\,,
\end{equation*}
while induction on $n$ yelds
\begin{subequations}
\begin{align}
E^n\az P_l^{\pm} &=(-q^{l+1/2})^n\frac{[l]!}{[l-n]!}\,P_{l-n}^{\pm}B^n
&&\mr{if}\;0\leq n\leq l\;, \label{eq:indJP} \\
E^n\az P_l^{\pm} &=0 &&\mr{if}\;n>l\geq 0\;. \label{eq:zero}
\end{align}
\end{subequations}
Using covariance of the action, together with (\ref{eq:coJn}),
we obtain:
\begin{align*}
E^{l+m}\az (P_l^-P_l^+) &=\sum_{k=0}^{l+m}\sqbn{l+m}{k}\,(E^{l+m-k}\az P_l^-)(E^k\az P_l^+) \\
 &=\sum_{\substack{0\leq p,r\leq l,\\ p+r=l-m}}\frac{[l+m]!}{[l-p]![l-r]!}
\,(E^{l-p}\az P_l^-)(E^{l-r}\az P_l^+)\;.
\end{align*}
In the second step we used the indexes $p=k-m$ and $r=l-k$, and noticed
that the terms with $p,r<0$ are zero, due to (\ref{eq:zero}).
In turn, using (\ref{eq:indJP}):
\begin{equation*}
E^{l+m}\az (P_l^-P_l^+)=(-q^{l+1/2})^{l+m}[l+m]!\sum_{p+r=l-m}
\sqbn{l}{p}\sqbn{l}{r}\,P_p^- B^{m+r}P_r^+ B^{l-r}\;.
\end{equation*}
Now, $AB^n=q^{2n}B^nA$ implies (when $m\geq 0$, so that $B^m$ is well-defined) that
\begin{equation}\label{eq:previous}
P_p^-B^m =B^m\left(\prod_{h=1}^p(q^{2(m+h)}A-1)\right) \;,\qquad
B^rP_r^+ =\left(\prod_{h=1}^r(q^{2(h-r)}A+s^2)\right)B^r\;.
\end{equation}
With this, equations (\ref{eq:previous}) and (\ref{eq:ins}) we get the
final result, summarized in a proposition.

\begin{prop}\label{prop:Laplace}
The spherical functions $\ket{l,m}$, defined by
\begin{align*}
\ket{l,m} &:=N_{l,m}\,B^mp_{l,m}   &\mr{if}\;m\geq 0\;, \\
\ket{l,m} &:=(-q)^m\ket{l,-m}^*    &\mr{if}\;m<0\;, \\
p_{l,m}   &:=\sum_{p+r=l-m}\sqbn{l}{p}\sqbn{l}{r}
  \prod_{k=1}^p(q^{2(m+k)}A-1)\prod_{h=1}^r(q^{2(h-r)}A+s^2)\;, \\
N_{l,m}   &:=c_l(-1)^{l+m}q^{-(l-m)/2}\sqrt{\frac{[l+m]!\,[l-m]!}{[2l]!}}\;,
\end{align*}
are an orthonormal basis of $\A(S^2_{qs})$ (with respect to the inner
product associated to the Haar functional). The left action of $U_q(su(2))$
on this basis is described by equations (\ref{eq:rep}).
\end{prop}

The normalization constants $c_l\in\R$ have been rescaled by $(-q)^l$
in order to have $\,\ket{l,l}=c_lB^l$. Their value will be determined
later on.
The degree $l-m$ polynomials $p_{l,m}$ are proportional to big $q$-Jacobi polynomials
(see~\cite[Section~4.5.3]{KS97} and~\cite{NM90}). 

Explicit expressions of some spherical harmonics are the following:
\begin{align*}
(-q)^l\ket{l,-l} &=c_l(B^*)^l\;. \\
(-q)^l\ket{l-1,-l+1} &=-qc_{l-1}(B^*)^{l-1}\;. \\
(-q)^l\ket{l,-l+1}   &=c_l[2l]^{-1/2}E\az (B^*)^l\;. \\
(-q)^l\ket{l+1,-l+1} &=c_{l+1}(-q)^{-1}\big([2l+1][2l+2][2]\big)^{-1/2}E^2\az (B^*)^{l+1}\;,
\end{align*}
where $E\az (B^*)^l$ and $E^2\az (B^*)^{l+1}$ are given by
\begin{align*}
E\az (B^*)^l &=q^{1/2}\bigl\{q^l[2l]A-[l](1-s^2)\bigr\}(B^*)^{l-1}\;, \\
E^2\az (B^*)^{l+1} &=\Big\{q^{2l+1}[2l+1][2l+2]A^2-q^{l+1}[l+1]([2l]+[2l+2])(1-s^2)A \\
 & \phantom{=}\;+q[l][l+1](1-s^2)^2-q[2l+2]s^2\Big\}(B^*)^{l-1}\;,
\end{align*}
again proved by induction on $l$.
Then, by direct computation one checks that the following equalities hold:
\begin{align*}
A\ket{l,-l+1} &=A^+_{l,-l+1}\ket{l+1,-l+1}+A^0_{l,-l+1}\ket{l,-l+1}+A^-_{l,-l+1}\ket{l-1,-l+1}\;, \\
B\ket{l,-l}   &=B^+_{l,-l}\ket{l+1,-l+1}+B^0_{l,-l}\ket{l,-l+1}+B^-_{l,-l}\ket{l-1,-l+1}\;,
\end{align*}
with coefficients
\begin{align*}
A^+_{l,-l+1} &=-\frac{c_l}{c_{l+1}}\,q^{-l+1/2}\sqrt{\frac{[2][2l]}{[2l+1][2l+2]}}\;,
& B^+_{l,-l} &=\frac{c_l}{c_{l+1}}q^{-2l}\sqrt{\frac{[2]}{[2l+1][2l+2]}}\;, \\
A^0_{l,-l+1} &=(1-s^2)q^{-l}\frac{[l][l+1][2l-2]}{[2l][2l+2][l-1]}\;,
& B^0_{l,-l} &=-(1-s^2)\frac{1-q^2}{q^{l+1/2}\sqrt{[2l]}}\,\frac{[l][l+1]}{[2l+2]}\;, \\
A^-_{l,-l+1} &=-\frac{c_l}{c_{l-1}}\,q^{-l-1/2}\,\frac{
               [2l]^2s^2+[l]^2(1-s^2)^2}{[2l][2l+1]\sqrt{[2l]}}\;,
& \qquad B^-_{l,-l} &=-\frac{c_l}{c_{l-1}}\,\frac{[2l]^2s^2+[l]^2(1-s^2)^2}{[2l][2l+1]}\;.
\end{align*}
These coefficients will be needed later.
Again, by direct computation we have
\begin{equation*}
A\ket{l-1,-l+1}=-\frac{c_{l-1}}{c_l}\,\frac{q^{-l+1/2}}{\sqrt{[2l]}}\ket{l,-l+1}
                +\textrm{orthogonal vectors}\;. \\
\end{equation*}
The representation is a $*$-representation (i.e.~$A=A^*$) only if last coefficient
equals $A^-_{l,-l+1}$. From this we deduce the condition
\begin{equation}\label{eq:cl}
\frac{c_l^2}{c_{l-1}^2}=q\,\frac{[2l][2l+1]}{[2l]^2s^2+[l]^2(1-s^2)^2}\;.
\end{equation}
When $q=1$, $c_0^{-2}$ is the area of the sphere; we fix $c_0^{-2}=4\pi$ for all $q$.
For $l\geq 1$, the coefficients $c_l$ are determined by (\ref{eq:cl}) up to a sign,
that we choose equal to $+1$:
$$
c_l=\sqrt{\frac{q^l[2l+1]!}{4\pi}}\,
\left(\,\prod_{k=1}^l\frac{1}{[2k]^2s^2+[k]^2(1-s^2)^2}\right)\!\!\,\rule{0pt}{18pt}^{\frac{1}{2}} \;.
$$

\noindent
Summarizing, the normalization constants $N_{l,m}$ of Proposition \ref{prop:Laplace} are given by
\begin{equation*}
N_{l,m}=(-1)^{l+m}\,\sqrt{\frac{q^m[2l+1]}{4\pi}}\,\sqrt{[l+m]!\,[l-m]!}
\left(\,\prod_{k=1}^l\frac{1}{[2k]^2s^2+[k]^2(1-s^2)^2}\right)\!\!\,\rule{0pt}{18pt}^{\frac{1}{2}} \;.
\end{equation*}

Using (\ref{eq:cl}) we derive the final expression for the coefficients $A_{l,-l+1}^{+,0}$
and $B_{l,-l}^{\pm,0}$:
\begin{subequations}\label{eq:bound}
\begin{align}
A^+_{l,-l+1} &=-q^{-l}\frac{\sqrt{\rule{0pt}{10pt}[2l+2]^2s^2+[l+1]^2(1-s^2)^2}}{[2l+2]}
               \sqrt{\frac{[2][2l]}{[2l+1][2l+3]}}\;, \\
A^0_{l,-l+1} &=(1-s^2)q^{-l}\frac{[l][l+1][2l-2]}{[2l][2l+2][l-1]}\;, \\
B^+_{l,-l} &=q^{-2l-1/2}\frac{\sqrt{\rule{0pt}{10pt}[2l+2]^2s^2+[l+1]^2(1-s^2)^2}}{[2l+2]}
             \sqrt{\frac{[2]}{[2l+1][2l+3]}}\;, \\
B^0_{l,-l} &=-(1-s^2)\frac{1-q^2}{q^{l+1/2}\sqrt{[2l]}}\,\frac{[l][l+1]}{[2l+2]}\;, \\
B^-_{l,-l} &=-q^{1/2}\sqrt{\frac{[2l]^2s^2+[l]^2(1-s^2)^2}{[2l][2l+1]}}\;.
\end{align}
\end{subequations}

\section{The left regular representation}\label{sec:reg}
Notice that, due to noncommutativity, left regular and right regular representations
of $\A(S^2_{qs})$ are not equivalent; indeed, they are related by the involution:
\begin{equation}\label{eq:ltor}
f\ket{l,m}=(-q)^m\{\ket{l,-m}\!f^*\}^*\;.
\end{equation}
In this section we compute the left action of the generators of
$\A(S^2_{qs})$ on spherical functions.
Since $hx=(h_{(1)}\az x)h_{(2)}$, for all $h\in U_q(su(2))$ and $x\in\A(S^2_{qs})$,
we deduce the commutation rules:
\begin{subequations}
\begin{align}
& KB=qBK\;,\qquad\quad\; KA=AK\;,                \label{eq:arA}\\
& EB=q^{-1}BE\;,\qquad EA=AE-q^{-1/2}BK\;. \label{eq:arB}
\end{align}
\end{subequations}
From (\ref{eq:arA}) we get that $B\ket{l,m}$ is in the subspace
spanned by $\ket{l',m+1}$ and $A\ket{l,m}$ is in the subspace spanned by
$\ket{l',m}$. Consider
\begin{equation*}
B\ket{l,m}=\sum_{l'\in\N}B^{l'}_{l,m}\ket{l',m+1}\;,\qquad
A\ket{l,m}=\sum_{l'\in\N}A^{l'}_{l,m}\ket{l',m}\;.
\end{equation*}
with $B^{l'}_{l,m}=0$ if $l'<m+1$ and $A^{l'}_{l,m}=0$ if $l'<m$.
From (\ref{eq:arB}) we get
\begin{align*}
E^{l-m+1}B\ket{l,m} &=q^{-l+m-1}BE^{l-m+1}\ket{l,m}=0\;, \\
\left<l',m+1\right|BE^{l'+m+2} &=q^{l'+m+2}\left<l',m+1\right|E^{l'+m+2}B=0\;,
\end{align*}
and in turn,
\begin{align*}
B_{l,m}^{l'}\sqrt{[l'-m-1]\ldots [l'-l-1]}\sqrt{[l'+m+2]\ldots [l'+l+2]} &=0\;, \\
\bar{B}_{l,m}^{l'}\sqrt{[l+m]\ldots [l-l'-1]}\sqrt{[l-m+1]\ldots [l+l']} &=0\;.
\end{align*}
The square roots in the first line are all different from zero if $l'>l+1$,
while the square roots in the second line are all different from zero if $l>l'+1$. Hence,
\begin{equation*}
B_{l,m}^{l'}=0\quad\forall\;l'\notin\{l-1,l,l+1\}\;.
\end{equation*}
A similar reasoning applies to $A\ket{l,m}$. From formula (\ref{eq:coJn}) for the
coproduct of $E^n$ we deduce:
\begin{equation*}
E^nA =\sum_{k=0}^n\sqbn{n}{k}(E^k\az A)E^{n-k}K^k
       =AE^n-[n]q^{-1/2}BE^{n-1}K\;.
\end{equation*}
This implies that $E^{l-m+2}A\ket{l,m}=0$ and
$\left<l',m\right|AE^{l'+m+2}=0$. Then
\begin{align*}
A_{l,m}^{l'}\sqrt{[l'-m]\ldots [l'-l-1]}\sqrt{[l'+m+1]\ldots [l'+l+2]}     &=0\;, \\
\bar{A}_{l,m}^{l'}\sqrt{[l+m]\ldots [l-l'-1]}\sqrt{[l-m+1]\ldots [l+l'+2]} &=0\;.
\end{align*}
As for $B$, these imply $A_{l,m}^{l'}=0$ if $l'\notin\{l,l\pm 1\}$. We can thus state:

\begin{prop}\label{prop:left}
The left regular representation of $\A(S^2_{qs})$ is defined by:
\begin{align*}
A\ket{l,m} &=A^+_{l,m}\ket{l+1,m}+A_{l,m}^0\ket{l,m}+\bar{A}^+_{l-1,m}\ket{l-1,m}\;, \\
B\ket{l,m} &=B^+_{l,m}\ket{l+1,m+1}+B_{l,m}^0\ket{l,m+1}+B^-_{l,m}\ket{l-1,m+1}\;,
\end{align*}
for suitable $A^{0,+}_{l,m},B^{0,\pm}_{l,m}\in\C$ and with $B^*$ the adjoint of $B$.
\end{prop}

The representation is equivariant by construction, and the equivariance condition
is sufficient to determine the coefficients. In fact, (\ref{eq:arB}) tells us that
the expressions
\begin{align*}
EB\ket{l,m}   &=B_{l,m}^+\sqrt{[l-m][l+m+3]}\ket{l+1,m+2}+ \\
  & \hspace{-5mm} +B^0_{l,m}\sqrt{[l-m-1][l+m+2]}\ket{l,m+2}
                  +B_{l,m}^-\sqrt{[l-m-2][l+m+1]}\ket{l-1,m+2}\;, \\
EA\ket{l,m}   &=A^+_{l,m}\sqrt{[l-m+1][l+m+2]}\ket{l+1,m+1}+ \\
  & \hspace{-5mm} +A^0_{l,m}\sqrt{[l-m][l+m+1]}\ket{l,m+1}
                  +\bar{A}^+_{l-1,m}\sqrt{[l-m-1][l+m]}\ket{l-1,m+1}\;,
\end{align*}
must be equal, respectively, to
\begin{align*}
q^{-1}BE\ket{l,m}
  &=q^{-1}\sqrt{[l-m][l+m+1]}\Big\{B^+_{l,m+1}\ket{l+1,m+2}+ \\
  & \hspace{-5mm} +B^0_{l,m+1}\ket{l,m+2}+B^-_{l,m+1}\ket{l-1,m+2}\Big\} \;, \\
(AE-q^{-1/2}BK)\ket{l,m}
  &=\Big\{\sqrt{[l-m][l+m+1]}\,A^+_{l,m+1}-q^{m-1/2}B^+_{l,m}\Big\}\ket{l+1,m+1}+ \\
  & \hspace{-5mm} +\Big\{\sqrt{[l-m][l+m+1]}\,A^0_{l,m+1}-q^{m-1/2}B^0_{l,m}\Big\}\ket{l,m+1}+ \\
  & \hspace{-5mm} +\Big\{\sqrt{[l-m][l+m+1]}\,\bar{A}^+_{l-1,m+1}-q^{m-1/2}B^-_{l,m}\Big\}\ket{l-1,m+1}\;.
\end{align*}
From the independence of the vectors $\ket{l,m}$, we derive the conditions
\begin{subequations}
\begin{align}
\frac{q^{-m}B^+_{l,m}}{\sqrt{[l+m+1][l+m+2]}} &=\frac{q^{-m-1}B^+_{l,m+1}}{\sqrt{[l+m+2][l+m+3]}} \\
\frac{q^{-m}B^0_{l,m}}{\sqrt{[l-m][l+m+1]}}   &=\frac{q^{-m-1}B^0_{l,m+1}}{\sqrt{[l-m-1][l+m+2]}} \\
\frac{q^{-m}B^-_{l,m}}{\sqrt{[l-m][l-m-1]}}   &=\frac{q^{-m-1}B^-_{l,m+1}}{\sqrt{[l-m-1][l-m-2]}} \\
\sqrt{[l-m][l+m+1]}\,A^+_{l,m+1}-q^{m-1/2}B^+_{l,m}         &=A^+_{l,m}\sqrt{[l-m+1][l+m+2]}  \label{eq:uno} \\
A^0_{l,m+1}-\frac{q^{m-1/2}B^0_{l,m}}{\sqrt{[l-m][l+m+1]}}  &=A^0_{l,m}   \label{eq:due} \\
\sqrt{[l-m+1][l+m+2]}\,\bar{A}^+_{l,m+1}-q^{m-1/2}B^-_{l+1,m} &=\bar{A}^+_{l,m}\sqrt{[l-m][l+m+1]}    \label{eq:tre}
\end{align}
\end{subequations}
The first three conditions are of the form $f_{l,m}=f_{l,m+1}$,
stating that they involve quantities $f_{l,m}$ which do not
depend on $m$. These can be rewritten as `$f_{l,m}=f_{l,-l}$',
that is:
\begin{align*}
\frac{q^{-m}B^+_{l,m}}{\sqrt{[l+m+1][l+m+2]}} &=\frac{q^lB^+_{l,-l}}{\sqrt{[2]}}  \;, \\
\frac{q^{-m}B^0_{l,m}}{\sqrt{[l-m][l+m+1]}}   &=\frac{q^lB^0_{l,-l}}{\sqrt{[2l]}} \;, \\
\frac{q^{-m}B^-_{l,m}}{\sqrt{[l-m][l-m-1]}}   &=\frac{q^lB^-_{l,-l}}{\sqrt{[2l][2l-1]}}\;.
\end{align*}
Using the boundary condition in (\ref{eq:bound}) we arrive at the final result:
\begin{subequations}\label{eq:B}
\begin{align}
B^+_{l,m} &=q^{-(l-m)-1/2}\frac{\sqrt{\rule{0pt}{10pt}[2l+2]^2s^2+[l+1]^2(1-s^2)^2}}{[2l+2]}
             \sqrt{\frac{[l+m+1][l+m+2]}{[2l+1][2l+3]}}\;, \\
B^0_{l,m} &=-(1-s^2)(1-q^2)q^{m-1/2}\,\frac{[l][l+1]}{[2l][2l+2]}\,\sqrt{[l-m][l+m+1]}\;, \\
B^-_{l,m} &=-q^{l+m+1/2}\frac{\sqrt{\rule{0pt}{10pt}[2l]^2s^2+[l]^2(1-s^2)^2}}{[2l]}
            \sqrt{\frac{[l-m][l-m-1]}{[2l-1][2l+1]}} \;.
\end{align}
\end{subequations}
Equations (\ref{eq:uno}-\ref{eq:due}) have solution
\begin{subequations}\label{eq:A}
\begin{align}
A^+_{l,m} &=-q^{m-1}\frac{\sqrt{\rule{0pt}{10pt}[2l+2]^2s^2+[l+1]^2(1-s^2)^2}}{[2l+2]}
             \sqrt{\frac{[l+m+1][l-m+1]}{[2l+1][2l+3]}}\;, \\
A^0_{l,m} &=(1-s^2)q^{-1}(1+q^{2m})\frac{[l][l+1]}{[2l][2l+2]}\;,
\end{align}
\end{subequations}
while equation (\ref{eq:tre}) is automatically satisfied.

The next step is to lift diagonally this representation on the
Hilbert space completion of $\A(S^2_{qs})^2$, we call it $\HH$,
and then decompose it into two subrepresentations on $\HH_{\pm}$,
the Hilbert space completion of the modules $\A(S^2_{qs})^2(1-e_s)$
and $\A(S^2_{qs})^2e_s$ for a suitable invariant idempotent $e_s$.
These will be the $q$-analogue of the chiral representations of
$\A(S^2)$, and elements in $\HH_{\pm}$ will be the $q$-analogue
of Weyl spinors.

\section{Chiral representations}
As explained in Section \ref{sec:3.2}, equivariant representations
correspond to invariant idempotents.
Let $\sigma:U_q(su(2))\to\mr{Mat}_2(\C)$ be the representation $V_{1/2}$,
written in matrix form
\begin{equation}\label{eq:qPaul}
\sigma(K)=\ma{14pt}{q^{1/2} & 0 \\ 0 & q^{-1/2}}\;,\qquad
\sigma(E)=\ma{14pt}{0 & 1 \\ 0 & 0}\;,\qquad
\sigma(F)=\ma{14pt}{0 & 0 \\ 1 & 0}\;.
\end{equation}
For any $s$ (and $q$) there is a $2\times 2$ idempotent which is
invariant, in the sense of (\ref{eq:cov}), and given by
\begin{equation}\label{eq:qBott}
e_s=\frac{1}{1+s^2}\ma{14pt}{s^2+A & qB \\ q^{-1}B^* & 1-q^2A}\;.
\end{equation}
It also satisfies $\kappa(e^*_s)=e_s$ with $\kappa(a):=K^2\az a$
the modular automorphism~\cite[Section 11.3.4]{KS97} thus,
by Lemma \ref{lemma:next} it corresponds to an orthogonal
projection via the formula (\ref{eq:pi}).
We notice that a projection (a hermitian idempotent) equivalent
to $e_s$ is given by the matrix $\sigma(K)^{-1}e_s\sigma(K)$;
this was introduced in \cite{BM98}, where the authors
compute its charge by pairing it with the cyclic $0$-cocycle
of \cite{MNW91}, and later generalized to projections of arbitrary
charge first for the standard sphere in \cite{HM99}, and then for
all Podle\'s spheres in \cite{BM00}.

Let $Y_{l,m}=\ket{l,m}$ be the spherical functions defined in Proposition
\ref{prop:Laplace} and identify $V_{1/2}$ with $\C^2$ via the map
$\ket{1/2,1/2}\mapsto (1,0)$ and $\ket{1/2,-1/2}\mapsto (0,1)$.
Using Clebsch-Gordan coefficients in equation (\ref{eq:Clebsch})
we define a basis of $\A(S^2_{qs})^2$,
\begin{subequations}\label{eq:ket12}
\begin{align}
\ket{l,m}_1 &=\left(
\sqrt{\tfrac{q^{l-m}[l+m]}{[2l]}}\,Y_{l-\frac{1}{2},m-\frac{1}{2}}
\;,\;
\sqrt{\tfrac{q^{-(l+m)}[l-m]}{[2l]}}\,Y_{l-\frac{1}{2},m+\frac{1}{2}}
\right)\;, \\
\ket{l,m}_2 &=\left(
-\sqrt{\tfrac{q^{-(l+m+1)}[l-m+1]}{[2l+2]}}\,Y_{l+\frac{1}{2},m-\frac{1}{2}}
\;,\;
\sqrt{\tfrac{q^{l-m+1}[l+m+1]}{[2l+2]}}\,Y_{l+\frac{1}{2},m+\frac{1}{2}}\right) \;,
\end{align}
\end{subequations}
(with $l\in\N+\frac{1}{2}$ and $m=-l,-l+1,\ldots,l$) which
corresponds to the decomposition
$$
\left(\bigoplus\nolimits_{l\in\N}V_l\right)
\otimes V_{1/2}\simeq\bigoplus\nolimits_{l\in\N+\frac{1}{2}}(V_l\oplus V_l)\;.
$$

Once fixed the eigenvalue of $K$ and the spin $l$, the corresponding eigenvectors
of $\A(S^2_{qs})^2$ span a two-dimensional space. Since $e_s$ is equivariant
we can search for a basis of this two-dimensional space made of eigenvectors of $e_s$.
In (\ref{eq:pi}) the idempotent $e_s$ multiplies vectors from the right, thus
to find eigenvectors of the corresponding projection we must work with the
right regular representation, related to the left regular one by (\ref{eq:ltor}),
which means
\begin{align*}
\ket{l,m}A   &=q^{-2m}A\ket{l,m}\;, \\
-\ket{l,m}q^{-1}B^* &=B^+_{l,-m}\ket{l+1,m-1}+B^0_{l,-m}\ket{l,m-1}
                     +B^-_{l,-m}\ket{l-1,m-1}\;, \\
-\ket{l,m}qB &=B^-_{l+1,-m-1}\ket{l+1,m+1}+B^0_{l,-m-1}\ket{l,m+1}
              +B^+_{l-1,-m-1}\ket{l-1,m+1}\;.
\end{align*}
A straightforward computation shows that the basis (\ref{eq:ket12})
satisfies
\begin{align*}
\ket{l,m}_1e_s =\; &
\frac{s^2q^{l+\frac{1}{2}}+q^{-l-\frac{1}{2}}}{1+s^2}\,\frac{[l+\tfrac{1}{2}]}{[2l+1]}
\ket{l,m}_1+c_{l+\frac{1}{2},s}\,\frac{1}{1+s^2}\,\ket{l,m}_2\;, \\
\ket{l,m}_2e_s =\; & \frac{q^{l+\frac{1}{2}}+s^2q^{-l-\frac{1}{2}}}{1+s^2}\,
\frac{[l+\frac{1}{2}]}{[2l+1]}\ket{l,m}_2
+c_{l+\frac{1}{2},s}\,\frac{1}{1+s^2}\,\ket{l,m}_1 \;,
\end{align*}
where $c_{l,s}$ is defined by
\begin{equation*}
c_{l,s}=\frac{\sqrt{\rule{0pt}{10pt}[2l]^2s^2+[l]^2(1-s^2)^2}}{[2l]}
       =\frac{[l]}{[2l]}\sqrt{(q^ls^2+q^{-l})(q^l+q^{-l}s^2)}\;.
\end{equation*}
Then the vectors $\ket{l,m}_{\pm}$, defined by
\begin{subequations}\label{eq:base}
\begin{align}
\sqrt{\frac{[2l+1]}{[l+\frac{1}{2}]}}\,\ket{l,m}_+ =\; &
\sqrt{\frac{q^{l+\frac{1}{2}}s^2+q^{-l-\frac{1}{2}}}{1+s^2}}\ket{l,m}_1+
\sqrt{\frac{q^{l+\frac{1}{2}}+q^{-l-\frac{1}{2}}s^2}{1+s^2}}\ket{l,m}_2\;, \\
\sqrt{\frac{[2l+1]}{[l+\frac{1}{2}]}}\,\ket{l,m}_- =\; &
\sqrt{\frac{q^{l+\frac{1}{2}}+q^{-l-\frac{1}{2}}s^2}{1+s^2}}\ket{l,m}_1-
\sqrt{\frac{q^{l+\frac{1}{2}}s^2+q^{-l-\frac{1}{2}}}{1+s^2}}\ket{l,m}_2\;,
\end{align}
\end{subequations}
satisfy the equations
\begin{equation*}
\ket{l,m}_+(1-e_s)=0\;,\quad\qquad\ket{l,m}_-e_s=0\;.
\end{equation*}
The `spinorial harmonics' in (\ref{eq:base}) form an
orthonormal basis of $\HH:=L^2(S^2_{qs})^2$.
We call $\HH_+$ the Hilbert space with basis $\ket{l,m}_+$,
and $\HH_-$ the Hilbert space with basis $\ket{l,m}_-$.
Both are finitely generated projective modules,
\begin{center}
$\HH_+=L^2(S^2_{qs})^2e_s\,\,$ and $\;\,\HH_-=L^2(S^2_{qs})^2(1-e_s)\;$.
\end{center}
They carry representations of $U_q(su(2))$ --- each equivalent to
$\bigoplus_{l\in\N+\frac{1}{2}}V_l$ --- and `chiral' representations
of $\A(S^2_{qs})$, which are equivariant with respect to the action
of $U_q(su(2))$.

Using (\ref{eq:ket12}) and (\ref{eq:base}) one can write $\ket{l,m}_\pm$
in terms of spherical functions. Then, from Proposition \ref{prop:left}
one computes the action of the generators $A$ and $B$ on
the vectors $\ket{l,m}_\pm$. The result is described in the following proposition.

\begin{prop}\label{prop:chiral}
The action of generators of $\A(S^2_{qs})$ on the basis $\ket{l,m}_\pm$
is given by:{\footnotesize
\begin{align*}
A\ket{l,m}_+=\; &
-q^{m-1}\,\sqrt{\frac{[l+\frac{1}{2}][l+\frac{3}{2}]}{[2l+1][2l+3]}}\,
\frac{\sqrt{[l+m+1][l-m+1]}}{[2l+2]}\,\sqrt{(q^{l+\frac{1}{2}}+q^{-l-\frac{1}{2}}s^2)
(q^{l+\frac{3}{2}}s^2+q^{-l-\frac{3}{2}})}\ket{l+1,m}_+ \\ &
-q^{m-1}\,\sqrt{\frac{[l-\frac{1}{2}][l+\frac{1}{2}]}{[2l-1][2l+1]}}\,\frac{
\sqrt{[l+m][l-m]}}{[2l]}\,\sqrt{(q^{l+\frac{1}{2}}s^2+q^{-l-\frac{1}{2}})
(q^{l-\frac{1}{2}}+q^{-l+\frac{1}{2}}s^2)} \ket{l-1,m}_+ \\ &
+\left((1-s^2)\frac{(1-q^2)[l-\tfrac{1}{2}][l+\tfrac{3}{2}]+1}{1+q^2}
-1\right)\frac{[l+m+1][l-m]-q^{-2}[l+m][l-m+1]}{[2l][2l+2]}\ket{l,m}_+ \\ &
+\frac{1-s^2}{1+q^2}\ket{l,m}_+\;, \\
A\ket{l,m}_- =\; &
-q^{m-1}\,\sqrt{\frac{[l+\frac{1}{2}][l+\frac{3}{2}]}{[2l+1][2l+3]}}\,
\frac{\sqrt{[l+m+1][l-m+1]}}{[2l+2]}\,\sqrt{(q^{l+\frac{1}{2}}s^2+q^{-l-\frac{1}{2}})
(q^{l+\frac{3}{2}}+q^{-l-\frac{3}{2}}s^2)}\ket{l+1,m}_- \\ &
-q^{m-1}\,\sqrt{\frac{[l-\frac{1}{2}][l+\frac{1}{2}]}{[2l-1][2l+1]}}\,
\frac{\sqrt{[l+m][l-m]}}{[2l]}\,\sqrt{(q^{l+\frac{1}{2}}+q^{-l-\frac{1}{2}}s^2)
(q^{l-\frac{1}{2}}s^2+q^{-l+\frac{1}{2}})}\ket{l-1,m}_- \\ &
+\left((1-s^2)\frac{(1-q^2)[l-\tfrac{1}{2}][l+\tfrac{3}{2}]+1}{1+q^2}
+s^2\right)\frac{[l+m+1][l-m]-q^{-2}[l+m][l-m+1]}{[2l][2l+2]}\ket{l,m}_- \\ &
+\frac{1-s^2}{1+q^2}\ket{l,m}_-\;, \\
B\ket{l,m}_+ =\; &
q^{-(l-m)-1/2}\,\sqrt{\frac{[l+\frac{1}{2}][l+\frac{3}{2}]}{[2l+1][2l+3]}}\,
\frac{\sqrt{[l+m+1][l+m+2]}}{[2l+2]}\,\sqrt{(q^{l+\frac{1}{2}}+q^{-l-\frac{1}{2}}s^2)
(q^{l+\frac{3}{2}}s^2+q^{-l-\frac{3}{2}})}\ket{l+1,m+1}_+ \\ &
-q^{l+m+1/2}\,\sqrt{\frac{[l-\frac{1}{2}][l+\frac{1}{2}]}{[2l-1][2l+1]}}
\,\frac{\sqrt{[l-m][l-m-1]}}{[2l]}\,\sqrt{(q^{l+\frac{1}{2}}s^2+q^{-l-\frac{1}{2}})
(q^{l-\frac{1}{2}}+q^{-l+\frac{1}{2}}s^2)}\ket{l-1,m+1}_+ \\ &
+q^{m-\frac{1}{2}}\frac{\sqrt{[l+m+1][l-m]}}{[2l][2l+2]}
\Big\{s^2+q^2-(1-s^2)(1-q^2)[l-\tfrac{1}{2}][l+\tfrac{3}{2}]\Big\}\ket{l,m+1}_+ \;, \\
B\ket{l,m}_- =\; &
q^{-(l-m)-1/2}\,\sqrt{\frac{[l+\frac{1}{2}][l+\frac{3}{2}]}{[2l+1][2l+3]}}
\,\frac{\sqrt{[l+m+1][l+m+2]}}{[2l+2]}\,\sqrt{(q^{l+\frac{1}{2}}s^2+q^{-l-\frac{1}{2}})
(q^{l+\frac{3}{2}}+q^{-l-\frac{3}{2}}s^2)}\ket{l+1,m+1}_- \\ &
-q^{l+m+1/2}\,\sqrt{\frac{[l-\frac{1}{2}][l+\frac{1}{2}]}{[2l-1][2l+1]}}
\,\frac{\sqrt{[l-m][l-m-1]}}{[2l]}\,\sqrt{(q^{l+\frac{1}{2}}+q^{-l-\frac{1}{2}}s^2)
(q^{l-\frac{1}{2}}s^2+q^{-l+\frac{1}{2}})}\ket{l-1,m+1}_- \\ &
-q^{m-\frac{1}{2}}\frac{\sqrt{[l+m+1][l-m]}}{[2l][2l+2]}
\Big\{s^2q^2+1+(1-s^2)(1-q^2)[l-\tfrac{1}{2}][l+\tfrac{3}{2}]\Big\}
\ket{l,m+1}_- \;. \\
\end{align*} }
\end{prop}

\vspace{-10pt}

\noindent
When $s=1$ these representations are equivalent to those
considered in~\cite{DLPS05} (the generators there being $b:=q^2A$ and $a:=B$),
and when $s=0$ they are equivalent to those considered in~\cite{DS03}
(the generators there being $A'=A$ and $B'=q^{-1}B$).

\section{Comparison with the literature}
Integrable representations of the crossed product algebra $\A(S^2_{qs})\rtimes U_q(su(2))$
were classified in \cite{SW03}. We rewrite them here.
Firstly, the following set of generators for $\A(S^2_{qs})$
\begin{equation}\label{eq:genchange}
x_1=(q[2])^{\frac{1}{2}}\,B\;,\qquad x_0=-q[2]A+(1-s^2)
\;,\qquad x_{-1}=-q^{-1}x_1^*\;,
\end{equation}
satisfy the commutation rules
\begin{subequations}\label{eq:DDLWtoD}
\begin{align}
\rule{0pt}{14pt}\label{eqq:uno}
x_1x_0-q^{-2}x_0x_1-t(1-q^{-2})x_1    &=0\;, \\
\rule{0pt}{14pt}\label{eqq:due}
x_{-1}x_0-q^2x_0x_{-1}-t(1-q^2)x_{-1} &=0\;, \\
\rule{0pt}{14pt}\label{eqq:tre}
-[2]x_{-1}x_1+(q^2x_0+t)(x_0-t)-[2]^2(1-t) &=0\;, \\
\rule{0pt}{14pt}\label{eqq:quattro}
-[2]x_1x_{-1}+(q^{-2}x_0+t)(x_0-t)-[2]^2(1-t) &=0\;.
\end{align}
\end{subequations}
with parameter $t:=1-s^2\in [0,1]$. These are exactly equations
(3.1) of \cite{DDLW07}.
We rewrite the representations in Proposition 4.1 of \cite{DDLW07} as follows.

\begin{prop}\label{prop:crossrep}
For fixed $N\in\frac{1}{2}\Z$, we call $\HH_N$ the Hilbert space
completion of
\begin{equation*}
\mr{Span}\big\{\ket{l,m;N}\,;\;l-|N|,l-|m|\in\N\big\}\;,
\end{equation*}
in which the vectors $\ket{l,m;N}$ are declared to be orthonormal.
A $*$-representation of $\A(S^2_{qs})$ on $\HH_N$ is given by
\begin{align*}
     x_1\ket{l,m;N} &=\tilde{B}^+_{l,m}\ket{l+1,m+1;N}+\tilde{B}_{l,m}^0\ket{l,m+1;N}+\tilde{B}^-_{l,m}\ket{l-1,m+1;N}\;, \\
     x_0\ket{l,m;N} &=\tilde{A}^+_{l,m}\ket{l+1,m;N}+\tilde{A}_{l,m}^0\ket{l,m;N}+\tilde{A}^+_{l-1,m}\ket{l-1,m;N}\;, \\
-qx_{-1}\ket{l,m;N} &=\tilde{B}^-_{l+1,m-1}\ket{l+1,m-1;N}+\tilde{B}_{l,m-1}^0\ket{l,m-1;N}+\tilde{B}^+_{l-1,m-1}\ket{l-1,m-1;N}\;,
\end{align*}
with coefficients
\begin{align*}
\tilde{A}^+_{l,m} &=q^m\sqrt{[2][l-m+1][l+m+1]}\,\alpha_{N,l+1}(s) \;,\\
\tilde{A}^0_{l,m} &=\bigl([2l]-q^{l+m+1}[2][l-m]\bigr)\,\beta_{N,l}(s) \;, \\
\tilde{B}^+_{l,m} &=q^{-l+m}\sqrt{[l+m+1][l+m+2]}\,\alpha_{N,l+1}(s) \;,\\
\tilde{B}^0_{l,m} &=-q^{m+2}\sqrt{[2][l-m][l+m+1]}\,\beta_{N,l}(s) \;,\\
\tilde{B}^-_{l,m} &=-q^{l+m+1}\sqrt{[l-m][l-m-1]}\,\alpha_{N,l}(s) \;,
\end{align*}
and
\begin{align*}
\alpha_{N,l}(s) &=\frac{1}{[2l]}\sqrt{\frac{[2][l+N][l-N]}{[2l-1][2l+1]}}\,
   \sqrt{(q^{l+N}+q^{-(l+N)}s^2)(q^{l-N}s^2+q^{-(l-N)})} \;, \\
\beta_{N,l}(s) &=\frac{\big([2l]-q^{l+N+1}[2][l-N]\big)
   (1-s^2)+q[2][2N]s^2}{q^2[2l][2l+2]} \;.
\end{align*}
\end{prop}

\bigskip

Now we want to check (i) that the representations are indeed correct,
for any $N\in\frac{1}{2}\Z$; (ii) that these are the same equations
as in Proposition 4.1 of \cite{DDLW07}; (iii) that for $N=\pm 1/2$
we get the same representations of Proposition \ref{prop:chiral}
by identifying $\ket{l,m;\pm\oh}$ with $\ket{l,m}_\mp$ and using
the change of generators (\ref{eq:genchange}).

Concerning point (ii), the non-trivial part is to show that
\begin{equation}\label{eq:checkab}
[2l-1]^{\frac{1}{2}}[2l]^{\frac{1}{2}}\,\alpha_{N,l}(s)=\alpha_N(l)
\;,\qquad
[2l]\beta_{N,l}(s)=\beta_N(l)\;,
\end{equation}
where $\alpha_N(l)$ and $\beta_N(l)$ are given by equations
(4.3) of \cite{DDLW07},
\begin{subequations}\label{eq:abNl}
\begin{align}
\alpha_N(l) &=\sqrt{\frac{[2][l+N][l-N][2l]}{[2l+1]}}
\sqrt{\frac{1-t}{[l]^2}+q^{-2N}\frac{(t-1+q^{2N})^2}{[2l]^2}} \;, \\
\beta_N(l)  &=\frac{[2N]\bigl([2]-q^{\mr{sign}\,N}t\bigr)
+t(q^{-1}-q)[l-|N|][l+|N|+1]}{q[2l+2]} \;. \label{eq:betaNl}
\end{align}
\end{subequations}
That (\ref{eq:betaNl}) and equation (4.3a) of \cite{DDLW07} are
the same follows from the identity
$$
[l][l+1]-[n][n+1]=[l-n][l+n+1] \;.
$$
To prove (\ref{eq:checkab}) is a straightforward computation,
based on the relations $[2l][l]^{-1}=q^l+q^{-l}$ and
$[2](1-t)=[2]-q^{\mr{sign}\,N}t-q^{-\mr{sign}\,N}t$,
from which one deduces the identities
\begin{gather*}
(q^{l+N}+q^{-(l+N)}s^2)(q^{l-N}s^2+q^{-(l-N)})=
[2l]^2[l]^{-2}s^2+q^{-2N}(q^{2N}-s^2)^2 \;, \\
[2l]-q^{l+N+1}[2][l-N]=q^{1-\mr{sign}\,N}[2N]+(1-q^2)[l-|N|][l+|N|+1] \;,
\end{gather*}
and on the definition $t=1-s^2$.
It is still not obvious that $\alpha^0_0(l,m;N)$ in (4.2b)
equals $\tilde{A}^0_{l,m}$. This follows from the identity
$$
[2l]-q^{l+m+1}[2][l-m]=[l+m][l-m+1]-q^2[l+m+1][l-m] \;,
$$
which concludes the check of point (ii). With previous identity,
also the verification of point (iii) is a long but simple
computation.

To control that the ones defined in Proposition \ref{prop:crossrep}
are representations --- i.e.~that equations (\ref{eq:DDLWtoD}) are
satisfied --- we make use of {\mate}.
Now, (\ref{eqq:due}) is the adjoint of (\ref{eqq:uno}). Moreover,
(\ref{eqq:tre}) and (\ref{eqq:quattro}) are identities between
hermitian elements: in matrix notation, we need to compute only
the upper-triangular elements. Thus, we need only to verify
that the following matrix elements are zero:{\small
\begin{align*}
\inner{l+2,m+1|\,(\ref{eqq:uno})\,|l,m} &=
  B^+_{l+1,m}A^+_{l,m}-q^{-2}A^+_{l+1,m+1}B^+_{l,m}  \;, \\
\inner{l+1,m+1|\,(\ref{eqq:uno})\,|l,m} &=
  B^0_{l+1,m}A^+_{l,m}+B^+_{l,m}(A^0_{l,m}-q^{-2}A^0_{l+1,m+1}+tq^{-2}-t)
  -q^{-2}A^+_{l,m+1}B^0_{l,m}  \;, \\
\inner{l,m+1|\,(\ref{eqq:uno})\,|l,m}   &=
  B^+_{l-1,m}A^+_{l-1,m}-q^{-2}A^+_{l-1,m+1}B^-_{l,m}+
  B^-_{l+1,m}A^+_{l,m}-q^{-2}A^+_{l,m+1}B^+_{l,m} \\ &\phantom{=}\;+
  B^0_{l,m}\big(  A^0_{l,m}-q^{-2}A^0_{l,m+1}+tq^{-2}-t  \big)  \;, \\
\inner{l-1,m+1|\,(\ref{eqq:uno})\,|l,m} &=
  B^0_{l-1,m}A^+_{l-1,m}+B^-_{l,m}(A^0_{l,m}-q^{-2}A^0_{l-1,m+1}+tq^{-2}-t)
  -q^{-2}A^+_{l-1,m+1}B^0_{l,m}  \;, \\
\inner{l-2,m+1|\,(\ref{eqq:uno})\,|l,m} &=
  B^-_{l-1,m}A^+_{l-1,m}-q^{-2}A^+_{l-2,m+1}B^-_{l,m}  \;, \\
\inner{l+2,m|\,(\ref{eqq:tre})\,|l,m}   &=
  q^{-1}[2]B^-_{l+2,m}B^+_{l,m}+q^2A^+_{l+1,m}A^+_{l,m}  \;, \\
\inner{l+1,m|\,(\ref{eqq:tre})\,|l,m}   &=
  q^{-1}[2](B^-_{l+1,m}B^0_{l,m}+B^0_{l+1,m}B^+_{l,m})+q^2A^+_{l,m}
  \big( A^0_{l,m}+A^0_{l+1,m}+tq^{-2}-t \big)  \;, \\
\inner{l,m|\,(\ref{eqq:tre})\,|l,m}     &=
  q^{-1}[2]\big\{(B^+_{l,m})^2+(B^0_{l,m})^2+(B^-_{l,m})^2\big\}+
  q^2(A^+_{l,m})^2+q^2(A^+_{l-1,m})^2 \\ &\phantom{=}\;+
  (q^2A^0_{l,m}+t)(A^0_{l,m}-t)-[2]^2(1-t)  \;, \\
\inner{l+2,m|\,(\ref{eqq:quattro})\,|l,m}   &=
  q^{-1}[2]B^+_{l+1,m-1}B^-_{l+1,m-1}+q^{-2}A^+_{l+1,m}A^+_{l,m}  \;, \\
\inner{l+1,m|\,(\ref{eqq:quattro})\,|l,m}   &=
  q^{-1}[2](B^+_{l,m-1}B^0_{l,m-1}+B^0_{l+1,m-1}B^-_{l+1,m-1})+q^{-2}A^+_{l,m}
  \big( A^0_{l,m}+A^0_{l+1,m}+tq^2-t \big)  \;, \\
\inner{l,m|\,(\ref{eqq:quattro})\,|l,m}     &=
  q^{-1}[2]\big\{(B^+_{l-1,m-1})^2+(B^0_{l,m-1})^2+(B^-_{l+1,m-1})^2\big\}+
  q^{-2}(A^+_{l,m})^2+q^{-2}(A^+_{l-1,m})^2 \\ &\phantom{=}\;+
  (q^{-2}A^0_{l,m}+t)(A^0_{l,m}-t)-[2]^2(1-t) \;,
\end{align*}
}where `tildes' have been removed to make the equations
clearer. The inputs and outputs of {\mate} are listed in Table
\ref{tab:math}: the zeroes in the outputs prove the statement.
Coefficients $\{\tilde{A}^+_{l,m},\tilde{A}^0_{l,m},
\tilde{B}^+_{l,m},\ldots\}$ are called respectively
$\{$\texttt{Apiu[l,m]},\texttt{Azero[l,m]},%
\texttt{Bpiu[l,m]}$,\ldots\}$.

\begin{table}
\begin{center}\begin{scriptsize}\begin{ttfamily}
\begin{longtable}{|p{15cm}|}
\hline
\rule{0pt}{20pt}\hspace*{1cm} {\small{\mate} code} \\
\rule{0pt}{20pt}\hspace*{-5mm}%
\includegraphics[width=\textwidth]{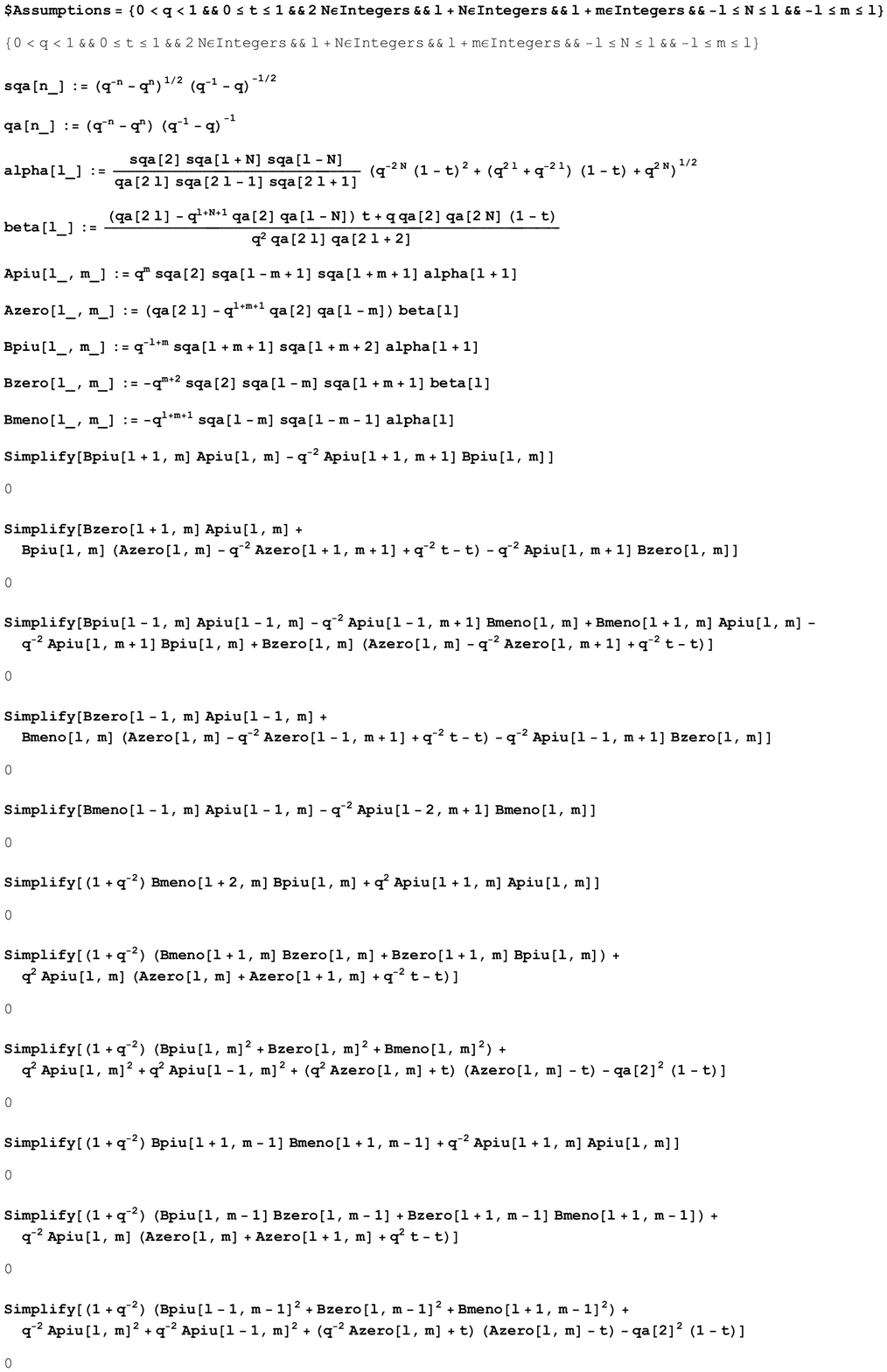} \\
\hline
\end{longtable}
\end{ttfamily}\end{scriptsize}

\caption{Check of the representations in Proposition \ref{prop:crossrep}.}\label{tab:math}

\end{center}
\end{table}

\newpage
\section{A family of equivariant spectral triples}\label{sec:6.6}
From now on, $N\in\frac{1}{2}\Z$ is fixed and different from $0$.
Remind that the Hilbert space $\HH_N$ has orthonormal basis $\ket{l,m;N}$ with
$l-|N|,l-|m|\in\N$, and carries the equivariant representation of $\A(S^2_{qs})$
in Proposition~\ref{prop:crossrep}. We wish to study equivariant even spectral
triples on $\HH:=\HH_{-N}\oplus\HH_N$, with natural grading $\gamma:=1\oplus -1$.

Let $F$ be the unitary operator\footnote{This operator $F$ should not be confused
with the generator of $U_q(su(2))$ denoted with the same symbol.}
$$
F\ket{l,m;\pm N}:=\ket{l,m;\mp N} \;.
$$
Given any divergent sequence of \emph{positive} numbers $\{\lambda_n\}_{n\in\N}$,
we define a self-adjoint operator $D_\lambda=|D_\lambda|F$ on
$\HH:=\HH_N\oplus\HH_{-N}$ by
$$
|D_\lambda|\ket{l,m;\pm N}:=\lambda_l\ket{l,m;\pm N}
$$
with domain the linear span of basis vectors.
Since $\lambda_n>0$ and $\lambda_n^{-1}\to 0$, $D_\lambda$ is invertible
and with compact inverse. It commutes with the action of $U_q(su(2))$,
being proportional to the identity in each irreducible representation.
So, $D_\lambda$ is a candidate for an equivariant Dirac operator.
For a fixed $N$, since all $D_\lambda$ have the same sign $F$, they correspond
to different `metric' structures that are `conformally' equivalent (c.f.~e.g.~\cite{Bar07}).

Among this family, we distinguish the operator
\begin{equation}\label{eq:isoD}
D\ket{l,m;\pm N}:=(l-|N|+1)\ket{l,m;\mp N} \;,
\end{equation}
which for $N=\pm\frac{1}{2}$ is isospectral to the Dirac operator of the
round sphere $S^2$ (cf.~e.g.~Chapter 2 of \cite{Var06}).

\begin{lemma}
The coefficients $\alpha_N(l)$ and $\beta_N(l)$ in (\ref{eq:abNl})
satisfy the inequalities
\begin{subequations}\label{eq:ineqAB}
\begin{align}
\bigl|\alpha_N(l)-s\sqrt{1+q^2}\bigr| &\leq C_{N,s}q^l \label{eq:ineqAa}\;,\\
\bigl|\beta_N(l)-(1-s^2)\bigr| &\leq C'_{N,s}q^l \label{eq:ineqB} \;,\\
\intertext{and in particular if $s=0$}
|\alpha_N(l)| &\leq C_{N,0}q^l \;, \label{eq:ineqAb}
\end{align}
\end{subequations}
with $C_{N,s}$ and $C'_{N,s}$ some positive constants.
\end{lemma}
\begin{prova}
Remind the definition of the parameter $t$, $t:=1-s^2$.
Equations (\ref{eq:ineqB}) and (\ref{eq:ineqAb}) are an immediate
consequence of the definition (\ref{eq:abNl}) and of the inequalities
$$
[n]\leq (q^{-1}-q)^{-1}q^{-n}\quad\mr{and}\quad
[n+1]^{-1}\leq q^n\;,\quad\textrm{for all}\;n\geq 0\;.
$$
By (\ref{eq:ineqAb}), we need to prove (\ref{eq:ineqAa}) only for $s\neq 0$,
i.e.~$t\neq 1$.
For $t\in[0,1)$ and $l\geq 1$, set 
$$
u_l:=\frac{[l+N][l-N][2l]}{q[2l+1][l]^2}\,
\frac{1-t+q^{-2N}[2l]^{-2}[l]^{2}(t-1+q^{2N})^2}{1-t}-1 \;.
$$
From $\,[l+N][l-N][2l]\leq q[2l+1][l]^2\,$ we obtain
\begin{equation*}
u_l\leq
\frac{q^{-2N}[2l]^{-2}[l]^2(t-1+q^{2N})^2}{1-t}=
C_N(1-q^2)^2[2l]^{-2}[l]^2\leq C_Nq^{2l} \;,
\end{equation*}
with $C_N:=\frac{\{q^N-q^{-N}(1-t)\}^2}{(1-q^2)^2(1-t)}$.
On the other hand, from the equalities $[l+N][l-N]=[l]^2-[N]^2$
and \mbox{$[2l]-q[2l+1]=-q^{2l+1}$} we get
$$
u_l\geq\frac{[2l][l+N][l-N]}{q[2l+1][l]^2}-1=
\frac{-[2l][N]^2-q^{2l+1}[l]^2}{q[2l+1][l]^2}\geq
-q^{2l-2}\frac{[N]^2}{1-q^2}-q^{4l-1} \;.
$$
Hence there exists a constant $C'_N\in\R^+$ such that $|u_l|\leq C'_Nq^{2l}$.
Moreover from the inequality $|\sqrt{1+u}-1|\leq |u|\,$, which holds 
for any $u\geq -1$, it follows that
$$
\left|\frac{\alpha_N(l)}{\sqrt{(1+q^2)(1-t)}}-1\right|=
\bigl|\sqrt{1+u_l}-1 \bigr|\leq |u_l|\leq C'_Nq^{2l}\;.
$$
This is exactly (\ref{eq:ineqAa}).
\end{prova}

\begin{lemma}\label{lemma:6.8}
The coefficients
$$
C^{j,k}_{i,\pm}(l,m)=\inner{l+j,m+k;\pm N|x_i|l,m;\pm N}
$$
satisfy (\ref{eq:ineqqn}), that is
$|C_{i,+}^{j,k}(l,m)-C_{i,-}^{j,k}(l,m)|\leq K_i^{j,k}q^l$
for suitable $K_i^{j,k}>0$. If $s=0$, they satisfy (\ref{eq:Kprime}) too,
that is $|C_{i,\pm}^{j,k}(l,m)|\leq K'^{j,k}_{i,\pm}\,q^l$ for suitable
$K'^{j,k}_{i,\pm}>0$.
\end{lemma}
\begin{prova}
By Proposition \ref{prop:crossrep} and equation (\ref{eq:checkab}) we
see that $C^{j,k}_{i,\pm}(l,m)$ is zero if $k\neq i$ and $j\notin\{0,\pm 1\}$,
while
\begin{align*}
C^{+1,i}_{i,\pm}(l,m) &=c^{+1}_i(l,m)\alpha_{\pm N}(l+1) \;,\\
C^{0,i}_{i,\pm}(l,m) &=c^{0}_i(l,m)\beta_{\pm N}(l) \;,\\
C^{-1,i}_{i,\pm}(l,m) &=c^{-1}_i(l,m)\alpha_{\pm N}(l) \;,
\end{align*}
where $c^{j}_i(l,m)$ are independent of $\pm N$ and bounded by a constant.
Thus, concerning (\ref{eq:ineqqn}) we just need to prove that
\begin{align*}
|\alpha_N(l)-\alpha_{-N}(l)|\lesssim q^l \;,\\
|\beta_N(l)-\beta_{-N}(l)|\lesssim q^l \;,
\end{align*}
which follow from (\ref{eq:ineqAa}, \ref{eq:ineqB}) and by
triangle inequality:
\begin{align*}
|\alpha_N(l)-\alpha_{-N}(l)|\leq
\big|\alpha_N(l)-s\sqrt{1+q^2}\big|+\big|\alpha_{-N}(l)-s\sqrt{1+q^2}\big|
\lesssim q^l \;,\\
|\beta_N(l)-\beta_{-N}(l)|\leq
\big|\beta_N(l)-(1-s^2)\big|+\big|\beta_{-N}(l)-(1-s^2)\big|
\lesssim q^l \;.
\end{align*}
If $s=0$ and for any $j\neq 0$, $C^{j,i}_{i,\pm}(l,m)$ is
the product of uniformly bounded coefficients $c^{j}_i(l,m)$
times either $\alpha_N(l)$ or $\alpha_N(l+1)$, so
by (\ref{eq:ineqAb}) it follows (\ref{eq:Kprime}).
\end{prova}

We are in the situation described in Section \ref{sec:2.8r}.
By previous Lemma, for all $s\in[0,1]$ we have Lemma \ref{lemma:5.1}
which we rewrite as a proposition,

\begin{prop}\label{prop:FredS2qs}
The data $(\A(S^2_{qs}),\HH,F,\gamma)$ is a even $1$-summable
Fredholm module.
\end{prop}

\noindent
and we have also Lemma \ref{lemma:5.2}, which claims that

\begin{prop}
If $\lambda_{n+1}-\lambda_n$ is a bounded sequence, the data
$(\A(S^2_{qs}),\HH,D_\lambda,\gamma)$ is a regular even spectral triple.
\end{prop}

\noindent
For $s=0$, condition (\ref{eq:Kprime}) is satisfied and previous result
can be improved by using Lemma \ref{lemma:5.4} as follows.

\begin{prop}
If $\{q^n\lambda_n\}$ is a bounded sequence, the data $(\A(S^2_{q0}),\HH,D_\lambda,\gamma)$
is a even Lipschitz-regular spectral triple.
The triple is regular if $\lambda_n$ is at most of polynomial growth.
\end{prop}

\begin{lemma}
The Fredholm module in Proposition \ref{prop:FredS2qs} is non-trivial.
In particular $$\mr{Index}(F_{e_s}^+)=2N\;,$$ with $e_s$ the idempotent
in (\ref{eq:qBott}).
\end{lemma}
\begin{prova}
Recall that $N\neq 0$.
We use (\ref{eq:ind}) with $m=0$ and find
\begin{align*}
\mr{Index}(F_{e_s}^+) &=\frac{1}{2}\,\tr_{\HH\otimes\C^2}(\gamma F[F,e_s])
  =\frac{1-q^2}{2(1+s^2)}\tr_{\HH}(\gamma F[F,A]) \\
 &=\frac{-1}{1+s^2}\,\frac{1-q^2}{1+q^2}\,\sum_{l,m}\frac{[2l]-q^{l+m+1}[2][l-m]}{[2l]}
   \big\{\beta_{-N}(l)-\beta_{N}(l)\big\} \\
 &=(q^{-2}-1)[2N]\,\sum_{l,m}\frac{[2l]-q^{l+m+1}[2][l-m]}{[2l][2l+2]} \;.
\end{align*}
Computing the sum in $m=-l,...,l$ and setting $n=l-|N|+1$, we obtain
$$
\mr{Index}(F_{e_s}^+)=2\,\mr{sign}(N)\,(1-x^{2|N|})\sum_{n=1}^\infty
\left(
\frac{n+|N|-1}{1-x^{2(n+|N|-1)}}x^{n-1}-
\frac{n+|N|}{1-x^{2(n+|N|)}}x^n
\right) \;,
$$
where $x=q^2\in\,]0,1[$. Last sum is a telescoping series, and the sum is
$$
\mr{Index}(F_{e_s}^+)=2\,\mr{sign}(N)\,(1-x^{2|N|})
\frac{n+|N|-1}{1-x^{2(n+|N|-1)}}x^{n-1}\Big|_{n=1}=2N\;.
$$
The non-triviality of the Fredholm module comes from 
the non-vanishing of the index of $F_{e_s}^+$.
\end{prova}

If $D_\lambda=D$ is the Dirac operator in (\ref{eq:isoD}),
the multiplicity of the $n$th eigenvalue $n+1$ is 
$4n+4|N|+2$ and so by Lemma \ref{lemma:5.3},

\begin{lemma}
The metric dimension of the linear Dirac operator $D$ is $2$.
\end{lemma}

For this Dirac operator we now compute the dimension spectrum and the
top residue. To this aim, we introduce
an auxiliary representation which is obtained from simpler operators.
Let $\hat{\HH}\simeq\ell^2(\Z\times\N)$ be the Hilbert space
having orthonormal basis $\kkett{n,k}$, with $n\in\frac{1}{2}\Z$ and $n+k\in\N$.
Consider the bounded operators $\alpha,\beta\in\B(\hat{\HH})$ defined by
$$
\alpha\kkett{n,k}=\sqrt{1-q^{2(n+k+1)}}\kkett{n+\oh,k+\oh} \;,\qquad
\beta\kkett{n,k}=-q^{n+k}\kkett{n+\oh,k-\oh}\;.
$$
These operators satisfy the defining relations of $\A(SU_q(2))$, as
can easily be seen. The morphism $\A(S^2_{qs})\hookrightarrow\A(SU_q(2))$
in (\ref{eq:SUtoS2}) leads to a $*$-representation
$\pi:\A(S^2_{qs})\to\B(\hat{\HH})$ given by
\begin{subequations}\label{eq:aprep}
\begin{align}
\pi(A)\kkett{n,k} =\; &
-s\,q^{n+k-1}\sqrt{1-q^{2(n+k)}}\kkett{n-1,k} \notag \\ &
+(1-s^2)q^{2(n+k)}\kkett{n,k} \notag \\ &
-s\,q^{n+k}\sqrt{1-q^{2(n+k+1)}}\kkett{n+1,k} \;, \\
\pi(B)\kkett{n,k} =\; &
-s\,q^{2(n+k)+1}\kkett{n-1,k+1} \notag \\ &
-(1-s^2)q^{n+k+1}\sqrt{1-q^{2(n+k+1)}}\kkett{n,k+1} \notag \\ &
+s\sqrt{1-q^{2(n+k+1)}}\sqrt{1-q^{2(n+k+2)}}\,\kkett{n+1,k+1} \;.
\end{align}
\end{subequations}
Our Hilbert space of spinors $\HH$ is identified with a subspace
of $\hat{\HH}\oplus\hat{\HH}$. That is, for each $N$ there is a
partial isometry $W_N:\HH_N\to\hat{\HH}$, given by
$$
W_N\ket{l,m;N}:=\kkett{l-|N|,m+|N|} \;,
$$
and with left inverse $W_N^*:\hat{\HH}\to\HH_N$, given by
$$
W_N^*\kkett{n,k}=\begin{cases}
\ket{n+|N|,k-|N|;N} &\mr{if}\;n\in\N\;\mr{and}\;k\leq n+2|N|\;, \\
0 &\mr{otherwise}\;.
\end{cases}
$$
We denote
$$
Q:=W_{-N}\oplus W_N\;,\qquad
P:=W_{-N}^*\oplus W_N^*\;.
$$
Since $W^*_NW_N=id_{\HH_N}$,
the Hilbert space $\HH$ is isomorphic to its image in
$\hat{\HH}\otimes\C^2\simeq\hat{\HH}\oplus\hat{\HH}$
via the map $Q$.
The Dirac operator $D$ in (\ref{eq:isoD}) is the ``restriction'' of
the self-adjoint operators $\hat{D}=\hat{F}|\hat{D}|$ on $\hat{\HH}\otimes\C^2$
determined by
$$
|\hat{D}|\kkett{n,k}=(n+1)\kkett{n,k}\;,\qquad
\hat{F}=\Big(\textrm{\footnotesize$\begin{array}{cc}
0 & 1 \\ 1 & 0 \end{array}$}\Big)\;,
$$
in the sense that $Q$ intertwines $D$ and $\hat{D}$
(i.e.~$\hat{D}Q=QD$ on $\HH$).
The same holds for $|\hat{D}|$ and $\hat{F}$.
For $a\in\A(S^2_{qs})$, the diagonal lift $\pi(a)\otimes id_{\C^2}$
of $\pi(a)$ to $\hat{\HH}\otimes\C^2$ will be denoted simply by $\pi(a)$.

\begin{lemma}\label{lem:ap}
The operator $a-P\pi(a)Q$ belongs to $\op(D)$
for all $a\in\A(S^2_{qs})$.
\end{lemma}
\begin{prova}
For $a\in\A(S^2_{qs})$, define two bounded operators $R_a:\HH\to\HH$
and $T_a:\hat{\HH}\to\HH$ by
$$
R_a:=a-P\pi(a)Q\;,\qquad
T_a:=P\pi(a)(1-QP) \;,
$$
and notice that if $a$ is a generator of $\A(S^2_{qs})$,
$T_a$ has finite rank. On the other hand, if $T_a$ and $T_b$ have
finite rank for some $a,b\in\A(S^2_{qs})$, then $T_{ab}$ has
finite rank too, due to the identity
$$
T_{ab}=T_a\,\pi(b)(1-QP)+P\pi(a)Q\,T_b \;.
$$
By the above considerations, we conclude that $T_a$ has
finite rank for all $a\in\A(S^2_{qs})$. Therefore,
$\op$ being an ideal in $\opz$, using the identity
$$
R_{ab}=aR_b+R_aP\pi(b)Q-T_a\pi(b)Q\;,
$$
we conclude that $R_{ab}$ is a smoothing operators whenever $R_a,R_b\in\op$.
It is then enough to prove this lemma for the generators
of the algebra. Using (\ref{eq:ineqAB}) and the inequality\linebreak
$0\leq 1-\sqrt{1-u}\leq u$, valid for $u\in [0,1]$, one proves that
$R_A$, $R_B$ and $R_{B^*}$ have coefficients bounded by $q^l$,
hence they are rapid decay matrices. This concludes the proof.
\end{prova}

The next proposition is the main result of this section.

\begin{prop}
The dimension spectrum is $\Sigma=\{1,2\}$.
\end{prop}
\begin{prova}
Let $\Psi^0$ be the algebra generated by
$\A(S^2_{qs})$, by $[D, a]$ for all $a\in\A(S^2_{qs})$ and by iterated applications
of the derivation $\delta$ (cf.~Section \ref{sec:2.2}). 
Let $\mf{C}$ be the $*$-algebra (of bounded operators on $\hat{\HH}$)
generated by $\alpha,\beta,\alpha^*,\beta^*$ and $F$. 
By Lemma~\ref{lem:ap}, $\A(S^2_{qt})\subset P\mf{C} Q+\op$. 
Note that 
\begin{align*}
[\hat{F},\alpha]&=0, & [\hat{D},\alpha]&=\tfrac{1}{2}\,\alpha \,\hat{F}, &
[|\hat{D}|,\alpha]&=\tfrac{1}{2}\,\alpha, \\
[\hat{F},\beta]&=0, &
[\hat{D},\beta ]&=\tfrac{1}{2}\,\beta  \,\hat{F}, &
[|\hat{D}|,\alpha]&=\tfrac{1}{2}\,\beta.
\end{align*}
Thus $P\mf{C}Q$ is invariant under application of $\delta$ and $[D,(\cdot)]$ and
hence $\Psi^0\subset P\mf{C} Q+\op$.

We shall compute the singularities of zeta functions associated to the monomials
$$
S:=P\alpha^i\beta^j(\beta^*)^kQ F\quad\mr{and}\quad
T:=P\alpha^i\beta^j(\beta^*)^kQ\;,
$$
where $i\in\Z$ and $j,k\in\N$ and we use the notation $\alpha^i:=(\alpha^*)^{|i|}$ for $i<0$.
From the commutation relations of $\alpha$ and $\beta$, it is clear that these monomials
span $P\mf{C}Q$.

Firstly, note that the zeta-function associated with a bounded off-diagonal operator is
identically zero in the half-plane $\mr{Re}\,z>2$
and so is its holomorphic extension to the entire complex plane.
This is the case for the monomials $S$ due to the presence of $F$. The other monomials 
$T$
shift the index $l$ by $(i+j-k)/2$ and the index $m$ by $(i-j+k)/2$ and 
therefore are also off-diagonal operators unless these shifts are zero, 
which happens when $i=0$ and $j=k$.
Hence only monomials $T=P\beta^k(\beta^*)^kQ=P(\beta\beta^*)^kQ$ contribute
to the dimension spectrum.

For $k=0$, $T=1$ and the corresponding zeta function is
$$
\zeta_{1}(z)=\sum_{n=1}^\infty (4n+4|N|-2)n^{-z}=4\zeta(z-1)+(4|N|-2)\zeta(z) \;,
$$
where $\zeta(z)$ is the Riemann zeta function, meromorphic in $\C$ with
a simple pole at $1$ and with residue $1$. 
Since $1\in\Psi^0$, this shows that $2\in\Sigma$. For $|N|\neq 1/2$, this
proves also that $1\in\Sigma$.

When $k>0$, $T= P(\beta\beta^*)^kQ$. So 
$T\ket{l,m}_\pm=q^{2k(l+m)}\ket{l,m}_\pm$
and the associated $\zeta$-function is
\begin{align*}
\zeta_T(z) &=\sum_{l\in\N+|N|}^\infty 2(l-|N|+1)^{-z}\sum_{l+m=0}^{2l}(q^{2k})^{l+m} 
=\sum_{l\in\N+|N|}^\infty 2(l-|N|+1)^{-z}\frac{1-q^{2k(2l+1)}}{1-q^{2k}} \\ &
=\tfrac{2}{1-q^{2k}}\,\zeta(z)+\text{holomorphic function}\;.
\end{align*}
Therefore $\Sigma$ may contain, besides 2, at most the additional point 1. 
We still have to check if $1\in\Sigma$ for $|N|=1/2$.
For this, we take $\;A^2=P\pi(A)^2Q\mod\op(D)\,$.
Then, using (\ref{eq:aprep}), we get (modulo holomorphic functions)
$$
\zeta_{A^2}(z)\sim 2\sum_{l+\frac{1}{2}=1}^\infty(l+\tfrac{1}{2})^{-z}
\sum_{l+m=0}^{2l}\Big\{(1-[2]^2s^2+s^4)q^{4(l+m)}+[2]s^2q^{2(l+m)-1}\Big\}
\sim 2\frac{1+s^4}{1-q^4}\,\zeta(z) \;,
$$
and
$\,\mr{Res}_{z=1}\zeta_{A^2}(z)=2\frac{1+s^4}{1-q^4}\neq 0$ 
for all $s\in[0,1]$. Thus $\Sigma=\{1,2\}$ for any
value of $N$.
\end{prova}

Let $\A(S^1)$ denote the polynomial $*$-algebra generated by the function $u$,
given by $u(\theta)=e^{i\theta}$. A $*$-algebra morphism $\sigma_s:\A(S^2_{qs})\to\A(S^1)$
is defined by
\begin{equation}\label{eq:mor}
\sigma_s(A)=0\;,\qquad\sigma_s(B)=su\;.
\end{equation}

\begin{prop}
The top residue of $\zeta_a(z)$, $a\in\A(S^2_{qs})$, is
\begin{equation}\label{eq:rem}
\nint a|D|^{-2}=\frac{2}{\pi}\int_{S^1}\sigma_s(a)(\theta)\de \theta \;.
\end{equation}
For $s=0$, $\sigma_0(a)\in\C$ is constant and equation (\ref{eq:rem})
simplifies to $\,4\,\sigma_0(a)$.
\end{prop}

\begin{prova}
It is sufficient to prove (\ref{eq:rem}) for the basis elements
$A^jB^k$, $j\in\N$ and $k\in\Z$, and then extend it to $\A(S^2_{qs})$
by linearity. We use the notation $B^k:=(B^*)^{|k|}$ if $k<0$.

Since $\sigma_s(A^jB^k)=\delta_{j0}(su)^k$, the right hand side of
equation (\ref{eq:rem}) is zero unless $j=k=0$. We next show that the
left hand side of (\ref{eq:rem}) also vanishes unless $j=k=0$.
When $j=k=0$, the relation $\zeta_{1}(z)=4\zeta(z-1)$ fixes the
normalization constant.

Now $A^jB^k$ is off-diagonal if $k\neq 0$ since it shifts the index $m$ by $k$.
It remains to prove that $\zeta_{A^j}(z)=
\zeta_{P\pi(A)^jQ}(z)\,+$
{\it holomorphic function} has no singularity in $z=2$.
For $j\neq 0$, $P\pi(A)^jQ$ satisfies the inequality
$$
\big|\inner{l,m,\pm N|P\pi(A)^jQ|l,m,\pm N}\big|\leq c_jq^{l+m}\;,
$$
for some positive constants $c_j$'s.
From this inequality, we deduce that $\zeta_{A^j}(z)$
is a convergent series for all $z$ with $\mr{Re}\,z>1$. In particular,
it is finite for $z=2$.
\end{prova}

Last point to discuss is the existence of a real structure $J$.
We can assume $N>0$, as the choice $N<0$ is equivalent to changing the sign
of the grading operator.
The operator $T$ defined by $T\ket{l,m;\pm N}=(-q)^{m+N}\ket{l,-m;\mp N}$
is equivariant --- one can repeat verbatim the proof of (\ref{eq:simi}).
The antiunitary part $J$ of $T$ is
\begin{equation}\label{eq:JS2qs}
J\ket{l,m;\pm N}=(-1)^{m+N}\ket{l,-m;\mp N} \;.
\end{equation}

\begin{prop}
The antilinear isometry $J$ defined by (\ref{eq:JS2qs}) is the
antiunitary part of an antilinear operator, it satisfies $J^2=(-1)^{2N}$,
$J\gamma=-\gamma J$, $JD=DJ$ and
$$
[a,Jb^*J^{-1}]\in\op\;,\qquad [[D,a],Jb^*J^{-1}]\in\op\;,
$$
for all $a,b\in\A(S^2_{qs})$.
It is then a (`weak') equivariant real structure with $K\!O$-dimension
$2$~\textup{mod}~$8$ (resp.~$6$~\textup{mod}~$8$) if $\,2N$ is odd (resp.~even).
\end{prop}

\begin{prova}
Since $[D,a]-\delta(a)F=|D|[F,a]$ and $(JF)b(JF)^{-1}-JbJ^{-1}=J[F,b]FJ^{-1}$
are smoothing operators for all $a,b$, we can replace $[D,\,.\,]$ by
$\delta$ and $J$ by $J_0:=JF$ in the verification. Therefore, we want to
show
$$
[a,J_0b^*J_0^{-1}]\in\op\;,\qquad [\delta(a),J_0b^*J_0^{-1}]\in\op\;.
$$
Now, let $\hat{J}$ be the antilinear operator on $\hat{\HH}$ given by
$$
\hat{J}\!\kkett{n,k}=i^{2k}\!\kkett{n,-k+2N} \;,
$$
and note that $P\hat{J}Q=J_0$. We can replace $a$ by $P\pi(a)Q$ since
the difference is a smoothing operator, thus Lemma is proved if we
prove that
$$
P[\pi(a),\hat{J}\pi(b^*)\hat{J}^{-1}]Q\in\op\;,\qquad
P[[|\hat{D}|,\pi(a))],\hat{J}\pi(b^*)\hat{J}^{-1}]Q\in\op\;.
$$
By Leibniz rule it is sufficient to prove the previous equations
taking instead of $\pi(a)$ and $\pi(b^*)$ two of the generators
$\{\alpha,\beta,\alpha^*,\beta^*\}$ of the enlarged algebra.
These are all eigenvectors for $[|\hat{D}|,\,.\,]$, thus
the first order condition follows from the commutant condition.
From
$$
[a,\hat{J}b\hat{J}^{-1}]=-[a^*,\hat{J}b^*\hat{J}^{-1}]=
\hat{J}[b,\hat{J}a\hat{J}^{-1}]\hat{J}^{-1}\;,
$$
once the commutant condition is proved for a fixed couple of
elements $(a,b)$, it is proved also for $(a^*,b^*)$ and $(b,a)$.
Therefore it is sufficient to
consider the cases $(a,b)=(\alpha,\alpha), (\alpha^*,\alpha),
(\beta,\alpha)$, $(\beta^*,\alpha), (\beta,\beta)$ and $(\beta^*,\beta)$.
By direct computation
$$
[\alpha,\hat{J}\alpha\hat{J}^{-1}]  =0 \;,\qquad
[\alpha^*,\hat{J}\alpha\hat{J}^{-1}]=0 \;,\qquad
[\beta,\hat{J}\beta\hat{J}^{-1}]    =0 \;,
$$
and
\begin{align*}
[\beta,\hat{J}\alpha\hat{J}^{-1}]\kkett{n,k}
   &=i\left(\sqrt{q^{2k}-q^{4(n+N+1)}}-
     \sqrt{q^{2k}-q^{4(n+N+1)-2}}\right)\!\kkett{n+1,k-1} \;,\\
[\beta^*,\hat{J}\alpha\hat{J}^{-1}]\kkett{n,k}
   &=i\left(\sqrt{q^{2k}-q^{4(n+N)}}-
     \sqrt{q^{2k}-q^{4(n+N+1)-2}}\right)\!\kkett{n+1,k-1} \;,\\
[\beta^*,\hat{J}\beta\hat{J}^{-1}]  &=i(q-q^{-1})q^{2(n+N)}\kkett{n,k+1} \;.
\end{align*}
Coefficients in the last equations are bounded by $q^n$,
hence projected to $\HH$ are bounded by $q^l$ and
this concludes the proof.
\end{prova}

\section{Spectral triples for pre-$C^*$-algebras}
Things become more complicated when one considers an algebra which is not
made of polynomials in a finite number of generators. We discuss an example
in this section, which is a generalization of the construction in Appendix B of \cite{DD06}.

We call $C(S^2_{qs})$ the universal $C^*$-algebra associated with $\A(S^2_{qs})$.
The relations defining $\A(S^2_{qs})$ can be rewritten as
\begin{equation*}
BB^*+\big(A-\tfrac{1-s^2}{2}\big)^2=\big(\tfrac{1+s^2}{2}\big)^2\;,\qquad
B^*\!B+\big(q^2A-\tfrac{1-s^2}{2}\big)^2=\big(\tfrac{1+s^2}{2}\big)^2\;.
\end{equation*}
The first equation implies the following inequality for the operator
norm in an arbitrary bounded $*$-representation of the algebra:
\begin{equation*}
||A-\tfrac{1-s^2}{2}||\leq \tfrac{1+s^2}{2} \;.
\end{equation*}
Hence, by triangle inequality
\begin{align*}
||q^2A-\tfrac{1-s^2}{2}|| &=q^2||A-\tfrac{1-s^2}{2}-(q^{-2}-1)\tfrac{1-s^2}{2}||
\leq q^2||A-\tfrac{1-s^2}{2}||+(1-q^2)\tfrac{1-s^2}{2} \\
&\rule{0pt}{14pt} \leq \tfrac{1+(2q^2-1)s^2}{2}<\tfrac{1+s^2}{2}\;.
\end{align*}
This implies (cf.~Lemma 2.2.4 of \cite{Lan98}) that the element
\begin{equation*}
B^*\! B=\big(\tfrac{1+s^2}{2}\big)^2-\big(q^2A-\tfrac{1-s^2}{2}\big)^2
\end{equation*}
is invertible in the $C^*$-algebra. We can then construct a projection $P\in C(S^2_{qs})$,
\begin{equation}
P=1-B(B^*\!B)^{-1}\! B^*\;,
\end{equation}
and the corresponding representation of the algebra on the module $C(S^2_{qs})P$.
For that we work first with $\A(S^2_{qs})P$, which is an $\A(S^2_{qs})$-module.

First of all, notice that if $[A,C]=0$ and $C^{-1}$ exists, then
$[A,C^{-1}]=-C^{-1}[A,C]C^{-1}=0$. This proves that $(B^*B)^{-1}$
commutes with $A$, and so in particular $AP=PA$.

From $B^*P=0$, $(A+s^2)(1-A)=BB^*$ and $AP=PA$ we deduce that
the elements:
\begin{equation}
P_+=\tfrac{1}{1+s^2}(A+s^2)P\;,\qquad
P_-=\tfrac{1}{1+s^2}(1-A)P\;,
\end{equation}
are projections, are mutually orthogonal ($P_+P_-=P_-P_+=0$) and their
sum is $P$. Then, $\A(S^2_{qs})P$ is $\Z_2$-graded and
decomposes as $\A(S^2_{qs})P_+\oplus\A(S^2_{qs})P_-$.

The (left) kernel of $P_+$ is the left ideal generated by $B^*$ (as $B^*P=0$)
and $1-A$. The (left) kernel of $P_-$ is the left ideal generated by $B^*$
and $A+s^2$. A linear basis of both modules is then $\{B^nP_\pm\}$.
For $s\neq 0$ we can normalize the vectors as follows.
We take as linear basis of $\A(S^2_{qs})P_{\pm}$ the vectors
\begin{equation}\label{eq:thevect}
\ket{n}_+=c^+_n\,B^n P_+\;,\qquad \ket{n}_-=c^-_n\,B^nP_-\;,
\qquad\;n\in\N\;,
\end{equation}
where the normalization constants are defined by $c_0^{\pm}:=1$ and
\begin{align*}
c^+_n &:=\left\{\prod\nolimits_{k=1}^n(1-q^{2k})(s^2+q^{2k})\right\}^{-1/2} \;, \\
c^-_n &:=s^{-n}\left\{\prod\nolimits_{k=1}^n(1-q^{2k})(1+q^{2k}s^2)\right\}^{-1/2} \;,
\end{align*}
for all $n\geq 1$. The basis (\ref{eq:thevect}) is orthonormal for the inner product
induced by the Haar functional, and we have a $\Z_2$-graded $*$-representation
of $C(S^2_{qs})$ on the Hilbert space completion, naturally isomorphic to
$\ell^2(\N)\oplus\ell^2(\N)$, defined by
\begin{subequations}\label{eq:murep}
\begin{align}\label{eq:murepA}
A\ket{n}_+ &=q^{2n}\ket{n}_+ &
B\ket{n}_+ &=\sqrt{(1-q^{2\smash[t]{(n+1)}})(s^2+q^{2\smash[t]{(n+1)}})}\ket{n+1}_+ \\
A\ket{n}_- &=-s^2q^{2n}\ket{n}_- & \;
B\ket{n}_- &=s\sqrt{(1-q^{2\smash[t]{(n+1)}})(1+q^{2\smash[t]{(n+1)}}s^2)}\ket{n+1}_-
\end{align}
\end{subequations}
with $B^*$ the adjoint of $B$. These are the well known unique
irreducible infinite--dimensional representations of $S^2_{qs}$.
They are inequivalent and not faithful (as follows from the observation
that $A$ has positive spectrum in the former representation and negative
in the latter), while their direct sum is known to be faithful \cite{She91}.
They are also equivariant with respect to the Hopf subalgebra of $U_q(su(2))$
generated by the element $K$; being $U(1)$ the Lie group associated to this
commutative and cocommutative Hopf algebra, these representations are usually
referred to as $U(1)$-equivariant~\cite{CP03b}.

\begin{figure}[h]\begin{center}

\includegraphics[width=14cm]{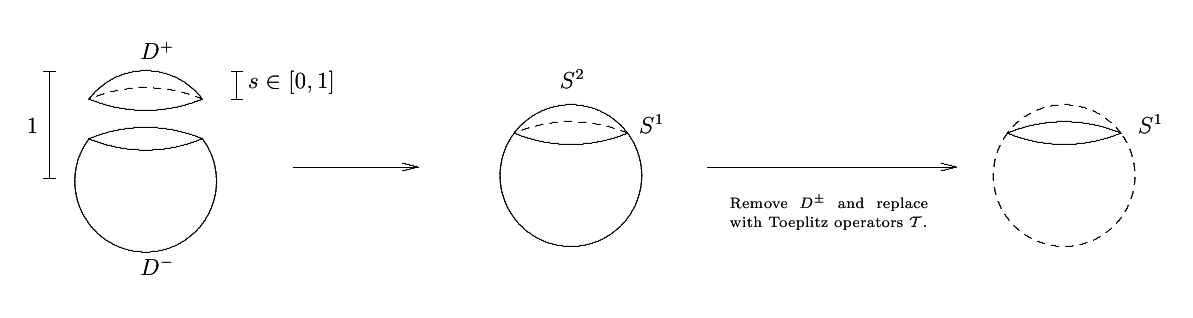}

\vspace{-15pt}

\caption{Geometric visualization of $S^2_{qs}$}\label{tab}

\end{center}\end{figure}

The representations (\ref{eq:murep}) were used in~\cite{She91} to prove that
the sequence
\begin{equation}\label{eq:seqA}
0\to\mc{K}\oplus\mc{K}\to C(S^2_{qs})\stackrel{\sigma_s}{\to}C(S^1)\to 0
\end{equation}
is exact (for $s\neq 0$). This is the noncommutative analogue of the sequence
\begin{equation*}
0\to C_0(\textrm{open disk})\oplus C_0(\textrm{open disk})\to C(S^2)\to C(S^1)\to 0\;,
\end{equation*}
and allows to interpret $C(S^2_{qs})$ as the union of two
noncommutative open disks $\mc{K}$ glued along a commutative circle,
as shown pictorially in Figure \ref{tab}.
Here $\sigma_s:C(S^2_{qs})\to C(S^1)$ is the completion of the $*$-algebra
morphism (\ref{eq:mor}).

When $s=0$, the constants $c_n^+$ are well defined and the `$+$'-representation
is still valid, while the `$-$'-representation becomes trivial. We'll consider
only the `$+$'-representation in this case.
Such a representation is faithful and irreducible, and has been used in~\cite{She91}
to prove that the following sequence is exact,
\begin{equation}\label{eq:seqB}
0\to\mc{K}\to C(S^2_{q0})\stackrel{\sigma_0}{\to}\C\to 0\;,
\end{equation}
with $\sigma_0:C(S^2_{q0})\to\C$ the $C^*$-algebra morphism
defined by $\sigma_0(A)=\sigma_0(B)=0$ (and of course $\sigma_0(1)=1$).
Equation (\ref{eq:seqB}) is the noncommutative analogue of the sequence
\begin{equation*}
0\to C_0(\textrm{open disk})\to C(S^2)\to\C\simeq C(\{*\})\to 0\;,
\end{equation*}
and allows us to interpret $C(S^2_{q0})$ as the one-point compactification
of a noncommutative disk.

Let us define the operators $N$ and $F$ as
\begin{equation*}
N\ket{n}_{\pm}=(n+1)\ket{n}_{\pm}\;,\qquad
F\ket{n}_{\pm}=\ket{n}_{\mp}\;,\qquad\forall\;n\in\N\;.
\end{equation*}
For $s\neq 0$, $(\A(S^2_{qs}),\ell^2(\N)\oplus\ell^2(\N), N\cdot F)$
is the $U(1)\times U(1)$-equivariant spectral triple
studied in~\cite{CP03b} (we shifted by $+1$ the spectrum of $N$
to make it invertible), while $(\A(S^2_{q0}),\ell^2(\N),N)$ is the
analogue one for the standard Podle\'s sphere, and has a positive Dirac
operator.

Since trace class operators do not contribute to the noncommutative integral,
one immediately proves that for $s\neq 0$
\begin{equation*}
\nint aN^{-1}=\frac{1}{\pi}\int_{S^1}\sigma_s(a)\de\theta\qquad\forall\;a\in\A(S^2_{qs})\;,
\end{equation*}
and for $s=0$
\begin{equation*}
\nint aN^{-1}=\sigma_0(a)\qquad\forall\;a\in\A(S^2_{q0})\;.
\end{equation*}
When considering only polynomials $\A(S^2_{qs})$, one can prove that
the triple is regular and the dimension spectrum is $\{1\}$.
If $s\neq 0$, this result is still valid if one considers the pre-$C^*$-algebra
$C^\infty(S^2_{qs})$ entering in the exact sequence
\begin{equation}\label{eq:preCa}
0\to\mc{S}\oplus\mc{S}\to C^\infty(S^2_{qs})\stackrel{\sigma_s}{\to}C^\infty(S^1)\to 0\;.
\end{equation}
where $\mc{S}$ is the pre-$C^*$-algebra of rapid decay matrices on $\ell^2(\N)$.

For $s=0$ the sequence defining the algebra of `smooth function' on the sphere is:
\begin{equation}\label{eq:preCb}
0\to\mc{S}\to C^\infty(S^2_{q0})\stackrel{\sigma_0}{\to}\C\to 0\;.
\end{equation}
In the remaining of the section we prove that the spectral triples above can
be extended to regular spectral triples over $C^\infty(S^2_{qs})$, and compute
their dimension spectrum.

\subsection{Spectral triples over the quantum disk}\label{sec:qDisk}
Let $\ket{n}$ be the canonical orthonormal basis of $\ell^2(\N)$, and
$N\ket{n}:=(n+1)\ket{n}$. The unilateral shift $w$, defined by
$w\ket{n}=\ket{n+1}$, generates the Toeplitz $C^*$-algebra $\mc{T}$,
that can be viewed as algebra of `continuous functions' on a quantum
version of the closed disk \cite{Con04}.

We define $\A^\infty\subset\mc{T}$ as the linear span of the elements:
\begin{equation*}
f=\smash[b]{\sum_{n\in\N}}(f_nw^n+f_{-n-1}(w^*)^{n+1})+\smash[b]{\sum_{j,k\in\N}}f_{jk}w^j
(1-ww^*)(w^*)^k\qquad\,\;\{f_n\}\in\mc{S}(\Z),\;\{f_{jk}\}\in\mc{S}\;,
\end{equation*}
\rule{0pt}{16pt}where $\mc{S}(\Z)$ indicates rapid decay sequences and $\mc{S}$
rapid decay matrices on $\ell^2(\N)$.

By direct calculation one checks that $\A^\infty$ is a $*$-algebra and the map
$\sigma:\A^\infty\to C^\infty(S^1)$, given by $f\mapsto\sigma(f)$, $\sigma(f)(\theta):=
\sum_{n\in\Z}f_ne^{in\theta}$, is a surjective $*$-algebra
morphism (this follows from the simple observation that, via Fourier series,
$\mc{S}(\Z)$ equipped with convolution product is isomorphic to $C^\infty(S^1)$).
Furthermore, $\ker\sigma$ is a two-sided $*$-ideal in $\A^\infty$ isomorphic
to $\mc{S}$. At the level of abstract algebras, this follows from the equality
\begin{equation*}
w^j(1-ww^*)(w^*)^kw^{j'}(1-ww^*)(w^*)^{k'}=\delta_{j'k}w^j(1-ww^*)(w^*)^{k'}
\end{equation*}
and becomes more evident in matrix form, as
\begin{equation*}
\inner{j|f|k}=f_{j,k}\;,\qquad\quad
\forall\;f=\sum\nolimits_{j,k\in\N}f_{jk}w^j(1-ww^*)(w^*)^k\in\ker\sigma\;.
\end{equation*}
Then, we have the short exact sequence
\begin{equation}\label{eq:seq}
0\to\mc{S}\to\A^\infty\stackrel{\sigma}{\to}C^\infty(S^1)\to 0\;.
\end{equation}
Notice that the action of $U(1)$ on polynomials in $\{w,w^*\}$, given by
$w\mapsto e^{i\theta}w$, defines a one-parameter group of automorphisms
implemented on $\ell^2(\N)$ by the unitary operators $e^{i\theta N}$
(i.e.~$x\mapsto e^{i\theta N}xe^{-i\theta N}$). Being implemented by
unitary operators, it extends to a (strongly continuous) action
of $U(1)$ on the $C^*$-algebra $\mc{T}$, and then on $\A^\infty$.

We conclude with the statement,
\begin{prop}
$\A^\infty$ is a Fr\'echet space and a pre-$C^*$-algebra.
\end{prop}
\begin{prova}
The direct sum of two Fr\'echet spaces is a Fr\'echet space too,
and the observation that $\A^\infty$ as a vector space is the direct
sum of $C^\infty(S^1)$ and $\mc{S}$ (the sequence (\ref{eq:seq}) splits) proves the first claim.
To prove that $\A^\infty$ is a pre-$C^*$-algebra it is sufficient to show
that it contains the inverse of each element $f\in\A^\infty$,
whenever $f$ is invertible in $\mc{T}$ (cf.~\cite{Bos90}). 
To reach this conclusion, one can easily adapt the proof
of Proposition 1 in~\cite{Con04} to the present case.
\end{prova}

We identify $\A^\infty\subset\mc{T}$ with `smooth functions' over the quantum disk.

\begin{prop}\label{prop:6.8}
The data $(\A^\infty,\ell^2(\N),N)$ is a regular spectral triple,
the dimension spectrum is $\{1\}$ and the top residue is
$$
\nint fN^{-1}=\frac{1}{2\pi}\int_{S^1}\sigma(f)\de\theta\;,
$$
for all $f\in\A^\infty$, where $\sigma$ is the map in (\ref{eq:seq}).
\end{prop}

\begin{prova}
The operators $w^k,(w^*)^k$ are eigenvectors of the derivation
$\delta=[N,\,.\,]\,$: $\,\delta(w^k)=kw^k$ and $\delta((w^*)^k)=-k(w^*)^k$.
If $\{f_n\}\in\mc{S}(\Z)$, $\{f_{jk}\}\in\mc{S}$ then
$\{nf_n\}\in\mc{S}(\Z)$ and $\{(j-k)f_{jk}\}\in\mc{S}$ too.
This proves that $\delta(f)\in\A^\infty$ for all $f\in\A^\infty$.
Being $\A^\infty$ bounded $\delta$-invariant operators, $\A^\infty
=\A^\infty\cup [N,\A^\infty]\subset\opz$ and the triple is regular.

The algebra $\Psi^0$ of order $\leq 0$ pseudo-differential operators
coincides with $\A^\infty$. By Weierstrass theorem $\tr_{\,\ell^2(\N)}(fN^{-z})$ is
a holomorphic function of $z\in\C$ for all $f\in\mc{S}=\ker\sigma$,
thus $\ker\sigma$ does not contribute to the dimension spectrum,
and we have to consider only elements of the form
$$
f=\sum\nolimits_{n\in\N}(f_nw^n+f_{-n-1}(w^*)^{n+1}) \;.
$$
All the terms in $f$ are off-diagonal, except the one proportional to
$f_0=\frac{1}{2\pi}\int_{S^1}\sigma(f)\de\theta$. The function
$\tr_{\ell^2(\N)}(fN^{-z})=f_0\sum_{n\in\N}(n+1)^{-z}$ has meromorphic
extension $f_0\,\zeta(z)$, where $\zeta(z)$ is the Riemann zeta-function.
This concludes the proof.
\end{prova}

\subsection{The $U(1)\times U(1)$-equivariant spectral triple over $S^2_{qs}$, $s\neq 0$}
In this section we consider $s\neq 0$, and introduce the symbols $\mu:=\mu_+\oplus\mu_-$ for
the representation of $\A(S^2_{qs})$ over $\ell^2(\N)\oplus\ell^2(\N)$ defined by
(\ref{eq:murep}). In terms of the operators $q^N$ and $w$ one has
\begin{align*}
\mu_+(A)&=q^{2N-2}\;, &
\mu_+(B)&=w\sqrt{(1-q^{2N+2})(s^2+q^{2N+2})} \;, \\
\mu_-(A)&=-s^2q^{2N-2}\;, &
\mu_-(B)&=sw\sqrt{(1-q^{2N+2})(1+q^{2N+2}s^2)} \;.
\end{align*}
We have 
$\,q^N\in\mc{S}\subset\A^\infty\,$, with $\A^\infty$ the algebra of the previous section.
Since $\mr{Sp}(1-q^{2N+2})\subset\R^+$ and the square root is holomorphic in an
angular sector containing $\R^+$, by stability under holomorphic functional calculus
we deduce that $\mu_{\pm}(A)$ and $\mu_{\pm}(B)$ belong to $\A^\infty$;
in particular, $\mu_{\pm}(A)\in\mc{S}$ and $\,\mu_{\pm}(B)=sw\!\mod\mc{S}$.
If we identify $\A(S^2_{qs})$ with its representation $\mu$, and call
\begin{equation*}
C^\infty(S^2_{qs}):=\{(x_+,x_-)\in\A^\infty\oplus\A^\infty\,|\,\sigma(x_+)=\sigma(x_-)\}\;,
\end{equation*}
then $\A(S^2_{qs})\subset C^\infty(S^2_{qs})\subset C(S^2_{qs})$ for all $s\neq 0$.
That is, `smooth functions' over $S^2_{qs}$ contain `polynomial' ones
and are contained in `continuous' ones.

From the direct sum of two copies of (\ref{eq:seq}), passing to the diagonal
($\sigma(x_+)=\sigma(x_-)$) we obtain the exact sequence (\ref{eq:preCa}).
The morphism $\sigma(x_+)=\sigma(x_-)\equiv\sigma_s(x)$ is the one defined by equation
(\ref{eq:mor}), for all $x=x_+\oplus x_-\in C(S^2_{qs})$ (the equality holds for the
generators $A,B$ and then for all the $C^*$-algebra,
all members of the equation being $C^*$-algebra morphisms).

The definition of $C^\infty(S^2_{qs})$ does not depend on $s$, so the
quantum spheres (with $s\neq 0$) are all isomorphic already at the smooth level.

Since the smooth algebra do not depend on $s$, we can apply
the analysis made in~\cite{DD06} for the equatorial Podle\'s sphere to prove that:

\begin{prop}
$(C^\infty(S^2_{qs}),\ell^2(\N)\oplus\ell^2(\N),N\cdot F)$ is a \emph{regular}
even spectral triple. The dimension spectrum is $\{1\}$. The noncommutative integral is
$$
\nint a\,N^{-1}=\frac{1}{\pi}\int_{S^1}\sigma_s(a)(\theta)\de\theta \;,
$$
for all $a\in C^\infty(S^2_{qs})$ and with $\sigma_s$ the $*$-algebra morphism
in equation (\ref{eq:mor}).
\end{prop}

\begin{prova}
Concerning the first claim, the non-trivial part is to prove regularity.
Let $\mf{B}$ be the algebra generated by $C^\infty(S^2_{qs})$ and by commutators
with $D=NF$, and let $\Psi^0$ the algebra generated by $\mf{B}$ and iterated
applications of $\delta=[N,\,.\,]$ to $\mf{B}$.

Notice that for any $x=x_+\oplus x_-\in C^\infty(S^2_{qs})$,
$x_+-x_-\in\ker\sigma=\mc{S}$. Since we have
$$
[N\cdot F,x]=\delta(x)F+N(x_+-x_-)F\gamma\;,
$$
with $\gamma=1\oplus -1$ the grading,
and $\,N\cdot\mc{S}\subset\mc{S}$, then
$$
\mf{B}\subset \mf{C}:=C^\infty(S^2_{qs})+C^\infty(S^2_{qs})\cdot F+\mc{S}\cdot F\gamma\;.
$$
Recall that $\delta(\A^\infty)\subset\A^\infty$: this means that
$\mf{C}$ is $\delta$-invariant and then belongs to the smooth domain $\opz$.
This proves regularity.
In particular, the algebra $\Psi^0$ is still contained in $\mf{C}$.

Any $T\in\Psi^0\subset\mf{C}$ can be written as $T=x_+\oplus x_-$ plus off-diagonal
terms, with $x_\pm\in\A^\infty$.
Therefore, for $z\in\C$ with sufficiently large real part,
the `zeta-type' function associated to $T$ is given by
$$
\zeta_T(z)=\tr_{\ell^2(\N)}\big((x_+ +x_-)N^{-z}\big) \;.
$$
By Proposition \ref{prop:6.8}, $\zeta_T(z)$ has meromorphic extension on
$\C\smallsetminus\{1\}$ and the residue in $z=1$ is
$$
\nint T\,N^{-1}=\frac{1}{2\pi}\int_{S^1}\{\sigma(x_+)(\theta)+\sigma(x_-)(\theta)\}
\de\theta \;.
$$
If $T=x_+\oplus x_-$ belongs to the algebra $C^\infty(S^2_{qs})$,
then $\sigma(x_+)+\sigma(x_-)=2\sigma_s(T)$, and this concludes the proof.
\end{prova}

\begin{prop}
$(C^\infty(S^2_{qs}),\ell^2(\N)\oplus\ell^2(\N),F)$ is a $1$-summable
non-trivial Fredholm module.
\end{prop}

\begin{prova}
For all $x=x_+\oplus x_-\in C^\infty(S^2_{qs})$, $[F,x]=(x_+-x_-)F\gamma$
is a rapid decay matrix, so in particular it is of trace class.
Non-triviality follows from the computation
$$
\mr{Index}(F_{e_s}^+)=\tfrac{1}{2}\tr(\gamma F[F,e_s])=
\tfrac{1-q^2}{1+s^2}\tr(\gamma\,\mu(A))
=(1-q^2)\sum\nolimits_{n\in\N}q^{2n}=+1\;,
$$
with $e_s$ the idempotent in (\ref{eq:qBott}).
The proof of the non-triviality of this Fredholm module already appear,
for $s\neq 0,1$, in~\cite{CP03b} at the $C^*$-algebra level and
using the projection $\ket{0}\!\left<0\right|\oplus 0$.
\end{prova}

\subsection{The $U(1)$-equivariant spectral triple over $S^2_{q0}$}
In this section we consider $s=0$ and call $\mu$ the representation of
$\A(S^2_{q0})$ over $\ell^2(\N)$ defined by (\ref{eq:murep}).
In terms of the operators $q^N$ and $w$ one has
\begin{equation*}
\mu(A)=q^{2N-2}\;,\qquad
\mu(B)=q^N\sqrt{1-q^{2N}}\,w\;.
\end{equation*}
Thus $\mu(A),\mu(B)$ belongs to $\A^\infty$, and in particular
they are in $\ker\sigma=\mc{S}$, so $\A(S^2_{q0})$ is a sub-algebra
of the unitization $\mc{S}^+$ of $\mc{S}$,
$$
\A(S^2_{q0})\subset C^\infty(S^2_{q0}):=\mc{S}^+\subset C(S^2_{q0})\;.
$$

By restricting Proposition \ref{prop:6.8} to the subalgebra $\mc{S}^+$
of $\A^\infty$ we obtain the following,

\begin{prop}
$(C^\infty(S^2_{q0}),\ell^2(\N),N)$ is a regular spectral triple,
the dimension spectrum is $\{1\}$ and the noncommutative integral is
$$
\nint aN^{-1}=\sigma_0(a)
$$
for all $f\in\A^\infty$, where $\sigma_0$ is the (restriction of
the) map in (\ref{eq:seqB}).
\end{prop}

Since the Dirac operator $N$ is positive, in the $s=0$
case one has no non-trivial Fredholm module associated with the triple.
The spectral triples discussed in this sections appear very different
for $s\neq 0$ and for $s=0$.

\section{The homological dimension}
Let $\kappa(a)=K^2\az a$ be the modular automorphism for $\A(S^2_{qs})$.
We now give a simple proof that $H\!H^\kappa_2(\A(S^2_{qs}))\neq 0$,
i.e.~that twisting the Hochschild homology `cures' the dimension drop.
The twisted Hochschild homology groups of $S^2_{qs}$ were computed in
\cite{Had07}.

Let $\gamma_s:=2e_s-1$ with $e_s$ the idempotent in (\ref{eq:qBott}), and
call
$$
\eta:=\sigma(K^2)=
\Big(\textrm{\footnotesize$\begin{array}{cc}
q & 0 \\ 0 & q^{-1} \end{array}$}\Big) \;,
$$
with $\sigma$ the representation in (\ref{eq:qPaul}).
The tensor
$$
\omega_2 :=\tfrac{1}{8}\tr_{\C^2}\big\{\eta\gamma_s\dot{\otimes}(1+\gamma_s)\dot{\otimes}(1-\gamma_s)\big\}
         =\tfrac{1}{8}\sum_{ijkl}\eta_{ij}(\gamma_s)_{jk}\otimes(1+\gamma_s)_{kl}
           \otimes(1-\gamma_s)_{li}
$$
is invariant.
Indeed, by (\ref{eq:cov}),
\begin{align*}
h\az\omega_2 &=\tfrac{1}{8}\tr_{\C^2}(\eta\sigma(h_{(1)})^t
\gamma_s\dot{\otimes}(1+\gamma_s)\dot{\otimes}(1-\gamma_s)\sigma(S^{-1}(h_{(2)}))^t) \\
&=\tfrac{1}{8}\tr_{\C^4}(\sigma(S^{-1}(h_{(2)}))^t\eta\sigma(h_{(1)})^t
\gamma_s\dot{\otimes}(1+\gamma_s)\dot{\otimes}(1-\gamma_s))\;,
\end{align*}
where in last step we used cyclicity of the trace.
But $\eta=\sigma(K^2)=\sigma(K^2)^t$ and so
$$
\sigma(S^{-1}(h_{(2)}))^t\eta\sigma(h_{(1)})^t=
\sigma\big( h_{(1)}K^2S^{-1}(h_{(2)}) \big)=
\sigma\big( h_{(1)}S(h_{(2)})K^2 \big)=\epsilon(h)\eta\;,
$$
where we used (\ref{eq:Squad}).
This proves the invariance of $\omega_2$, i.e.~$h\az\omega_2=\epsilon(h)\omega_2$.

By using $\kappa(\gamma_s)\eta=\eta\gamma_s$ and $\gamma_s^2=1$
we compute
\begin{align*}
b_\kappa\omega_2 &=\tfrac{1}{4}\,\tr_{\C^2}(\eta\gamma_s\otimes 1-1\otimes\eta\gamma_s) \\
&=\tfrac{1}{4}\bigl\{\tr_{\C^2}(\eta\gamma_s)\otimes 1-1\otimes\tr_{\C^2}(\eta\gamma_s)\bigr\}\;.
\end{align*}
with $b_\kappa$ the boundary operator defined in Section \ref{sec:2.8}.
Since $\tr_{\C^2}(\eta\gamma_s)=(q^{-1}-q)\frac{1-s^2}{1+s^2}$ is a scalar,
we get $b_\kappa\omega_2=0$ and so:

\begin{prop}
The tensor $\omega_2$ is an invariant twisted Hochschild $2$-cycle.
\end{prop}

For $q=1$, the image of $\omega_2$ under the isomorphism (\ref{eq:HKR})
is the volume form on $S^2$. This justify the interpretation of $\omega_2$
as `twisted orientation'.

We now apply Lemma \ref{lemma:ch} to construct twisted cocycles.
Consider the Fredholm module $(\A(S^2_{qs}),\HH,F,\gamma)$ of
Section \ref{sec:6.6}.
From the crossed product relation, the automorphism $\kappa$ is implemented
on $\HH$ by $K^{-2}$.

\begin{lemma}
The multilinear map
$$
\mr{ch}^{F,\kappa}_n(a_0,\ldots,a_n):=\tfrac{1}{2}\,\tr_{\HH}(K^{-2}\gamma F[F,a_0][F,a_1]\ldots[F,a_n])
$$
is a $U_q(su(2))$-invariant twisted Hochschild $n$-cocycle on $\A(S^2_{qs})$,
for all even $n\geq 2$.
\end{lemma}
\begin{prova}
Invariance follows from (\ref{eq:Squad}), that is $S^{-1}(h)K^{-2}=K^{-2}S(h)$,
and from equivariance of the representation of $\A(S^2_{qs})$ on $\HH$.
By Lemma \ref{lemma:ch}, to prove that
$\mr{ch}^{F,\kappa}_n$ is a cocycle we need just to check that it is
well defined; that is, the operators
$$
K^{-2}[F,a_0][F,a_1]\ldots[F,a_n]
$$
must be of trace class for all $n\geq 2$ and for all $a_0,\ldots,a_n\in\A(S^2_{qs})$.
Since trace-class operators are a two-sided ideal in $\B(\HH)$, it is sufficient
to do the computation for $n=2$. By Lemma \ref{lemma:6.8} the matrix coefficients
of $[F,a_j]$ are bounded by $q^l$, thus coefficients of
$K^{-2}[F,a_0][F,a_1][F,a_2]$ are bounded by $q^{3l-2m}\leq q^l$ and the
operator is of trace class.
\end{prova}

That $H\!H^\kappa_2(\A(S^2_{qs}))$ is not zero is a corollary of the
following three lemmas.

\begin{lemma}\label{lemma:Qindex}
Let $\HH'=\HH'_+\oplus\HH'_-$ be a (separable) Hilbert space and
$P:\HH'_+\to\HH'_-$ be a Fredholm operator with quasi-inverse $P^*$.
Let $Q$ be an even selfadjoint operator with dense domain in $\HH'$
and discrete spectrum.
Let $n\in\N$ such that $Q(1-PP^*)^n$ and $Q(1-P^*P)^n$ are of trace
class, and assume $[P,Q]=0$ on $\HH'_+$. Then
$$
Q\textrm{-}\mr{Index}(P)=\tr_{\HH'_+}Q(1-P^*P)^n-\tr_{\HH'_-}Q(1-PP^*)^n \;,
$$
where the `$Q$-index' is defined by
$$
Q\textrm{-}\mr{Index}(P):=\tr_{\ker P}(Q)-\tr_{\mr{coker}\,P}(Q) \;.
$$
\end{lemma}

\begin{prova}
Since $Q$ commutes with $P^*P$ (resp.~with $PP^*$), we can find an orthonormal
basis of $\HH'_+$ (resp.~of $\HH'_-$) made of joint eigenvectors of
$Q$ and $P^*P$ (resp.~$Q$ and $PP^*$). Let
\begin{align*}
E_0(\lambda,\mu)&=\big\{v\in\HH'_+\,\big|\,P^*Pv=\lambda v,\;Qv=\mu v\big\}\;, \\
E_1(\lambda,\mu)&=\big\{v\in\HH'_-\,\big|\,PP^*v=\lambda v,\;Qv=\mu v\big\}\;.
\end{align*}
Since $P^*P$ (resp.~$PP^*$) is (the identity plus) a positive compact operator,
it can be diagonalized and its eigenspaces are all finite dimensional.
Since $P(P^*P)=(PP^*)P$ and $Q$ commutes with $P$, the operator $P$ maps
$E_0(\lambda,\mu)$ into $E_1(\lambda,\mu)$. Similarly $P^*:E_1(\lambda,\mu)
\mapsto E_0(\lambda,\mu)$. For $\lambda\neq 0$ we have the linear isomorphisms
\begin{align*}
\lambda&=P^*P|_{E_0(\lambda,\mu)}:E_0(\lambda,\mu)\mapsto E_1(\lambda,\mu)\mapsto E_0(\lambda,\mu)\;, \\
\lambda&=PP^*|_{E_1(\lambda,\mu)}:E_1(\lambda,\mu)\mapsto E_0(\lambda,\mu)\mapsto E_1(\lambda,\mu)\;,
\end{align*}
thus $P$ and $P^*$ are both bijective linear maps and $\dim E_0(\lambda,\mu)=\dim E_1(\lambda,\mu)$
for all $\lambda\neq 0$. Using this equality we compute
\begin{align*}
\tr_{\HH'_+}Q(1-P^*P)^n&-\tr_{\HH'_-}Q(1-PP^*)^n=\sum_{\lambda,\mu}\mu(1-\lambda)^n
             \bigl\{\dim E_0(\lambda,\mu)-\dim E_1(\lambda,\mu)\bigr\} \\
             &=\sum_{\mu}\mu\bigl\{\dim E_0(0,\mu)-\dim E_1(0,\mu)\bigr\}
             =\tr_{\ker P}(Q)-\tr_{\mr{coker}\,P}(Q) \;.
\end{align*}
This concludes the proof.
\end{prova}

\begin{lemma}\label{lemma:6.26}
The pairing $\mr{ch}^{F,\kappa}_2(\omega_2)$ is given by the `$q$-index'
of the operator $F^+_{e_s}$ in (\ref{eq:Fpe}),
$$
\mr{ch}^{F,\kappa}_2(\omega_2)=q\textrm{-}\mr{Index}(F^+_{e_s}):=
\tr_{\ker F^+_{e_s}}(K^{-2}\eta)-\tr_{\mr{coker}\,F^+_{e_s}}(K^{-2}\eta) \;.
$$
\end{lemma}
\begin{prova}
By definition
$$
\mr{ch}^{F,\kappa}_2(\omega_2)=-\tfrac{1}{2}\,\tr_{\HH\otimes\C^2}(K^{-2}\eta\gamma F[F,e_s]^3) \;.
$$
Since $\gamma F$ and $K^{-2}\eta$ both commute with $[F,e_s]$, with
some algebraic manipulation and using $[F,e_s]=e_s[F,e_s]+[F,e_s]e_s$
we get
\begin{align*}
\mr{ch}^{F,\kappa}_2(\omega_2) &=\tr_{\HH\otimes\C^2}(K^{-2}\eta\gamma(e_s-e_sFe_sFe_s)^2) \\
&=\tr_{\HH'_+}(K^{-2}\eta(1-F^-_{e_s}F^+_{e_s})^2)
-\tr_{\HH'_-}(K^{-2}\eta(1-F^+_{e_s}F^-_{e_s})^2)  \;.
\end{align*}
where $\HH'_\pm=(1\pm\gamma)e_s(\HH\otimes\C^2)$,
$F^+_{e_s}:\HH'_+\to\HH'_-$ is given by (\ref{eq:Fpe})
and $F^-_{e_s}:=(F^+_{e_s})^*$. Now we apply Lemma \ref{lemma:Qindex}
to $Q=K^{-2}\eta$, $P=F^+_{e_s}$ and $\HH'=\HH'_+\oplus\HH'_-$,
and obtain the desired result (since $Q=K^{-2}\eta$ is given in
diagonal form in the basis of harmonic spinors, it is immediate
to extend it to a self-adjoint operator).
\end{prova}

\begin{lemma}\label{lemma:6.27}
We have
$$
q\textrm{-}\mr{Index}(F^+_{e_s})=[2N]\;.
$$
In particular, the class of $\omega_2$ in $H\!H^\kappa_2(\A(S^2_{qs}))$
is not zero.
\end{lemma}

\begin{prova}
The idempotent $e_s$ satisfies (\ref{eq:cov}) for all $h\in U_q(su(2))$,
with $\sigma$ the spin $1/2$ representation in (\ref{eq:qPaul}). As a
consequence, if we define the action of $h\in U_q(su(2))$ on $w\in\HH\otimes\C^2$
as in (\ref{eq:RmodAz}),
$$
h\aaz w:=\sigma(S^{-1}(h_{(1)}))^t\,h_{(2)}\az w\;,
$$
then
$$
h\aaz (e_sw)=e_s(h\aaz w) \;.
$$
Using the Clebsch-Gordan rules (\ref{eq:copClebsch}), with
$\ket{+}=q^{-1/2}\binom{0}{1}$ and $\ket{-}=-q^{1/2}\binom{1}{0}$,
we identify the following orthonormal basis of $\HH\otimes\C^2$:
\begin{align*}
v^{\pm,\uparrow}_{l,m} &:=
\left(\begin{array}{r}
-\sqrt{\tfrac{q^{l+m+1}[l-m]}{[2l]}}\ket{l-\oh,m+\oh;\mp N} \\
\sqrt{\tfrac{q^{-(l-m+1)}[l+m]}{[2l]}}\ket{l-\oh,m-\oh;\mp N}
\end{array}\right)
 \;, & l=|N|+\tfrac{1}{2},|N|+\tfrac{3}{2},\ldots \;, \\
v^{\pm,\downarrow}_{l,m} &:=
\left(\begin{array}{r}
\sqrt{\tfrac{q^{-(l-m)}[l+m+1]}{[2l+2]}}\ket{l+\oh,m+\oh;\mp N} \\
\sqrt{\tfrac{q^{l+m}[l-m+1]}{[2l+2]}}\ket{l+\oh,m-\oh;\mp N}
\end{array}\right) \;,
& l=|N|-\tfrac{1}{2},|N|+\tfrac{1}{2},\ldots \;,
\end{align*}
which transform under the `black' action $\aaz$ as in (\ref{eq:rep}) (with $\omega=1$).
We set $v^{\pm,\uparrow}_{|N|-\frac{1}{2},m}:=0$ and start counting from $l=|N|-\frac{1}{2}$
for both $v^\uparrow$ and $v^\downarrow$.
Left multiplication of these vectors by $e_s$ gives
\begin{subequations}\label{eq:esv}
\begin{align}
e_sv^{\pm,\uparrow}_{l,m}
&=a^\pm_{l,s}v^{\pm,\uparrow}_{l,m}+b^\pm_{l,s}v^{\pm,\downarrow}_{l,m} \;,\\
e_sv^{\pm,\downarrow}_{l,m}
&=b^\pm_{l,s}v^{\pm,\uparrow}_{l,m}+(1-a^\pm_{l,s})v^{\pm,\downarrow}_{l,m} \;,
\end{align}
\end{subequations}
where for all $l\neq|N|-\tfrac{1}{2}$,
$$
a^\pm_{l,s}=\frac{1+q^2s^2+q^2\beta_{\mp N}(l-\tfrac{1}{2})}{q[2](1+s^2)} \;,\qquad
b^\pm_{l,s}=\frac{1}{1+s^2}\sqrt{\frac{[2l+2]}{[2][2l+1]}}\,\alpha_{\mp N}(l+\tfrac{1}{2}) \;,
$$
with $\alpha_N(l)$ and $\beta_N(l)$ the coefficients in (\ref{eq:abNl}),
while for $l=|N|-\tfrac{1}{2}$,
$$
a^\pm_{|N|-\frac{1}{2},s}=\delta_{\,\mr{sign}N,\mp 1} \;,\qquad
b^\pm_{|N|-\frac{1}{2},s}=0 \;.
$$
The map $N\to -N$ changes sign to the $q$-index (it is
equivalent to changing the sign of the grading $\gamma$), thus
it is sufficient to consider the case \fbox{$N>0$}\ .

Since $\alpha_N(l)=0$ if and only if $l=|N|$, we have
$b^\pm_{l,s}\neq 0$ for all $l\neq |N|-\frac{1}{2}$
and the sum $(a^\pm_{l,s})^2+(b^\pm_{l,s})^2$
is invertible. Therefore, the following orthonormal
basis of $\HH\otimes\C^2$ is well defined:
\begin{align*}
w^{\pm,||}_{l,m} &:=
\frac{1}{(a^\pm_{l,s})^2+(b^\pm_{l,s})^2}
\big(a^\pm_{l,s}\,v^{\pm,\uparrow}_{l,m}
+b^\pm_{l,s}\,v^{\pm,\downarrow}_{l,m}\big)\;,
&&\forall\;l\geq|N|+\tfrac{1}{2}\;, \\
w^{\pm,\perp}_{l,m} &:=
\frac{1}{(a^\pm_{l,s})^2+(b^\pm_{l,s})^2}
\big(b^\pm_{l,s}\,v^{\pm,\uparrow}_{l,m}
-a^\pm_{l,s}\,v^{\pm,\downarrow}_{l,m}\big) \;,
&&\forall\;l\geq|N|+\tfrac{1}{2}\;, \\
v^{||}_m &:=v^{+,\downarrow}_{|N|-\frac{1}{2},m} \;, \\
v^\perp_m &:=v^{-,\downarrow}_{|N|-\frac{1}{2},m} \;.
\end{align*}
From (\ref{eq:esv}) it follows that this basis is made of eigenvectors of $e_s$:
$$
e_sw^{\pm,||}_{l,m}=w^{\pm,||}_{l,m}\;,\qquad
e_sv^{||}_m=v^{||}_m\;,\qquad
e_sw^{\pm,\perp}_{l,m}=e_sv^{\perp}_m=0\;.
$$
By definition, the kernel of $F^+_{e_s}$ is the span of the vectors
$\{w^{+,||}_{l,m},v^{||}_m\}$ which are annihilated by $e_sF$,
while the cokernel is the span of the vectors
$\{w^{-,||}_{l,m}\}$ which are annihilated by $e_sF$.

Since $F$ sends $v^\pm$ into $v^\mp$, we have:
$$
e_sFw^{\pm,||}_{l,m}=
\frac{(a^\mp_{l,s})^2+(b^\mp_{l,s})^2}{(a^\pm_{l,s})^2+(b^\pm_{l,s})^2}
(a^+_{l,s}a^-_{l,s}+b^+_{l,s}b^-_{l,s})
w^{\mp,||}_{l,m} \;.
$$
The fraction in last equation is always different from zero.
If $a^+_{l,s}a^-_{l,s}+b^+_{l,s}b^-_{l,s}=0$, then 
$w^{+,||}_{l,m}$ is in the kernel of $F^+_{e_s}$, $w^{-,||}_{l,m}$ is
in the cokernel and by
$$
\big<w^{+,||}_{l,m},K^{-2}\eta\,w^{+,||}_{l,m}\big>=
\big<w^{-,||}_{l,m},K^{-2}\eta\,w^{-,||}_{l,m}\big>=q^{-2m}
$$
their contributions to $q$-Index$(F^+_{e_s})$ cancel each other.
Since
$$
e_sFv^{||}_m=e_sv^\perp_m=0 \;,
$$
the vectors $v^{||}_m=v^{+,\downarrow}_{|N|-\frac{1}{2},m}$
are in the kernel of $F^+_{e_s}$ and we have
$$
q\textrm{-}\mr{Index}(F^+_{e_s})=
\sum_{m=-N+1/2}^{N-1/2}\big<
v^{+,\downarrow}_{|N|-\frac{1}{2},m},K^{-2}\eta\,v^{+,\downarrow}_{|N|-\frac{1}{2},m}
\big>=\sum_{m=-N+1/2}^{N-1/2}q^{-2m}=[2N]\;.
$$
This concludes the proof.
\end{prova}

As a corollary, since $H\!H^\kappa_2(\A(S^2_{qs}))$ has (at least)
one non-zero element, it is not zero and so there is no dimension
drop.

A result similar to Lemma \ref{lemma:6.27} was obtained by \cite{Wag07}
for the standard sphere: he computed the $q$-index of a fixed $F$
(coming from the Fredholm module with $N=1/2$) twisted with the
generic idempotent with charge $n\in\Z$; here the idempotent
is fixed (is the one with charge $1$) while the Fredholm module
is generic and depends on the half integer $N$.

We conclude with some remarks. We have constructed a non-trivial
element $[\omega_2]\in H\!H^\kappa_2$, with $\kappa$
the modular automorphism, given on generators by $\kappa(A)=A$,
$\kappa(B)=q^2B$, $\kappa(B^*)=q^{-2}B^*$. From Theorem 4.6 of
\cite{Had07} we know that $H\!H^\sigma_2\simeq\C$ if $\sigma$
is a positive power of the modular automorphism $\sigma_{\mr{mod}}$,
and $H\!H^\sigma_2\simeq 0$ for all other automorphisms. From
equation (5) of \cite{Had07} we deduce that his generators are related
to ours by a replacement $B\leftrightarrow B^*$, and from the definition
of $\sigma_{\mr{mod}}$ \cite[eq.~(17)]{Had07} we deduce that $\kappa$
and $\sigma_{\mr{mod}}$ are indeed the same automorphism. Thus, we can
conclude that $H\!H^\kappa_2\simeq\C$ and that $[\omega_2]$ is a
generator.

Dually, we proved that $\mr{ch}^{F,\kappa}_2$ has non-trivial
class in $H\!H_\kappa^2$. We want to compare this result
with \cite{SW04}. In that paper the authors define a twisted
$2$-cocycle for the standard Podle\'s sphere,
with twist $\sigma$ given by the \emph{inverse} of the modular
automorphism \cite[equation (29)]{SW04}:
$\sigma(a)=K^{-2}\az a=\kappa^{-1}(a)=\sigma_{\mr{mod}}^{-1}(a)$,
for all $a\in\A(S^2_{q0})$. 
It is then consistent that their cocycle, as explained in Section 6
of \cite{Had07}, has a non-trivial class in cyclic cohomology but is
trivial in Hochschild cohomology.


\typeout{Capitolo 4}

\chapter{The $4$-dimensional quantum orthogonal sphere}\label{chap:S4q}

\noindent
In this chapter we construct an equivariant regular
even spectral triple on $S^4_q$ following \cite{DDL06}. 
This spectral triple is a deformation of the canonical one associated
to the round sphere $S^4$, and the Dirac operator has undeformed spectrum.
We discuss analytical properties, construct a weak real structure,
and prove that the noncommutative integral is the integral over the space of
characters of the algebra. We give also an explicit expression of
the Haar state in a basis of monomials in $x_i$'s, and construct
a twisted Hochschild $4$-cycle which may be interpreted as the
`orientation'.

We repeat here Definition \ref{def:S4q} of the coordinate algebra
$\A(S^4_q)$ of the quantum orthogonal $4$-sphere. This is the $*$-algebra generated
by five elements $x_0=x_0^*,x_i,x_i^*$ (with $i=1,2$), subject to the
relations:
\begin{align*}
&x_ix_j     = q^2x_jx_i   \;, \qquad\quad\forall\;\;0\leq i<j\leq 2 \;, \\
&x_i^*x_j   =q^2x_jx_i^*  \;, \qquad\quad\forall\;\;i\neq j \;, \\
&[x_1^*,x_1]=(1-q^4)x_0^2 \;, \\
&[x_2^*,x_2]=x_1^*x_1-q^4x_1x_1^* \;, \\
&x_0^2+x_1x_1^*+x_2x_2^*=1 \;.
\end{align*}
Recall that $0<q<1$.

\section{The symmetry Hopf algebra $U_q(so(5))$}\label{sec:7.1}
We call $U_q(so(5))$ the real form of the Drinfeld-Jimbo deformation
of $so_{\C}(5)$, corresponding to the Euclidean signature $(+,+,+,+,+)$;
it is a real form of the Hopf algebra called $\breve{U}_q(so(5))$
in~\cite[Section 6.1.2]{KS97}. As a $*$-algebra, 
$U_q(so(5))$ is generated by the elements $K_i=K_i^*,K_i^{-1},E_i,F_i:=E_i^*$,
$i=1,2$ ($i\to 3-i$ with respect to the notations of~\cite{KS97}), with relations:
\begin{equation*}\begin{array}{c}
[K_1,K_2]=0\;,\qquad
K_iK_i^{-1}=K_i^{-1}K_i=1\;,
\\ \rule{0pt}{20pt}
[E_i,F_j]=\delta_{ij}\frac{K_j^2-K_j^{-2}}{q^j-q^{-j}}\;,
\\ \rule{0pt}{20pt}
K_iE_iK_i^{-1}=q^iE_i\;,\qquad
K_iE_jK_i^{-1}=q^{-1}E_j\;\;\mr{if}\;i\neq j\;,
\end{array}\end{equation*}
and with the ones obtained by conjugation and Serre relations.
These are explicitly given by
\begin{subequations}\label{eq:Serre}
\begin{align}
E_1E_2^2-(q^2+q^{-2})E_2E_1E_2+E_2^2E_1 &=0\;, \\*
E_1^3E_2-(q^2+1+q^{-2})(E_1^2E_2E_1-E_1E_2E_1^2)-E_2E_1^3 &=0\;,
\end{align}
\end{subequations}
together with their adjoints. Serre relations can be written in a more compact
form by defining $[a,b]_q:=q^2ab-ba$. Then, (\ref{eq:Serre}) are equivalent to
$$
[E_2,[E_1,E_2]_q]_q=0\;,\qquad
[E_1,[E_1,[E_2,E_1]_q]_q]=0\;.
$$
The Hopf algebra structure $(\Delta,\epsilon,S)$ of $U_q(so(5))$ is given by:
\begin{equation*}
\begin{array}{c}
\Delta K_i=K_i\otimes K_i\;,\quad \Delta E_i=E_i\otimes K_i+K_i^{-1}\otimes E_i\;, 
\\ \rule{0pt}{3.5ex}\epsilon(K_i)=1\;,\quad \epsilon(E_i)=0\;,\\ \rule{0pt}{3.5ex}
S(K_i)=K_i^{-1}\;,\quad S(E_i)=-q^iE_i\;.
\end{array}
\end{equation*}
For any non-negative $n_1,n_2$ such that $n_2\in\frac{1}{2}\N$ and $n_2-n_1\in\N$ there
is an irreducible representation of $U_q(so(5))$ whose representation space we denote
by $V_{(n_1,n_2)}$~\cite{KS97}.
We call it ``the representation with highest weight $(n_1,n_2)$'', since
the highest weight vector is an eigenvector of $K_1$ and $K_1K_2$
with eigenvalues $q^{n_1}$ and $q^{n_2}$, respectively.

Irreducible representations with highest weight $(0,l)$ and $(\frac{1}{2},l)$
(the ones that we need explicitly) can be found in~\cite{Cha94b} and are recalled 
presently. 
We use the shorthand notation $V_l:=V_{(0,l)}$ if $l\in\N$ and $V_l:=V_{(\frac{1}{2},l)}$
if $l\in\N+\frac{1}{2}$.
The vector space $V_l$, for all $l\in\frac{1}{2}\N$, has orthonormal basis $\ket{l,m_1,m_2;j}$,
with the labels $(m_1,m_2,j)$ satisfying the following constraints.
For $l\in\N$:
$$
j=0,1,\ldots,l\;,\qquad
j-|m_1|\in\N\;,\qquad
l-j-|m_2|\in 2\N\;,
$$
while for $l\in\N+\frac{1}{2}$:
$$
j=\tfrac{1}{2},\tfrac{3}{2},\ldots,l-1,l\;,\qquad
j-|m_1|\in\N\;,\qquad
l+\tfrac{1}{2}-j-|m_2|\in\N\;.
$$
Notice that for any admissible $(l,m_1,m_2,j)$ there exists a unique $\epsilon\in\{0,\pm\frac{1}{2}\}$
such that \mbox{$l+\epsilon-j-m_2\in 2\N$} (that is, $\epsilon=0$ if $l\in\N$ and
$\epsilon=\frac{1}{2}(-1)^{l+\frac{1}{2}-j-m_2}$ if $l\in\N+\frac{1}{2}$). We shall need the
coefficients,
\begin{subequations}
\begin{align}
a_l(j,m_2) &=\frac{1}{[2]}\sqrt{\frac{[l-j-m_2+\epsilon][l+j+m_2+3+\epsilon]}
  {[2(j+|\epsilon|)+1][2(j-|\epsilon|)+3]}} \;, \\
b_l(j,m_2) &=2|\epsilon|\,\frac{\sqrt{[l-\epsilon(2j+1)-m_2+1][l-\epsilon(2j+1)+m_2+2]}}{[2j][2j+2]}
 \;,\label{eq:abc} \\
c_l(j,m_2) &=\frac{(-1)^{2\epsilon}}{[2]}\sqrt{\frac{[l-j+m_2+2-\epsilon][l+j-m_2+1-\epsilon]}
  {[2(j+|\epsilon|)-1][2(j-|\epsilon|)+1]}} \;,
\end{align}
\end{subequations}
where, as usual, $[z]:=(q^z-q^{-z})/(q-q^{-1})$ denotes the $q$-analogue of
a number $z$.

The $*$-representation $\sigma_l:U_q(so(5))\to\mr{End}(V_l)$ is defined by the
assignments,
\begin{align*}
\sigma_l(K_1)\ket{l,m_1,m_2;j} &=q^{m_1}\ket{l,m_1,m_2;j} \;, \\
\sigma_l(K_2)\ket{l,m_1,m_2;j} &=q^{m_2-m_1}\ket{l,m_1,m_2;j} \;, \\
\sigma_l(E_1)\ket{l,m_1,m_2;j} &=\sqrt{[j-m_1][j+m_1+1]}\ket{l,m_1+1,m_2;j} \;, \\
\sigma_l(E_2)\ket{l,m_1,m_2;j} &=
\sqrt{[j-m_1+1][j-m_1+2]}\,a_l(j,m_2)\ket{l,m_1-1,m_2+1;j+1} \\ &+
\sqrt{[j+m_1][j-m_1+1]}\,b_l(j,m_2)\ket{l,m_1-1,m_2+1;j} \\ &+
\sqrt{[j+m_1][j+m_1-1]}\,c_l(j,m_2)\ket{l,m_1-1,m_2+1;j-1} \;.
\end{align*}
The representation symbol $\sigma_l$ will be suppressed, when there is no risk of ambiguity.

For $l\in\N$ the representation $\sigma_l$ is \emph{real}. That is, there is an
antilinear map $C:V_l\to V_l$, which satisfies $C^2=1$ and $C\sigma_l(h)C=\sigma_l(S(h)^*)$.
It is explicitly given by 
\begin{equation}\label{eq:C}
C\ket{l,m_1,m_2;j}:=(-q)^{m_1}q^{3m_2}\ket{l,-m_1,-m_2;j}  \;.
\end{equation}

The operator
\begin{equation}\label{eq:Cdef}
\mc{C}_1:=q^{-1}K_1^2+qK_1^{-2}+(q-q^{-1})^2E_1F_1\;,
\end{equation}
is a Casimir for the subalgebra generated by $(K_1,K_1^{-1},E_1,F_1)$.
For future reference, we note the action of $\mc{C}_1$ on
a vector of $V_l$, with $l\in\frac{1}{2}\N$; it is 
\begin{equation}\label{eq:Caz}
\mc{C}_1\ket{l,m_1,m_2;j}=(q^{2j+1}+q^{-2j-1})\ket{l,m_1,m_2;j} \;.
\end{equation}

\section{General properties of the quantum $4$-sphere}\label{sec:7.2}

In the next Propositions we summarize some well known facts~\cite{KS97}.

\begin{prop}
The algebra $\A(S^4_q)$ is an $U_q(so(5))$-module $*$-algebra for the action
given by
\begin{gather*}
K_1\az x_1=qx_1 \;, \quad
K_2\az x_2=qx_2 \;, \quad
K_2\az x_1=q^{-1}x_1\;, \\
E_1\az x_0=q^{-1/2}x_1\;,\qquad
E_2\az x_1=x_2\;, \quad \\
F_1\az x_1=q^{1/2}[2]x_0 \;,\quad
F_1\az x_0=-q^{-3/2}x_1^*\;\;\quad
F_2\az x_2=x_1\;,
\end{gather*}
while the action is trivial in all other cases:
\begin{gather*}
K_1\az x_0=x_0 \;,\qquad K_2\az x_0=x_0 \;,\qquad K_1\az x_2=x_2 \;, \\
E_2\az x_0=E_1\az x_1=E_1\az x_2=E_2\az x_2=0 \;, \qquad
F_2\az x_0=F_2\az x_1=F_1\az x_2=0 \;.
\end{gather*}
\end{prop}

\noindent
The action on the $x_i^*$'s is determined by compatibility with the involution:
$$
K_i\az a^*=\{K_i^{-1}\az a\}^* \;,\qquad
E_1\az a^*=\{-qF_1\az a\}^* \;,\qquad
E_2\az a^*=\{-q^2F_2\az a\}^* \;.
$$

\begin{prova}
The algebra $\A(S^4_q)$ is the quotient of the $*$-algebra
$\C\!\inner{x}$ --- the free $*$-algebra
generated by the $x_i$'s --- by the ideal generated by
suitable degree $\leq 2$ polynomials.
The bijective linear map from the linear span of $\{x_0,x_1,x_1^*,x_2,x_2^*\}$
to the representation space $V_1$ defined (modulo a global
proportionality constant) by 
$$
x_2\mapsto\ket{0,1;0}\,,\;
x_1\mapsto\ket{1,0;1}\,,\;
x_0\mapsto (q[2])^{-1/2}\ket{0,0;1}\,,\;
x_1^*\mapsto -q\ket{-1,0;1}\,,\;
x_2^*\mapsto q^3\ket{0,-1;0}\,,
$$
is a unitary equivalence of $U_q(so(5))$-modules (
the real structure $C$ on $V_1$ being implemented by the $*$ operation
on $x_i$'s). This guarantee that $\C\!\inner{x}$
is an $U_q(so(5))$-module $*$-algebra.

It may help to draw the weight diagram of $V_{(0,1)}$, with each
generator of the $\A(S^4_q)$ in the position corresponding to its
weight.

\begin{footnotesize}\begin{center}\begin{tabular}{ccc}
\begin{tabular}{c}
\begindc{\commdiag}[20]
 \obj(1,3)[A]{$\bullet$}
 \obj(3,3)[B]{$\bullet$}
 \obj(5,3)[C]{$\bullet$}
 \obj(3,1)[D]{$\bullet$}
 \obj(3,5)[E]{$(0,1)$}
 \mor{A}{B}{}
 \mor{B}{C}{}
 \mor{D}{A}{}[\atright,\dasharrow]
 \mor{C}{E}{}[\atright,\dasharrow]
\enddc
\end{tabular}
&~~~~&
\begin{tabular}{c}
\begindc{\commdiag}[20]
 \obj(1,3)[A]{$x_1^*$}
 \obj(3,3)[B]{$x_0$}
 \obj(5,3)[C]{$x_1$}
 \obj(3,1)[D]{$x_2^*$}
 \obj(3,5)[E]{$x_2$}
 \mor{A}{B}{}
 \mor{B}{C}{}
 \mor{D}{A}{}[\atright,\dasharrow]
 \mor{C}{E}{}[\atright,\dasharrow]
\enddc
\end{tabular}
\end{tabular}\end{center}\end{footnotesize}

\noindent A solid arrow indicates points that can be joined by applying $E_1$ (the reverse
arrow corresponds to $F_1$), a dashed arrow indicates points that can be joined by applying
$E_2$ (the reverse arrow corresponds to $F_2$). In the left diagram, the highest weight
vector is denoted by its weight $(n_1,n_2)$. A bullet indicates a weight with multiplicity $1$,
an empty circle a weight with multiplicity $2$.

The polynomials generating the ideal which
defines $\A(S^4_q)$, that is
\begin{align*}
v_{ij}&:=x_ix_j-q^2x_jx_i      \;, \qquad 0\leq i<j\leq 2 \;, \\
w_{ij}&:=x_i^*x_j-q^2x_jx_i^*  \;, \qquad i\neq j \;, \\
u_1&:=[x_1^*,x_1]-(1-q^4)x_0^2 \;, \\
u_2&:=[x_2^*,x_2]-x_1^*x_1-q^4x_1x_1^* \;, \\
c&:=x_0^2+x_1x_1^*+x_2x_2^*-1 \;,
\end{align*}
span the real representations $V_0$ and
$V_{(1,1)}$ inside the tensor product $V_1\otimes V_1$.
More precisely, $c$ is invariant and spans $V_0$ while
the remaining polynomials are drawn in the weight diagram
of the adjoint representation $V_{(1,1)}$ in the position
corresponding to their weight.

\begin{footnotesize}\begin{center}\begin{tabular}{ccc}
\begin{tabular}{c}
\begindc{\commdiag}[20]
 \obj(3,3)[B]{$\bullet$}
 \obj(5,3)[C]{$\bullet$}
 \obj(7,3)[D]{$\bullet$}
 \obj(3,5)[F]{$\bullet$}
 \obj(5,5)[G]{$\circ$}
 \obj(7,5)[H]{$\bullet$}
 \obj(3,7)[L]{$\bullet$}
 \obj(5,7)[M]{$\bullet$}
 \obj(7,7)[N]{$(1,1)$}
 \mor{B}{C}{}
 \mor{C}{D}{}
 \mor{F}{G}{}
 \mor{G}{H}{}
 \mor{L}{M}{}
 \mor{M}{N}{}
 \mor{C}{F}{}[\atright,\dasharrow]
 \mor{D}{G}{}[\atright,\dasharrow]
 \mor{G}{L}{}[\atright,\dasharrow]
 \mor{H}{M}{}[\atright,\dasharrow]
\enddc
\end{tabular}
&~~~&
\begin{tabular}{c}
\begindc{\commdiag}[20]
 \obj(3,3)[B]{$v_{12}^*$}
 \obj(5,3)[C]{$v_{02}^*$}
 \obj(7,3)[D]{$w_{21}$}
 \obj(3,5)[F]{$v_{01}^*$}
 \obj(5,5)[G]{$\{u_1,u_2\}$}
 \obj(7,5)[H]{$v_{01}$}
 \obj(3,7)[L]{$w_{12}$}
 \obj(5,7)[M]{$v_{02}$}
 \obj(7,7)[N]{$v_{12}$}
 \mor{B}{C}{}
 \mor{C}{D}{}
 \mor{F}{G}{}
 \mor{G}{H}{}
 \mor{L}{M}{}
 \mor{M}{N}{}
 \mor{C}{F}{}[\atright,\dasharrow]
 \mor{D}{G}{}[\atright,\dasharrow]
 \mor{G}{L}{}[\atright,\dasharrow]
 \mor{H}{M}{}[\atright,\dasharrow]
\enddc
\end{tabular}
\end{tabular}\end{center}\end{footnotesize}

As before, one moves in the direction of solid (resp.~dashed) arrows by applying $E_1$
(resp.~$E_2$).
So, for example, one checks that $E_1\az v_{01}=0$ and
$F_1\az v_{01}=-q^{-1/2}u_1$.

Since $\A(S^4_q)$ is the quotient of the $U_q(so(5))$-module $*$-algebra
$\C\!\inner{x}$ by an ideal which is invariant under $U_q(so(5))$, 
it is itself an $U_q(so(5))$-module $*$-algebra.
\end{prova}

\begin{prop}\label{prop:S4q}
There is an isomorphism $\,\A(S^4_q)\simeq\bigoplus_{l\in\N}V_l\,$ of left $U_q(so(5))$-modules.
\end{prop}
\begin{prova}
A linear basis for $\A(S^4_q)$ is made of monomials $x_0^{n_0}x_1^{n_1}(x_1^*)^{n_2}x_2^{n_3}$
with $n_0,n_1,n_2\in\N$, $n_3\in\Z$ and with the notation $x_2^{n_3}:=(x_2^*)^{|n_3|}$
if $n_3<0$.
Such a monomial is an eigenvector of $K_1$ (resp.~$K_1K_2$)
with eigenvalue $n_1-n_2$ (resp.~$n_3$). Thus, the generic
eigenvector $v$ of $K_1$ and $K_2$ is a sum of
monomials with $n_1-n_2$ and $n_3$ fixed. Since $E_2$
sends $x_1$ to $x_0$, being the basis monomials
linearly independent,
we deduce that $v$ is annihilated by $E_2$ if{}f it contains
only monomials with $n_1=0$:
$$
v=\sum\nolimits_{n_0}c_{n_0}x_0^{n_0}(x_1^*)^{n_2}x_2^{n_3} \;.
$$
Finally, the condition $E_1\az v=0$ 
implies that $v$ must be independent of $x_0$ and $x_1^*$.
We conclude that highest weight vectors are of the form $x_2^l$ ---
corresponding to the weight $(0,l)$, $l\in\N$ ---
and $\A(S^4_q)$ decomposes as a multiplicity free
direct sum of highest weight representations with weights $(0,l)$.
\end{prova}

The algebra $\A(S^4_q)$ has two inequivalent irreducible infinite dimensional
representations. The representation space is the Hilbert space $\ell^2(\N^2)$
and the representations are given by  
\begin{align}\label{4irreps}
x_0\ket{k_1,k_2}_\pm &:=\pm q^{2(k_1+k_2)}\ket{k_1,k_2}_\pm \;, \nonumber \\
x_1\ket{k_1,k_2}_\pm &:=q^{2k_2}\sqrt{1-q^{4(k_1+1)}}\ket{k_1+1,k_2}_\pm \;,\\
x_2\ket{k_1,k_2}_\pm &:=\sqrt{1-q^{4(k_2+1)}}\ket{k_1,k_2+1}_\pm \nonumber \;.
\end{align}
The direct sum of these representations, with obvious grading $\gamma$ and  
operator $F$ given by $F\ket{k_1,k_2}_\pm:=\ket{k_1,k_2}_\mp$, constitutes
a $1$-summable Fredholm module over $\A(S^4_q)$. 

\bigskip

Recall the Definition \ref{def:SUq2} of $SU_q(2)$ and the Definition
\ref{def:S2q} of the equatorial Podle\'s sphere $S^2_{q^2}=S^2_{q1}$ with
parameter $q^2$.

\begin{prop}
There is a $*$-algebra morphism $\varphi:\A(S^4_q)\to\A(SU_q(2))\otimes\A(S^2_{q^2})$
defined by:
\begin{align}
\varphi(x_0) &=-(\alpha\beta+\beta^*\alpha^*)\otimes A\;, \nonumber \\
\varphi(x_1) &=\bigl(-\alpha^2+q\,(\beta^*)^2\bigr)\otimes A\;, \label{eq:embed} \\
\varphi(x_2) &=1\otimes B\;. \nonumber
\end{align}
\end{prop}
\begin{prova}
One proves by direct computation that the five $\varphi(x_i),\varphi(x_i)^*$ satisfy
all the defining relations of $\A(S^4_q)$.
\end{prova}

\section{The modules of chiral spinors}\label{sec:7.3}

We apply the general theory of Section~\ref{sec:3.2}, 
to the case $\A=\A(S^4_q)$ and $\mc{U}=U_q(so(5))$. In the present case
$\kappa(a)=K_1^8K_2^6\az a$ is the modular automorphism \cite[Section 11.3.4]{KS97}.
We shall use the notations of Section~\ref{sec:7.1} for the irreducible
representations $(V_l,\sigma_l)$ of $U_q(so(5))$. 

By Proposition \ref{prop:S4q} we have the equivalence
\mbox{$\A(S^4_q)\simeq\bigoplus_{l\in\N}V_l$} as left $U_q(so(5))$-modules.
Using Lemma~\ref{lemma:3.2} for $N=1$, we deduce that on the vector space
$\bigoplus_{l\in\N}V_l$ there exists at least one $*$-representation of
the crossed product $\A(S^4_q)\rtimes U_q(so(5))$ that extends the $*$-representation
$\,\bigoplus_{l\in\N}\sigma_l\,$ of $U_q(so(5))$.

Let $e\in\mr{Mat}_4(\A(S^4_q))$ be the following idempotent:
\begin{equation}\label{eq:P}
e:=\frac{1}{2}\left(\begin{array}{cccc}
   1+x_0       &  q^3x_2   & -qx_1        & 0 \\
   q^{-3}x_2^* &  1-q^2x_0 &  0           & q^3x_1 \\
  -q^{-1}x_1^* &  0        &  1-q^2x_0    & q^3x_2 \\
   0           &  qx_1^*   &  q^{-3}x_2^* & 1+q^4x_0
  \end{array}\right)\;.
\end{equation}
By direct computation one proves that $K_1^8K_2^6\az e^*=e=e^2$ and then,
by Lemma \ref{lemma:next},
$e$ defines an orthogonal projection $\pi$, by equation~\eqref{eq:pi}, on the linear space
$\A(S^4_q)^4$ with inner product~\eqref{eq:in}.
Next, let $\sigma:U_q(so(5))\to\mr{Mat}_4(\C)$ be the $*$-representation defined by
\begin{subequations}\label{eq:spin}
\begin{equation}
\sigma(K_1)=\ma{35pt}{
  q^{1/2} & 0 & 0 & 0 \\
  0 & q^{1/2} & 0 & 0 \\
  0 & 0 & q^{-1/2} & 0 \\
  0 & 0 & 0 & q^{-1/2}
}\;,\qquad
\sigma(K_2)=\ma{35pt}{
  1 & 0 & 0 & 0 \\
  0 & q^{-1} & 0 & 0 \\
  0 & 0 & q & 0 \\
  0 & 0 & 0 & 1
}\;,
\end{equation}
\begin{equation}
\sigma(E_1)=\ma{35pt}{
  0 & 0 & 1 & 0 \\
  0 & 0 & 0 & 1 \\
  0 & 0 & 0 & 0 \\
  0 & 0 & 0 & 0
}\;,\qquad
\sigma(E_2)=\ma{35pt}{
  0 & 0 & 0 & 0 \\
  0 & 0 & 0 & 0 \\
  0 & 1 & 0 & 0 \\
  0 & 0 & 0 & 0
}\;.
\end{equation}
\end{subequations}
Again, by direct computation one proves that:
\begin{subequations}\label{eq:sub}
\begin{align}
K_i\az e &=\sigma(K_i)\,e\,\sigma(K_i)^{-1} \;,\\
E_i\az e &=\sigma(F_i)\,e\,\sigma(K_i)^{-1}-q^{-i}\sigma(K_i)^{-1}e\,\sigma(F_i) \;.
\end{align}
\end{subequations}
Since $\sigma(K_i)=\sigma(K_i)^t$ and $\sigma(F_i)=\sigma(E_i)^t$, we conclude
that condition (\ref{eq:cov}) is satisfied and that $\pi$ and $\pi^\perp=1-\pi$ project
$\A(S^4_q)^4$ onto $*$-subrepresentations of $\A(S^4_q)\rtimes U_q(so(5))$.

We state the main Proposition of this section.

\begin{prop}\label{prop:main}
There exist two inequivalent representations of the crossed product
algebra $\A(S^4_q)\rtimes U_q(so(5))$ on $\bigoplus_{l\in\N+\frac{1}{2}}V_l$
that extend the representation $\,\bigoplus_{l\in\N+\frac{1}{2}}\sigma_l\,$
of $U_q(so(5))$.
\end{prop}

The proof is in two steps. We first prove (in Lemma \ref{lemma:ptwo})
that $\pi(\A(S^4_q)^4)$ and $\pi^\perp(\A(S^4_q)^4)$ are not equivalent
as representations of the algebra $\A(S^4_q)$. Then we prove (in Lemma
\ref{lemma:pone}) that as $U_q(so(5))$ representations they are both equivalent
to $\bigoplus_{l\in\N+\frac{1}{2}}V_l$.

\begin{lemma}\label{lemma:ptwo}
The idempotent $e$ in (\ref{eq:P}) splits $\A(S^4_q)^4$ into two \emph{inequivalent}
$*$-representations of the crossed product algebra $\A(S^4_q)\rtimes U_q(so(5))$.
\end{lemma}

\begin{prova}
To prove the statement we apply Lemma~\ref{lemma:inv}.
We use the Fredholm module associated to the representation on
$\ell^2(\N)\oplus\ell^2(\N)$, and defined by equation (\ref{4irreps}).
From (\ref{eq:ind}) we have
\begin{align*}
\mr{Index}(F_e^+)
&=\tfrac{1}{2}\tr_{\ell^2(\N)\otimes\C^8}(\gamma F[F,e]) \\
\rule{0pt}{14pt}
&=\tfrac{1}{4}(1-q^2)^2\tr_{\ell^2(\N)\otimes\C^2}(\gamma F[F,x_0]) \\
\rule{0pt}{14pt}
&=(1-q^2)^2\sum\nolimits_{k_1,k_2\in\N}q^{2(k_1+k_2)}=1\;.
\end{align*}
By Lemma~\ref{lemma:inv} this proves that the two representations of $\A(S^4_q)$
are inequivalent.
The statement of the Proposition~\ref{prop:main} follows from the obvious
observation that if the two representations of the crossed product algebra
were equivalent, their restrictions to representations of $\A(S^4_q)$
would be equivalent too.
\end{prova}

\begin{lemma}\label{lemma:pone}
$\pi(\A(S^4_q)^4)\simeq\pi^\perp(\A(S^4_q)^4)\simeq\bigoplus_{l\in\N+\frac{1}{2}}V_l$
as $U_q(so(5))$ representations.
\end{lemma}
\begin{prova}
In this proof, `$\simeq$' always means equivalence of representations of $U_q(so(5))$.

Since $\sigma$ in \eqref{eq:spin} is just the spin representation $V_{1/2}$ in matrix form,
the representation of $U_q(so(5))$ on $\A(S^4_q)^4$ is the Hopf tensor product
of the representation over $\A(S^4_q)$ with the representation $V_{1/2}$.
From $\A(S^4_q)\simeq\bigoplus_{l\in\N}V_l$ and from the decomposition
$V_l\otimes V_{1/2}\simeq V_{l-\frac{1}{2}}\oplus V_{l+\frac{1}{2}}$
for all $l\in\{1,2,3,\ldots\}$, we deduce that
$
\A(S^4_q)^4\simeq\bigoplus\nolimits_{l\in\N+\frac{1}{2}}(V_l\oplus V_l)
$
and then,
$$
\pi(\A(S^4_q)^4)\simeq\bigoplus\nolimits_{l\in\N+\frac{1}{2}}m^+_lV_l\;,\qquad
\pi^\perp(\A(S^4_q)^4)\simeq\bigoplus\nolimits_{l\in\N+\frac{1}{2}}m^-_lV_l\;,
$$
with multiplicities $m^\pm_l$ to be determined, such that $m_l^++m^-_l=2$.
For $l\in\N+\frac{1}{2}$, the vectors 
$$
v_l^\pm:=x_2^{l-\frac{1}{2}}(1\pm x_0,\pm q^3x_2,\mp qx_1,0) \;,
$$
are highest weight vectors, being annihilated by both $E_1$ and $E_2$, and
have weight $(\frac{1}{2},l)$. Furthermore, $v^+_l(1-e)=v^-_le=0$. Thus,
$v^+_l\in\pi(\A(S^4_q)^4)$ and $v^-_l\in\pi^\perp(\A(S^4_q)^4)$.

Then in both modules $\pi(\A(S^4_q)^4)$ and $\pi^\perp(\A(S^4_q)^4)$
each representation $V_l$, $l\in\N+\frac{1}{2}$, appears with multiplicity
$m_l^\pm\geq 1$. Since $m_l^++m^-_l=2$, we deduce that $m_l^\pm=1$ for all
$l\in\N+\frac{1}{2}$.
\end{prova}

\section{Equivariant representations}\label{sec:7.4}
Next, we construct $U_q(so(5))$-equivariant representations of $\A(S^4_q)$ which
classically correspond to the left regular and chiral spinor representations.
The representation spaces will be (the closure of) $\bigoplus_{l\in\N}V_l$
and $\bigoplus_{l\in\N+\frac{1}{2}}V_l$, respectively.

Equivariance of a representation means that it is a representation of the crossed
product algebra $\A(S^4_q)\rtimes U_q(so(5))$. The latter is defined by the crossed relations
$ha=(h_{(1)}\az a)h_{(2)}$ for all $a\in\A(S^4_q)$ and $h\in U_q(so(5))$. Explicitly, the
relations among the generators read:
\begin{equation}\label{eq:crossA}\begin{array}{ccc}
[K_1,x_0]=0\;, & K_1x_1=qx_1K_1\;, & K_1x_2=x_2K_1\;, \\
\rule{0pt}{18pt}
[K_2,x_0]=0\;, & K_2x_1=q^{-1}x_1K_2\;, & K_2x_2=qx_2K_2\;, \\
\rule{0pt}{18pt}
[E_1,x_0]=q^{-1/2}x_1K_1\;, & E_1x_1=q^{-1}x_1E_1\;, & E_1x_2=x_2E_1\;, \\
\rule{0pt}{18pt}
[F_1,x_0]=-q^{-1/2}K_1x_1^*\;, &\hspace{3mm} F_1x_1=q^{-1}x_1F_1+q^{1/2}[2]x_0K_1\;,\hspace{3mm} & F_1x_2=x_2F_1\;, \\
\rule{0pt}{18pt}
[E_2,x_0]=0\;, & E_2x_1=qx_1E_2+x_2K_2\;, & E_2x_2=q^{-1}x_2E_2\;, \\
\rule{0pt}{18pt}
[F_2,x_0]=0\;, & F_2x_1=qx_1F_2\;, & \hspace{-5mm}F_2x_2=q^{-1}x_2F_2+x_1K_2\;.\hspace{-5mm}
\end{array}\end{equation}
In previous Section we proved that on $\bigoplus_{l\in\N}V_l$ there is at
least one equivariant representation, the left regular one, and that on
$\bigoplus_{l\in\N+\frac{1}{2}}V_l$ there are at least two equivariant
representations, corresponding to the projective modules $\A(S^4_q)^4e$
and $\A(S^4_q)^4(1-e)$. In this section we prove that on such spaces
there are no other equivariant representations (modulo unitary equivalences)
besides the ones just mentioned.

Let us denote with $\ket{l,m_1,m_2;j}$ the basis of the representation
space $V_l$ of $U_q(so(5))$ as discussed in Section \ref{sec:7.1}.
From the first two lines of (\ref{eq:crossA}) we deduce that 
\begin{subequations}\label{eq:coeff}
\begin{align}
x_0\ket{l,m_1,m_2;j} &=\sum_{l',j'}A_{j,j',l,l'}^{m_1,m_2}\ket{l',m_1,m_2;j'} \;, \\
x_1\ket{l,m_1,m_2;j} &=\sum_{l',j'}B_{j,j',l,l'}^{m_1,m_2}\ket{l',m_1+1,m_2;j'} \;, \\
x_2\ket{l,m_1,m_2;j} &=\sum_{l',j'}C_{j,j',l,l'}^{m_1,m_2}\ket{l',m_1,m_2+1;j'} \;,\label{eq:coeffC}
\end{align}
\end{subequations}
with coefficients to be determined. Notice that from the crossed relations
\begin{align*}
x_1\ket{l,m_1,m_2;j} &=(F_2x_2-q^{-1}x_2F_2)K_2^{-1}\ket{l,m_1,m_2;j} \;,\\
x_0\ket{l,m_1,m_2;j} &=q^{-1/2}[2]^{-1}(F_1x_1-q^{-1}x_1F_1)K_1^{-1}\ket{l,m_1,m_2;j} \;,
\end{align*}
the matrix coefficients of $x_0$ and $x_1$ can be expressed in term
of the coefficients of $x_2$.

\begin{lemma}
Let $k\in\N$. The following formul{\ae} hold:
\begin{subequations}\label{eq:nul}
\begin{align}
F_1^k\ket{l,m_1,m_2;j} &=
  \left\{\begin{array}{ll}
  =0 &\quad\mr{if}\;k>j+m_1 \\
  \neq 0&\quad\mr{if}\;k\leq j+m_1
  \end{array}\right. \label{eq:nulA} \\
E_1^k\ket{l,m_1,m_2;j} &=
  \left\{\begin{array}{ll}
  =0 &\quad\mr{if}\;k>j-m_1 \\
  \neq 0&\quad\mr{if}\;k\leq j-m_1
  \end{array}\right. \;\; . \label{eq:nulB}
\end{align} 
\end{subequations} 
\end{lemma}

\begin{prova}
By direct computation:
\begin{align*}
F_1^k\ket{l,m_1,m_2;j} &=\sqrt{[j+m_1][j+m_1-1]\ldots [j+m_1-k+1]}\,\times \\ 
& \qquad \times\sqrt{[j-m_1+1][j-m_1+2]\ldots [j-m_1+k]}\,\ket{l,m_1-k,m_2;j}\;.
\end{align*}
The second square root is always different from zero since the
$q$-analogues are in increasing order and $j-m_1+1\geq 1$. In
the first square root $q$-analogues are in decreasing order and
are all different from zero if and only if $j+m_1-k+1\geq 1$.
This proves equation (\ref{eq:nulA}).

In the same way one establishes (\ref{eq:nulB}), by using the equation
\begin{align*}
E_1^k\ket{l,m_1,m_2;j} &=\sqrt{[j-m_1][j-m_1-1]\ldots [j-m_1-k+1]}\,\times \\ 
& \qquad
\times\sqrt{[j+m_1+1][j+m_1+2]\ldots [j+m_1+k]}\,\ket{l,m_1+k,m_2;j}\;. \qedhere
\end{align*}
\end{prova}

\begin{lemma}\label{lemma:A}
The coefficients in (\ref{eq:coeff}) satisfy:
$$
A^{m_1,m_2}_{j,j',l,l'}=B^{m_1,m_2}_{j,j',l,l'}=0\;\;\mr{if}\;\;|j-j'|>1\;,
\qquad C^{m_1,m_2}_{j,j',l,l'}=0\;\;\mr{if}\;\;j'\neq j\;.
$$
\end{lemma}

\begin{prova}
From (\ref{eq:crossA}), (\ref{eq:nulA}) and (\ref{eq:nulB}) we derive:
\begin{align*}
E_1^{j-m_1+1}x_1\ket{l,m_1,m_2;j} &=q^{-j+m_1-1}x_1E_1^{j-m_1+1}\ket{l,m_1,m_2;j}=0\;, \\
F_1^{j'+m_1+2}x_1^*\ket{l',m_1+1,m_2;j'} &=q^{j'+m_1+2}x_1^*F_1^{j'+m_1+2}\ket{l',m_1+1,m_2;j'}=0\;.
\end{align*}
We expand the left hand sides and use the independence of the vectors
$E_1^{j-m_1+1}\ket{l',m_1+1,m_2;j'}$ and
$F_1^{j'+m_1+2}\ket{l,m_1,m_2;j}$ to arrive at the conditions:
\begin{align*}
B_{j,j',l,l'}^{m_1,m_2}\Big\{E_1^{j-m_1+1}\ket{l',m_1+1,m_2;j'}\Big\}  &=0\;, \\
\bar{B}_{j,j',l,l'}^{m_1,m_2}\Big\{F_1^{j'+m_1+2}\ket{l,m_1,m_2;j}\Big\} &=0\;.
\end{align*}
By (\ref{eq:nulB}) the curly brackets in the first line are different from zero
if $j-m_1+1\leq j'-m_1-1$, i.e.~$B_{j,j',l,l'}^{m_1,m_2}$ must be zero
if $j'\geq j+2$. By (\ref{eq:nulA}) the curly brackets in the second line is different
from zero if $j'+m_1+2\leq j+m_1$, i.e.~$\bar{B}_{j,j',l,l'}^{m_1,m_2}$ must be zero
if $j'\leq j-2$. This proves $1/3$ of the statement
\begin{equation*}
B_{j,j',l,l'}^{m_1,m_2}=0\quad\forall\;j'\notin\{j-1,j,j+1\}\;.
\end{equation*}
A similar argument applies to $x_0$. From the coproduct of $E_1^n$ we deduce:
\begin{equation*}
E_1^nx_0 =\sum_{k=0}^n\sqbn{n}{k}(E_1^k\az x_0)E_1^{n-k}K_1^k
       =x_0E_1^n-[n]q^{-1/2}x_1E_1^{n-1}K_1\;.
\end{equation*}
This implies that $E_1^{j-m_1+2}x_0\ket{l,m_1,m_2;j}=0$ and
$F_1^{j'+m_1+2}x_0\ket{l',m_1,m_2;j}=0$. From these conditions
we deduce that also $x_0$ shift $j$ by $\{0,\pm 1\}$ only.

Finally, let $\mc{C}_1$ be the Casimir element in equation (\ref{eq:Cdef}). Then
$[\mc{C}_1,x_2]=0$ and from (\ref{eq:Caz}) we deduce that $x_2$ is diagonal
on the index $j$.
\end{prova}

\begin{lemma}\label{le:B}
The coefficients in (\ref{eq:coeffC}) satisfy
$$
C^{m_1,m_2}_{j,j',l,l'}=0\;\;\mr{if}\;\;|l-l'|>1\;\;\textup{or if}\;\;
|l-l'|=0\;\;\mr{and}\;\;l\in\N\;.
$$
\end{lemma}
\begin{prova}
The elements $\{x_i,x_i^*\}$ are a basis of the irreducible representation $V_1$.
Covariance of the action tells that $x_i\ket{l,m_1,m_2;j}$ and
$x_i^*\ket{l,m_1,m_2;j}$ are a basis of the tensor representation
$V_1\otimes V_l$. Equations (14--15) in Chapter 7 of~\cite{KS97} tells
that $V_1\otimes V_l\simeq V_{l-1}\oplus V_{l+1}$ if $l\in\N$
and that $V_1\otimes V_l\simeq V_{l-1}\oplus V_l\oplus V_{l+1}$ if $l\in\N+\frac{1}{2}$
(with $V_{l-1}$ omitted if $l-1<0$). This Clebsh-Gordan decomposition
tells that $x_2\ket{l,m_1,m_2;j}$ is in the linear span of the basis
vectors $\ket{l',m_1,m_2+1;j}$ with $l'-l=\pm 1$ if $l\in\N$ or
with $l'-l=0,\pm 1$ if $l\in\N+\frac{1}{2}$. This concludes the proof of the Lemma.
\end{prova}

\subsection{Computing the coefficients of $x_2$}
From Lemma~\ref{le:B}, we have to consider only the cases $j'=j$, $|l'-l|\leq 1$ if $l\in\N+\frac{1}{2}$
or $|l'-l|=1$ if $l\in\N$.
The condition $[E_1,x_2]=0$ implies that $C_{j,j,l,l'}^{m_1,m_2}=:C_{j,l,l'}^{m_2}$ is
independent of $m_1$.
Equations $(E_2x_2-q^{-1}x_2E_2)\ket{l,-j,m_2;j}=0$ and
$(F_2x_2^*-qx_2^*F_2)\ket{l',j,m_2+1;j}=0$ imply, respectively:{\small
\begin{gather*}
C^{m_2}_{j,l,l'}\sqrt{[l'-j-m_2-1+\epsilon'][l'+j+m_2+4+\epsilon']}\,=
C^{m_2+1}_{j+1,l,l'}q^{-1}\sqrt{[l-j-m_2+\epsilon][l+j+m_2+3+\epsilon]}\;, \\
C^{m_2}_{j,l,l'}\sqrt{[l+j-m_2+3-\epsilon][l-j+m_2-\epsilon]}\,=
C^{m_2-1}_{j+1,l,l'}q\sqrt{[l'+j-m_2+2-\epsilon'][l'-j+m_2+1-\epsilon']}\;,
\end{gather*}
}with $\epsilon,\epsilon'\in\{0,\pm\frac{1}{2}\}$ determined by the conditions
$l+\epsilon-j-m_2\in 2\N$ and $l'-\epsilon'-j-m_2\in 2\N$.
Notice that if $l'-l\in 2\N+1$ then $\epsilon'=\epsilon$, while if
$l'-l\in 2\N$ then $\epsilon'=-\epsilon$.
Looking at the cases $l'-l=\pm 1$, we deduce that
$$
\frac{ q^{-\frac{1}{2}(j+m_2)} }{\sqrt{[l+j+m_2+3+\epsilon]}}\, C_{j,l,l+1}^{m_2}
\qquad\mr{and}\qquad
\frac{ q^{-\frac{1}{2}(j+m_2)} }{\sqrt{[l-j-m_2+\epsilon]}}\, C_{j,l,l-1}^{m_2} 
$$
depend on $j+m_2$ only through their parity (i.e.~they depend only on the value of $\epsilon$). Similarly,
$$
\frac{ q^{\frac{1}{2}(j-m_2)} }{\sqrt{[l-j+m_2+2-\epsilon]}}\, C_{j,l,l+1}^{m_2} 
\qquad\mr{and}\qquad
\frac{ q^{\frac{1}{2}(j-m_2)} }{\sqrt{[l+j-m_2+1-\epsilon]}}\, C_{j,l,l-1}^{m_2} 
$$
depend on $j-m_2$ only through their parity.
Combining these informations, we deduce that the following elements do not depend on
the exact value of $j$, $m_2$, but only on the value of $\epsilon$,
\begin{align*}
\frac{ q^{-m_2} }{\sqrt{[l+j+m_2+3+\epsilon][l-j+m_2+2-\epsilon]}}\, C_{j,l,l+1}^{m_2} &=:C_{l,l+1}(\epsilon) \;, \\
\frac{ q^{-m_2} }{\sqrt{[l-j-m_2+\epsilon][l+j-m_2+1-\epsilon]}}\, C_{j,l,l-1}^{m_2} 
 &=:C_{l,l-1}(\epsilon) \;.
\end{align*}
If $l\in\N$ there are no other coefficients $C_{j,l,l'}^{m_2}$ to compute.
If $l\notin\N$, we have to compute also $C_{j,l,l}^{m_2}$.
In this case $\epsilon'=-\epsilon$ and we get{\small
\begin{gather*}
C^{m_2}_{j,l,l}\sqrt{[l-j-m_2-1-\epsilon][l+j+m_2+4-\epsilon]}\,=
C^{m_2+1}_{j+1,l,l}q^{-1}\sqrt{[l-j-m_2+\epsilon][l+j+m_2+3+\epsilon]}\;, \\
C^{m_2}_{j,l,l}\sqrt{[l+j-m_2+3-\epsilon][l-j+m_2-\epsilon]}\,=
C^{m_2-1}_{j+1,l,l}q\sqrt{[l+j-m_2+2+\epsilon][l-j+m_2+1+\epsilon]}\;.
\end{gather*}
}Again, looking at the two cases $\epsilon=\pm\frac{1}{2}$ we deduce that
$$
\frac{ q^{-\frac{1}{2}(j+m_2)} }{\sqrt{[l+\frac{1}{2}-j-m_2]}}\, C_{j,l,l}^{m_2} 
\;\;\;\mr{if}\;\epsilon=\tfrac{1}{2}
\quad\mr{and}\quad
\frac{ q^{-\frac{1}{2}(j+m_2)} }{\sqrt{[l+\frac{1}{2}+j+m_2+2]}}\, C_{j,l,l}^{m_2} 
\;\;\;\mr{if}\;\epsilon=-\tfrac{1}{2}
$$
do not depend on $j+m_2$ (this time $\epsilon$ is fixed, then the parity of $j+m_2$ is fixed). Similarly,
$$
\frac{ q^{\frac{1}{2}(j-m_2)} }{\sqrt{[l+\frac{1}{2}-j+m_2+1]}}\, C_{j,l,l}^{m_2} 
\;\;\;\mr{if}\;\epsilon=\tfrac{1}{2}
\quad\mr{and}\quad
\frac{ q^{\frac{1}{2}(j-m_2)} }{\sqrt{[l+\frac{1}{2}+j-m_2+1]}}\, C_{j,l,l}^{m_2} 
\;\;\;\mr{if}\;\epsilon=-\tfrac{1}{2}
$$
do not depend on $j-m_2$.
Combining these informations, we deduce that the following element do not depend on
the exact value of $j$, $m_2$, but only on the value of $\epsilon$:
$$
\frac{ q^{-m_2} }{\sqrt{[l-2\epsilon j-m_2+1-\epsilon]
[l-2\epsilon j+m_2+2-\epsilon]}}\, C_{j,l,l}^{m_2} =:C_{l,l}(\epsilon) \;.
$$
The denominator of the left hand side is just $[2j][2j+2]b_l(j,m_2)$
with $b_l$ the coefficient in equation (\ref{eq:abc}).
The formula $C_{j,l,l}^{m_2}=q^{m_2}[2j][2j+2]b_l(j,m_2)C_{l,l}(\epsilon)$ is
valid for all $l$, since $b_l(j,m_2)$ vanish for $l$ integer.

Summarizing, we find that
\begin{subequations}\label{eq:Cjllm}
\begin{align}
C_{j,l,l+1}^{m_2} &=q^{m_2}\sqrt{[l+j+m_2+3+\epsilon][l-j+m_2+2-\epsilon]} \, C_{l,l+1}(\epsilon) \;, \\
C_{j,l,l}^{m_2} &=q^{m_2}[2j][2j+2]b_l(j,m_2) \, C_{l,l}(\epsilon) \;, \\
C_{j,l,l-1}^{m_2} &=q^{m_2}\sqrt{[l-j-m_2+\epsilon][l+j-m_2+1-\epsilon]} \, 
C_{l,l-1}(\epsilon)\;,
\end{align}
\end{subequations}
with coefficients $C_{l,l'}(\epsilon)$ to be determined.

\subsection{Computing the coefficients of $x_1$}
From Lemma~\ref{lemma:A}, we have to consider only the three cases 
$j'=j,j\pm 1$.
Using equation $E_1x_1=q^{-1}x_1E_1$ we get, 
\begin{align*}
\frac{ q^{-m_1} }{\sqrt{[j+m_1+1][j+m_1+2]}}\, B_{j,j+1,l,l'}^{m_1,m_2}  &=
\frac{ q^{-m_1-1} }{\sqrt{[j+m_1+2][j+m_1+3]}}\, B_{j,j+1,l,l'}^{m_1+1,m_2} \;, \\
\frac{ q^{-m_1} }{\sqrt{[j-m_1][j+m_1+1]}}\, B_{j,j,l,l'}^{m_1,m_2}  &=
\frac{ q^{-m_1-1} }{\sqrt{[j-m_1-1][j+m_1+2]}}\, B_{j,j,l,l'}^{m_1+1,m_2} \;, \\
\frac{ q^{-m_1} }{\sqrt{[j-m_1][j-m_1-1]}}\, B_{j,j-1,l,l'}^{m_1,m_2} &=
\frac{ q^{-m_1-1} }{\sqrt{[j-m_1-1][j-m_1-2]}}\, B_{j,j-1,l,l'}^{m_1+1,m_2} \;.
\end{align*}
We see that the left hand sides of these three equations are independent of $m_1$, and call:
\begin{subequations}\label{eq:Bm}
\begin{align}
B_{j,j+1,l,l'}^{m_1,m_2} &=:q^{m_1}\sqrt{[j+m_1+1][j+m_1+2]}\,B_{j,j+1,l,l'}^{m_2} \;, \\
B_{j,j,l,l'}^{m_1,m_2}    &=:q^{m_1}\sqrt{[j-m_1][j+m_1+1]}\,B_{j,j,l,l'}^{m_2} \;, \\
B_{j,j-1,l,l'}^{m_1,m_2} &=:q^{m_1}\sqrt{[j-m_1][j-m_1-1]}\,B_{j,j-1,l,l'}^{m_2} \;.
\end{align}
\end{subequations}
Imposing the condition $x_1K_2=F_2x_2-q^{-1}x_2F_2$ on the subspace
spanned by $\ket{l,j,m_2;j}$ (so $m_1=j$ and $B_{j,j,l,l'}^{m_1,m_2}=B_{j,j-1,l,l'}^{m_1,m_2}=0$
on this subspace) we get:
$$
q^{m_2}B^{m_2}_{j,j+1,l,l'}=c_{l'}(j+1,m_2)C^{m_2}_{j,l,l'}-q^{-1}c_l(j+1,m_2-1)C^{m_2-1}_{j+1,l,l'}\;.
$$
From this we deduce that, that since coefficients $C_{j,l,l'}^{m_2}$ vanishes
for $|l-l'|>1$, also $B_{j,j+1,l,l'}^{m_2}$ is zero in these cases.
In the remaining three cases $l'=l,l\pm 1$, using equation
(\ref{eq:Cjllm}) we get:
\begin{subequations}\label{eq:Bjp}
\begin{align}
B_{j,j+1,l,l+1}^{m_2} &=(-1)^{2\epsilon} q^{l-j+m_2-\epsilon}
  \sqrt{\frac{[l+j+m_2+3+\epsilon][l+j-m_2+3-\epsilon]}{[2(j+|\epsilon|)+1]
  [2(j-|\epsilon|)+3]}}\, C_{l,l+1}(\epsilon) \;, \\
B_{j,j+1,l,l}^{m_2} &=-2\epsilon q^{2\epsilon\,l-j+m_2-2+3\epsilon}
\frac{\sqrt{[l+\frac{1}{2}+j-2\epsilon\,m_2+2][l+\frac{1}{2}-j-2\epsilon\,m_2]}}{[2j+2]}
  \, C_{l,l}(\epsilon) \;, \\
B_{j,j+1,l,l-1}^{m_2} &=(-1)^{2\epsilon+1} q^{-l-j+m_2-3+\epsilon}
  \sqrt{\frac{[l-j+m_2-\epsilon][l-j-m_2+\epsilon]}{[2(j+|\epsilon|)+1]
 [2(j-|\epsilon|)+3]}}\, C_{l,l-1}(\epsilon) \;.
\end{align}
\end{subequations}
Imposing $qx_1^*K_2=x_2^*E_2-q^{-1}E_2x_2^*$ on the subspace spanned by
$\ket{l',-j+1,m_2;j-1}$  (so $B_{j,j,l,l'}^{m_1,m_2}=B_{j,j+1,l,l'}^{m_1,m_2}=0$
on this subspace) we get:
$$
q^{m_2}B^{m_2}_{j,j-1,l,l'}=a_{l'}(j-1,m_2)C^{m_2}_{j,l,l'}-q^{-1}a_l(j-1,m_2-1)C^{m_2-1}_{j-1,l,l'}\;.
$$
We deduce that $B_{j,j-1,l,l'}^{m_2}$ vanishes if $|l-l'|>1$, while in the three remaining
cases $l'=l,l\pm 1$ using equation (\ref{eq:Cjllm}) we get:
\begin{subequations}\label{eq:Bjm}
\begin{align}
B_{j,j-1,l,l+1}^{m_2} &=q^{l+j+m_2+1+\epsilon}
  \sqrt{\frac{[l-j-m_2+2+\epsilon][l-j+m_2+2-\epsilon]}{[2(j+|\epsilon|)-1][2(j-|\epsilon|)+1]}}\, C_{l,l+1}(\epsilon) \;, \\
B_{j,j-1,l,l}^{m_2} &=-2\epsilon q^{-2\epsilon\,l+j+m_2-1-3\epsilon}
\frac{\sqrt{[l+\frac{1}{2}+j+2\epsilon\,m_2+1][l+\frac{1}{2}-j+2\epsilon\,m_2+1]}}{[2j]}\, C_{l,l}(\epsilon)  \;, \\
B_{j,j-1,l,l-1}^{m_2} &=-q^{-l+j+m_2-2-\epsilon}
  \sqrt{\frac{[l+j+m_2+1+\epsilon][l+j-m_2+1-\epsilon]}{[2(j+|\epsilon|)-1][2(j-|\epsilon|)+1]}}\, C_{l,l-1}(\epsilon)\;.
\end{align}
\end{subequations}
Moreover, the condition 
$\inner{l',j,m_2;j|x_1K_2+q^{-1}x_2F_2-F_2x_2|l,j-1,m_2;j}=0$
implies that
\begin{equation}\label{eq:xxx}
q^{m_2}B^{m_2}_{j,j,l,l'}=b_{l'}(j,m_2)C^{m_2}_{j,l,l'}-q^{-1}b_l(j,m_2-1)C^{m_2-1}_{j,l,l'}\;.
\end{equation}
A further elaboration on these coefficients is postponed to after the following section.

\subsection{Computing the coefficients of $x_0$}
The condition $q^{1/2}[2]x_0K_1=F_1x_1-q^{-1}x_1F_1$ implies:
$$
q^{m_1+\frac{1}{2}}[2]A_{j,j',l,l'}^{m_1,m_2}=\sqrt{[j'-m_1][j'+m_1+1]}\, B_{j,j',l,l'}^{m_1,m_2}-q^{-1}\sqrt{[j+m_1][j-m_1+1]}\, B_{j,j',l,l'}^{m_1-1,m_2}\;.
$$
In the three non-trivial cases $j'-j=1,0,-1$, using (\ref{eq:Bm}), we get:
\begin{subequations}\label{eq:Am}
\begin{align}
A_{j,j+1,l,l'}^{m_1,m_2} &=q^{j+m_1-\frac{1}{2}}\sqrt{[j+m_1+1][j-m_1+1]}\,B^{m_2}_{j,j+1,l,l'} \;, \\
A_{j,j,l,l'}^{m_1,m_2} &=[2]^{-1}q^{-2-\frac{1}{2}}\big(q^{j+m_1+1}[2][j-m_1]-[2j]\big)B_{j,j,l,l'}^{m_2}\;, \\
A_{j,j-1,l,l'}^{m_1,m_2} &=-q^{-j+m_1-1-\frac{1}{2}}\sqrt{[j+m_1][j-m_1]}\,B^{m_2}_{j,j-1,l,l'} \;.
\end{align}
\end{subequations}
The hermiticity condition $x_0=x_0^*$ means that $A^{m_1,m_2}_{j,j+1,l,l'}=\bar{A}^{m_1,m_2}_{j+1,j,l'l}$
and $A^{m_1,m_2}_{j,j,l,l'}=\bar{A}^{m_1,m_2}_{j,j,l'l}$. Thus, from
(\ref{eq:Am}) follows that:
$$
B^{m_2}_{j+1,j,l'l}=-q^{2j+2}\bar{B}^{m_2}_{j,j+1,l,l'} \;,\qquad
B^{m_2}_{j,j,l,l'}=\bar{B}^{m_2}_{j,j,l'l} \;.
$$
Using (\ref{eq:Bjp}), the first equation turns out to be equivalent to
the following conditions:
\begin{subequations}\label{eq:sol}
\begin{equation}
C_{l+1,l}(\epsilon)=(-1)^{2\epsilon}q^{2l+4}\bar{C}_{l,l+1}(\epsilon)\;,\qquad
C_{l,l}(\epsilon)=\bar{C}_{l,l}(-\epsilon)\;.
\end{equation}
The second of equation together with (\ref{eq:xxx}) implies:
$$
b_{l'}(j,m_2)C^{m_2}_{j,l,l'}-q^{-1}b_l(j,m_2-1)C^{m_2-1}_{j,l,l'}=
b_l(j,m_2)C^{m_2}_{j,l'l}-q^{-1}b_{l'}(j,m_2-1)C^{m_2-1}_{j,l'l}\;.
$$
That is, using (\ref{eq:Cjllm}):
\begin{equation}
C_{l,l+1}(\epsilon)=C_{l,l+1}(-\epsilon)\;,\qquad
C_{l,l}(\epsilon)=\bar{C}_{l,l}(\epsilon)\;.
\end{equation}
\end{subequations}

\subsection{Again the coefficients of $x_1$}
Now, using (\ref{eq:sol}) together with (\ref{eq:xxx}) we are able to compute
the last coefficients. Notice that from (\ref{eq:abc}) the coefficients $b_l$ vanish
if $\epsilon=0$ (i.e.~in the left regular representation), and then from (\ref{eq:xxx})
$B^{m_2}_{j,j,l,l'}$ vanish too if $\epsilon=0$.
Moreover, from Lemma \ref{lemma:A} $B^{m_2}_{j,j,l,l'}$ vanish also if $|l-l'|>1$.
In the three cases $l'=l,l\pm 1$, using equation (\ref{eq:Cjllm}) we get:
\begin{align*}
B_{j,j,l,l+1}^{m_2} &=2|\epsilon|[2]q^{l+m_2+1+\epsilon(2j+1)}
  \frac{\sqrt{[l+2\epsilon\,j-m_2+2+\epsilon][l-2\epsilon\,j+m_2+2-\epsilon]}}{[2j][2j+2]}\, C_{l,l+1}(\epsilon) \;, \\
B_{j,j,l,l}^{m_2} &=\frac{2|\epsilon|}{[2j][2j+2]}\Big\{
                  [l-\epsilon(2j+1)-m_2+1][l-\epsilon(2j+1)+m_2+2]+ \\
&\phantom{=\frac{2|\epsilon|\,C_{l,l}(\epsilon)}{[2j][2j+2]}\, C_{l,l}(\epsilon)\Big\{}
                 -q^{-2}[l+\epsilon(2j+1)-m_2+2][l+\epsilon(2j+1)+m_2+1]\Big\} \\
         &=-\frac{2|\epsilon|}{[2j][2j+2]}\,
           \frac{q^{-2\epsilon(2j+1)}[2l+4]-q^{2\epsilon(2j+1)}[2l+2]-[2]q^{2m_2}}{1-q^2}\, C_{l,l}(\epsilon) \;, \\
B_{j,j,l,l-1}^{m_2} &=-2|\epsilon|[2]q^{-l+m_2-2-\epsilon(2j+1)}
  \frac{\sqrt{[l+2\epsilon\,j+m_2+1+\epsilon][l-2\epsilon\,j-m_2+1-\epsilon]}}{[2j][2j+2]}\, C_{l,l-1}(\epsilon) \;.
\end{align*}
We have inserted the factor $2|\epsilon|$, so that the expressions remain valid also
when $\epsilon=0$.

\subsection{The condition on the radius}
Orbits for $SO(5)$ are spheres of arbitrary radius, equivariance
alone  not imposing constraints on the radius. Similarly, for the quantum spheres 
one has to impose a constraint on the radius to determine the coefficients of the representation.
In fact, this will determine $C_{l,l+1}(0)$, $C_{l,l+1}(\frac{1}{2})$ and
$C_{l,l}(\frac{1}{2})$ only up to a phase.
Different choices of the phases correspond to unitary equivalent
representations and without losing generality we choose
$C_{l,l'}(\epsilon)\in\R$. A possible expression for the radius is 
$q^8x_0^2+q^4x_1^*x_1+x_2^*x_2$, and because of the defining relations
of $\A(S^4_q)$ this must be equal to 1. Let then, 
$$
r(l,m_1,m_2;j):=\inner{l,m_1,m_2;j|q^8x_0^2+q^4x_1^*x_1+x_2^*x_2|l,m_1,m_2;j}\;.
$$
All these matrix coefficients must be $1$. In particular,
for $l\in\N$ the condition $\,r(l,0,l;0)=1\,$ implies (up to a phase) that
\begin{equation}\label{eq:lint}
C_{l,l+1}(0)=\frac{q^{-l-3/2}}{\sqrt{[2l+3][2l+5]}}\;.
\end{equation}
For $l\in\N+\frac{1}{2}$ we first require that
$\,r(l,\frac{1}{2},l;\frac{1}{2})=r(l,-\frac{1}{2},l;\frac{1}{2})\,$
obtaining two possibilities:
$$
C_{l,l}(\tfrac{1}{2})=\pm\frac{[2]q^{l+2}}{[2l+2]}C_{l,l+1}(\tfrac{1}{2})\;.
$$
Then imposing $\,r(l,\tfrac{1}{2},l;\tfrac{1}{2})=1\,$, yields (up to a phase)
\begin{equation}
C_{l,l+1}(\tfrac{1}{2})=\frac{q^{-l-3/2}}{[2l+4]}\;, 
\end{equation}
hence,
\begin{equation}\label{eq:twocases}
C_{l,l}(\tfrac{1}{2})=\pm\frac{q^{1/2}[2]}{[2l+2][2l+4]}\;.
\end{equation}
With these, all the coefficients are completely determined.

\subsection{Explicit form of the representations}
Let us recall what we know on the equivariant representations of
the algebra $\A(S^4_q)$.

By the decomposition $\A(S^4_q)\simeq\bigoplus_{l\in\N}V_l$ into
irreducible representations of $U_q(so(5))$, there exists (at least)
one representation of $\A(S^4_q)\rtimes U_q(so(5))$ on the vector
space $\bigoplus_{l\in\N}V_l$ extending the representation
$\bigoplus_{l\in\N}\sigma_l$ of $U_q(so(5))$.
As we computed above, the equivariance uniquely
determines (for $l\in\N$, up to unitary equivalence) the matrix
coefficients of the representation, whose expression is
characterized by (\ref{eq:lint}).
On the other hand, by Proposition \ref{prop:main} there are
(at least) two inequivalent representations of $\A(S^4_q)\rtimes U_q(so(5))$ on the vector
space $\bigoplus_{l\in\N+\frac{1}{2}}V_l$ extending the representation
$\bigoplus_{l\in\N+\frac{1}{2}}\sigma_l$ of $U_q(so(5))$. These correspond,
by Lemma \ref{lemma:ptwo}, to the projective modules $\A(S^4_q)^4e$ and
$\A(S^4_q)^4(1-e)$, with $e$ the idempotent in equation (\ref{eq:P}).
The computation above (for $l\in\N+\frac{1}{2}$), which culminates
in equation (\ref{eq:twocases}), tells us that there are only
two possibilities for the matrix coefficients (up to unitary equivalence).
Therefore, the two possible choices in (\ref{eq:twocases}) must correspond
to the inequivalent representations associated with the projective
modules $\A(S^4_q)^4e$ and $\A(S^4_q)^4(1-e)$.

Let us summarize these results in the following
two theorems, which correspond to the scalar (i.e.~left regular)
and chiral spinor representations, respectively.

\begin{thm}\label{thm:lr}
The vector space
$\A(S^4_q)$ has orthonormal basis $\ket{l,m_1,m_2;j}$ with,
$$
l\in\N\;,\qquad 
j=0,1,\ldots,l\;,\qquad
j-|m_1|\in\N\;,\qquad
l-j-|m_2|\in 2\N\;.
$$
We call $L^2(S^4_q)$ the Hilbert space completion of $\A(S^4_q)$. Modulo a
unitary equivalence, the left regular representation is given by
\begin{align*}
x_0\ket{l,m_1,m_2;j} &=A_{j,m_1}C^+_{l,j,m_2}\ket{l+1,m_1,m_2;j+1} \\
                     &+A_{j,m_1}C^-_{l,j,m_2}\ket{l-1,m_1,m_2;j+1} \\
                     &+A_{j-1,m_1}C^-_{l+1,j-1,m_2}\ket{l+1,m_1,m_2;j-1} \\
                     &+A_{j-1,m_1}C^+_{l-1,j-1,m_2}\ket{l-1,m_1,m_2;j-1} \;, \\
x_1\ket{l,m_1,m_2;j} &=B^+_{j,m_1}C^+_{l,j,m_2}\ket{l+1,m_1+1,m_2;j+1} \\
                     &+B^+_{j,m_1}C^-_{l,j,m_2}\ket{l-1,m_1+1,m_2;j+1} \\
                     &+B^-_{j,m_1}C^-_{l+1,j-1,m_2}\ket{l+1,m_1+1,m_2;j-1} \\
                     &+B^-_{j,m_1}C^+_{l-1,j-1,m_2}\ket{l-1,m_1+1,m_2;j-1} \;, \\
x_2\ket{l,m_1,m_2;j} &=D_{l,j,m_2}^+\ket{l+1,m_1,m_2+1;j} \\
                     &+D_{l,j,m_2}^-\ket{l-1,m_1,m_2+1;j} \;,
\end{align*}
with coefficients
\begin{align*}
A_{j,m_1} &=q^{m_1-1}\sqrt{\frac{[j+m_1+1][j-m_1+1]}{[2j+1][2j+3]}} \;, \\
B^+_{j,m_1} &=q^{-j+m_1-1/2}\sqrt{\frac{[j+m_1+1][j+m_1+2]}{[2j+1][2j+3]}} \;, \\
B^-_{j,m_1} &=-q^{j+m_1+1/2}\sqrt{\frac{[j-m_1][j-m_1-1]}{[2j-1][2j+1]}} \;.
\end{align*}
and
\begin{align*}
C^+_{l,j,m_2} &=q^{m_2-1}\sqrt{\frac{[l+j+m_2+3][l+j-m_2+3]}{[2l+3][2l+5]}} \;, \\
C^-_{l,j,m_2} &=-q^{m_2-1}\sqrt{\frac{[l-j+m_2][l-j-m_2]}{[2l+1][2l+3]}} \;, \\
D_{l,j,m_2}^+ &=q^{-l+m_2-3/2}\sqrt{\frac{[l+j+m_2+3][l-j+m_2+2]}{[2l+3][2l+5]}} \;, \\
D_{l,j,m_2}^- &=q^{l+m_2+3/2}\sqrt{\frac{[l-j-m_2][l+j-m_2+1]}{[2l+1][2l+3]}} \;.
\end{align*}
\end{thm}

The two chiral spinorial representations (corresponding to the
sign $\pm$ in equation (\ref{eq:twocases})) are described in
the following Theorem.

\begin{thm}\label{thm:cr}
Let $\HH_\pm$ be two Hilbert spaces with orthonormal basis $\ket{l,m_1,m_2;j}_\pm$, where
$$
l\in\N+\tfrac{1}{2}\;,\qquad
j=\tfrac{1}{2},\tfrac{3}{2},\ldots,l\;,\qquad
j-|m_1|\in\N\;,\qquad
l+\tfrac{1}{2}-j-|m_2|\in\N\;.
$$
Let $\epsilon=\pm\frac{1}{2}$ be defined by $l+\epsilon-j-m_2\in 2\N$. On each space $\HH_\pm$
there is an equivariant $*$-representation of $\A(S^4_q)$ defined by:
\begin{align*}
x_0\ket{l,m_1,m_2;j}_\pm &=A^+_{j,m_1}C^+_{l,j,m_2}\ket{l+1,m_1,m_2;j+1}_\pm \\
                     &\mp A^+_{j,m_1}C^0_{l,j,m_2}\ket{l,m_1,m_2;j+1}_\pm \\
                     &+A^+_{j,m_1}C^-_{l,j,m_2}\ket{l-1,m_1,m_2;j+1}_\pm \\
                     &+A^0_{j,m_1}H^+_{l,j,m_2}\ket{l+1,m_1,m_2;j}_\pm \\
                     &\pm A^0_{j,m_1}H^0_{l,j,m_2}\ket{l,m_1,m_2;j}_\pm \\
                     &+A^0_{j,m_1}H^+_{l-1,j,m_2}\ket{l-1,m_1,m_2;j}_\pm \\
                     &+A^+_{j-1,m_1}C^-_{l+1,j-1,m_2}\ket{l+1,m_1,m_2;j-1}_\pm \\
                     &\mp A^+_{j-1,m_1}C^0_{l,j-1,m_2}\ket{l,m_1,m_2;j-1}_\pm \\
                     &+A^+_{j-1,m_1}C^+_{l-1,j-1,m_2}\ket{l-1,m_1,m_2;j-1}_\pm \;, \\
x_1\ket{l,m_1,m_2;j}_\pm &=B^+_{j,m_1}C^+_{l,j,m_2}\ket{l+1,m_1+1,m_2;j+1}_\pm \\
                     &\mp B^+_{j,m_1}C^0_{l,j,m_2}\ket{l,m_1+1,m_2;j+1}_\pm \\
                     &+B^+_{j,m_1}C^-_{l,j,m_2}\ket{l-1,m_1+1,m_2;j+1}_\pm \\
                     &+B^0_{j,m_1}H^+_{l,j,m_2}\ket{l+1,m_1+1,m_2;j}_\pm \\
                     &\pm B^0_{j,m_1}H^0_{l,j,m_2}\ket{l,m_1+1,m_2;j}_\pm \\
                     &+B^0_{j,m_1}H^+_{l-1,j,m_2}\ket{l-1,m_1+1,m_2;j}_\pm \\
                     &+B^-_{j,m_1}C^-_{l+1,j-1,m_2}\ket{l+1,m_1+1,m_2;j-1}_\pm \\
                     &\mp B^-_{j,m_1}C^0_{l,j-1,m_2}\ket{l,m_1+1,m_2;j-1}_\pm \\
                     &+B^-_{j,m_1}C^+_{l-1,j-1,m_2}\ket{l-1,m_1+1,m_2;j-1}_\pm \;, \\
x_2\ket{l,m_1,m_2;j}_\pm &=D_{l,j,m_2}^+\ket{l+1,m_1,m_2+1;j}_\pm \\
                     &\pm D_{l,j,m_2}^0\ket{l,m_1,m_2+1;j}_\pm \\
                     &+D_{l,j,m_2}^-\ket{l-1,m_1,m_2+1;j}_\pm \;,
\end{align*}
with coefficients
\begin{align*}
A^+_{j,m_1} &=q^{m_1-1}\frac{\sqrt{[j+m_1+1][j-m_1+1]}}{[2j+2]} \;, \\
A^0_{j,m_1} &=q^{-2}\frac{q^{j+m_1+1}[2][j-m_1]-[2j]}{[2j][2j+2]} \;, \\
B^+_{j,m_1} &=q^{-j+m_1-1/2}\frac{\sqrt{[j+m_1+1][j+m_1+2]}}{[2j+2]} \;, \\
B^0_{j,m_1} &=(1+q^2)q^{m_1-1/2}\frac{\sqrt{[j-m_1][j+m_1+1]}}{[2j][2j+2]} \;, \\
B^-_{j,m_1} &=-q^{j+m_1+1/2}\frac{\sqrt{[j-m_1][j-m_1-1]}}{[2j]} \;.
\end{align*}
and
\begin{align*}
C^+_{l,j,m_2} &=-q^{m_2-1-\epsilon}\frac{\sqrt{[l+j+m_2+3+\epsilon][l+j-m_2+3-\epsilon]}}{[2l+4]} \;, \\
C^0_{l,j,m_2} &=[4\epsilon]\,q^{2\epsilon\,l+m_2-1+3\epsilon}
    \frac{\sqrt{[l+\frac{1}{2}+j-2\epsilon\,m_2+2][l+\frac{1}{2}-j-2\epsilon\,m_2]}}
    {[2l+2][2l+4]} \;, \\
C^-_{l,j,m_2} &=-q^{m_2-1+\epsilon}\frac{\sqrt{[l-j+m_2-\epsilon][l-j-m_2+\epsilon]}}{[2l+2]} \;, \\
H^+_{l,j,m_2} &=q^{m_2-1+\epsilon(2j+1)}\frac{\sqrt{[l+2\epsilon j-m_2+2+\epsilon]
    [l-2\epsilon j+m_2+2-\epsilon]}}{[2l+4]} \;, \\
H^0_{l,j,m_2} &=\tfrac{[l-\epsilon(2j+1)-m_2+1][l-\epsilon(2j+1)+m_2+2]
       -q^{-2}[l+\epsilon(2j+1)-m_2+2][l+\epsilon(2j+1)+m_2+1]}{[2l+2][2l+4]} \;, \\
D_{l,j,m_2}^+ &=q^{-l+m_2-3/2}\frac{\sqrt{[l+j+m_2+3+\epsilon][l-j+m_2+2-\epsilon]}}{[2l+4]} \;, \\
D_{l,j,m_2}^0 &=[2]q^{m_2+1/2}\frac{\sqrt{[l-2\epsilon j-m_2+1-\epsilon]
    [l-2\epsilon j+m_2+2-\epsilon]}}{[2l+2][2l+4]} \;, \\
D_{l,j,m_2}^- &=-q^{l+m_2+3/2}\frac{\sqrt{[l-j-m_2+\epsilon][l+j-m_2+1-\epsilon]}}{[2l+2]} \;.
\end{align*}
These two representations are inequivalent and correspond to the
projective modules $\A(S^4_q)^4e$ and $\A(S^4_q)^4(1-e)$, with $e$ the
idempotent in equation (\ref{eq:P}).
\end{thm}

\section{The Dirac operator on the quantum orthogonal $4$-sphere}\label{sec:7.5}
We start by constructing a non-trivial Fredholm module on the
quantum orthogonal sphere. With different representations a non-trivial
Fredholm module was already constructed in~\cite{HL04}.

\begin{prop}\label{prop:notr}
Consider the representations of $\A(S^4_q)$ on $\HH_\pm$ given in Theorem \ref{thm:cr}.
Then, the datum $(\A(S^4_q),\HH,F,\gamma)$ is a $1$-summable even Fredholm module,
where $\HH:=\HH_+\oplus\HH_-$, $\gamma$ is the natural grading and $F\in\B(\HH)$ is defined by
$$
F\ket{l,m_1,m_2;j}_\pm:=\ket{l,m_1,m_2;j}_\mp\;.
$$
This Fredholm module is non-trivial. In particular,
\begin{equation}\label{eq:pairS4q}
\mr{Index}(F_e^+)=
\tfrac{1}{2}\tr_{\HH\otimes\C^4}(\gamma F[F,e])=1\;,
\end{equation}
with $e$ the idempotent defined by equation (\ref{eq:P}).
\end{prop}
\begin{prova}
That $F=F^*$, $F^2=1$ and $\gamma F+F\gamma=0$ is obvious. Then, it is enough to show
that $[F,x_i]\in\mc{L}^1(\HH)$ for $i=0,1,2$. From this and Leibniz rule
it follows that $[F,a]$ is trace class, and then compact, for all $a\in\A(S^4_q)$.

Now, notice that
\begin{align}
[F,x_0]\ket{l,m_1,m_2;j}_\pm=\;&\mp 2A^+_{j,m_1}C^0_{l,j,m_2} \ket{l,m_1,m_2;j+1}_\mp \nonumber \\
                               &\pm 2A^0_{j,m_1}H^0_{l,j,m_2} \ket{l,m_1,m_2;j}_\mp \label{eq:Fx0} \\
                               &\mp 2A^+_{j-1,m_1}C^0_{l,j-1,m_2} \ket{l,m_1,m_2;j-1}_\mp \nonumber \;, \\
[F,x_1]\ket{l,m_1,m_2;j}_\pm=\;&\mp 2B^+_{j,m_1}C^0_{l,j,m_2} \ket{l,m_1+1,m_2;j+1}_\mp \nonumber \\
                               &\pm 2B^0_{j,m_1}H^0_{l,j,m_2} \ket{l,m_1+1,m_2;j}_\mp \nonumber \\
                               &\mp 2B^-_{j,m_1}C^0_{l,j-1,m_2} \ket{l,m_1+1,m_2;j-1}_\mp \;, \nonumber \\
[F,x_2]\ket{l,m_1,m_2;j}_\pm=\;&\pm 2D^0_{l,j,m_2}\ket{l,m_1,m_2+1;j}_\mp \;. \nonumber
\end{align}
All the coefficients appearing in these equations are bounded by $q^{2l}$.
Thus the commutators are trace class and this concludes the first part of the proof.

To prove non-triviality it is enough to prove \eqref{eq:pairS4q}.
Substituting \eqref{eq:P} into \eqref{eq:pairS4q} yields
$$
\mr{Index}(F_e^+)
=\tfrac{(1-q^2)^2}{4}\tr_{\HH}(\gamma F[F,x_0])\;.
$$
and in turn, using equation (\ref{eq:Fx0}),
$$
\mr{Index}(F_e^+)
=(1-q^2)^2\sum_{l,j,m_1,m_2}A^0_{j,m_1}H^0_{l,j,m_2} \;.
$$
Summing over $m_1$ from $-j$ to $j$ we obtain that
\begin{align*}
\mr{Index}(F_e^+)
        &=q^{-3}(1-q^2)^2\sum_{l,j,m_2}\frac{[l+\epsilon(2j+1)-m_2+2]
        [l+\epsilon(2j+1)+m_2+1]}{[2l+2][2l+4][2j][2j+2]}\times \\ &
        \qquad\qquad \qquad\qquad \qquad\qquad\qquad\qquad\qquad
        \times
        \sum_{m_1}\big\{q^{2j+2}+q^{-2j}-[2]q^{2m_1+1}\big\} \\
        &=\sum_{l,j}\frac{(2j+1)(q^{2j+1}+q^{-2j-1})-[2][2j+1]}
        {[2l+2][2l+4][2j][2j+2]} \times \\ &
        \qquad\qquad \qquad\qquad
        \times\sum_{m_2}\Big\{q^{2l+2\epsilon(2j+1)+3}+q^{-2l-2\epsilon(2j+1)-3}
        -q^{2m_2-1}-q^{-2m_2+1}\Big\}\;.
\end{align*}

The sum over $m_2$ requires additional care.
For $\epsilon$ fixed, $l-\epsilon-j+m_2=0,2,4,\ldots,2(l-j)$. If we call
$2i:=l-\epsilon-j+m_2$ and sum first over $i=0,1,\ldots,l-j$ and then over
$\epsilon=\pm 1/2$ we get:
\begin{align*}
& \mr{Index}(F_e^+)
           =\sum_{l,j}\frac{(2j+1)(q^{2j+1}+q^{-2j-1})-[2][2j+1]}
           {[2l+2][2l+4][2j][2j+2]} \times \\ & 
           \qquad\times
           \sum_{2\epsilon=\pm 1}\Big\{(l-j+1)(q^{2l+2\epsilon(2j+1)+3}+q^{-2l-2\epsilon(2j+1)-3})
           -(q^{2\epsilon-1}+q^{-2\epsilon+1})[2]^{-1}[2(l-j+1)]\Big\} \\
           &\qquad=\sum_{l,j}\frac{(2j+1)(q^{2j+1}+q^{-2j-1})-[2][2j+1]}
           {[2l+2][2l+4][2j][2j+2]} \times \\ &
           \qquad\qquad\qquad\qquad\times
           \Big\{(l-j+1)(q^{2l+3}+q^{-2l-3})(q^{2j+1}+q^{-2j-1})
           -[2][2(l-j+1)]\Big\} \\
           &\qquad=:\sum_{l,j}f_{lj}(q)=:f(q) \;.
\end{align*}
We call $f_{lj}(q)$ the generic term of last series, explicitly written as
\begin{align*}
&f_{lj}(q)=(1-q^2)^4\,\frac{(2j+1)(1+q^{4j+2})-\frac{1+q^2}{1-q^2}(1-q^{4j+2})}
           {(1-q^{4l+4})(1-q^{4l+8})(1-q^{4j})(1-q^{4j+4})}
  \times \\ & \qquad\qquad\qquad\qquad\times
           q^{2l-1}\Big\{(l-j+1)(1+q^{4l+6})(1+q^{4j+2})
           -\tfrac{1+q^2}{1-q^2}\,q^2(q^{4j}-q^{4l+4})\Big\}\;,
\end{align*}
and consider it as a function of $q\in [\hspace{1pt}0,1[\,$. Notice that each
$f_{lj}(q)$ is a $C^\infty$ function of $q$ (they are rational functions whose
denominators never vanish for $0\leq q<1$). From the inequality
$$
0\leq f_{lj}(q)\leq 4(2j+1)q^{2l-1}
$$
we deduce (using the Weierstrass M-test~\cite[pag.~301]{Arf85}) that the series
is absolutely (hence uniformly) convergent in each interval $[0,q_0]\subset
[\hspace{1pt}0,1[\,$. Then, it converges to a function $f(q)$ which is continuous
in $[\hspace{1pt}0,1[\,$. 
Being the index of a Fredholm operator, $f(q)$ is integer valued in
$]\hspace{1pt}0,1[$; by continuity it is constant and can be computed in the limit
$q\to 0$. In this limit we have $f_{lj}(q)=2j(l-j+1)q^{2l-1}+O(q^{2l})$.
Thus, $f_{lj}(0)=\delta_{l,1/2}\delta_{j,1/2}$ and $\mr{Index}(F_e^+)=f(0)=1$.
\end{prova}

The next step is to define a spectral triple whose Fredholm module
is the one described in Proposition~\ref{prop:notr}.

\begin{prop}\label{prop:st}
Let $D$ be the (unbounded) operator on $\HH:=\HH_+\oplus\HH_-$ defined by
$$
D\ket{l,m_1,m_2;j}_\pm:=(l+\tfrac{3}{2})\ket{l,m_1,m_2;j}_\mp\;.
$$
Then, the datum $(\A(S^4_q),\HH,D,\gamma)$ is a $U_q(so(5))$-equivariant regular even
spectral triple of \mbox{metric} dimension $4$.
\end{prop}

\noindent{\bf Remark:} The operator $D$ is isospectral to the classical Dirac
operator on $S^4$ (whose spectrum has been computed in~\cite{Tra93,CH96}).
When $q=1$, this spectral triple becomes the canonical one associated to the spin
structure of $S^4$.

\begin{prova}
Clearly the representation of the algebra is even, $D$ is odd, with compact
resolvent and $4^+$-summable (being isospectral to the classical Dirac operator
on $S^4$). To prove that $D$ defines a regular spectral triple, one repeats
verbatim the proof of Lemma \ref{lemma:5.2}.
Finally, since $D$ is proportional to the identity in any irreducible subrepresentation
$V_l$ of $U_q(so(5))$, it commutes with all $h\in U_q(so(5))$ and it is equivariant.
\end{prova}

As a preparation for the study of the dimension spectrum in Section~\ref{sec:7.6},
let us explicitly verify the $4$-summability of $D$. As one can easily check, the
dimension of $V_l$ is~\cite{CH96}
$$
\dim V_l=\tfrac{2}{3}(l+\tfrac{5}{2})(l+\tfrac{3}{2})(l+\tfrac{1}{2}) \;.
$$
From this we get
$$
\tr(|D|^{-s})=\sum\nolimits_{l\in\N+\frac{1}{2}}2(l+\tfrac{3}{2})^{-s}\dim V_l
             =\tfrac{4}{3}\sum\nolimits_{n=1}^\infty (n^2-1)n^{-s+1} \;,
$$
where $n=l+\frac{3}{2}$ (and we added the term with $n=1$ since it is
identically zero). The above series is convergent in the right half-plane
$\{s\in\C\,|\,\mr{Re}\,s>4\}$, thus $D$ has metric dimension $4$.

Notice that $\tr(|D|^{-s})$ has meromorphic extension on $\C$ given by
\begin{equation}\label{eq:zeta}
\zeta_1(s)=\tfrac{4}{3}\big\{\zeta(s-3)-\zeta(s-1)\big\}\;,
\end{equation}
where $\zeta(s)$ is the Riemann zeta-function.
We recall that $\zeta(s)$ has a simple pole in $s=1$ as unique singularity
and that $\,\mr{Res}_{s=1}\zeta(s)=1$. 

\section{The dimension spectrum and residues}\label{sec:7.6}

Let $\mc{J}_2$ be the two-sided ideal in $\B(\HH)$ defined in (\ref{eq:Jr}), with $r=2$.
We recall here the definition:
$$
\mc{J}_2:=\big\{T\in\opz\,\big|\,T|D|^{-p}\in\mc{L}^1(\HH)\;\forall\;p>2\big\}\;.
$$
Clearly, if $T$ is of trace class, so is $T|D|^{-p}$ for any positive $p$,
and $\mc{L}^1(\HH)\subset\mc{J}_2$. Since $\op\subset\mc{L}^1(\HH)$,
smoothing operators belongs to $\mc{J}_2$ as well. On the other hand, $\mc{J}_2$
is strictly bigger than $\mc{L}^1(\HH)$; indeed, 
the operator $L_q\in\B(\HH)$, given by
$$
L_q\ket{l,m_1,m_2;j}_\pm:=q^{j+\frac{1}{2}}\ket{l,m_1,m_2;j}_\pm\;,
$$
is not of trace class but belongs to $\mc{J}_2$, by the following proposition.
\begin{prop}\label{prop:jei}
For any $s\in\C$ with $\mr{Re}\,s>2$ one has that
$$
\zeta_{L_q}(s):=\sum_{l,j,m_1,m_2}(l+\tfrac{3}{2})^{-s}q^{j+\frac{1}{2}}=\frac{4q}{(1-q)^2}\,
\Big(\zeta(s-1)-\frac{1+q}{1-q}\,\zeta(s)\Big)+\textup{holomorphic function}\;,
$$
where $\zeta(s)$ is the Riemann zeta-function.
In particular, this means that $L_q$ belongs to the ideal $\mc{J}_2$.
Furthermore, since the series $\zeta_{L_q}(0)$ is divergent,
$L_q$ is not of trace class.
\end{prop}
\begin{prova}
Calling $n:=l+\frac{3}{2}$, $k:=j+\frac{1}{2}$, we have
$$
\zeta_{L_q}(s)=4\sum_{n=2}^\infty n^{-s}\sum_{k=1}^{n-1}k(n-k)q^k\;.
$$
We can sum starting
from $n=1$ and for $k=0,\ldots,n$ (we simply add zero terms) to get
$$
\zeta_{L_q}(s)=4\sum_{n=1}^\infty n^{-s}\big\{nq\partial_q-(q\partial_q)^2\big\}\sum_{k=0}^n q^k
    =4\sum_{n=1}^\infty n^{-s}\big\{nq\partial_q-(q\partial_q)^2\big\}\frac{1-q^{n+1}}{1-q}\;.
$$
Terms decaying as $q^n$ give a holomorphic function of $s$, thus modulo
holomorphic functions,
$$
\zeta_{L_q}(s)\sim 4\sum_{n=1}^\infty n^{-s}\left\{n\,\tfrac{q}{(1-q)^2}
-\tfrac{q(1+q)}{(1-q)^3}\right\} \;.
$$
Last series is summable for all $s$ with $\mr{Re}\,s>2$, and its sum
can be written in terms of the Riemann zeta-function as in the statement
of the proposition.
\end{prova}

\subsection{An approximated representation}\label{sec:HH}
Let $\hat{\HH}$ be a Hilbert space with orthonormal basis $\kkett{l,m_1,m_2;j}_\pm$
labelled by,
$$
l\in\tfrac{1}{2}\Z\;,\qquad
l+j\in\Z\;,\qquad
j+m_1\in\N\;,\qquad
l+\tfrac{1}{2}-j+m_2\in\N\;.
$$
Let $I$ be the labelling set of the Hilbert space $\HH_\pm$ as in Theorem
\ref{thm:cr}, and given by
$$
I:=\big\{(l,m_1,m_2,j)\;\big|\;l\in\N+\tfrac{1}{2}\;,\;\;j=\tfrac{1}{2},\tfrac{3}{2},...,l\;,
\;\;j-|m_1|\in\N\;,\;\;l+\tfrac{1}{2}-j-|m_2|\in\N\big\}\;.
$$
Notice that $I$ is the subset of labels of $\hat{\HH}$ satisfying
$l\in\N+\tfrac{1}{2}$, $m_1\leq j\leq l$ and $m_2\leq l+\tfrac{1}{2}-j$.
Define the inclusion $Q:\HH\to\hat{\HH}$ and the adjoint projection 
$P:\hat{\HH}\to\HH$ by,
\begin{align*}
Q\ket{l,m_1,m_2;j}_\pm &:=\kkett{l,m_1,m_2;j}_\pm 
\;\quad\textrm{for all}\;(l,m_1,m_2,j)\in I\;, \\
P\kkett{l,m_1,m_2;j}_\pm &:=\bigg\{\begin{array}{ll}
\ket{l,m_1,m_2;j}_\pm &\mr{if}\;(l,m_1,m_2,j)\in I\;, \\
0 &\mr{otherwise}\;.
\end{array}
\end{align*}
Clearly, $PQ=id_{\HH}$.
The Hilbert space $\hat{\HH}$ carries a bounded $*$-representation of the algebra $\A(SU_q(2))\otimes\A(S^2_{q^2})$
defined by,
\begin{align*}
\alpha\kkett{l,m_1,m_2;j}_\pm &=\sqrt{1-q^{2(j+m_1+1)}}\kkett{l+\oh,m_1+\oh,m_2;j+\oh}_\pm \;,\\
\beta\kkett{l,m_1,m_2;j}_\pm  &=q^{j+m_1}\kkett{l+\oh,m_1-\oh,m_2;j+\oh}_\pm \;,\\
A\kkett{l,m_1,m_2;j}_\pm &=q^{l-j+m_2-\epsilon}\kkett{l,m_1,m_2;j}_\pm \;,\\
B\kkett{l,m_1,m_2;j}_\pm &=\sqrt{1-q^{2(l-j+m_2+2-\epsilon)}}\kkett{l+1,m_1,m_2+1;j}_\pm \;,
\end{align*}
where, as before, $\epsilon:=\frac{1}{2}(-1)^{l+\frac{1}{2}-j-m_2}$.
Composition of such a representation with the algebra embedding
$\A(S^4_q)\hookrightarrow\A(SU_q(2))\otimes\A(S^2_{q^2})$ given in equation
(\ref{eq:embed}) results into a $*$-representation $\pi:\A(S^4_q)\to\B(\hat{\HH})$.
The sandwich $\tilde{\pi}(a):=P\pi(a)Q$ defines a $*$-linear map
$\tilde{\pi}:\A(S^4_q)\to\B(\HH)$.

\begin{prop}\label{prop:app}
With $\mc{J}_2$ the class of operators defined in equation (\ref{eq:Jr}), one has that
the difference $a-\tilde{\pi}(a)\in\mc{J}_2$ for all $a\in\A(S^4_q)$.
\end{prop}
\begin{prova}
Define $\hat{\mc{J}}_2$ as the collection of bounded operators $T:\HH\to\hat{\HH}$ such that
$T|D|^{-p}$ is trace class for all $p>2$. Since trace class operators are a two sided ideal
in bounded operators, the space $\hat{\mc{J}}_2$ is stable when multiplied from the left
by bounded operators: $T_1\in\B(\hat{\HH})$ and $T_2\in\hat{\mc{J}}_2\;\Rightarrow\;T_1T_2\in\hat{\mc{J}}_2$.

Next, suppose that $a,b$ satisfy $a-\tilde{\pi}(a)\in\mc{J}_2$ and $b-\tilde{\pi}(b)\in\mc{J}_2$
and consider the following algebraic identity
$$
ab-\tilde{\pi}(ab)=\bigl\{a-\tilde{\pi}(a)\bigr\}b+P\pi(a)\bigl\{Qb-\pi(b)Q\bigr\}\;.
$$
Since $\mc{J}_2$ is a two-sided ideal in $\opz$, the first summand
is in $\mc{J}_2$. The stability of $\hat{\mc{J}}_2$ discussed above implies
that $\pi(a)\bigl\{Qb-\pi(b)Q\bigr\}\in\hat{\mc{J}}_2$,
but if $T\in\hat{\mc{J}}_2$ clearly $PT\in\mc{J}_2$. Hence the second
summand lies in $\mc{J}_2$ too. Thus, $ab-\tilde{\pi}(ab)\in\mc{J}_2$ whenever
this property holds for each of the operators $a,b$.
We conclude that it is enough to show that $a-\tilde{\pi}(a)\in\mc{J}_2$
when $a$ is a generator of $\A(S^4_q)$. By Proposition \ref{prop:jei},
this amounts to prove that the matrix elements of $a-\tilde{\pi}(a)$ 
are bounded in modulus by $q^{j+\frac{1}{2}}$.

Let us have a close look at the coefficients of $a\in\{x_i,x_i^*\}$ in Theorem \ref{thm:cr}.
Firstly, $A^+_{j,m_1}$, $B^+_{j,m_1}$, $B^-_{j,m_1}$, $q^{-2j}A^0_{j,m_1}$
and $q^{-2j}B^0_{j,m_1}$ are uniformly bounded by a constant, as one can see
by writing explicitly the $q$-analogues in their expressions so getting:
\begin{align*}
A^+_{j,m_1} &=q^{j+m_1}\frac{\sqrt{(1-q^{2(j+m_1+1)})(1-q^{2(j-m_1+1)})}}{1-q^{4j+4}} \;, \\
q^{-2j}A^0_{j,m_1} &=(1-q^2)
\frac{[2]q^{2(j+m_1)}-q^{4j+1}-q^{-1}}{(1-q^{4j})(1-q^{4j+4})} \;, \\
B^+_{j,m_1} &=\frac{\sqrt{(1-q^{2(j+m_1+1)})(1-q^{2(j+m_1+2)})}}{1-q^{4j+4}}\;, \\
q^{-2j}B^0_{j,m_1} &=(1+q^2)q^{j+m_1+1}\frac{\sqrt{(1-q^{2(j-m_1)})
  (1-q^{2(j+m_1+1)})}}{(1-q^{4j})(1-q^{4j+4})} \;, \\
B^-_{j,m_1} &=-q^{2(j+m_1)+1}\frac{\sqrt{(1-q^{2(j-m_1)})(1-q^{2(j-m_1-1)})}}{1-q^{4j}} \;.
\end{align*}
Analogously, the coefficients
$q^{2j}H^0_{l,j,m_2}$, $C_{l,j,m_2}^0$ and $D_{l,j,m_2}^0$ are seen to be bounded by $q^l$ . Thus, modulo rapid decay matrices (i.e.~smoothing
operators),
\begin{subequations}\label{eq:apprep}
\begin{align}
x_0\ket{l,m_1,m_2;j} &\simeq A^+_{j,m_1}C^+_{l,j,m_2}\ket{l+1,m_1,m_2;j+1} \notag\\
                     &+A^+_{j,m_1}C^-_{l,j,m_2}\ket{l-1,m_1,m_2;j+1} \notag\\
                     &+A^0_{j,m_1}H^+_{l,j,m_2}\ket{l+1,m_1,m_2;j} \notag\\
                     &+A^0_{j,m_1}H^+_{l-1,j,m_2}\ket{l-1,m_1,m_2;j} \notag\\
                     &+A^+_{j-1,m_1}C^-_{l+1,j-1,m_2}\ket{l+1,m_1,m_2;j-1} \notag\\
                     &+A^+_{j-1,m_1}C^+_{l-1,j-1,m_2}\ket{l-1,m_1,m_2;j-1} \;, \\
x_1\ket{l,m_1,m_2;j} &\simeq B^+_{j,m_1}C^+_{l,j,m_2}\ket{l+1,m_1+1,m_2;j+1} \notag\\
                     &+B^+_{j,m_1}C^-_{l,j,m_2}\ket{l-1,m_1+1,m_2;j+1} \notag\\
                     &+B^0_{j,m_1}H^+_{l,j,m_2}\ket{l+1,m_1+1,m_2;j} \notag\\
                     &+B^0_{j,m_1}H^+_{l-1,j,m_2}\ket{l-1,m_1+1,m_2;j} \notag\\
                     &+B^-_{j,m_1}C^-_{l+1,j-1,m_2}\ket{l+1,m_1+1,m_2;j-1} \notag\\
                     &+B^-_{j,m_1}C^+_{l-1,j-1,m_2}\ket{l-1,m_1+1,m_2;j-1} \;, \\
x_2\ket{l,m_1,m_2;j} &\simeq D_{l,j,m_2}^+\ket{l+1,m_1,m_2+1;j} \notag\\
                     &+D_{l,j,m_2}^-\ket{l-1,m_1,m_2+1;j} \;.
\end{align}
\end{subequations}
Since modulo smoothing operators the representations are the same we are omitting the label `$\pm$' in the vector basis.
Furthermore, using the inequalities
\begin{equation}\label{eq:ineqJ}
0\leq (1-qu)^{-1}-1\leq q (1-q)^{-1} \,u\;,\qquad
0\leq 1-(1-u)^{\frac{1}{2}}\leq u\;,
\end{equation}
which are valid when $0\leq u\leq 1$,
we prove that modulo terms bounded by $q^l$, one has
\begin{subequations}\label{eq:coefsim}
\begin{align}
C^+_{l,j,m_2} &\simeq -q^{l-j+m_2-\epsilon}\sqrt{1-q^{2(l+j+m_2+3+\epsilon)}} \;, \\
C^-_{l,j,m_2} &\simeq -q^{l+j+m_2+1+\epsilon}\sqrt{1-q^{2(l-j+m_2-\epsilon)}} \;, \\
H^+_{l,j,m_2} &\simeq q^{l+m_2+1}\sqrt{q^{2\epsilon(2j+1)}-q^{2(l+m_2+2)}} \;, \\
D_{l,j,m_2}^+ &\simeq \sqrt{1-q^{2(l+j+m_2+3+\epsilon)}}\sqrt{1-q^{2(l-j+m_2+2-\epsilon)}} \;, \\
D_{l,j,m_2}^- &\simeq -q^{2(l+m_2)+3} \;.
\end{align}
\end{subequations}
Up to now, we neglected only smoothing contributions (the above approximation
will be needed when dealing with the real structure later on). We use again
(\ref{eq:ineqJ}) to get a rougher approximation by neglecting terms bounded by $q^j$. This yields
\begin{subequations}\label{eq:get}
\begin{align}
A^+_{j,m_1} &\simeq \tilde{A}^+_{j,m_1}:=q^{j+m_1}\sqrt{1-q^{2(j+m_1+1)}} \;, \\
A^0_{j,m_1}H^+_{l,j,m_2}  &\simeq 0 \;, \\
B^+_{j,m_1} &\simeq \tilde{B}^+_{j,m_1}:=\sqrt{(1-q^{2\smash[t]{(j+m_1+1)}})(1-q^{2\smash[t]{(j+m_1+2)}})} \;, \\
B^0_{j,m_1}H^+_{l,j,m_2}  &\simeq 0 \;, \\
B^-_{j,m_1} &\simeq \tilde{B}^-_{j,m_1}:=-q^{2(j+m_1)+1} \;,\\
C^+_{l,j,m_2} &\simeq \tilde{C}_{l,j,m_2}:=-q^{l-j+m_2-\epsilon} \;, \\
C^-_{l,j,m_2} &\simeq 0 \;, \\
D_{l,j,m_2}^+ &\simeq \tilde{D}_{l,j,m_2}:=\sqrt{1-q^{2(l-j+m_2+2-\epsilon)}} \;, \\
D_{l,j,m_2}^- &\simeq 0 \;.
\end{align}
\end{subequations}
Plugging these coefficients in the equations for the $x_i$'s we see that,
modulo operators in the ideal $\mc{J}_2$, we get
\begin{align*}
x_0\ket{l,m_1,m_2;j} &\simeq \tilde{A}^+_{j,m_1}\tilde{C}_{l,j,m_2}\ket{l+1,m_1,m_2;j+1} \\
                     &\quad +\tilde{A}^+_{j-1,m_1}\tilde{C}_{l-1,j-1,m_2}\ket{l-1,m_1,m_2;j-1} \\
                     &=-P(\alpha\beta A)Q\ket{l,m_1,m_2;j} - P(\beta^*\alpha^* A)Q\ket{l,m_1,m_2;j} \;, \\
x_1\ket{l,m_1,m_2;j} &\simeq \tilde{B}^+_{j,m_1}\tilde{C}_{l,j,m_2}\ket{l+1,m_1+1,m_2;j+1} \\
                     &\quad +\tilde{B}^-_{j,m_1}\tilde{C}_{l-1,j-1,m_2}\ket{l-1,m_1+1,m_2;j-1} \\
                     &=-P(\alpha^2 A)Q\ket{l,m_1,m_2;j} + P(q(\beta^*)^2 A)Q\ket{l,m_1,m_2;j} \;, \\
x_2\ket{l,m_1,m_2;j} &\simeq \tilde{D}_{l,j,m_2}\ket{l+1,m_1,m_2+1;j} \\
                     &=PB\,Q \ket{l+1,m_1,m_2+1;j} \;.
\end{align*}
The observation that
$$
-P(\alpha\beta+\beta^*\alpha^*)AQ=\tilde{\pi}(x_0)\;,\qquad
P(-\alpha^2+q(\beta^*)^2)AQ=\tilde{\pi}(x_1)\;,\qquad
PBQ=\tilde{\pi}(x_2)\;,
$$
concludes the proof.
\end{prova}

\subsection{The dimension spectrum and the top residue}
The approximation modulo $\mc{J}_2$ allows considerable simplifications when getting information on the part of
the dimension spectrum contained in the half plane $\mr{Re}\,s>2$.
To study the part of the dimension spectrum in the left half plane
$\mr{Re}\,s\leq 2$ would require a less drastic approximation which we are lacking at the moment.

\begin{prop}
In the region $\mr{Re}\,s>2$ the dimension spectrum $\Sigma$ of the spectral triple $(\A(S^4_q),\HH,D,\gamma)$ given in Proposition~\ref{prop:st}
consists of the
two points $\{3,4\}$, which are simple poles of the zeta-functions. The top residue
coincides with the integral on the subspace of classical points of $S^4_q$, that is
\begin{equation}\label{eq:topres}
\nint a|D|^{-4}=\frac{2}{3\pi}\int_0^{2\pi}\!\!\sigma(a)(\theta)\de\theta\;,
\end{equation}
with $\sigma:\A(S^4_q)\to\A(S^1)$ the $*$-algebra morphism defined by
$\sigma(x_0)=\sigma(x_1)=0$ and $\sigma(x_2)=u$, where $u$, given by
$u(\theta):=e^{i\theta}$, is the unitary generator of $\A(S^1)$.
\end{prop}
\begin{prova}
Let $\Psi^0$ be the $*$-algebra generated by $\A(S^4_q)$, by $[D,\A(S^4_q)]$
and by iterated applications of the derivation $\delta$. Let $\mathfrak{A}\subset
\A(SU_q(2))\otimes\A(S^2_{q^2})\otimes\mr{Mat}_2(\C)$ be the $*$-algebra generated
by $\alpha$, $\beta$, $\alpha^*$, $\beta^*$, $A$, $B$, $B^*$ and $F$.
By Proposition~\ref{prop:app} there is an inclusion $\A(S^4_q)\subset P\mathfrak{A}Q+\mc{J}_2$.

A linear basis for $\mathfrak{A}$ is given by,
\begin{equation}\label{eq:linb}
T:=\alpha^{k_1}\beta^{n_1}(\beta^*)^{n_2}A^{n_3}B^{k_2}F^h\;,
\end{equation}
where $h\in\{0,1\}$, $n_i\in\N$, $k_i\in\Z$ and
with the notation $\alpha^{k_1}:=(\alpha^*)^{-k_1}$ if $k_1<0$ and
$B^{k_2}:=(B^*)^{-k_2}$ if $k_2<0$. For this operator,
$$
\delta(PTQ)=\big(\tfrac{1}{2}(k_1+n_1-n_2)+k_2\big)PTQ\qquad\mr{and}\qquad
[D,PTQ]=\delta(PTQ)F\;.
$$
Thus, $P\mathfrak{A}Q$ is invariant under application of $\delta$ and $[D,\,.\,]$
and hence $\Psi^0\subset P\mathfrak{A}Q+\mc{J}_2$.

For the part of the dimension spectrum in the right half
plane $\mr{Re}\,s>2$, we can neglect $\mc{J}_2$ and consider only the singularities
of zeta-functions associated with elements in $P\mathfrak{A}Q$. By linearity of
the zeta-functions, it is enough to consider the generic basis element in
equation (\ref{eq:linb}).

Such a $T$ shifts $l$ by $\frac{1}{2}(k_1+n_1-n_2)+k_2$, $m_1$ by $\frac{1}{2}(k_1-n_1+n_2)$,
$m_2$ by $k_2$, $j$ by $\frac{1}{2}(k_1+n_1-n_2)$, and flips the chirality if $h=1$.
Thus it is off-diagonal unless $h=k_i=0$ and $n_1=n_2$. The zeta-function associated with
a bounded off-diagonal operator is identically zero in the half-plane $\mr{Re}\,z>4$,
and so is its holomorphic extension to the entire complex plane. It remains to consider the cases $T=P(\beta\beta^*)^kA^nQ$, with $n,k\in\N$.

If $n$ and $k$ are both different from zero, one finds
$$
\zeta_T(s)=2\sum_{l,j,m_1,m_2} (l+\tfrac{3}{2})^{-s}q^{n(l-j+m_2-\epsilon)+2k(j+m_1)}
 =2\sum_{l,j,m_2} (l+\tfrac{3}{2})^{-s}q^{n(l-j+m_2-\epsilon)}\frac{1-q^{2k(2j+1)}}{1-q^{2k}}\;.
$$
For $\epsilon$ fixed, set $2i:=l-\epsilon-j+m_2=0,2,\ldots,2(l-j)$. Then,
\begin{align*}
\zeta_T(s)&=2\sum_{l,j} (l+\tfrac{3}{2})^{-s}\frac{1-q^{2k(2j+1)}}{1-q^{2k}}
\sum_{\epsilon=\pm 1/2}\sum_{i=0}^{l-j}q^{2ni}=
4\sum_{l,j} (l+\tfrac{3}{2})^{-s}\frac{1-q^{2k(2j+1)}}{1-q^{2k}}
\frac{1-q^{2n(l-j+1)}}{1-q^{2n}} \\
&=4\zeta(s-1)-4\,\frac{1+(1-q^{4k})^{-1}+(1-q^{2n})^{-1}}{(1-q^{2k})(1-q^{2n})}\,\zeta(s)
  +\textrm{holomorphic function}\;, 
\end{align*}
which has meromorphic extension on $\C$ with simple pole in $s=\{1,2\}$.

If $n=0$ and $k\neq 0$,
\begin{align*}
\zeta_T(s) &=4\sum_{l,j}(l+\tfrac{3}{2})^{-s}(l-j+1)\frac{1-q^{2k(2j+1)}}{1-q^{2k}} \\
  &=\tfrac{4}{1-q^{2k}}\left(\tfrac{1}{2}\,\zeta(s-2)-\big(\tfrac{1}{2}+\tfrac{1}{1-q^{4k}}\big)
  \zeta(s-1)+\tfrac{q^{4k}}{(1-q^{4k})^2\log q^{4k}}\,\zeta(s)\right)
  +\textrm{hol.~function}\;,
\end{align*}
which has meromorphic extension on $\C$ with simple pole in $s=\{1,2,3\}$.

If $n\neq 0$ and $k=0$,
\begin{align*}
\zeta_T(s)&=4\sum_{l,j} (l+\tfrac{3}{2})^{-s}(2j+1)\frac{1-q^{2n(l-j+1)}}{1-q^{2n}} \\
  &=\tfrac{4}{1-q^{2n}}\left\{\zeta(s-2)-\Big(1+\tfrac{2q^{2n}}{1-q^{2n}}\Big)\zeta(s-1)
  +\tfrac{2q^{2n}}{1-q^{2n}}\Big(1+\tfrac{q^{2n}}{(1-q^{2n})\log q^{2n}}\Big)\zeta(s)\right\}
  +\textrm{hol.~fun.}\;, 
\end{align*}
which has meromorphic extension on $\C$ with simple pole in $s=\{1,2,3\}$.

Finally, if both $n$ and $k$ are zero we get (cf.~equation (\ref{eq:zeta})),
$
\zeta_T(s)=\tfrac{4}{3}\big\{\zeta(s-3)-\zeta(s-1)\big\}
$,
and this is meromorphic with simple poles in $\{2,4\}$. Thus, the part of the
dimension spectrum in the region $\mr{Re}\,s>2$ consists at
most of the two points $\{3,4\}$ and both are simple poles.

Since we have considered the enlarged algebra $P\mathfrak{A}Q+\mc{J}_2$, it suffices
to prove that there exists an $a\in\Psi^0$ whose zeta-function is
singular in both points $s=3$ and $s=4$. We take $a=x_2x_2^*$. From the definition
$$
\tilde{\pi}(x_2x_2^*)\ket{l,m_1,m_2;j}_\pm=
(1-q^{2(l-\epsilon-j+m_2)})\ket{l,m_1,m_2;j}_\pm\;.
$$
Then, modulo functions that are holomorphic when $\mr{Re}\,s>2$, we have
$$
\zeta_{x_2x_2^*}(s)\sim\zeta_{\tilde{\pi}(x_2x_2^*)}(s)= \zeta_1(s)
-2\sum_{l,j,m_1,m_2}(l+\tfrac{3}{2})^{-s}q^{2(l-\epsilon-j+m_2)}
\sim\tfrac{4}{3}\,\zeta(s-3)-\tfrac{4}{1-q^4}\,\zeta(s-2)\;.
$$
This proves the first part of the proposition,
that is $\Sigma\cap\{\mr{Re}\,s>2\}=\{3,4\}$.

The proof of equation (\ref{eq:topres}) is based on the observation that the residue
in $s=4$ of $\zeta_T$, for $T$ a basis element of $P\mathfrak{A}Q$, is zero unless
$T=1$. That is, it depends only on the image of $T$ under the map sending $\beta,A$ and $F$ to $0$ while $\alpha\mapsto e^{i\phi}$ and
$B\mapsto e^{i\theta}$. Composing this map with $\tilde{\pi}$ we get the morphism $\sigma:\A(S^4_q)\to \A(S^1)$ of the proposition and that
$$
\nint a|D|^{-4}\propto\int_0^{2\pi}\!\!\sigma(a)\de\theta\;.
$$
The equality $\,\int\mkern-16mu- \,|D|^{-4}=\frac{4}{3}\,$ fixes the proportionality
constant.
\end{prova}

\section{The real structure}\label{sec:7.7}
Classically, if $(\A(M),\HH,D,\gamma)$ is the canonical spectral triple associated with a
$4$-dimensional spin manifold $M$, there exists an antilinear isometry $J$ on $\HH$,
named the \emph{real structure}, satisfying the following compatibility condition
\begin{equation}\label{eq:noA}
J^2=-1\;,\quad\qquad
J\gamma=\gamma J\;,\quad\qquad
JD=DJ\;.
\end{equation}
There are also two additional conditions involving the coordinate algebra $\A(M)$:
\begin{equation}\label{eq:further}
[a,JbJ^{-1}]=0\;,\qquad[[D,a],JbJ^{-1}]=0\;,\qquad\forall\;a,b\in\A(M)\;.
\end{equation}
The real structure on $S^4$ is equivariant and equivariance is sufficient to
determine $J$. 

In the deformed situation one has to be careful on how to implement equivariance. Let us start with the working hypothesis that equivariance for $J$ is the requirement that it satisfies $Jh=S(h)^*J$ for all $h\in U_q(so(5))$. 
Then, consider the Casimir operator $\mc{C}_1$ given in equation
(\ref{eq:Cdef}). This operator commutes with $J$ since $S(\mc{C}_1)^*=\mc{C}_1$ and from its expression, $\mc{C}_1\ket{l,m_1,m_2;j}=(q^{2j+1}+q^{-2j-1})\ket{l,m_1,m_2;j}$,
we conclude that $J$ leaves the index $j$ invariant.
Compatibility with $\gamma$ and $D$ in equation (\ref{eq:noA}) and equivariance
with respect to $h=K_1$ and $h'=K_2$ yields
$$
J\ket{l,m_1,m_2;j}_\pm=c_\pm(l,m_1,m_2;j)\ket{l,-m_1,-m_2;j}_\pm\;,
$$
with some constants  $c_\pm$  to be determined. Equivariance with respect to $h=E_1$
implies
$$
c_\pm(l,m_1,m_2;j)=(-1)^{m_1+1/2}q^{m_1}c_\pm(l,m_2;j)\;.
$$
For $h=E_2$, looking at the piece diagonal in $j$ we deduce that the dependence on $m_2$
is through a factor $q^{3m_2}$; and looking at the piece shifting $j$ by $\pm 1$ we conclude that
$$
c_\pm(l,m_1,m_2;j)=(-1)^{j+m_1}q^{m_1+3m_2}c_\pm(l)\;.
$$
Such an operator $J$  cannot be antiunitary unless $q=1$. At $q=1$ the antiunitarity
condition requires that $c_\pm(l)\in U(1)$ and modulo a unitary equivalence we can choose
$c_\pm(l)=i^{2l+1}$. In conclusion for $q=1$ the operator
\begin{equation}\label{eq:J}
J\ket{l,m_1,m_2;j}_\pm=i^{2l+1}(-1)^{j+m_1}\ket{l,-m_1,-m_2;j}_\pm\;,
\end{equation}
is \emph{the} real structure on $S^4$ (modulo a unitary equivalence). 

For $q\neq 1$ we
keep (\ref{eq:J}) as the real structure and notice that conditions \eqref{eq:noA} are
satisfied, but $J$ no longer satisfies the requirement $Jh=S(h)^*J$ for all $h\in U_q(so(5))$. 
Nevertheless, $J$ is the
antiunitary part of an  antilinear operator $T$ that has this property. The antilinear operator $T$ defined by 
$$
T\ket{l,m_1,m_2;j}_\pm=i^{2l+1}(-1)^{j+m_1}q^{m_1+3m_2}\ket{l,-m_1,-m_2;j}_\pm\;,
$$
has the $J$ in \eqref{eq:J} as antiunitary part and it is \emph{equivariant}, i.e. it is such that $Th=S(h)^*T$ for all $h\in U_q(so(5))$.

Next, we turn to the conditions (\ref{eq:further}). 
In parallel with the cases of the manifold of $SU_q(2)$ in \cite{DLS05} and
of Podle{\'s} spheres in \cite{DLPS05,DDLW07}, once again we need to modify
them.
For instance, the commutator $[x_2,Jx_2J]$ is not zero, as one can see
by computing the matrix element
\begin{align*}
f(l,j,m_2)&:=\pm \inner{l+1,m_1,m_2;j\big|[x_2,Jx_2J]\big|l,m_1,m_2;j}_\pm \\
&=D_{l+1,j,m_2-1}^0D_{l,j,-m_2}^+-D_{l,j,m_2-1}^+D_{l,j,-m_2}^0
+D_{l+1,j,-m_2-1}^0D_{l,j,m_2}^+-D_{l,j,-m_2-1}^+D_{l,j,m_2}^0 \;,
\end{align*}
which for $j=\frac{1}{2}$ and $m_2=l$ is
$$
f(l,\tfrac{1}{2},l)=-q^{-l-4}(1-q^2)^2[2]\frac{(q^{l-1}+q^{-l+1})
\sqrt{[2l+3]}\,[l+1][l+2][l+3]}{[2l+2][2l+4]^2[2l+6]}\neq 0\;.
$$
It is relatively easy to prove that the two conditions are satisfied
modulo the ideal $\mc{J}_2$. It is much more cumbersome computationally to
show that they are in fact satisfied modulo the smaller ideal of
smoothing operators.

\begin{prop}\label{prop:real}
Let $J$ be the antilinear isometry given by (\ref{eq:J}). Then,
$$
[a,JbJ]\in\mc{J}_2\;,\quad[[D,a],JbJ]\in\mc{J}_2\;,\qquad\forall\;a,b\in\A(S^4_q)\;.
$$
\end{prop}
\begin{prova}
We lift $J$ and $D$ to the Hilbert space $\hat{\HH}$
defined in Section \ref{sec:HH}, as follows:
\begin{align*}
\hat{J}\!\kkett{l,m_1,m_2;j}_\pm &=i^{2l+1}(-1)^{j+m_1}\!\kkett{l,-m_1,-m_2;j}_\pm \;, \\
\hat{D}\!\kkett{l,m_1,m_2;j}_\pm &=(l+\tfrac{3}{2})\!\kkett{l,m_1,m_2;j}_\mp \;.
\end{align*}
Notice that $\hat{J}^2=-1$ on $\hat{\HH}$ (thanks to the phase $i^{2l+1}$ that is
irrelevant when restricted to $\HH$).

Let now $\{\alpha,\beta,\alpha^*,\beta^*,A,B,B^*\}$ be the operators defined in Section
\ref{sec:HH}, generators of the algebra $\A(SU_q(2))\otimes\A(S^2_q)$.
Due to Proposition~\ref{prop:app} it is enough to prove that for all
pairs $(a,b)$ of such generators,
the commutators $[a,\hat{J}b\hat{J}]$ and $[[\hat{D},a],\hat{J}b\hat{J}]$
are weighted shifts with weight which are bounded by $q^{2j}$. From 
$$
[\hat{D},\alpha]=\tfrac{1}{2}\,\alpha\hat{F}\;,\quad
[\hat{D},\beta]=\tfrac{1}{2}\,\beta\hat{F}\;,\quad
[\hat{D},A]=0\;,\quad
[\hat{D},B]=B\hat{F}\;,
$$
the condition on $[[\hat{D},a],\hat{J}b\hat{J}]$ follows
from the same condition on $[a,\hat{J}b\hat{J}]$,
and we have to compute only the latter commutators.

Since $[a,\hat{J}b^*\hat{J}]=-[a^*,\hat{J}b\hat{J}]^*$ and 
$[b,\hat{J}a\hat{J}]=\hat{J}[a,\hat{J}b\hat{J}]\hat{J}$,
we have to check the $16$ combinations in the following table.
\begin{center}\begin{tabular}{c|ccccccc}
\hline\hline $b \backslash a$ &
           $\alpha$ & $\alpha^*$ & $\beta$ & $\beta^*$ & $A$ & $B$ & $B^*$ \\
\hline
$\alpha$   & $\bullet$ & $\bullet$ & $\times$ & $\times$ & $\bullet$ & $\bullet$ & $\bullet$ \\
$\beta$    &        &          & $\bullet$ & $\times$ & $\bullet$ & $\bullet$ & $\bullet$ \\
$B$        &        &          &       &         & $\bullet$ & $\bullet$ & $\bullet$ \\
$A$        &        &          &       &         & $\bullet$ &   & \\
\hline\hline
\end{tabular}
\end{center}
By direct computations one shows that bullets in the table correspond to
vanishing commutators. On the other hand, the commutators corresponding to
the crosses in the table are given, on the subspace with $j-|m_1|\in\N$, by
\begin{align*}
[\beta^*,\hat{J}\alpha\hat{J}]\kkett{l,m_1,m_2;j}_\pm &=q^{j+m_1}\Big\{\sqrt{1-q^{2(j-m_1+1)}}-
\sqrt{1-q^{2(j-m_1)}}\Big\}\kkett{l,m_1,m_2;j}_\pm \\
[\beta,\hat{J}\alpha\hat{J}]\kkett{l,m_1,m_2;j}_\pm &=-[\beta^*,\hat{J}\alpha\hat{J}]\kkett{l+1,m_1-1,m_2;j+1}_\pm\;,\\
[\beta^*,\hat{J}\beta\hat{J}]\kkett{l,m_1,m_2;j}_\pm &=-[2]q^{2j}\kkett{l,m_1+1,m_2;j}_\pm \;.
\end{align*}
Since $1-u\leq\sqrt{1-u}\leq 1$ for all $u\in[0,1]$, we have that
$$
0\leq q^{j+m_1}\Big\{\sqrt{1-q^{2(j-m_1+1)}}-\sqrt{1-q^{2(j-m_1)}}\Big\}
\leq q^{j+m_1}(1-1+q^{2(j-m_1)})\leq q^{2j}\;.
$$
Then, all three non-zero commutators are weighted shifts with weights bounded by
$q^{2j}$. 
\end{prova}

\begin{prop}
Let $J$ be the antilinear isometry given by (\ref{eq:J}). Then,
$$
[a,JbJ]\in\op\;,\quad[[D,a],JbJ]\in\op\;,\qquad\forall\;a,b\in\A(S^4_q)\;.
$$
\end{prop}

\begin{prova}
By Leibniz rule, it is sufficient to prove the statement
when $a$ and $b$ are generators of the algebra. The operators
$[D,a]-\delta(a)F=|D|[F,a]$ are smoothing.
Thus,
it is enough to show that
\begin{equation}\label{eq:aJbJ}
[a,JbJ]\in\op\;,\qquad
[\delta(a),JbJ]\in\op\;,
\end{equation}
for any pair $(a,b)$ of generators.
From
$$
[b,JaJ]=J[a,JbJ]J \;,\qquad
[\delta(b),JaJ]=-J[\delta(a),JbJ]J+\delta([b,JaJ]) \;,
$$
it follows that if (\ref{eq:aJbJ}) is satisfied for a particular
pair $(a,b)$, then it is satisfied for $(b,a)$ too. From
\begin{equation}\label{eq:symstar}
[a^*,Jb^*J]=-[a,JbJ]^*\;,\qquad
[\delta(a^*),b^*]=[\delta(a),b]^* \;,
\end{equation}
we see that if (\ref{eq:aJbJ}) is satisfied for a
pair $(a,b)$, then it is satisfied for $(a^*,b^*)$ too. 
With these symmetries we need to check only the following
$9$ cases out of $25$:
\begin{center}\begin{tabular}{c|ccccc}
\hline\hline $b \backslash a$ &
           $x_0$ & $x_1$ & $x_1^*$ & $x_2$ & $x_2^*$ \\
\hline
$x_0$   & $\bullet$ & $\bullet$ && $\bullet$ \\
$x_1$  && $\bullet$ & $\bullet$ & $\bullet$ & $\bullet$ \\
$x_2$&&&& $\bullet$ & $\bullet$ \\
\hline\hline
\end{tabular}
\end{center}

\noindent
From equations (\ref{eq:apprep}) and (\ref{eq:coefsim}) we see that modulo smoothing operators
\begin{align*}
\oh\bigl\{x_0+\delta(x_0)\bigr\}\ket{l,m_1,m_2;j} 
        &\simeq \boxed{ A^+_{j,m_1}\hat{C}^+_{l,j,m_2}\ket{l+1,m_1,m_2;j+1} } \\
        &+q^{-2j}A^0_{j,m_1}\hat{H}^+_{l,j,m_2}\ket{l+1,m_1,m_2;j} \\
        &+A^+_{j-1,m_1}\hat{C}^-_{l+1,j-1,m_2}\ket{l+1,m_1,m_2;j-1} \;,\\
\oh\bigl\{x_0-\delta(x_0)\bigr\}\ket{l,m_1,m_2;j}
        &\simeq A^+_{j,m_1}\hat{C}^-_{l,j,m_2}\ket{l-1,m_1,m_2;j+1} \\
        &+q^{-2j}A^0_{j,m_1}\hat{H}^+_{l-1,j,m_2}\ket{l-1,m_1,m_2;j} \\
        &+\boxed{ A^+_{j-1,m_1}\hat{C}^+_{l-1,j-1,m_2}\ket{l-1,m_1,m_2;j-1} } \;, \\
\oh\bigl\{x_1+\delta(x_1)\bigr\}\ket{l,m_1,m_2;j}
        &\simeq \boxed{ B^+_{j,m_1}\hat{C}^+_{l,j,m_2}\ket{l+1,m_1+1,m_2;j+1} } \\
        &+q^{-2j}B^0_{j,m_1}\hat{H}^+_{l,j,m_2}\ket{l+1,m_1+1,m_2;j} \\
        &+B^-_{j,m_1}\hat{C}^-_{l+1,j-1,m_2}\ket{l+1,m_1+1,m_2;j-1} \;,\\
\oh\bigl\{x_1-\delta(x_1)\bigr\}\ket{l,m_1,m_2;j}
        &\simeq B^+_{j,m_1}\hat{C}^-_{l,j,m_2}\ket{l-1,m_1+1,m_2;j+1} \\
        &+q^{-2j}B^0_{j,m_1}\hat{H}^+_{l-1,j,m_2}\ket{l-1,m_1+1,m_2;j} \\
        &+\boxed{ B^-_{j,m_1}\hat{C}^+_{l-1,j-1,m_2}\ket{l-1,m_1+1,m_2;j-1} } \;, \\
\oh\bigl\{x_2+\delta(x_2)\bigr\}\ket{l,m_1,m_2;j}
        &\simeq \boxed{\boxed{ \hat{D}_{l,j,m_2}^+\ket{l+1,m_1,m_2+1;j} }} \;,\\
\oh\bigl\{x_2-\delta(x_2)\bigr\}\ket{l,m_1,m_2;j}
        &\simeq\hat{D}_{l,j,m_2}^-\ket{l-1,m_1,m_2+1;j} \;.
\end{align*}
where
\begin{subequations}\label{eq:Dhat}
\begin{align}
\hat{C}^+_{l,j,m_2} &=-q^{l-j+m_2-\epsilon}\sqrt{1-q^{2(l+j+m_2+3+\epsilon)}} \;, \\
\hat{C}^-_{l,j,m_2} &=-q^{l+j+m_2+1+\epsilon}\sqrt{1-q^{2(l-j+m_2-\epsilon)}} \;, \\
\hat{H}^+_{l,j,m_2} &=q^{2j}q^{l+m_2+1}\sqrt{q^{2\epsilon(2j+1)}-q^{2(l+m_2+2)}} \;, \\
\hat{D}_{l,j,m_2}^+ &=\sqrt{1-q^{2(l+j+m_2+3+\epsilon)}}\sqrt{1-q^{2(l-j+m_2+2-\epsilon)}} \;,\\
\hat{D}_{l,j,m_2}^- &=-q^{2(l+m_2)+3} \;.
\end{align}
\end{subequations}
We have divided the terms in three classes, which need to be analyzed
separately.

All terms $T$ which are not `boxed' have coefficients which are uniformly
bounded by $q^{l+m_2}$; since the conjugation with $J$ changes sign of
the labels $m_1,m_2$, for such $T$'s, the coefficients of $JTJ$ are
uniformly bounded by $q^{l-m_2}$.
They give products (and so commutators) with coefficients bounded by
$q^{l+m_2}q^{l-m_2}=q^{2l}$, and so (these products) are smoothing operators.

Analogously, the coefficients of single-boxed terms are bounded by $q^{l-j+m_2}$,
and becomes smoothing when multiplied by the $J$-conjugated of non-boxed
terms (as $q^{l-j+m_2}q^{l-m_2}\leq q^l$), and vice versa for the
product of a non-boxed term with the $J$-conjugated of a single-boxed
one ($q^{l+m_2}q^{l-j-m_2}\leq q^l$).

Next we consider pairs of single-boxed terms.
A closer look at the single-boxed terms in $x_0\pm\delta(x_0)$ and $x_1-\delta(x_1)$
(and then $x_1^*+\delta(x_1^*)$) shows that they have coefficients bounded by $q^{l+m_1+m_2}$,
and become smoothing when multiplied by the $J$-conjugated of one of them
($q^{l+m_1+m_2}q^{l-m_1-m_2}=q^{2l}$).
Last single-boxed term is the one in $x_1+\delta(x_1)$ (and $x_1^*-\delta(x_1^*)$).
The relevant terms for the commutators involving them are
\begin{align*}
&\oh[x_1+\delta(x_1),Jx_0J]\ket{l,m_1,m_2;j} \\
&\qquad \simeq\bigl\{A^+_{j+1,-m_1-1}\hat{C}^+_{l+1,j+1,-m_2}B^+_{j,m_1}\hat{C}^+_{l,j,m_2} + \\ &\qquad\qquad\qquad\qquad
         -A^+_{j,-m_1}\hat{C}^+_{l,j,-m_2}B^+_{j+1,m_1}\hat{C}^+_{l+1,j+1,m_2}\bigr\}\ket{l+2,m_1+1,m_2;j+2} \\
&\qquad\qquad +\bigl\{A^+_{j,-m_1-1}\hat{C}^+_{l,j,-m_2}B^+_{j,m_1}\hat{C}^+_{l,j,m_2} + \\ &\qquad \qquad\qquad\qquad
         -A^+_{j-1,-m_1}\hat{C}^+_{l-1,j-1,-m_2}B^+_{j-1,m_1}\hat{C}^+_{l-1,j-1,m_2}\bigr\}\ket{l,m_1+1,m_2;j} \;, \\
&\oh[x_1+\delta(x_1),Jx_1J]\ket{l,m_1,m_2;j} \\
&\qquad \simeq\bigl\{B^+_{j+1,m_1-1}\hat{C}^+_{l+1,j+1,m_2}B^+_{j,-m_1}\hat{C}^+_{l,j,-m_2} + \\ &\qquad\qquad\qquad\qquad
         -B^+_{j+1,-m_1-1}\hat{C}^+_{l+1,j+1,-m_2}B^+_{j,m_1}\hat{C}^+_{l,j,m_2}\bigr\}\ket{l+2,m_1,m_2;j+2} \\
&\qquad\qquad +\bigl\{B^+_{j-1,m_1-1}\hat{C}^+_{l-1,j-1,m_2}B^-_{j,-m_1}\hat{C}^+_{l-1,j-1,-m_2} + \\ &\qquad\qquad\qquad\qquad
         -B^-_{j+1,-m_1-1}\hat{C}^+_{l,j,-m_2}B^+_{j,m_1}\hat{C}^+_{l,j,m_2}\bigr\}\ket{l,m_1,m_2;j} \;,\\
&\oh[x_1-\delta(x_1),Jx_1J]\ket{l,m_1,m_2;j} \\
&\qquad \simeq\bigl\{B^-_{j+1,m_1-1}\hat{C}^+_{l,j,m_2}B^+_{j,-m_1}\hat{C}^+_{l,j,-m_2} + \\ &\qquad\qquad\qquad\qquad
         -B^-_{j,m_1}\hat{C}^+_{l-1,j-1,m_2}B^+_{j-1,-m_1-1}\hat{C}^+_{l-1,j-1,-m_2}\bigr\}\ket{l,m_1,m_2;j} \;,\\
&\oh[x_1^*+\delta(x_1^*),Jx_1J]\ket{l,m_1,m_2;j} \\
&\qquad\simeq \bigl\{B^-_{j+2,m_1-2}\hat{C}^+_{l+1,j+1,m_2}B^+_{j,-m_1}\hat{C}^+_{l,j,-m_2} + \\ &\qquad\qquad\qquad\qquad
         -B^+_{j+1,-m_1+1}\hat{C}^+_{l+1,j+1,-m_2}B^-_{j+1,m_1-1}\hat{C}^+_{l,j,m_2}\bigr\}\ket{l+2,m_1-2,m_2;j+2} \;,\\
&\oh[x_1^*-\delta(x_1^*),Jx_1J]\ket{l,m_1,m_2;j}  \\
&\qquad\simeq \bigl\{B^+_{j,m_1-2}\hat{C}^+_{l,j,m_2}B^+_{j,-m_1}\hat{C}^+_{l,j,-m_2} + \\ &\qquad\qquad\qquad\qquad
          -B^+_{j-1,-m_1+1}\hat{C}^+_{l-1,j-1,-m_2}B^+_{j-1,m_1-1}\hat{C}^+_{l-1,j-1,m_2}\bigr\}\ket{l,m_1-2,m_2;j} \\
&\qquad\qquad +\bigl\{B^+_{j-2,m_1-2}\hat{C}^+_{l-2,j-2,m_2}B^-_{j,-m_1}\hat{C}^+_{l-1,j-1,-m_2} + \\ &\qquad\qquad\qquad\qquad
          -B^-_{j-1,-m_1+1}\hat{C}^+_{l-2,j-2,-m_2}B^+_{j-1,m_1-1}\hat{C}^+_{l-1,j-1,m_2}\bigr\}\ket{l-2,m_1-2,m_2;j-2} \;.
\end{align*}
We need to estimate products of the form $\hat{C}^+_{l,j,-m_2}\hat{C}^+_{l+i,j+i,m_2}$,
for which, modulo smoothing operators, we find
\begin{align*}
\hat{C}^+_{l,j,-m_2}\hat{C}^+_{l+i,j+i,m_2}
&\simeq\sqrt{q^{2(l-j)}-q^{2l}q^{2(l-m_2+3-\epsilon)}}\sqrt{q^{2(l-j)}-q^{2l}q^{2(l+m_2+3+2i+\epsilon)}} \\
&\simeq \sqrt{q^{2(l-j)}}\sqrt{q^{2(l-j)}}=q^{2(l-j)} \;.
\end{align*}
Using this we get
\begin{align*}
&\oh[x_1+\delta(x_1),Jx_0J]\ket{l,m_1,m_2;j} \\
&\qquad \simeq q^{2(l-j)}\bigl\{A^+_{j+1,-m_1-1}B^+_{j,m_1}
         -A^+_{j,-m_1}B^+_{j+1,m_1}\bigr\}\ket{l+2,m_1+1,m_2;j+2} \\
&\qquad +q^{2(l-j)}\bigl\{A^+_{j,-m_1-1}B^+_{j,m_1}
         -A^+_{j-1,-m_1}B^+_{j-1,m_1}\bigr\}\ket{l,m_1+1,m_2;j} \;, \\
&\oh[x_1+\delta(x_1),Jx_1J]\ket{l,m_1,m_2;j} \\
&\qquad \simeq q^{2(l-j)}\bigl\{B^+_{j+1,m_1-1}B^+_{j,-m_1}
         -B^+_{j+1,-m_1-1}B^+_{j,m_1}\bigr\}\ket{l+2,m_1,m_2;j+2} \\
&\qquad +q^{2(l-j)}\bigl\{B^+_{j-1,m_1-1}B^-_{j,-m_1}
         -B^-_{j+1,-m_1-1}B^+_{j,m_1}\bigr\}\ket{l,m_1,m_2;j} \;,\\
&\oh[x_1-\delta(x_1),Jx_1J]\ket{l,m_1,m_2;j} \\
&\qquad \simeq q^{2(l-j)}\bigl\{B^-_{j+1,m_1-1}B^+_{j,-m_1}-B^-_{j,m_1}B^+_{j-1,-m_1-1}\bigr\}\ket{l,m_1,m_2;j} \;,\\
&\oh[x_1^*+\delta(x_1^*),Jx_1J]\ket{l,m_1,m_2;j} \\
&\qquad\simeq q^{2(l-j)}\bigl\{B^-_{j+2,m_1-2}B^+_{j,-m_1}-B^+_{j+1,-m_1+1}B^-_{j+1,m_1-1}\bigr\}\ket{l+2,m_1-2,m_2;j+2} \;,\\
&\oh[x_1^*-\delta(x_1^*),Jx_1J]\ket{l,m_1,m_2;j}  \\
&\qquad\simeq q^{2(l-j)}\bigl\{B^+_{j,m_1-2}B^+_{j,-m_1}-B^+_{j-1,-m_1+1}B^+_{j-1,m_1-1}\bigr\}\ket{l,m_1-2,m_2;j} \\
&\qquad +q^{2(l-j)}\bigl\{B^+_{j-2,m_1-2}B^-_{j,-m_1}-B^-_{j-1,-m_1+1}B^+_{j-1,m_1-1}\bigr\}\ket{l-2,m_1-2,m_2;j-2} \;.
\end{align*}
To prove that these commutators are smoothing we still need to check
that the terms in braces are bounded by $q^j$ (since $q^{2(l-j)}q^j\leq q^l$
is of rapid decay). This is done by using equations (\ref{eq:get}).
For example the first two braces are identically
zero, while the third one is
\begin{align*}
B^+_{j+1,m_1-1}B^+_{j,-m_1}-B^+_{j+1,-m_1-1}B^+_{j,m_1}
&=\tilde{B}^+_{j+1,m_1-1}\tilde{B}^+_{j,-m_1}-\tilde{B}^+_{j+1,-m_1-1}\tilde{B}^+_{j,m_1}
+O(q^j) \\ &=0+O(q^j) \;.
\end{align*}

What remains to control are the commutators $[x_2+\delta(x_2),JbJ]$
for $b=x_0,x_1,x_2$ and the commutators $[x_2^*-\delta(x_2^*),JbJ]$
for $b=x_1,x_2$ (which involve the `doubly-boxed' term).

The operators $x_2$ and $\delta(x_2)$ do not shift $m_1,j$ and have
coefficients independent of $m_1$. Thus, any operator acting only on
the label $m_1$ and with coefficients depending only on $m_1,j$,
commutes with $x_2$ and $\delta(x_2)$ and so can be neglected.
In particular, $x_0$ and $x_1$ can be written as sums of products of
operators of this kind (commuting with $x_2$ and $\delta(x_2)$) by
operators $y_i$'s,
\begin{align*}
y_1\ket{l,m_1,m_2;j} & :=\hat{C}^+_{l,j,m_2}\ket{l+1,m_1,m_2;j+1} \;,\\
y_2\ket{l,m_1,m_2;j} & :=\hat{H}^+_{l,j,m_2}\ket{l+1,m_1,m_2;j} \;,\\
y_3\ket{l,m_1,m_2;j} & :=\hat{C}^-_{l,j,m_2}\ket{l-1,m_1,m_2;j+1} \;,
\end{align*}
and their adjoints.
To prove that the commutators
$[x_2+\delta(x_2),JbJ]$,
for $b=x_0,x_1,x_2$, and $[x_2^*-\delta(x_2^*),JbJ]$,
for $b=x_1,x_2$, are smoothing, is sufficient to establish
the same for $b=y_1,y_2,y_3$.
For these operators we have
\begin{align*}
\oh[x_2+\delta(x_2),&\,Jy_1J]\ket{l,m_1,m_2;j} \\ & \simeq
        \bigl\{\hat{C}^+_{l+1,j,-m_2-1}\hat{D}_{l,j,m_2}^+
        -\hat{D}_{l+1,j+1,m_2}^+\hat{C}^+_{l,j,-m_2}\bigr\}
        \ket{l+2,m_1,m_2+1;j+1} \\ &=
        \hat{C}^+_{l,j,-m_2}\bigl\{\hat{D}_{l,j,m_2}^+
        -\hat{D}_{l+1,j+1,m_2}^+\bigr\}\ket{l+2,m_1,m_2+1;j+1} \;,\\
\oh[x_2+\delta(x_2),&\,Jy_2J]\ket{l,m_1,m_2;j} \\ & \simeq
        \bigl\{\hat{H}^+_{l+1,j,-m_2-1}\hat{D}_{l,j,m_2}^+
        -\hat{D}_{l+1,j,m_2}^+\hat{H}^+_{l,j,-m_2}\bigr\}
        \ket{l+2,m_1,m_2+1;j} \\ &=
        \hat{H}^+_{l,j,-m_2}\bigl\{\hat{D}_{l,j,m_2}^+
        -\hat{D}_{l+1,j,m_2}^+\bigr\}\ket{l+2,m_1,m_2+1;j} \;,\\
\oh[x_2+\delta(x_2),&\,Jy_3J]\ket{l,m_1,m_2;j} \\ & \simeq
        \bigl\{\hat{C}^-_{l+1,j,-m_2-1}\hat{D}_{l,j,m_2}^+
        -\hat{D}_{l-1,j+1,m_2}^+\hat{C}^-_{l,j,-m_2}\bigr\}
        \ket{l,m_1,m_2+1;j+1} \\ &=
        \hat{C}^-_{l,j,-m_2}\bigl\{\hat{D}_{l,j,m_2}^+
        -\hat{D}_{l-1,j+1,m_2}^+\bigr\}\ket{l,m_1,m_2+1;j+1} \;,\\
\oh[x_2^*-\delta(x_2^*),&\,Jy_1J]\ket{l,m_1,m_2;j} \\ & \simeq
        \bigl\{\hat{C}^+_{l-1,j,-m_2+1}\hat{D}_{l-1,j,m_2-1}^+
        -\hat{D}_{l,j+1,m_2-1}^+\hat{C}^+_{l,j,-m_2}\bigr\}
        \ket{l,m_1,m_2-1;j+1} \\ &=
        \hat{C}^+_{l,j,-m_2}\bigl\{\hat{D}_{l-1,j,m_2-1}^+
        -\hat{D}_{l,j+1,m_2-1}^+\bigr\}\ket{l,m_1,m_2-1;j+1} \;,\\
\oh[x_2^*-\delta(x_2^*),&\,Jy_2J]\ket{l,m_1,m_2;j} \\ & \simeq
        \bigl\{\hat{H}^+_{l-1,j,-m_2+1}\hat{D}_{l-1,j,m_2-1}^+
        -\hat{D}_{l,j,m_2-1}^+\hat{H}^+_{l,j,-m_2}\bigr\}
        \ket{l,m_1,m_2-1;j} \\ &=
        \hat{H}^+_{l,j,-m_2}\bigl\{\hat{D}_{l-1,j,m_2-1}^+
        -\hat{D}_{l,j,m_2-1}^+\bigr\}\ket{l,m_1,m_2-1;j} \;,\\
\oh[x_2^*-\delta(x_2^*),&\,Jy_3J]\ket{l,m_1,m_2;j} \\ & \simeq
        \bigl\{\hat{C}^-_{l-1,j,-m_2+1}\hat{D}_{l-1,j,m_2-1}^+
        -\hat{D}_{l-2,j+1,m_2-1}^+\hat{C}^-_{l,j,-m_2}\bigr\}
        \ket{l-2,m_1,m_2+1;j+1} \\ &=
        \hat{C}^-_{l,j,-m_2}\bigl\{\hat{D}_{l-1,j,m_2-1}^+
        -\hat{D}_{l-2,j+1,m_2-1}^+\bigr\}\ket{l-2,m_1,m_2+1;j+1} \;.
\end{align*}
Now $\hat{D}_{l\pm 1,j+1,m_2}^+-\hat{D}_{l,j,m_2}^+$ is bounded
by $q^{l\pm j+m_2}$, and $q^{l\pm j+m_2}\hat{C}^\pm_{l,j,-m_2}$
is bounded by $q^{2l}$. Next, $\hat{H}^+_{l,j,-m_2}$ is bounded by
$q^{l+j-m_2}$, and $q^{l+j-m_2}\hat{D}_{l,j,m_2}^+\simeq 1$.
This proves that all previous commutators are smoothing.
For the $y_i^*$'s the same statement follows from the symmetry (\ref{eq:symstar}).

We arrived at last two commutators.
Modulo smoothing operators, the first one is
\begin{align*}
\oh[x_2+\delta(x_2),&\,Jx_2J]\ket{l,m_1,m_2;j} \\ & \simeq
\bigl\{\hat{D}_{l+1,j,-m_2-1}^+\hat{D}_{l,j,m_2}^+
-\hat{D}_{l+1,j,m_2-1}^+\hat{D}_{l,j,-m_2}^+\bigr\}\ket{l+2,m_1,m_2;j} \\
&\qquad\qquad +\bigl\{\hat{D}_{l+1,j,-m_2-1}^-\hat{D}_{l,j,m_2}^+
-\hat{D}_{l-1,j,m_2-1}^+\hat{D}_{l,j,-m_2}^-\bigr\}\ket{l,m_1,m_2;j} \\
&=\hat{D}_{l,j,-m_2}^-\bigl\{\hat{D}_{l,j,m_2}^+
-\hat{D}_{l-1,j,m_2-1}^+\bigr\}\ket{l,m_1,m_2;j} \;, \\
\intertext{where the second equality follows from the fact
that both $\hat{D}_{l,j,m_2}^+$ and $\hat{D}_{l,j,m_2}^-$ in
(\ref{eq:Dhat}) depend on $l$ and $m_2$ only through their sum.
For the same reason we have also that}
\oh[x_2^*-\delta(x_2^*),&\,Jx_2J]\ket{l,m_1,m_2;j} \\ & \simeq
\bigl\{\hat{D}_{l-1,j,-m_2+1}^+\hat{D}_{l-1,j,m_2-1}^+
-\hat{D}_{l,j,m_2-2}^+\hat{D}_{l,j,-m_2}^+\bigr\}\ket{l,m_1,m_2-2;j} \\
&\qquad\qquad +\bigl\{\hat{D}_{l-1,j,-m_2+1}^-\hat{D}_{l-1,j,m_2-1}^+
-\hat{D}_{l-2,j,m_2-2}^+\hat{D}_{l,j,-m_2}^-\bigr\}\ket{l-2,m_1,m_2-2;j} \\
&=\hat{D}_{l,j,-m_2}^-\bigl\{\hat{D}_{l-1,j,m_2-1}^+
-\hat{D}_{l-2,j,m_2-2}^+\bigr\}\ket{l-2,m_1,m_2-2;j} \;.
\end{align*}
The final observation that $\hat{D}_{l,j,-m_2}^-\hat{D}_{l-i,j,m_2-i}^+\simeq\hat{D}_{l,j,-m_2}^-$,
for $i=0,1,2$, gives that these commutators vanish modulo smoothing operators.
\end{prova}

\section{An application: the Haar functional on $S^4_q$}
The representation of $\A(S^4_q)$ on $L^2(S^4_q)$ is the GNS representation
associated to the Haar functional $\varphi$. Then, we can recover the
Haar functional from the cyclic vector $1$. First, we have to identify
$1$ among the basis of `harmonic functions' $\ket{l,m_1,m_2;j}$.
It is the highest weight vector with weight $(0,0)$, then
$1\propto\ket{0,0,0;0}$. We'll denote it $\ket{0}$ for short.
Then the Haar functional can be expressed as
$$
\varphi(a)=\inner{0|a|0}
$$
for any $a\in\A(S^4_q)$. The $\mc{U}_q(so(5))$-invariance of right hand side
is evident, since $h\ket{0}=\epsilon(h)\ket{0}$ is the
trivial representation and by crossed product relations
$$
\varphi(h\az a)=\inner{0\smash[b]{|h_{(1)}aS(h_{(2)})|}0}=
\inner{0\smash[b]{|a\epsilon(h_{(1)})S(h_{(2)})|}0}=
\inner{0|ah|0}=\epsilon(h)\varphi(a)\;.
$$
Let us compute explicitly $\varphi$ in the basis 
of monomials $x_0^{n_0}x_1^{n_1}(x_1^*)^{n_2}x_2^{n_3}$,
with $n_0,n_1,n_2\in\N$, $n_3\in\Z$ and with the notation
$x_2^{n_3}:=(x_2^*)^{|n_3|}$ if $n_3<0$. 
Invariance of $\varphi$ implies that $\varphi(a)=0$ unless
$K_i\az a=a$ for all $i=1,2$. Applying this condition to the generic monomial
we get
$$
\varphi(x_0^{n_0}x_1^{n_1}(x_1^*)^{n_2}x_2^{n_3})=0
\qquad\mr{if}\;\,n_1\neq n_2\;\,\mr{or}\;\,n_3\neq 0\;.
$$
Furthermore, for $n_1=n_2=:k$, $n_3=0$ and for odd $n_0$, from
$(x_1^*)^k\ket{0}\propto\ket{k,-k,0;k}$ and using the commutation
rule $x_0x_1=q^2x_1x_0$ we get
$$
\varphi(x_0^{n_0}x_1^k(x_1^*)^k)=q^{2kn_0}
\big<0\big|x_1^kx_0^{n_0}(x_1^*)^k\big|0\big>=q^{2kn_0}
\inner{k,-k,0;k|x_0^{n_0}|k,-k,0;k}
$$
which is zero since $x_0^{n_0}$ is off-diagonal for odd $n_0$.
For $n=2j$ even we use `twisted' cyclicity.
Since $\kappa^{-1}(x_1^*)=q^2x_1^*$, we have
\begin{equation}\label{eq:twist}
\varphi(x_0^{2j-2}x_1^{k+1}(x_1^*)^{k+1})=
q^2\varphi(x_1^*x_0^{2j-2}x_1^{k+1}(x_1^*)^k)=
q^{4j-2}\varphi(x_0^{2j-2}x_1^*x_1^{k+1}(x_1^*)^k)\;.
\end{equation}
By induction one proves
$$
x_1^*x_1^{k+1}=x_1^{k+1}x_1^*+(q^{-4k}-q^4)x_0^2x_1^k\;,
$$
so equation (\ref{eq:twist}) can be rewritten as
\begin{align}
\varphi(x_0^{2j}x_1^k(x_1^*)^k)&=
q^{2(j+k)-3}\frac{[2j-1]}{[2k+2]}\varphi(x_0^{2j-2}x_1^{k+1}(x_1^*)^{k+1})=
\ldots = \nonumber \\ &=q^{2j(j+k)-3j}\frac{[2j-1]!!\,[2k]!!}{[2(j+k)]!!}
\;\varphi(x_1^{j+k}(x_1^*)^{j+k}) \;. \label{eq:subst}
\end{align}
As usual we denote with $[n]!!$ the $q$-analogue of the double factorial, given by
\begin{equation*}
[n]!!=\begin{cases}
[n][n-2][n-4]\ldots [2] & \mr{if}\;n\geq 2\;\textrm{is even}, \\
[n][n-2][n-4]\ldots [1] & \mr{if}\;n\geq 1\;\textrm{is odd}, \\
1 & \mr{if}\;n=0,-1.
\end{cases}
\end{equation*}
Finally, from the equality
$\,x_1^*\ket{l,0,-l;l}=-q^{-1/2}\sqrt{\frac{[2l+2]}{[2l+5]}}\ket{l+1,-l-1,0;l+1}\,$
we deduce that $\,(x_1^*)^{j+k}\ket{0}=(-1)^{j+k}q^{-(j+k)/2}\sqrt{[3]
\frac{[2(j+k)]!!}{[2(j+k)+3]!!}}\ket{j+k,-(j+k),0;j+k}\;$
and then
$$
\varphi(x_1^{j+k}(x_1^*)^{j+k})=q^{-(j+k)}[3]\frac{[2(j+k)]!!}{[2(j+k)+3]!!}\;.
$$
Substituting this into equation (\ref{eq:subst}) we get the following result.

\begin{prop}\label{prop:HaarS4q}
The Haar functional on $\A(S^4_q)$ is given by:
\begin{equation*}
\varphi(x_0^{n_0}x_1^{n_1}(x_1^*)^{n_2}x_2^{n_3})=
\begin{cases}
q^{(2j-1)(j+k)-3j}[3]\frac{[2j-1]!!\,[2k]!!}{[2(j+k)+3]!!} \quad
  & \mr{if}\;n_0=2j\in 2\N,\;n_1=n_2=k,\;n_3=0, \\
0 & \mr{otherwise.}
\end{cases}
\end{equation*}
\end{prop}

\noindent
It is instructive to compare the situation with the case $q=1$.
Cartesian coordinates on $S^4$ are expressed in spherical coordinates as
$$
x_0=\cos\theta\;,\quad
x_1=e^{i\alpha}\sin\theta\cos\phi\;,\quad
x_2=e^{i\beta}\sin\theta\sin\phi\;,
$$
with $\alpha,\beta\in [0,2\pi]$, $\theta\in[0,\pi]$ and $\phi\in[0,\pi/2]$.
The Riemann integral is
$$
\int_{S^4}\de\mr{Vol}\,S^4=
\int_{-1}^1(\sin\theta)^2\,\de\!\cos\theta
\int_0^1\cos\phi\,\de\!\cos\phi
\int_0^{2\pi}\de\alpha
\int_0^{2\pi}\de\beta\;,
$$
and it is normalized to $8\pi^2/3$.
The Haar measure on a monomial, in the non-trivial cases, yields:
\begin{align*}
\varphi(x_0^{2j}|x_1|^{2k})
&=\frac{3}{8\pi^2}\int_{S^4}x_0^{2j}|x_1|^{2k}\de\mr{Vol}\,S^4 \\
&=\frac{3}{2}\int_{-1}^1(\cos\theta)^{2j}(\sin\theta)^{2k+2}\,\de\!\cos\theta
\int_0^1(\cos\phi)^{2k+1}\,\de\!\cos\phi \\
&=\frac{3}{2k+2}\int_0^1\xi^{2j}(1-\xi^2)^{k+1}\,\de\xi \;.
\end{align*}
Integrating by parts
$$
\int_0^1\xi^{2j}(1-\xi^2)^{k+1}\,\de\xi=
\int_0^1(1-\xi^2)^{k+1}\,\de\frac{\xi^{2j+1}}{2j+1}
=\frac{2k+2}{2j+1}\int_0^1\xi^{2j+2}(1-\xi^2)^k\de\xi \;,
$$
and by induction
$$
\int_0^1\xi^{2j}(1-\xi^2)^{k+1}\,\de\xi=
\frac{(2k+2)!!\,(2j-1)!!}{(2j+2k+1)!!}\int_0^1\xi^{2j+2k+2}\de\xi
=\frac{(2k+2)!!\,(2j-1)!!}{(2j+2k+3)!!}\;.
$$
We arrive then at the final result
$$
\varphi(x_0^{2j}|x_1|^{2k})=3\frac{(2k)!!\,(2j-1)!!}{(2j+2k+3)!!}\;,
$$
which is the $q\to 1$ limit of Proposition \ref{prop:HaarS4q}.

\section{The homological dimension}
With $\kappa(a):=K_1^8K_2^6\az a$ the modular automorphism, we are
able to construct a twisted Hochschild $4$-cycle as follows.
First of all, note that $\kappa(x_0)=x_0$, $\kappa(x_1)=q^2x_1$
and $\kappa(x_2)=q^6x_2$. Thus, $\kappa$ is a scaling automorphism.
The element $K_1^8K_2^6$ implements the square of the antipode, that is
\begin{equation}\label{eq:Ssquare}
(K_1^8K_2^6)h(K_1^8K_2^6)^{-1}=S^2(h)\;,\qquad\forall\;h\in\mc{U}_q(so(5))\,,
\end{equation}
as one can check on the generators $E_i$
(a necessary and sufficient condition).

Let $\sigma:\mc{U}_q(so(5))\to\mr{Mat}_4(\C)$ be the representation
in (\ref{eq:spin}) and $e$ the idempotent in (\ref{eq:P}). By (\ref{eq:cov})
we have $\kappa(e)=K_1^8K_2^6\az e=\eta e\eta^{-1}$ where
$$
\eta:=\sigma(K_1^8K_2^6)=\mr{diag}(q^4,q^{-2},q^2,q^{-4})\;.
$$
Note that if we call $\gamma_i\in\mr{Mat}_N(\C)$ the following
deformation of Dirac's gamma matrices
$$
\gamma_0=\sma{
   1 &  0   &  0   & 0 \\
   0 & -q^2 &  0   & 0 \\
   0 &  0   & -q^2 & 0 \\
   0 &  0   &  0   & q^4}\;,
\qquad
\gamma_1=\sma{
   0 & 0 & -q & 0 \\
   0 & 0 & 0 & q^3 \\
   0 & 0 & 0 & 0 \\
   0 & 0 & 0 & 0}\;,
\qquad
\gamma_{-1}=\sma{
   0 & 0 & 0 & 0 \\
   0 & 0 & 0 & 0 \\
   -q^{-1} & 0 & 0 & 0 \\
   0 & q & 0 & 0}\;,
$$
$$
\gamma_2=\sma{
   0 & q^3 & 0 & 0 \\
   0 & 0 & 0 & 0 \\
   0 & 0 & 0 & q^3 \\
   0 & 0 & 0 & 0}\;,
\qquad
\gamma_{-2}=\sma{
   0 & 0 & 0 & 0 \\
   q^{-3} & 0 & 0 & 0 \\
   0 & 0 & 0 & 0 \\
   0 & 0 & q^{-3} & 0}\;,
$$
then the idempotent $e$ can be written as
$$
e=\tfrac{1}{2}(1+\gamma_0x_0+\gamma_1x_1+\gamma_2x_2+\gamma_{-1}x_1^*+\gamma_{-2}x_2^*)\;.
$$
\begin{lemma}\label{lemma:gamma}
$\tr_{\C^4}(\eta\gamma_i\gamma_j\gamma_k)=0\,$ for all $\,i,j,k\in\{0,\pm 1,\pm 2\}$.
\end{lemma}
\begin{prova}
The proof is a long and tedious computation (one has to compute $125$ traces),
in contrast with the $q=1$ case where one simply uses cyclicity of the trace
and the anticommutation rules of the Clifford algebra.
We quote the code for {\mate} which allows to perform the check with a personal computer.
\begin{quote}\begin{footnotesize}\begin{verbatim}
gamma[0]:={{1,0,0,0},{0,-q^2,0,0},{0,0,-q^2,0},{0,0,0,q^4}}
gamma[1]:={{0,0,-q,0},{0,0,0,q^3},{0,0,0,0},{0,0,0,0}}
gamma[-1]:={{0,0,0,0},{0,0,0,0},{-1/q,0,0,0},{0,q,0,0}}
gamma[2]:={{0,q^3,0,0},{0,0,0,0},{0,0,0,q^3},{0,0,0,0}}
gamma[-2]:={{0,0,0,0},{1/q^3,0,0,0},{0,0,0,0},{0,0,1/q^3,0}}
eta:={{q^4,0,0,0},{0,1/q^2,0,0},{0,0,q^2,0},{0,0,0,1/q^4}}
f[i_,j_,k_]:= Tr[eta.gamma[i].gamma[j].gamma[k]]
For[
  counter=0;i=-2,i<3,
  For[
    j=-2,j<3,
    For[
      k=-2,k<3,
      If[f[i,j,k]==0,Null,counter++,counter++];
      k++
      ];
    j++;Clear[k]
    ];
  i++;Clear[j,k]
  ];Print[counter]
\end{verbatim}\end{footnotesize}\end{quote}
The counter `\texttt{counter}' counts how many times the expression
\texttt{f$[$i,j,k$]$} ($=\tr_{\C^4}(\eta\gamma_i\gamma_j\gamma_k)$)
fails to be zero. Running this code one gets $0$, confirming that
the Lemma is true.
\end{prova}

\pagebreak

\begin{prop}
The element
$$
\omega_4:=\tr_{\C^4}(\eta(e-\tfrac{1}{2})^{\dot{\otimes}5})
=\sum\nolimits_{k_i=1}^4\eta_{k_0,k_1}(e_{k_1,k_2}-\tfrac{1}{2}\delta_{k_1,k_2})
\otimes\ldots\otimes (e_{k_5,k_0}-\tfrac{1}{2}\delta_{k_5,k_0})
$$
is the representative of an $\mc{U}_q(so(5))$-invariant twisted Hochschild $4$-cycle.
\end{prop}

\noindent
Note that $\omega_4$ is the sum of $4^5=1024$ terms, so it is very long
to write it down explicitly.

\begin{prova}
Using $\kappa(e)\eta=\eta e$ we write the boundary of $\omega_4$ as
\begin{align*}
b_\kappa(\omega_4) &=\sum_{j=0}^3(-1)^j(1+\delta_{j,0})\tr_{\C^4}(\eta (e-\tfrac{1}{2})^{\dot{\otimes}j}
\otimes 1\otimes (e-\tfrac{1}{2})^{\dot{\otimes}(3-j)})
=\frac{1}{8}\sum_{i,j,k}\tr_{\C^4}(\eta\gamma_i\gamma_j\gamma_k)\times \\
&\times\bigl(2\otimes x_i\otimes x_j\otimes x_k-x_i\otimes 1\otimes x_j\otimes x_k
+x_i\otimes x_j\otimes 1\otimes x_k-x_i\otimes x_j\otimes x_k\otimes 1\bigr) \;,
\end{align*}
and this is zero by Lemma \ref{lemma:gamma}.
By equation (\ref{eq:cov})
$$
h\az\omega_4=\tr_{\C^4}(\eta\sigma(h_{(1)})^t(e-\tfrac{1}{2})\dot{\otimes}
e^{\dot{\otimes}4}\sigma(S^{-1}(h_{(2)}))^t)
=\tr_{\C^4}(\sigma(S^{-1}(h_{(2)}))^t\eta\sigma(h_{(1)})^t
(e-\tfrac{1}{2})\dot{\otimes}e^{\dot{\otimes}4})\;,
$$
where in last step we used cyclicity of the trace.
But $\eta=\sigma(K_1^8K_2^6)=\sigma(K_1^8K_2^6)^t$ and by (\ref{eq:Ssquare})
$$
\sigma(S^{-1}(h_{(2)}))^t\eta\sigma(h_{(1)})^t=
\sigma\big( h_{(1)}K_1^8K_2^6S^{-1}(h_{(2)}) \big)=
\sigma\big( h_{(1)}S(h_{(2)})K_1^8K_2^6 \big)=\epsilon(h)\eta\;.
$$
From this it follows the invariance of $\omega_4$, that is $h\az\omega_4=\epsilon(h)\omega_4$.
\end{prova}

\noindent
For $q=1$ the image of $\omega_4$ under the isomorphism (\ref{eq:HKR}) is
the volume form of $S^4$. For $q\neq 1$ without the factor of $\eta$ Lemma
\ref{lemma:gamma} is no more valid, so the twist is essential to get a
coboundary. Non triviality of $\omega_4$ would prove that the `twisted'
Hochschild dimension is at least $4$, i.e.~there is no dimension drop.
To prove non-triviality one must pair $\omega_4$ with a $4$-cocycle
and prove that the pairing is not zero. A $4$-cocycle is constructed
in the following.
Recall that $K_1^8K_2^6\ket{l,m_1,m_2;j}_\pm=q^{2m_1+6m_2}\ket{l,m_1,m_2;j}_\pm$
(cf.~Section \ref{sec:7.1}).

\begin{lemma}\label{lemma:fin}
The operator
$$
K_1^{-8}K_2^{-6}[F,a_0]\ldots [F,a_n]\;,\qquad
\mr{with}\;\,a_0,\ldots,a_n\in\A(S^4_q)\,,
$$
is of trace class on $\HH$ for all $n\geq 3$.
\end{lemma}
\begin{prova}
Let $L_q$ be the operator defined by $L_q\ket{l,m_1,m_2;j}_\pm:=
q^{-\frac{2}{3}m_1-2m_2}\ket{l,m_1,m_2;j}_\pm$. Notice that $L_qF=FL_q$
and that the map $a\mapsto L_qaL_q^{-1}$ is an automorphism of $\A(S^4_q)$;
indeed $L_qx_0L_q^{-1}=x_0$, $L_qx_1L_q^{-1}=q^{-2/3}x_1$
and $L_qx_2L_q^{-1}=q^{-2}x_2$.

Since $K_1^{-8}K_2^{-6}=L_q^3$, we have to prove that $L_q^3[F,a_0]\ldots [F,a_n]$
is of trace class for all $n\geq 3$. Since $[F,a_j]$ is of trace class,
it is sufficient to prove that $L_q^3[F,a_0][F,a_1][F,a_2]$ is a bounded
operator. This operator can be written as:
$$
L_q[F,L_q^2a_0L_q^{-2}]L_q[F,L_qa_1L_q^{-1}]L_q[F,a_2]=L_q[F,a_0']
L_q[F,a_1']L_q[F,a_2']\;.
$$
Thus, it is sufficient to prove that $L_q[F,a]$ is bounded for all
$a\in\A(S^4_q)$. By Leibniz rule and $[F,a^*]=-[F,a]^*$, it is
sufficient to consider the cases $a=x_0,x_1,x_2$. We know
that $[F,x_i]$ is a sum of independent weighted shifts with weights
bounded by $q^{2l}$.
Then, $L_q[F,x_i]$ is a sum of independent weighted shifts with weights
bounded by
$$
q^{2l}q^{-\frac{2}{3}m_1-2m_2}=q^{2(l+\frac{1}{2}-j-m_2)}q^{\frac{2}{3}(j-m_1)}
q^{\frac{4}{3}j-1}\leq q^{\frac{4}{3}j-1}\leq q^{-1/3}
$$
and this concludes the proof.
\end{prova}

We are in the hypothesis of Lemma \ref{lemma:ch}, with $K=K_1^8K_2^6$
and $p=4$. Thus, we can construct the `twisted' analogue of Chern-Connes
cocycle.

\begin{prop}
The multilinear map
$$
\mr{ch}^{F,\kappa}_n(a_0,\ldots,a_n):=
\tr_{\HH}(K_1^{-8}K_2^{-6}\gamma F[F,a_0]\ldots [F,a_n])
$$
is the representative on an $\mc{U}_q(so(5))$-invariant twisted Hochschild
$n$-cocycle for all even $n\geq 4$.
\end{prop}
\begin{prova}
From Lemma \ref{lemma:ch}, $\mr{ch}^{F,\kappa}_n$ is a twisted Hochschild
$n$-cocycle. We have only to prove invariance.

Recall that the action of $\mc{U}_q(so(5))$ on multilinear maps is
$$
(h\az\phi)(a_0,a_1,\ldots,a_n):=\phi(S^{-1}(h_{(n+1)})\az
a_0,S^{-1}(h_{(n)})\az a_1,\ldots,S^{-1}(h_{(1)})a_n)\;.
$$
The rules of the crossed product and the
invariance of $F$ and $\gamma$ ($Fh=hF$, $F\gamma=\gamma F$) imply:
$$
\gamma F[F,S^{-1}(h_{(n+1)})\az a_0]\ldots [F,S^{-1}(h_{(1)})\az a_n]
=S^{-1}(h_{(2)})(\gamma F[F,a_0]\ldots [F,a_n])h_{(1)}\;.
$$
Recall that $K_1^8K_2^6$ is the element implementing the square
of the antipode, $(K_1^8K_2^6)h(K_1^8K_2^6)^{-1}=S^2(h)$. Thus,
\begin{align*}
\bigl\{h\az\mr{ch}^{F,\kappa}_n\bigr\}(a_0,\ldots,a_n) &=
\tr\bigl(K_1^{-8}K_2^{-6}S^{-1}(h_{(2)})\gamma F[F,a_0]\ldots [F,a_n]h_{(1)}\bigr) \\
&= \tr\bigl(S(h_{(2)})K_1^{-8}K_2^{-6}\gamma F[F,a_0]\ldots [F,a_n]h_{(1)}\bigr)\;.
\end{align*}
Cyclicity of the trace and the property $h_{(1)}S(h_{(2)})=\epsilon(h)$
imply that $h\az\mr{ch}^{F,\kappa}_n=\epsilon(h)\mr{ch}^{F,\kappa}_n$
for all $h\in\mc{U}_q(so(5))$.
\end{prova}

\noindent
Note that while the Fredholm module $(\A(S^4_q),\HH,F,\gamma)$ is
$1$-summable (cf.~Proposition \ref{prop:notr}), and one has a
cocycle $\mr{ch}^F_n$ for any even $n\geq 0$, the twisted cocycle
$\mr{ch}^{F,\kappa}_n$ is defined only for even~$n$, greater than $4$.

Now, we prove that the Hochschild homology class of $\omega_4$
is not zero. To this end, it is enough to prove the
non-vanishing of the pairing
\begin{equation}\label{eq:comp}
\mr{ch}^{F,\kappa}_4(\omega_4)=
\tr_{\HH\otimes\C^4}(\eta K_1^{-8}K_2^{-6}\gamma F[F,e]^5)\;,
\end{equation}
which depends only on the class of $\omega_4$ and $\mr{ch}^{F,\kappa}_4$
in Hochschild (co)homology.

Along the lines of Lemma \ref{lemma:6.26}, one proves that the pairing
(\ref{eq:comp}) is $2$ times the $q$-index of the operator $F^+_e$ in (\ref{eq:Fpe}),
$$
\tfrac{1}{2}\,\mr{ch}_4^{F,\kappa}(\omega_4)=\tr_{\ker F^+_e}(\eta K_1^{-8}K_2^{-6})
-\tr_{\mr{coker}\,F^+_e}(\eta K_1^{-8}K_2^{-6}) \;.
$$
Using representation theory, it is possible to compute explicitly the index.

\begin{lemma}
The class $[\omega_4]\in H\!H^\kappa_4$ is non-trivial.
In particular, $\mr{ch}_4^{F,\kappa}(\omega_4)=2$.
\end{lemma}
\begin{prova}
We lift the representation of $U_q(so(5))$ from $\HH_\pm$ to
$v\in\HH_\pm\otimes\C^4$ by taking the Hopf tensor product
with $V_{1/2}$, as in (\ref{eq:RmodAz}):
$$
h\aaz v=\sigma\bigl(S^{-1}(h_{(1)})\bigr)^th_{(2)}\az v\;,
$$
with $\sigma$ the representation (\ref{eq:spin}). Since $e$
satisfies (\ref{eq:cov}), it follows that $h\aaz (ew)=e(h\aaz w)$,
and the subspaces
$$
\HH^{||}_\pm:=e(\HH_\pm\otimes\C^4)\;,\qquad
\HH^\perp_\pm:=(1-e)(\HH_\pm\otimes\C^4)\;,
$$
carries sub-representations of $U_q(so(5))$.
As $\HH_\pm\simeq\bigoplus_{l\in\N+\frac{1}{2}}V_l$,
Clebsh-Gordan decomposition
\begin{align*}
V_{1/2}\otimes V_{1/2} &\simeq V_0\oplus V_1\oplus V_{(1,1)} \;,\\
V_{1/2}\otimes V_l &\simeq V_{l-\frac{1}{2}}\oplus V_{l+\frac{1}{2}}
\oplus V_{(1,l-\frac{1}{2})}\oplus V_{(1,l+\frac{1}{2})}
\quad\mr{if}\;l=\tfrac{3}{2},\tfrac{5}{2},\ldots
\end{align*}
tells us that
$$
\HH^{||}_\pm\simeq\bigoplus_{l\in\N}(m_l^{\pm,||}V_l\oplus n_l^{\pm,||}V_{(1,l+1)}) \;,\qquad
\HH^{\perp}_\pm\simeq\bigoplus_{l\in\N}(m_l^{\pm,\perp}V_l\oplus n_l^{\pm,\perp}V_{(1,l+1)}) \;,
$$
with multiplicities satisfying
$$
m_0^{\pm,||}+m_0^{\pm,\perp}=1\;,\qquad
m_l^{\pm,||}+m_l^{\pm,\perp}=2\;\forall\;l\geq 1\;,\qquad
n_l^{\pm,||}+n_l^{\pm,\perp}=2\;\forall\;l\geq 0\;.
$$
Consider the vectors ($l\geq 1$)
\begin{align*}
v^\pm_{(1,l)} &:=\ket{l-\oh,\oh,l-\oh;\oh}_\pm\otimes(0,0,0,1)^t \;,\\
v^\pm_{(0,l)} &:=\ket{l-\oh,\oh,l-\oh;\oh}_\pm\otimes(0,1,0,0)^t
  +\ket{l-\oh,-\oh,l-\oh;\oh}_\pm\otimes(0,0,0,1)^t \;,\\
v^\pm_{(0,0)} &:=\ket{\oh,\oh,\oh;\oh}_\pm\otimes(1,0,0,0)^t
  +\ket{\oh,\oh,-\oh;\oh}_\pm\otimes(0,1,0,0)^t \\
&\qquad +\ket{\oh,-\oh,\oh;\oh}_\pm\otimes(0,0,1,0)^t
  +\ket{\oh,-\oh,-\oh;\oh}_\pm\otimes(0,0,0,1)^t \;.
\end{align*}
These are highest weight vectors with weight $(1,l)$ and $(0,l)$ respectively.
They remain highest weight vectors when multiplied by $e$. So,
\begin{equation}\label{eq:qIndHwv}
v^{\pm,||}_{(n_1,n_2)} :=e v^\pm_{(n_1,n_2)} \qquad\mr{and}\qquad
v^{\pm,\perp}_{(n_1,n_2)} :=(1-e)v^\pm_{(n_1,n_2)}
\end{equation}
are highest weight vectors with weight $(n_1,n_2)$, for $(n_1,n_2)=(1,l)$
and $(0,l)$. By using the explicit form of the chiral representations,
one finds that $v^{-,||}_{(0,0)}=ev^-_{(0,0)}=0$ and
$v^{+,||}_{(0,0)}=(1-e)v^+_{(0,0)}=0$,
while $v^{\pm,||}_{(n_1,n_2)}\neq 0$ and $v^{\pm,\perp}_{(n_1,n_2)}\neq 0$
in all other cases.

This proves that $m_l^{\pm,||}=m_l^{\pm,\perp}=n_l^{\pm,||}=n_l^{\pm,\perp}=1$
for all $l\geq 1$, while $m_0^{+,||}=m_0^{-,\perp}=1$ and $m_0^{-,||}=m_0^{+,\perp}=0$.
That is,
$$
\HH^{||}_s\simeq (\delta_{s,+}V_0)\oplus\bigoplus\nolimits_{l\geq 1}(V_l\oplus V_{(1,l)}) \;,\qquad
\HH^{\perp}_s\simeq (\delta_{s,-}V_0)\oplus\bigoplus\nolimits_{l\geq 1}(V_l\oplus V_{(1,l)}) \;,
$$
with $s=\pm 1$. By equivariance, $F$ intertwines representations with the same highest
weight, while flipping the label `$\pm$'.
The kernel and cokernel of $F^+_e$ are still subrepresentations of $U_q(so(5))$,
thus it is enough to determine which of the highest weight vectors
(\ref{eq:qIndHwv}) are in the kernel and which in the cokernel.
It turns out that $\,\ker F^+_e=V_0\,$ is the subspace of $\HH_+$ spanned by
$v^+_{(0,0)}$ ($eFev^+_{(0,0)}=eFv^+_{(0,0)}=ev^-_{(0,0)}=0$),
and the cokernel is empty. So,
\begin{equation*}
\rule{0pt}{20pt}
\tfrac{1}{2}\,\mr{ch}_4^{F,\kappa}(\omega_4)=
\tfrac{(v^+_{(0,0)}|\eta K_1^{-8}K_2^{-6}|v^+_{(0,0)})}{(v^+_{(0,0)}|v^+_{(0,0)})}=1 \;.
\qedhere
\end{equation*}
\end{prova}

\medskip

\section{Check with {\mate}}
As a further control that the left regular and chiral representations of
the algebra $\A(S^4_q)$ are correct, we list the {\mate} code which allows
to check the commutation rules of the algebra.
Before that, let me explain how the computer program works.

\subsection*{Left regular representation}
We call $X_i(a,b,l,m,k,j)$, $i=0,1$, the coefficients such that
$$
x_i\ket{l,m_1,m_2;j}=\sum\nolimits_{a,b}X_i(a,b,l,m_1,m_2,j)\ket{l+a,m_1+\delta_{i1},m_2;j+b}
$$
and notice that they are zero if $a,b\neq \pm 1$, while we
call $X_2(a,l,m_2,j)$ the coefficients in
$$
x_2\ket{l,m_1,m_2;j}=\sum\nolimits_{a}X_2(a,l,m_2,j)\ket{l+a,m_1,m_2+1;j}
$$
and notice that they vanish for $a\neq\pm 1$
(they are called~\texttt{Xzero[a,b,l,m,k,j]}, \texttt{Xuno[a,b,l,m,k,j]}
and \texttt{Xdue[a,l,k,j]} in {\mate}). We then define the seven polynomials
\begin{align*}
\texttt{Pol}_1&=x_0x_2-q^2x_2x_0 \;,\\
\texttt{Pol}_2&=x_1x_2-q^2x_2x_1 \;,\\
\texttt{Pol}_3&=x_2^*x_1-q^2x_1x_2^* \;,\\
\texttt{Pol}_4&=x_0x_1-q^2x_1x_0 \;,\\
\texttt{Pol}_5&=[x_1,x_1^*]+(1-q^4)x_0^2 \;,\\
\texttt{Pol}_6&=[x_2,x_2^*]+x_1^*x_1-q^4x_1x_1^* \;,\\
\texttt{Pol}_7&=x_0^2+x_1x_1^*+x_2x_2^*-1 \;.
\end{align*}
The vanishing of these polynomials guarantees that the commutation
relations of $\A(S^4_q)$ are satisfied (the remaining commutation
relations are obtained by conjugation).
In the code listed at the end of the chapter we call
$$
\texttt{Comm[i,x,y]}:=
\inner{l+x,m',k',j+y| \texttt{Pol}_i |l,m,k,j} \;,
$$
where $m',k'$ are the unique values of the coefficients
for which the matrix element is different from zero
($m',k'$ depend on which polynomial one is considering).
For $i=1,2,3$ we have to check that $\texttt{Comm[i,x,y]}=0$
in the cases $x=-2,0,2$ and $y=-1,+1$ (they are automatically
zero in all other cases). This is done by the code
\begin{footnotesize}\begin{verbatim}
  For[counter = 0; i = 1, i <= 3,
    For[x = -2, x <= 2,
      For[y = -1, y <= 1,
        If[Simplify[Comm[i, x, y]] == 0, Null, counter++, counter++];
        Print["i=", i, ", x=", x, ", y=", y, ", Comm=", Simplify[Comm[i, x, y]]];
      y++; y++];
    x++; x++];
  i++]; Print["errors=", counter]
\end{verbatim}\end{footnotesize}
For each of the values of $i,x,y$ we are interested in, it prints a line
with the value of $\texttt{Comm[i,x,y]}$ (which is
always zero). The counter `\texttt{counter}' counts cumulatively how many
times these coefficients fail to be zero, and is printed at the end of the
cycle (it turns out to be zero).

For $i=4,5,6,7$ we do the same, but now we are interested in the cases
$x=-2,0,2$ and $y=-2,0,2$ (in the other cases coefficients are automatically
zero).

\subsection*{Chiral representations}
Let
$$
\ket{l,j,m_1,m_2,\epsilon}=
\delta_{2\epsilon,(-1)^{l+\frac{1}{2}-j+m_2}}\ket{l,m_1,m_2;j}_\pm \;.
$$
We first do the computation as $\epsilon\in\{\pm\oh\}$ were an independent
variable; then we should evaluate the result on the actual value
$\epsilon=\oh(-1)^{l+\frac{1}{2}-j+m_2}$, but this is not necessary
since the result is zero even for the `wrong' value of $\epsilon$.
In {\mate} code we call \texttt{Xzero[a,b,c,l,j,m,k,e]} the coefficient
$X_0^{a,b,c}(l,j,m,k,\epsilon)$, etc., in the equation
$$
x_i\ket{l,j,m,k,\epsilon}=\sum_{c=\pm\frac{1}{2}}\sum_{a,b=0,\pm 1}
X_i^{a,b,c}(l,j,m,k,\epsilon)\ket{l+a,j+b,m+\delta_{i1},k+\delta_{i2},c} \,.
$$
They are given by
$$
X_i^{a,b,c}(l,j,m,k,\epsilon)=s(\epsilon,a-b+\delta_{i2},c)Y_i^{a,b}(l,j,m,k,\epsilon)\;,
$$
with
$$
s(\epsilon,t,c)=\begin{cases}
1 & \mathrm{if}\;c=(-1)^t \epsilon \;,\\
0 & \mathrm{otherwise}\;,
\end{cases}
$$
and
\begin{align*}
Y_0^{a,1}(l,j,m,k,\epsilon)  &={A}^+_{j,m}{C}^a_{l,j,k}(\epsilon) \;,\\
Y_0^{a,0}(l,j,m,k,\epsilon)  &={A}^0_{j,m}{H}^a_{l,j,k}(\epsilon) \;,\\
Y_0^{a,-1}(l,j,m,k,\epsilon) &=Y_0^{-a,1}(l+a,j-1,m,k,\epsilon) \;, \\
Y_1^{a,1}(l,j,m,k,\epsilon)  &={B}^+_{j,m}{C}^a_{l,j,k}(\epsilon) \;,\\
Y_1^{a,0}(l,j,m,k,\epsilon)  &={B}^0_{j,m}{H}^a_{l,j,k}(\epsilon) \;,\\
Y_1^{a,-1}(l,j,m,k,\epsilon) &={B}^-_{j,m}{C}^{-a}_{l+a,p-1,k}(\epsilon) \;,\\
Y_2^{a,b}(l,j,m,k,\epsilon)  &=\delta_{b,0}{D}_{l,j,k}^a(\epsilon) \;.
\end{align*}
The factor $s(\epsilon,a-b+\delta_{i2},c)$ takes into account the fact that
when $l-j+m_2$ is shifted by $a-b+\delta_{i2}$, the sign of $\epsilon$ changes
or not, depending on the parity of that quantity.

With $\texttt{Pol}_i$ given by the same polynomials as in the previous section, we call
$$
\texttt{Comm[i,x,y,e]}:=
\inner{l+x,j+y,m',k',e'| \texttt{Pol}_i |l,j,m,k,e} \;,
$$
where $m',k',e'$ are the unique values of the coefficients
for which the matrix elements is different from zero
(they depend on which polynomial one is considering).
For $i=1,2,3$ we have to check that $\texttt{Comm[i,x,y,e]}=0$
in the cases $x=-2,-1,0,1,2$, $y=-1,0,+1$ and $e=\pm\oh$ (they are automatically
zero in all other cases). This is done by the code
\begin{footnotesize}\begin{verbatim}
  For[counter = 0; i = 1, i <= 3,
   For[e = -1/2, e <= 1/2,
    For[x = -2, x <= 2,
     For[y = -1, y <= 1,
      If[Simplify[Comm[i, x, y, e]] == 0, Null, counter++, counter++];
     y++];
    x++];
   e++];
  i++]; Print["errors=", counter]
\end{verbatim}\end{footnotesize}
As for the left regular representation, the value of each coefficient
$\texttt{Comm[i,x,y,e]}$ could be printed, but we don't do it to avoid
wasting paper.
Instead, only the number of coefficients failing to be zero is printed at the end.
For $i=4,5,6,7$ we do the same, but now we are interested in the cases
$x,y=-2,-1,0,1,2$ and $e=\pm\oh$ (in the other cases coefficients are
trivially zero). As one can see by running the code, the result
is always zero, confirming that the representation is correct.

\pagebreak

\begin{center}
\begin{tabular}{|p{\tabella}|}
\hline
\rule{0pt}{25pt}\hspace*{1.5cm}%
\texttt{{\mate} code -- left regular representation} \\
\hspace*{-5mm}%
\includegraphics[width=\textwidth]{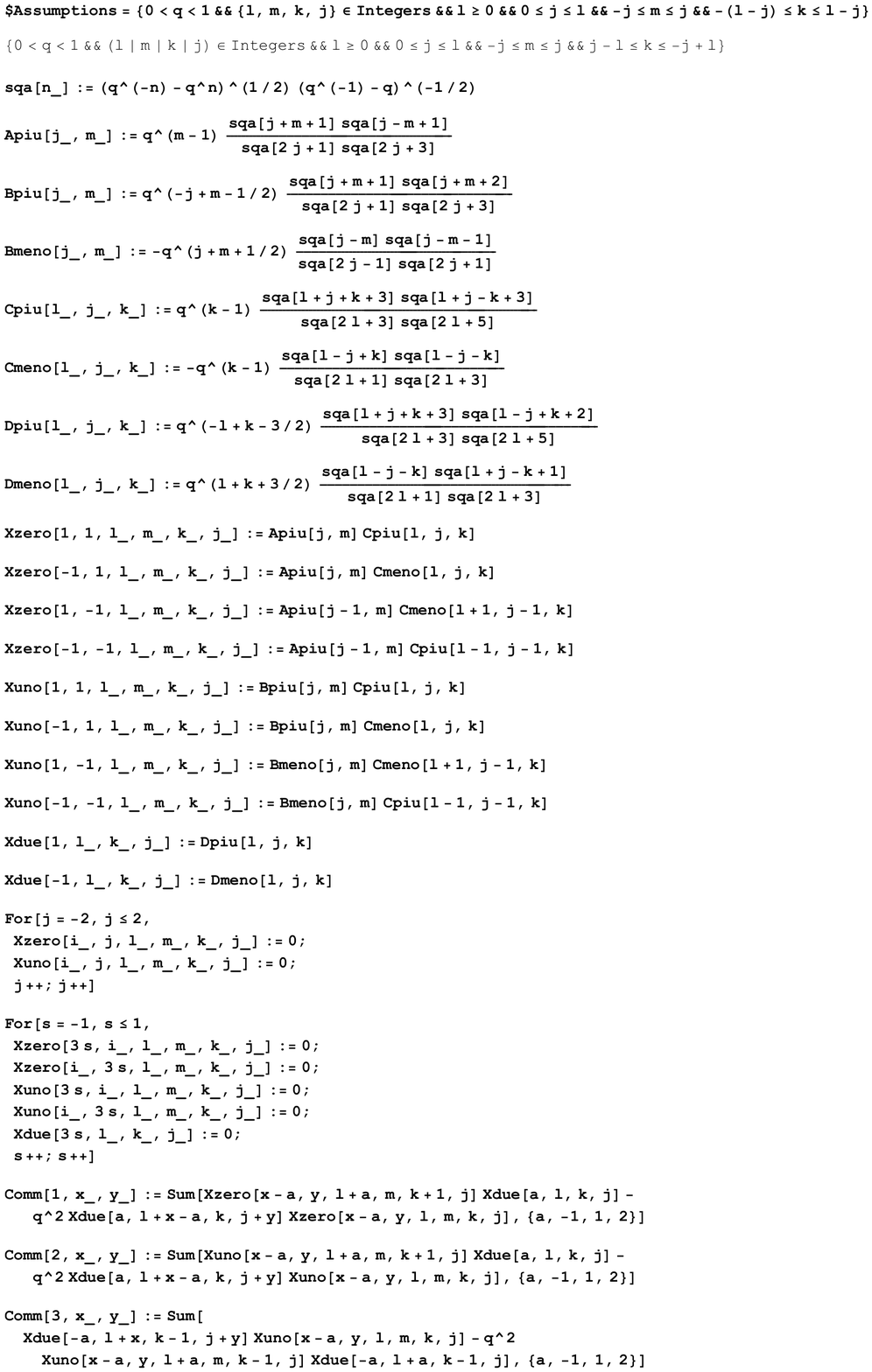} \\
\hline
\end{tabular}

\begin{tabular}{|p{\tabella}|}
\hline
\hspace*{-5mm}%
\includegraphics[width=\textwidth]{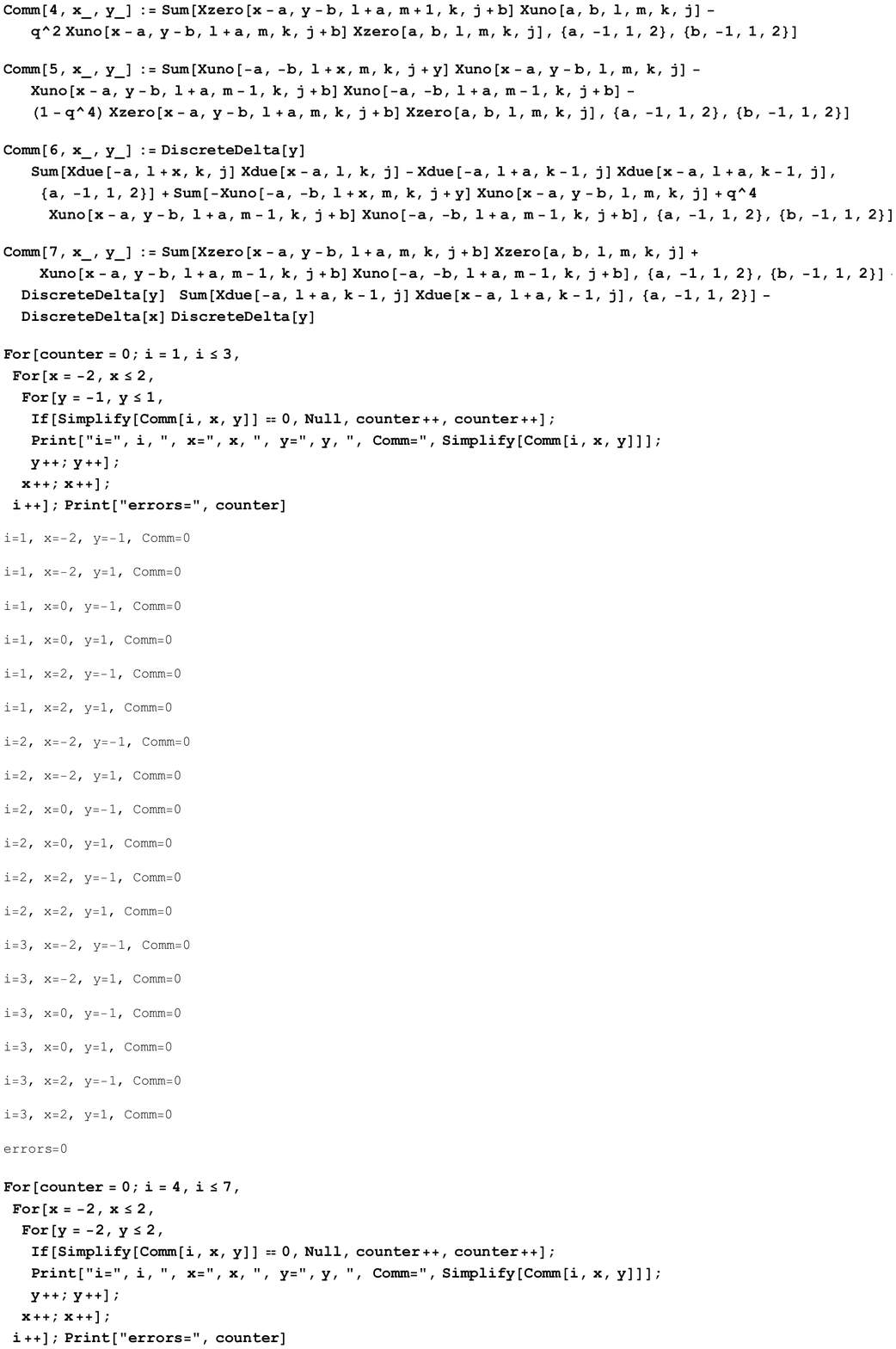} \\
\hline
\end{tabular}

\begin{tabular}{|p{\tabella}|}
\hline
\hspace*{-5mm}%
\includegraphics[width=\textwidth]{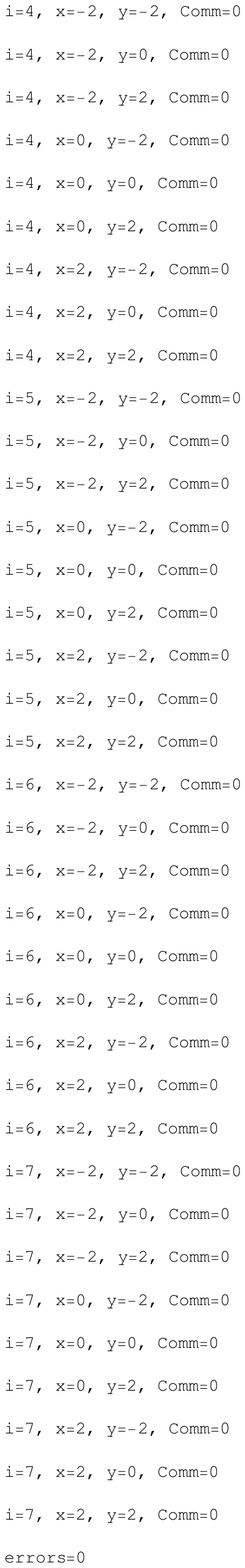} \\
\hline
\end{tabular}
\end{center}

\begin{center}
\begin{tabular}{|p{\tabella}|}
\hline
\rule{0pt}{25pt}\hspace*{1.5cm}%
\texttt{{\mate} code -- chiral representations} \\
\hspace*{-5mm}%
\includegraphics[width=\textwidth]{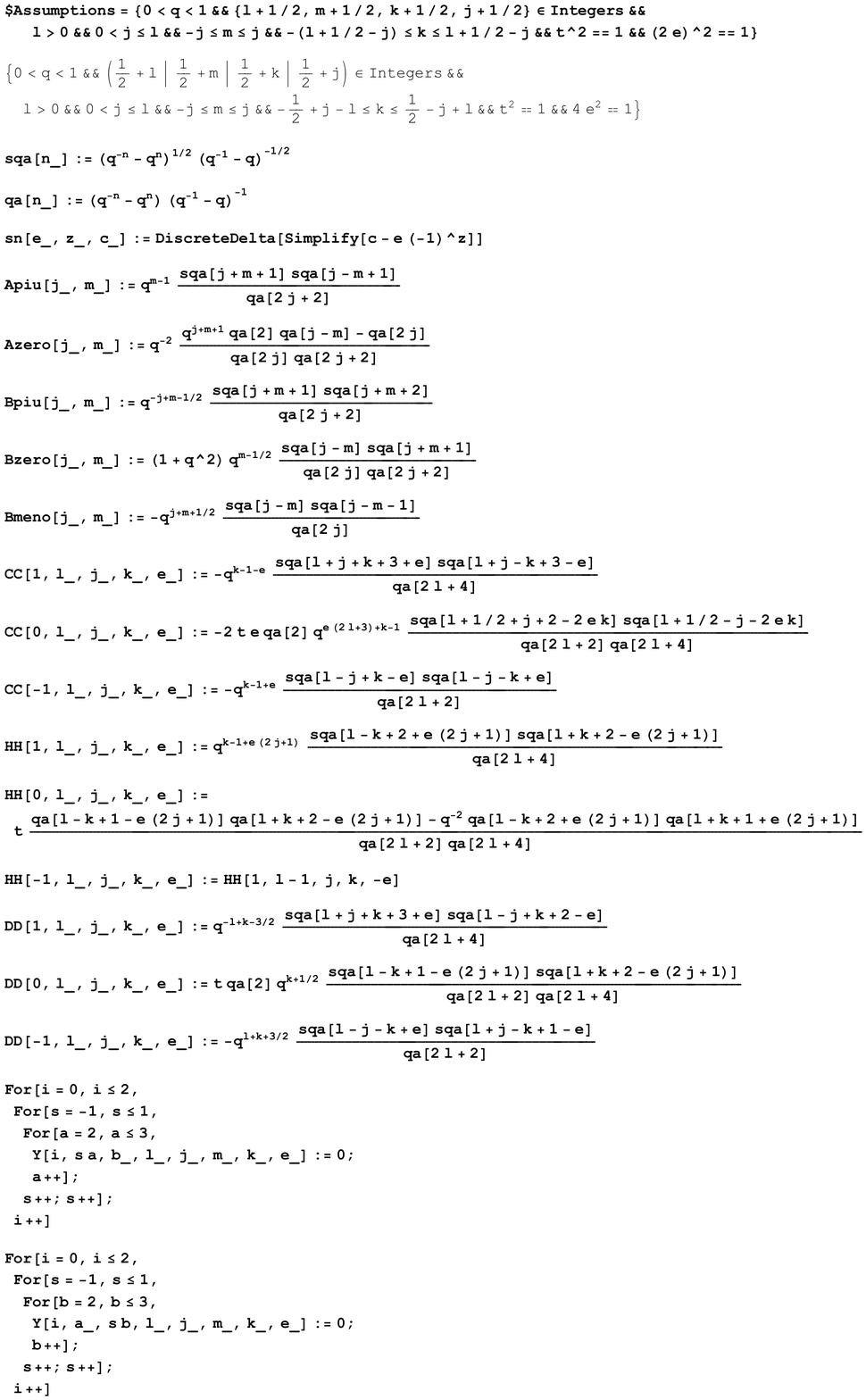} \\
\hline
\end{tabular}

\begin{tabular}{|p{\tabella}|}
\hline
\hspace*{-5mm}%
\includegraphics[width=\textwidth]{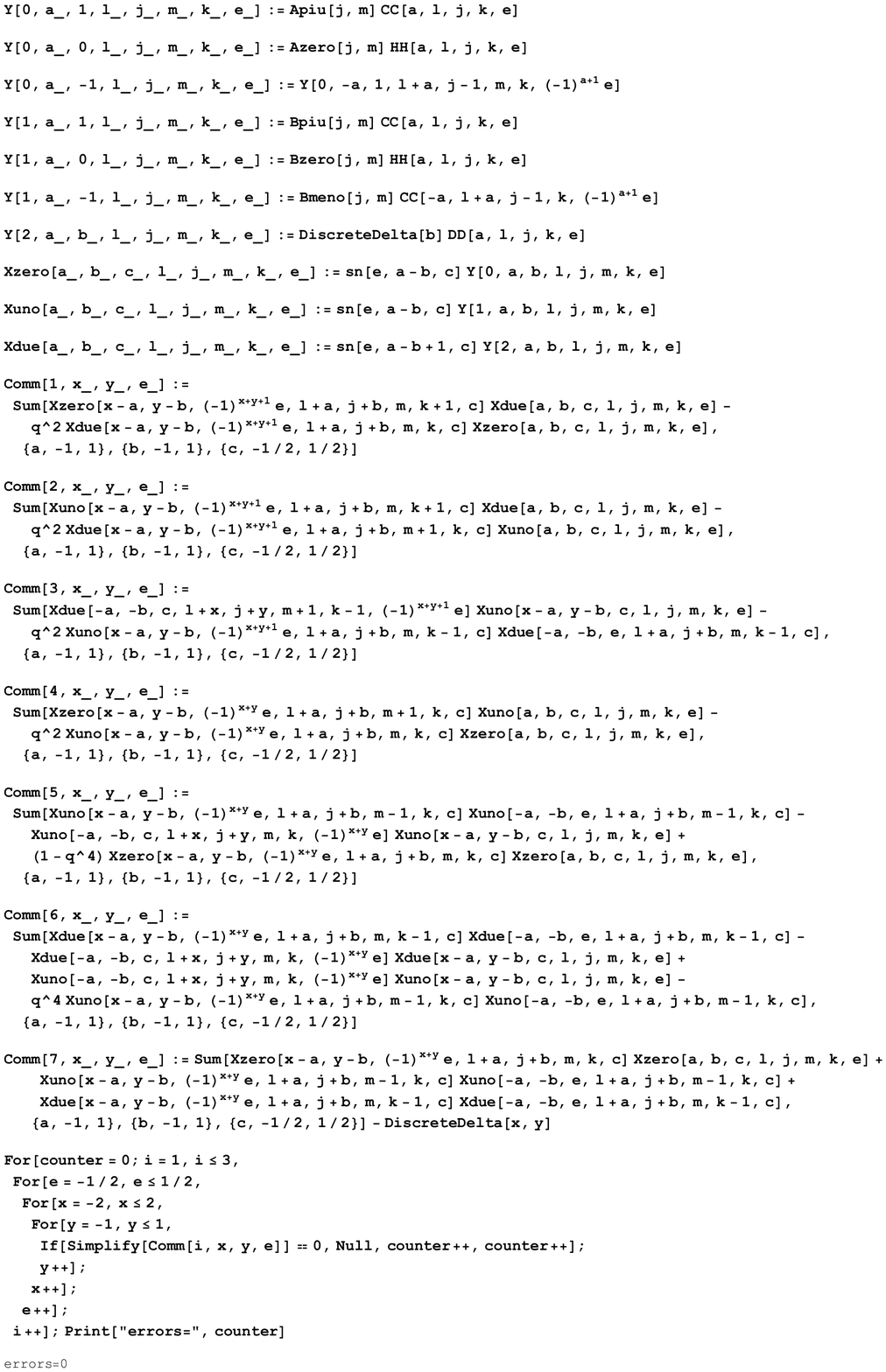} \\
\hline
\end{tabular}

\begin{tabular}{|p{\tabella}|}
\hline
\rule{0pt}{4.3cm}\hspace*{-5mm}%
\smash[t]{\includegraphics[width=\textwidth]{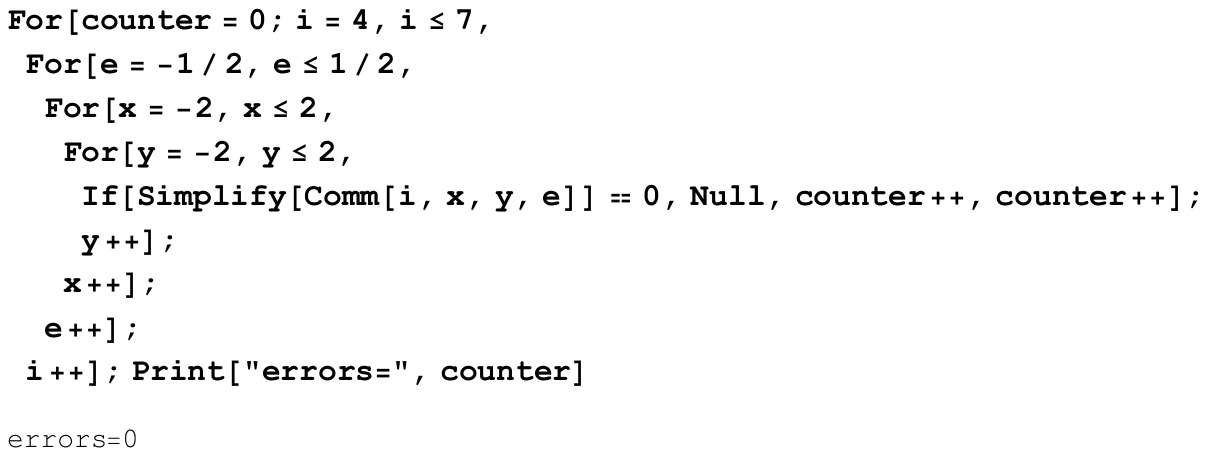}} \\
\hline
\end{tabular}
\end{center}


\typeout{Capitolo 5}

\chapter{Some remarks on odd-dimensional quantum spheres}\label{chap:CPlq}

\afo{5.5cm}{One should always generalize}
           {Carl Jacobi~\cite{PH81}.}

\noindent
The purpose of this chapter is to extend to odd-dimensional quantum spheres
some of the results which hold in the case of $SU_q(2)$~\cite{Con04,vSDL05},
Podle\'s spheres~\cite{DDLW07} and of the quantum $4$-sphere \cite{DDL06}.
We recall here the definition of the odd-dimensional quantum sphere
$S^{2\ell+1}_q$ of~\cite{VS91},
which for $\ell>1$ and $0<q<1$ is described by the $*$-algebra
$\A(S^{2\ell+1}_q)$ generated by $\{z_i,z_i^*\}_{i=1,\ldots,\ell+1}$ with relations
\begin{align*}
z_iz_j &=qz_jz_i &&\forall\;i<j \;,\\
z_i^*z_j &=qz_jz_i^* &&\forall\;i\neq j \;,\\
[z_1^*,z_1] &=0 \;,\\
[z_{i+1}^*,z_{i+1}] &=(1-q^2)\sum\nolimits_{j=1}^i z_jz_j^* &&\forall\;i=1,\ldots,\ell \;,\\
z_1z_1^*+z_2z_2^* &+\ldots+z_{\ell+1}z_{\ell+1}^*=1 \;.
\end{align*}
In the next section we present the left regular representation
of $\A(S^{2\ell+1}_q)$ and the `optimal' Dirac operator $D$ constructed
in \cite{CP07}.
Then, in Section \ref{sec:9.2} we define  a representation which
generalizes to all odd-dimensional sphere the `symbol map' of
$S^3_q=SU_q(2)$ of \cite{Con04,vSDL05} and which
differs from the left regular representation only by infinitesimals,
that in the present case means operators whose noncommutative integral
is zero. 
Finally, in Section \ref{sec:9.3} we use the symbol map to compute the
noncommutative integral associated with $D$. Section \ref{app:uno}
contains the proof of the main result of the chapter, that is Theorem \ref{thm:main}.
We argue that these results may be used to construct spectral
triples on quantum projective spaces, which are the `virtual'
spaces underlying the $*$-algebra $\A(\CP^\ell_q)\subset\A(S^{2\ell+1}_q)$
generated by the elements $p_{ij}:=z_i^*z_j$.

\section{The left regular representation}\label{sec:9.1}
We use, with minor changes, the notations of \cite{CP07}.
Irreducible unitary corepresentations of $\A(SU_q(\ell+1))$ are labelled by Young
tableaux $n=(n_1,\ldots,n_\ell)$, that are $\ell$-tuple of non-negative integers satisfying
$n_1\geq n_2\geq\ldots\geq n_\ell\geq 0$. We denote
$V_n$ the corresponding representation space.
Basis vectors of $V_n$ can be indexed by Gelfand-Tsetlin tableaux
(GT tableaux, for short), which are arrays of integers of the form
\begin{equation}\label{eq:GTT}
\GT{r}=\left(\begin{array}{cccccc}
r_{1,1} & r_{1,2} & \ldots & r_{1,\ell-1} & r_{1,\ell} & 0 \\
r_{2,1} & r_{2,2} & \ldots & r_{2,\ell-1} & r_{2,\ell} \\
r_{3,1} & r_{3,2} & \ldots & r_{3,\ell-1} \\
\ldots & \ldots & \ldots \\
r_{\ell,1} & r_{\ell,2} \\
r_{\ell+1,1}
\end{array}\right)
\end{equation}
with $r_{1,i}=n_i$ and $r_{i,j}\geq r_{i+1,j}\geq r_{i,j+1}$ for all
$i,j$ (from now on it is understood that $r_{ij}=0$ if $i,j$ are out of
range).
Usually GT tableaux are defined modulo a global rescaling
(two arrays are equivalent if they differ by a constant);
here for each equivalence class of GT tableaux, we choose the representative
which has zero in the top-right corner, $r_{1,\ell+1}=0$.

If $\HH$ is the Hilbert space completion of $\A(S^{2\ell+1}_q)$ with respect
to the inner product induced by the Haar state, an orthonormal basis
$\ket{n,h;\GT{r}}$ of $\HH$ is indexed by $n,h\in\N$ and
$\GT{r}$ is a GT tableau with Young tableau
$(r_{1,1},\ldots,r_{1,\ell})=(n+h,h,\ldots,h)$.

For $k=1,\ldots,\ell+1$ we denote
$$
\Gamma^{(k)}:=\big\{(m_1,\ldots,m_k)\in\N^k\,\big|\,1\leq m_i\leq \ell+2-i\;\forall\;i=1,\ldots,k\big\}\;.
$$
For $m\in\Gamma^{(k)}$ and $\GT{r}$ a GT tableau, we denote $\GT{r}^{(m)}$
the array with elements
$$
r^{(m)}_{ij}:=-\delta_{m_1,\ell+1}+\begin{cases}
r_{ij}+1 & \mr{if}\;j=m_i,\;1\leq i\leq k,\\
r_{ij}   & \mr{otherwise}.
\end{cases}
$$
The factor $-\delta_{m_1,\ell+1}$ is needed to select the representative of
the GT tableau equivalence class which has $r_{1,\ell+1}=0$.

The left regular representation $\pi:\A(S^{2\ell+1}_q)\to\B(\HH)$ can be
determined by using the inclusion into $\A(SU_q(\ell+1))$ defined by
$z_{\ell+2-i}^*\mapsto q^{-i+1}u^1_i$, with $u^i_j$ the generators of
the quantum $SU(\ell+1)$ group (cf.~\cite{CP07}). It is defined on generators by
%
\begin{align}\label{eq:Lrep}
\pi(z_{\ell+2-i}^*)\ket{n,h;\GT{r}}&=
q^{-\frac{1}{2}(\ell+h)}\frac{[n+1]^{\frac{1}{2}}}{[n+h+\ell+1]^{\frac{1}{2}}}
\sum_{m\in\Gamma^{(i)}:\,m_1=1}q^{\ell+h}q^{-i+1}
C_q(i,\GT{r},m)\ket{n+1,h;\GT{r}^{(m)}} \nonumber \\ &
\, +q^{-\frac{1}{2}n}\frac{[h+\ell-1]^{\frac{1}{2}}}{[n+h+\ell-1]^{\frac{1}{2}}}
\sum_{m\in\Gamma^{(i)}:\,m_1=\ell+1}q^{-i+1}
C_q(i,\GT{r},m)\ket{n,h-1;\GT{r}^{(m)}}
\end{align}
where we set to zero the summands for which $\GT{r}^{(m)}$ is \emph{not}
a GT tableau.
The coefficients $C_q(i,\GT{r},\GT{r}^{(m)})$ are the
Clebsch-Gordan coefficients in equation (5.3) of \cite{CP07} (their sign
can be found at page 220 of \cite{KS97}). Their explicit value is
\begin{align}\label{eq:coe} &
C_q(i,\GT{r},m):=
q^{\frac{1}{2}\left(1+m_1-2m_i-r_{1,m_1}+2r_{i,m_i}-r_{i,\ell+2-i}
+\sum_{k=1}^{\ell+1-i}(r_{i+1,k}-r_{i,k})\right)} \;
\left\{\prod_{k=1}^{i-1}\mr{sign}(m_k-m_{k+1})\right\}
\times\nonumber \\ & \;\,\times
\left\{\prod_{k=1}^{i-1}
\prod_{m_k\neq j=1}^{\ell+2-k}\frac{[r_{k,j}-r_{k+1,m_{k+1}}-j+m_{k+1}]}{[r_{k,j}-r_{k,m_k}-j+m_k]}
\prod_{m_{k+1}\neq j=1}^{\ell+1-k}\frac{[r_{k+1,j}-r_{k,m_k}-j+m_k-1]}{[r_{k+1,j}-r_{k+1,m_{k+1}}-j+m_{k+1}-1]}
\right\}\rule{-4pt}{18pt}^{1/2} \! \times\nonumber \\ & \;\,\times
\left(\frac{\prod_{k=1}^{\ell+1-i}[r_{i+1,k}-r_{i,m_i}-k+m_i-1]}
{\prod_{m_i\neq k=1}^{\ell+2-i}[r_{i,k}-r_{i,m_i}-k+m_i]}
\right)\rule{-4pt}{18pt}^{1/2} \! \;,
\end{align}
where we set $\,\mr{sign}(0):=1\,$, empty sums are set to $0$
and empty products are set to $1$.
This representation together with the `optimal' Dirac
operator $D$ of \cite{CP07}, given by\footnote{We shifted the spectrum of $|D|$ in
\cite{CP07} by $+1$ to get an invertible Dirac operator.}
\begin{equation}\label{eq:optD}
D\ket{n,h;\GT{r}}=\begin{cases}
-(h+1)\ket{n,h;\GT{r}} & \mr{if}\;n=0\,,\\
(n+h+1)\ket{n,h;\GT{r}} & \mr{if}\;n>0\,,
\end{cases}
\end{equation}
defines a spectral triple on $S^{2\ell+1}_q$. A peculiarity of this
Dirac operator is that for any other equivariant Dirac operator $D'$
(on the same Hilbert space) there exists two positive constants $a,b$
such that $|D'|\leq a|D|+b$, which explains the name `optimal'.
This essentially means that $D$ is the equivariant Dirac operator
corresponding to the minimal metric dimension.

From an analysis of the case $\ell=1$, where the bounded commutator
condition fails for $q=1$ (cf.~remarks 5.1 and 5.5 in \cite{CP03a}),
we argue that these spectral triples constructed with `scalar functions'
instead of spinors are purely quantum objects.

Notice that the coefficients in front of the sums in (\ref{eq:Lrep}) are
never greater than $1$. Moreover as proved in \cite[eq.~(4.10)]{CP07},
$C_q(i,\GT{r},m)$ are bounded by a constant (which is independent of $n,h,\GT{r}$).
Then, the coefficients in the first line of (\ref{eq:Lrep}) are bounded by $q^h$.
If we neglect terms bounded by $q^h$, the first line of
(\ref{eq:Lrep}) disappears.
We'll show in Section \ref{sec:9.3} that this `approximation modulo $q^h$'
yields a remarkable simplification of the left regular representation.

\section{The symbol map}\label{sec:9.2}
Although it is a very simple example, the equatorial Podle\'s sphere
seems to be the building block for the theory of $q$-deformed spaces, when doing
local computations. In many examples (cf.~\cite{Con04,vSDL05,DDLW07,DDL06},
the last two being discussed also in Chapters \ref{chap:S2qs} and \ref{chap:S4q} 
of this dissertation) one finds that representations of quantum spaces are described,
modulo `infinitesimals', by suitable tensor products of representations of the
coordinate algebra $\A(S^2_q)$ of the equatorial Podle\'s sphere. 
This is indeed a general feature of all odd-dimensional quantum
spheres, as we shall see later on.

Recall that $\A(S^2_q)$ is the $*$-algebra generated by $A=A^*$, $B$ and
$B^*$ with relations:
$$
AB=qBA \;,\qquad
BB^*+A^2=1 \;,\qquad
B^*B+q^2A^2=1 \;.
$$

\begin{prop}
Let $u$ be the unitary generator of $\A(S^1)$.
There is a $*$-algebra morphism $\varphi:\A(S^{2\ell+1}_q)\to\A(S^1)\otimes\A(S^2_q)^{\otimes\ell}$,
given by ($2\leq i\leq\ell$) \vspace{-5pt}
\begin{align*}
\varphi(z_1)&:=u\,\otimes\stackrel{\ell\;\mr{times}}{\overbrace{A\otimes\ldots\otimes A}} \;, \\
\varphi(z_i)&:=\;\stackrel{i-1\;\mr{times}}{\overbrace{1\otimes\ldots\otimes 1}}
\otimes\,B\,\otimes\stackrel{\ell+1-i\;\mr{times}}{\overbrace{A\otimes\ldots\otimes A}}\, \;, \\
\varphi(z_{\ell+1})&:=\;\stackrel{\ell\;\mr{times}}{\overbrace{1\otimes\ldots\otimes 1}}\otimes\, B \;.
\end{align*}
\end{prop}
\begin{prova}
One proves by direct computation that the elements $\{\varphi(z_j),\varphi(z_j)^*\}_{j=1,\ldots,\ell+1}$
satisfy all the defining relations of $\A(S^{2\ell+1}_q)$.
\end{prova}

Let $\Lambda$ be the set of $(n,h;\mathbf{a},\mathbf{b})$ with
$n,h\in\N$, $\mathbf{a}=(a_1,\ldots,a_{\ell-1})\in\Z^{\ell-1}$ and
$\mathbf{b}=(b_0,\ldots,b_{\ell-1})\in\Z^\ell$
satisfying the constraints
\begin{subequations}
\begin{gather}
0\leq a_1\leq a_2\leq\ldots\leq a_{\ell-1}\leq n\;,\\
b_0\leq b_1\leq b_2\leq\ldots\leq b_{\ell-1}\leq h\;.
\end{gather}
\end{subequations}
We call $\hat{\mc{H}}$ the Hilbert space with orthonormal basis
$\kkett{n,h;\mathbf{a},\mathbf{b}}$ labelled by
$(n,h;\mathbf{a},\mathbf{b})\in\Lambda$.

For $i=1,\ldots,\ell+1$, we denote by $\underline{i}\in\{0,1\}^\ell$ the array
$$
\underline{i}:=
(\,\,\stackrel{i-1\;\mr{times}}{\overbrace{0,0,\ldots,0}}\,,
\stackrel{\ell+1-i\;\mr{times}}{\overbrace{1,1,\ldots,1}}) \;.
$$
The proof of the following proposition is a simple computation.

\begin{prop}
The formula
$$
\rho_1(u)\kkett{n,h;\mathbf{a},\mathbf{b}}=\kkett{n,h+1;
\mathbf{a},\mathbf{b}+\underline{1}\,} \;,
$$
defines a $*$-representation $\rho_1:\A(S^1)\to\B(\hat{\mc{H}})$.
For $2\leq i\leq\ell+1$, the formul{\ae} (here $b_\ell:=h$)
\begin{align*}
\rho_i(A)\kkett{n,h;\mathbf{a},\mathbf{b}}
 &=q^{b_{i-1}-b_{i-2}}\kkett{n,h;\mathbf{a},\mathbf{b}} \;,\\
\rho_i(B)\kkett{n,h;\mathbf{a},\mathbf{b}}
 &=\sqrt{1-q^{\smash[t]{2(b_{i-1}-b_{i-2}+1)}}}
\kkett{n,h+1;\mathbf{a},\mathbf{b}+\underline{i}\,} \;,
\end{align*}
define a $*$-representation $\rho_i:\A(S^2_q)\to\B(\hat{\mc{H}})$.
Distinct representations 
are mutually commuting.
\end{prop}

Modulo multiplicities, each representation $\rho_i$ ($i\neq 1$) is
equivalent to the 
representation of
$\A(S^2_q)=\A(S^2_{q^{1/2}1})$ given by (\ref{eq:murepA}).
By applying $\rho_1\otimes\rho_2\otimes\ldots\otimes\rho_{\ell+1}$
to $\varphi(a)$, $a\in\A(S^{2\ell+1}_q)$, and multiplying all the
factors in the tensor product we get a $*$-representation
$\rho:\A(S^{2\ell+1}_q)\to\B(\hat{\mc{H}})$.
We want to compare this representation with the left regular representation $\pi$.
As one can easily guess, they agree modulo `infinitesimals'.

Notice that since the spectrum of $\rho_i(A)$ is positive, the image of each
$\rho_i$ ($i\neq 1$) is a noncommutative disk $\A(D_q^2)$. Thus the map
$\rho$ has image in $\A\bigl(S^1\times (D^2_q)^{\times\ell}\bigr)$.
Its extension to order zero pseudo-differential operators would be the
higher-dimensional analogue of the \emph{symbol map} of \cite{Con04} (cf.~also
\cite{DLS05}). The image of such extension would generalize to arbitrary $\ell$
the \emph{cosphere bundle} of $S^3_q=SU_q(2)$ worked out in \cite{Con04,DLS05}.

\section{The noncommutative integral}\label{sec:9.3}
We have a representation $\pi:\A(S^{2\ell+1}_q)\to\B(\HH)$, where
$\HH\simeq L^2(S^{2\ell+1}_q)$ has an orthonormal basis $\ket{n,h;\GT{r}}$
labelled by GT tableau $\GT{r}$ with Young tableau $(r_{1,1},\ldots,r_{1,\ell})=(n+h,h,\ldots,h)$.
In this basis, the operators $\pi(z_i)$ are given by (\ref{eq:Lrep}).
The space $\HH$ can be included in $\hat{\HH}$ as follows.

A careful look at the relations defining a GT tableau shows that
any GT tableau with Young tableau $(r_{1,1},\ldots,r_{1,\ell})=(n+h,h,\ldots,h)$
has $h=r_{1k}\geq r_{jk}\geq r_{1,j+k-1}=h$, for each $k\geq 2$
and $j+k\leq\ell+1$, i.e.~it is of the form
\begin{equation}\label{eq:form}
\GT{r}=\left(\begin{array}{cccccc}
n+h     & h       & \ldots & h            & h          & 0 \\
r_{2,1} & h       & \ldots & h            & r_{2,\ell} \\
r_{3,1} & h       & \ldots & r_{3,\ell-1} \\
\ldots  & \ldots  & \ldots \\
r_{\ell-1,1} & h  & r_{\ell-1,3} \\
r_{\ell,1} & r_{\ell,2} \\
r_{\ell+1,1}
\end{array}\right) \;,
\end{equation}
and we have $\bigl(n,h;\mathbf{a}(\GT{r}),\mathbf{b}(\GT{r})\bigr)\in\Lambda$,
where we called
$$
a_j(\GT{r}):=r_{\ell+1-j,1}-h \;,\qquad
b_k(\GT{r}):=h-r_{\ell+1-k,k+1} \;.
$$
On the other hand, if $(n,h;\mathbf{a},\mathbf{b})\in\Lambda$ satisfy
the additional constraints
\begin{equation}\label{eq:add}
b_1\geq 0\quad\mr{and}\quad b_0\geq -a_1\;,
\end{equation}
then the array
\begin{equation}\label{eq:lbchange}
\GT{r}(n,h;\mathbf{a},\mathbf{b}):=\left(\begin{array}{cccccc}
n+h          & h       & \ldots & h            & h          & 0 \\
a_{\ell-1}+h & h       & \ldots & h            & h-b_{\ell-1} \\
a_{\ell-2}+h & h       & \ldots & h-b_{\ell-2} \\
\ldots       & \ldots  & \ldots \\
a_2+h        & h       & h-b_2 \\
a_1+h        & h-b_1 \\
h-b_0
\end{array}\right)
\end{equation}
is a GT tableau. We define a partial isometry $W:\HH\to\hat{\mc{H}}$ and its adjoint
$W^*:\hat{\mc{H}}\to\HH$ by
\begin{align*}
W\ket{n,h;\GT{r}} &=\kkett{n,h;\mathbf{a}(\GT{r}),\mathbf{b}(\GT{r})} \\
W^*\!\kkett{n,h;\mathbf{a},\mathbf{b}} &=\bigg\{\begin{array}{ll}
\ket{n,h;\GT{r}(n,h;\mathbf{a},\mathbf{b})}\quad &\mr{if}\;(\mathbf{a},\mathbf{b})
\;\textrm{satisfy (\ref{eq:add})}, \\
0 &\mr{otherwise}\;.
\end{array}
\end{align*}
Clearly, $W^*W=id_{\HH}$.

Next theorem is the main result of this chapter. In this chapter,
we call \emph{shift operator} an operator on $\HH$ which, in matrix form in
the basis $\ket{n,h;\GT{r}}$, has at most one matrix element different from zero
in each row, and at most one matrix element different from zero in each column.
Its matrix elements will be called \emph{weights} of the shift.

\begin{thm}\label{thm:main}
For all $i=1,\ldots,\ell+1$, the operator $\pi(z_i)-W^*\rho(z_i)W$ is
a finite sum of shift operators on $\HH$ with weights bounded by $q^h$.
\end{thm}

\noindent
Since the proof is quite long, we relegate it to Section \ref{app:uno}.

For $\ell>1$, terms which we are neglecting (shifts with weights bounded by $q^h$)
are not smoothing operators for the `optimal' Dirac operator $D$ in (\ref{eq:optD}).
Hence this approximation does not allow us to compute the dimension spectrum
and the local index formula on $S^{2\ell+1}_q$.

The situation is different if we consider the subalgebra $\A(\CP^\ell_q)$,
generated by the elements $p_{ij}:=z_i^*z_j$.
For any fixed $N\in\N$, let $\HH_N$ be the Hilbert subspace of $\HH$ which is the
closure of the linear span
$$
\mr{Span}\big\{ \ket{n+N,n;\GT{r}}\,\big|\,n\geq 0 \big\} \;.
$$
As one can easily check on generators, $\HH_N\subset\HH$ is an invariant subspace
for $\A(\CP^\ell_q)$. If we call $\pi_N(a)$ the restriction of $\pi(a)$ to $\HH_N$,
for all $a\in\A(\CP^\ell_q)$, we get a bounded $*$-representation $\pi_N:\A(\CP^\ell_q)
\to\B(\HH_N)$. A corollary of Theorem \ref{thm:main} is the following.

\begin{cor}\label{cor}
For any generator $p_{ij}$ of $\A(\CP^\ell_q)$, $\pi_N(p_{ij})-W^*\rho(p_{ij})W$
is a finite sum of shift operators on $\HH_N$ with weights bounded by $q^n$
(which is a rapid decay matrix, and then a smoothing operator for $D$).
\end{cor}

With this, it is possible (but very involved) to construct Fredholm
modules and spectral triples for quantum projective spaces.

We now use Theorem \ref{thm:main} to compute the noncommutative integral
on $S^{2\ell+1}_q$.
As a first step, we show that $|D|^{-2\ell-1}$ is in the
Dixmier ideal $\mc{L}^{1,\infty}(\HH)$ and has non-vanishing noncommutative integral;
this proves that the metric dimension is $2\ell+1$ (as argued in \cite{CP07}).

\begin{prop}\label{prop:Dix}
The spectral triple $(\A(S^{2\ell+1}),\HH,D)$ has metric dimension $2\ell+1$,
that is $|D|^{-2\ell-1}\in\mc{L}^{1,\infty}(\HH)$ and $\tr_\omega(|D|^{-2\ell-1})=C_\ell\neq 0$.
The normalization constant $C_\ell$ is given by
\begin{equation}\label{eq:Cl}
C_\ell:=\frac{1}{\ell!}\sum_{i=0}^{\ell-1}\frac{(-1)^i}{i!\,(\ell+i)!} \;.
\end{equation}
\end{prop}

\begin{prova}
From Weyl's character formula \cite{IN66} we know that the dimension of the
irreducible corepresentation of $SU_q(\ell+1)$ with Young tableau
$(n_1,\ldots,n_\ell)$ is
\begin{equation}\label{eq:Weyl}
\frac{\prod_{1\leq j\leq k\leq\ell}\left(k-j+1
+\sum\nolimits_{i=j}^k\lambda_i\right)}{\prod\nolimits_{k=1}^\ell k!} \;,
\end{equation}
where $(\lambda_1,\ldots,\lambda_\ell)$ is the highest weight of the representation,
given by $\lambda_i=n_i-n_{i+1}$.
The multiplicity $\mu_k$ of the eigenvalue $k:=n+h+1$ of $|D|$ is given by
previous expression computed for $(n_1,n_2,\ldots,n_\ell)
=(k-1,h,\ldots,h)$, i.e.~$\lambda_i=(k-h-1)\delta_{i,1}+h\delta_{i,\ell}$,
and summed for $h=0,\ldots,k-1$, that is
\begin{align*}
\mu_k &=\sum_{h=0}^{k-1}\frac{\prod_{2\leq i\leq j\leq\ell-1}(j-i+1)}
  {\prod\nolimits_{j=1}^\ell j!}(\ell+k-1)
  \prod_{j=1}^{\ell-1}\big\{j+k-h-1\big\}
  \prod_{j=2}^\ell\big\{ \ell+1-j+h \big\} \\
&=\frac{k+\ell-1}{\ell!(\ell-1)!}\sum_{h=0}^{k-1}
  \prod_{j=1}^{\ell-1}(j+k-h-1)(j+h) \\
&=\frac{k+\ell-1}{\ell!(\ell-1)!}\sum_{m=1}^k
  \prod_{j=1}^{\ell-1}\bigl\{ (j+k)(j-1)+(k+1)m-m^2 \bigr\} \\
&\sim\frac{k+\ell-1}{\ell!(\ell-1)!}\sum_{m=1}^k
  (k+1-m)^{\ell-1}m^{\ell-1} \\
&=\frac{k+\ell-1}{\ell!(\ell-1)!}\sum_{m=1}^k\sum_{i=0}^{\ell-1}(-1)^i\binom{\ell-1}{i}
  (k+1)^{\ell-1-i}m^{i+\ell-1} \\
&=\frac{k+\ell-1}{\ell!(\ell-1)!}\sum_{i=0}^{\ell-1}(-1)^i\binom{\ell-1}{i}
  (k+1)^{\ell-1-i}H_k^{(-r)}\big|_{r=i+\ell-1} \;,
\end{align*}
with $H_k^{(-r)}:=\sum\nolimits_{m=1}^k m^r$ and `$\sim$' means modulo lower
order in $k$. For $r\in\N$ by Faulhaber's formula we have
$$
H_k^{(-r)}=\frac{1}{r+1}\sum_{j=0}^r\binom{r+1}{j}B_j\,k^{r+1-j}
$$
where $B_j$ are Bernoulli numbers; in particular $B_0=1$.
Thus,
\begin{align*}
\mu_k &\sim\frac{k}{\ell!}\sum_{i=0}^{\ell-1}\sum_{j=0}^{i+\ell-1}
  \frac{(-1)^i}{i!\,j!\,(\ell+i-j)!}B_j\,k^{2\ell-j-1} \\
  &\sim\frac{k^{2\ell}}{\ell!}\sum_{i=0}^{\ell-1}
  \frac{(-1)^i}{i!\,(\ell+i)!} \;.
\end{align*}
From this we deduce
$$
\tr(|D|^{-s})=\sum_{k=1}^\infty \mu_k k^{-s}=C_\ell\,\zeta(s-2\ell)+\ldots
$$
where $C_\ell$ is the constant in (\ref{eq:Cl}),
$\zeta$ is the Riemann zeta-function and the dots represent some function
that is holomorphic on the half plane $\,\mr{Re}\,s>2\ell$.
This in particular means that $|D|^{-(2\ell+1)}$ is in the Dixmier ideal.
Its Dixmier trace coincides with the residue of the previous function in $s=2\ell+1$,
which is given by
$$
\tr_\omega(|D|^{-(2\ell+1)})=\mr{Res}_{s=2\ell+1}\tr(|D|^{-s})
=C_\ell\,\mr{Res}_{s=2\ell+1}\zeta(s-2\ell)=C_\ell\;.
$$
Let $\ell=2n$ resp. $2n+1$. In both cases:
\begin{align*}
\ell!\,C_\ell &\geq\sum_{i=0}^{2n-1}\frac{(-1)^i}{i!\,(\ell+i)!}
=\sum_{i=1}^n\left(\frac{1}{(2i-2)!\,(\ell+2i-2)!}-\frac{1}{(2i-1)!\,(\ell+2i-1)!}\right) \\
&=\sum_{i=1}^n\frac{(2i-1)(\ell+2i-1)-1}{(2i-1)!\,(\ell+2i-1)!} >0 \;,
\end{align*}
being a sum of positive terms.
The observation that $C_\ell>0$ for any $\ell$ concludes the proof.
\end{prova}


Consider now the bounded operators $L_{j,k}$ defined by
\begin{equation}\label{eq:Ljk}
L_{j,k}\ket{n,h;\GT{r}}=q^{r_{j,k}}\ket{n,h;\GT{r}} \;,
\end{equation}
and the two-sided ideal $\mc{J}_{2\ell}\subset\B(\HH)$
defined in (\ref{eq:Jr}).

\begin{lemma}\label{lemma:5.6}
For all $j,k=1,\ldots,\ell+1$, we have $L_{j,k}\in\mc{J}_{2\ell}$.
\end{lemma}
\begin{prova}
Since $q^{r_{j,k}}\leq q^{r_{2,\ell}}$ for all $j,k$, it is enough to
prove the statement for $L_{2,\ell}$.
We have
$$
\tr(L_{2,\ell}|D|^{-s})=\sum_{k=1}^\infty k^{-s}\sum_{r_{2,1}=0}^{k-1}
\sum_{h=0}^{r_{2,1}}\sum_{r_{2,\ell}=0}^h q^{r_{2,\ell}} c(k,r_{2,1},h,r_{2,\ell})\;,
$$
where $c(k,r_{2,1},h,r_{2,\ell})$ is the number of $\GT{r}$ of
the form (\ref{eq:form}) with $(k,r_{2,1},h,r_{2,\ell})$ fixed.
Of course, $c(k,r_{2,1},h,r_{2,\ell})\leq c(k,r_{2,1},h,0)$ and
we have ($\sum_{r_{2,\ell}=0}^hq^{r_{2,\ell}}<(1-q)^{-1}$)
$$
\tr(L_{2,\ell}|D|^{-s})<\frac{1}{1-q}\sum_{k=1}^\infty k^{-s}
\sum_{i=0}^{k-1}\sum_{h=0}^i c(k,i,h,0)\;.
$$
Now, looking at (\ref{eq:form}) one realizes that $c(k,i,h,0)$ is the
dimension of the irreducible representation of $SU(\ell)$ with Young
tableau $(i,h,\ldots,h)$, so $\sum_{h=0}^ic(k,i,h,0)=\mu_{i+1}$ where
$\mu_{i+1}$ are the multiplicities appearing in the proof of Proposition
\ref{prop:Dix}, but for a replacement $\ell+1\to\ell$. We know that
$\mu_{i+1}$ is a polynomial in $i+1$ of order $2\ell-2$,
and so $\sum_{i=0}^{k-1}\mu_{i+1}$ is a polynomial in $k$ of order $2\ell-1$.
The observation that
$$
\sum\nolimits_{k=1}^\infty k^{-s+2\ell-1}
$$
converges for $s>2\ell$ concludes the proof.
\end{prova}

By Theorem \ref{thm:main} the difference $\pi(z_i)-W^*\rho(z_i)W$
is the product of a bounded operator by $L_{1,\ell}$ ($q^{r_{1,\ell}}=q^h$),
hence by Lemma \ref{lemma:5.6} it gives an element of $\mc{J}_{2\ell}$.
Repeating verbatim the proof of Lemma \ref{lem:ap} one proves
that $\pi(a)-W^*\rho(a)W\in\mc{J}_{2\ell}$ for all $a\in\A(S^{2\ell+1}_q)$,
which in particular means
$$
\nint a|D|^{-2\ell-1}:=
\mr{Res}_{s=2\ell+1}\tr\bigl(\pi(a)|D|^{-s}\bigr)
=\mr{Res}_{s=2\ell+1}\tr\bigl(W^*\rho(a)W|D|^{-s}\bigr) \;.
$$
With this, it's easy to compute the noncommutative integral.

\begin{prop}
Let $\sigma:\A(S^{2\ell+1}_q)\to\A(S^1)$ be the $*$-algebra morphism
defined by $\sigma(z_i):=0$ for $i\neq\ell+1$, $\sigma(z_{\ell+1})(\theta):=
e^{i\theta}$. For all $a\in\A(S^{2\ell+1}_q)$ we have
\begin{equation}\label{eq:NCint}
\nint a|D|^{-2\ell-1}=\frac{C_\ell}{2\pi}\int_0^{2\pi}\sigma(a)(\theta)\de\theta\;,
\end{equation}
with $C_\ell$ given by (\ref{eq:Cl}).
\end{prop}
\begin{prova}
It is enough to prove (\ref{eq:NCint}) for
$$
a=z_1^{a_1}(z_1^*)^{a'_1}\ldots z_\ell^{a_\ell}(z_\ell^*)^{a'_\ell}z_{\ell+1}^j
$$
with $a_i,b_i\in\N$, $j\in\Z$ and the notation $z_{\ell+1}^j=(z_{\ell+1}^*)^{-|j|}$
if $j<0$, since these monomials are a linear basis for the algebra. We have
$$
\frac{C_\ell}{2\pi}\int_0^{2\pi}\sigma(a)(\theta)\de\theta=
C_\ell\,\delta_{a_1,0}\delta_{a'_1,0}\ldots\delta_{a_\ell,0}\delta_{a'_\ell,0}\delta_{j,0}\;.
$$
Thus if $a=1$ the equality (\ref{eq:NCint}) follows from Proposition \ref{prop:Dix},
while for $a\neq 1$ we have to prove that also the left hand side of (\ref{eq:NCint}) is zero.
By previous Lemma we have
$$
\nint a|D|^{-2\ell-1}=\mr{Res}_{s=2\ell+1}\tr\bigl(W^*\rho(a)W|D|^{-s}\bigr) \;.
$$
For $i\neq\ell+1$, coefficients of $\rho(z_i)$ are bounded by $q^{r_{\ell+2-i,i}}$
(cf.~equation (\ref{eq:Ctilde})). Thus $W^*\rho(a)W\in\mc{J}_{2\ell}$
and $\int\mkern-16mu-\, a|D|^{-2\ell-1}=0$, unless $a_1=a_1'=\ldots=a_\ell=a_\ell'=0$.
It remains to prove that $\int\mkern-16mu-\, z^j_{\ell+1}|D|^{-2\ell-1}=0$ if $j\neq 0$,
but this is trivial since $\rho(z_{\ell+1})$ is an off-diagonal operator.
\end{prova}

\section[Proof of the main theorem]{On the behavior of the coefficients (\protect{\ref{eq:coe}})}\label{app:uno}
The purpose of this Section is to prove Theorem \ref{thm:main}.
Recall (\ref{eq:Lrep}) and (\ref{eq:coe}).

The coefficients in front of the sums in (\ref{eq:Lrep}) are bounded by
$1$. Then, only the behavior of the coefficients inside the sums
matters.
Since $C_q(i,\GT{r},m)$ are bounded by a constant, we immediately deduce
that coefficients in the first line of (\ref{eq:Lrep}) are bounded by $q^h$
and can be neglected.
We now focus on the second line, i.e. on the case $m_1=\ell+1$.

By equation (4.10) of \cite{CP07} there exists a positive constant $C$ such that
\begin{equation}\label{eq:CPineq}
|C_q(i,\GT{r},m)|\leq Cq^{K(i,\GT{r},m)+\sum_{j=m_i}^{\ell+1-i}(r_{i+1,j}-r_{i,j+1})}
\leq Cq^{K(i,\GT{r},m)}
\end{equation}
where
$$
K(i,\GT{r},m):=\sum_{k=1}^{i-1}\left(
\sum_{j=\min(m_k,m_{k+1})}^{\max(m_k,m_{k+1})-1}(r_{k+1,j}-r_{k,j+1})
+2\sum_{j=m_{k+1}+1}^{m_k-1}(r_{k,j}-r_{k+1,j})\right) \;.
$$
All terms with $K(i,\GT{r},m)\geq h$ will be neglected, since in this case
by (\ref{eq:CPineq}) we have $|C_q(i,\GT{r},m)|\leq Cq^h$.
Recall that we are interested in the case $(r_{1,1},\ldots,r_{1,\ell})=(n+h,h,\ldots,h)$
and $m_1=\ell+1$.
The first step is to prove the following Lemma.

\begin{lemma}
Let $m_1=\ell+1$ and $(r_{1,1},\ldots,r_{1,\ell})=(n+h,h,\ldots,h)$.
For any (admissible) $\,(m_2,\ldots,m_i)\neq(\ell,\ell-1,\ldots,\ell+2-i)$
we have $K(i,\GT{r},m)\geq h$.
\end{lemma}
\begin{prova}
We prove it by induction on $i$ starting from $i=2$ (for $i=1$ it is trivial since
$m=(m_1)=(\ell+1)$ is fixed). For $i=2$ recalling that $m_i\leq\ell+2-i$, we get
\begin{align*}
K(2,\GT{r},m)&=
\sum_{j=\min(m_1,m_2)}^{\max(m_1,m_2)-1}(r_{2,j}-r_{1,j+1})
+2\sum_{j=m_2+1}^{m_1-1}(r_{1,j}-r_{2,j}) \\ &\geq
(r_{2,\ell}-r_{1,\ell+1})+2\sum_{j=m_2+1}^{\ell}(r_{1,j}-r_{2,j}) \;.
\end{align*}
If $m_2=\ell$ last sum is empty, otherwise (i.e.~when $m_2<\ell$) it has at least the
term with $j=\ell$ and so
$$
K(2,\GT{r},m)\geq
(r_{2,\ell}-r_{1,\ell+1})+2(r_{1,\ell}-r_{2,\ell})=r_{1,\ell}+(r_{1,\ell}-r_{2,\ell})
\geq r_{1,\ell}=h \;.
$$
Hence $m_2\neq\ell$ implies $K(2,\GT{r},m)\geq h$.

Now we prove the inductive step. Assume the Lemma is true for $i$.
Let $m=(m_1,\ldots,m_i)$ and $m'=(m,m_{i+1})$, and notice that
$$
K(i+1,\GT{r},m')\geq K(i,\GT{r},m) \;.
$$
If $m\neq (\ell+1,\ldots,\ell+2-i)$ we have $K(i,\GT{r},m)\geq h$ (inductive
hypothesis) and so $K(i+1,\GT{r},m')\geq h$. It remains to prove that even
for $m=(\ell+1,\ldots,\ell+2-i)$, if $m_{i+1}\neq\ell+1-i$ we still have 
$K(i+1,\GT{r},m')\geq h$.

Take then $m$ as above and $m_{i+1}<\ell+1-i=m_i-1$.
Notice that in this case
$$
K(i,\GT{r},m)=\sum_{k=1}^{i-1}(r_{k+1,\ell+1-k}-r_{k,\ell+2-k})=r_{i,\ell+2-i} \;,
$$
(being a telescopic sum) and so
\begin{align*}
K(i+1,\GT{r},m') &=K(i,\GT{r},m)
+\left(\sum_{j=m_{i+1}}^{m_i-1}(r_{i+1,j}-r_{i,j+1})
+2\sum_{j=m_{i+1}+1}^{m_i-1}(r_{i,j}-r_{i+1,j})\right) \\ &\geq
r_{i,\ell+2-i}+(r_{i+1,\ell+1-i}-r_{i,\ell+2-i})
+2(r_{i,\ell+1-i}-r_{i+1,\ell+1-i}) \\ &=
r_{i,\ell+1-i}+(r_{i,\ell+1-i}-r_{i+1,\ell+1-i})
\geq r_{i,\ell+1-i}\geq r_{1,\ell}=h \;.
\end{align*}
In the first inequality we neglect all the terms in the sums with $j\neq m_i-1=\ell+1-i$
(both sums are non-empty since by hypothesis $m_{i+1}<\ell+1-i$).
In the second inequality we neglect the positive term inside the parenthesis.
Last inequality is due to the definition of GT tableaux.
This concludes the proof.
\end{prova}

By previous Lemma all terms in $\pi(z_{\ell+2-i}^*)$ have weights bounded
by $q^h$, but for the one corresponding to $m=m^i:=(\ell+1,\ell,\ldots,\ell+2-i)$.
Thus,
\begin{align*}
\pi(z_{\ell+2-i}^*)\ket{n,h;\GT{r}} &=
q^{-\frac{1}{2}n}\frac{[h+\ell-1]^{\frac{1}{2}}}{[n+h+\ell-1]^{\frac{1}{2}}}\,
q^{-i+1}C_q(i,\GT{r},m^i)\ket{n,h-1;\GT{r}^{\smash[t]{(m^i)}}} +\ldots \\ &=
q^{-i+1}C_q(i,\GT{r},m^i)\ket{n,h-1;\GT{r}^{\smash[t]{(m^i)}}} +\ldots \;,
\end{align*}
where dots denote the terms which we are neglecting.
The next step is to simplify $C_q(i,\GT{r},m^i)$ by neglecting terms bounded
by $q^h$.
Plugging $m_k=\ell+2-k$ ($k\leq i$) and $(r_{1,1},\ldots,r_{1,\ell})=(n+h,h,\ldots,h)$
into equation (\ref{eq:coe}) we get, after some algebraic manipulations,
\begin{align}
q^{-i+1}&C_q(i,\GT{r},m^i)=q^{r_{i,m_i}}\;
\left\{\prod_{k=1}^{\ell+1-i}
\frac{\qan{r_{i+1,k}-r_{i,m_i}-k+m_i-1}}{\qan{r_{i,k}-r_{i,m_i}-k+m_i}}
\right\}\rule{-4pt}{18pt}^{1/2} \! \times \nonumber\\ & \times
\left\{\prod_{k=1}^{i-1}
\prod_{j=1}^{\ell+1-k}\frac{\qan{r_{k,j}-r_{k+1,m_{k+1}}-j+m_{k+1}}}{\qan{r_{k,j}-r_{k,m_k}-j+m_k}}
\prod_{j=1}^{\ell-k}\frac{\qan{r_{k+1,j}-r_{k,m_k}-j+m_k-1}}{\qan{r_{k+1,j}-r_{k+1,m_{k+1}}-j+m_{k+1}-1}}
\right\}\rule{-4pt}{18pt}^{1/2} \;. \label{eq:ctmp}
\end{align}
We notice that (by definition of GT tableaux) all arguments of the exponentials
are non-negative (positive if in the denominator).
In fact, last denominator vanishes if $j=\ell-k$ and $r_{k+1,\ell-k}=r_{k+1,\ell-k+1}$,
but this `singularity' is just apparent. For $r_{k+1,\ell-k}=r_{k+1,\ell-k+1}$ the
move $\GT{r}\to\GT{r}^{(m^i)}=:\GT{r}'$ shifts $r_{k+1,\ell-k}$ by $-1$ and leaves $r_{k+1,\ell-k+1}$
invariant: hence the result is that $r_{k+1,\ell-k}'<r_{k+1,\ell-k+1}'$ and $\GT{r}'$
is not a GT tableaux. In conclusion, the term with $r_{k+1,\ell-k}=r_{k+1,\ell-k+1}$
(for which the denominator vanishes) never appears.

We use the notation $\,a\lesssim b\,$ if $\,a\leq Cb\,$ and $C$ is a positive
constant depending only on~$q$.

\begin{lemma}
For any $0<u\leq 1$ we have
\begin{equation}\label{eq:ineq}
0\leq 1-(1-u)^{\frac{1}{2}}\leq u \;,\qquad
0\leq (1-qu)^{-\frac{1}{2}}-1\lesssim u \;.
\end{equation}
\end{lemma}
\begin{prova}
That the differences in (\ref{eq:ineq}) are non-negative is obvious.
Since $0\leq 1-u<1$ we have
$$
1-u\leq\sqrt{1-u} \;,
$$
and this proves the first inequality in (\ref{eq:ineq}). We have also
\begin{align*}
\frac{1}{\sqrt{1-qu}}-1
&=\left(\frac{1}{\sqrt{1-qu}}+1\right)^{-1}
\left(\frac{1}{1-qu}-1\right) \\
&=\left(\frac{1}{\sqrt{1-qu}}+1\right)^{-1}
\frac{qu}{(1-qu)} \\
&\leq 2^{-1}\frac{qu}{1-q}\lesssim u \;,
\end{align*}
and this proves the second inequality.
\end{prova}
\begin{lemma}
For any $0<u,v\leq 1$ we have
\begin{equation}\label{eq:ineqb}
\left|1-\sqrt{\frac{1-u}{1-qv}}\right|\lesssim\max(u,v)\;.
\end{equation}
\end{lemma}
\begin{prova}
The algebraic identity
$$
1-\sqrt{\frac{1-u}{1-qv}}=
\frac{1}{2}\left(1-\frac{1}{\sqrt{1-qv}}\right)\left(1+\sqrt{1-u}\right)+
\frac{1}{2}\left(1+\frac{1}{\sqrt{1-qv}}\right)\left(1-\sqrt{1-u}\right) \;,
$$
together with the triangle inequality and the previous Lemma implies
\begin{align*}
\left|1-\sqrt{\frac{1-u}{1-qv}}\right| &\lesssim
\frac{1}{2}v(1+\sqrt{1-u})+
\frac{1}{2}\left(1+\frac{1}{\sqrt{1-qv}}\right)u \\
& \leq v+(1-q)^{-\frac{1}{2}}u\leq 
(1-q)^{-\frac{1}{2}}(u+v)\leq
2(1-q)^{-\frac{1}{2}}\max(u,v) \;.
\end{align*}
This concludes the proof.
\end{prova}

Applying (\ref{eq:ineqb}) to the coefficients in (\ref{eq:ctmp}) we get
\begin{align*}
q^{r_{i,m_i}}\left|1-\left\{\frac{\qan{r_{i+1,j}-r_{i,m_i}-j+m_i-1}}{\qan{r_{i,j}-r_{i,m_i}-j+m_i}}
\right\}\rule{-4pt}{18pt}^{1/2} \right| &\lesssim q^{r_{i+1,j}}\leq q^{h}
 \;, && \forall\;1\leq j\leq\ell-i\,,
\\
q^{r_{i,m_i}}\left|1-\left\{\frac{\qan{r_{k,j}-r_{k+1,m_{k+1}}-j+m_{k+1}}}
{\qan{r_{k,j}-r_{k,m_k}-j+m_k}}\right\}\rule{-4pt}{18pt}^{1/2}\right| &\lesssim
q^{r_{k,j}}\leq q^{h} \;, && \forall\;1\leq j\leq\ell+1-k\,,
\\
q^{r_{i,m_i}}\left|1-\left\{\frac{\qan{r_{k+1,j}-r_{k,m_k}-j+m_k-1}}
{\qan{r_{k+1,j}-r_{k+1,m_{k+1}}-j+m_{k+1}-1}}\right\}\rule{-4pt}{18pt}^{1/2}\right| &\lesssim
q^{r_{k+1,j}}\leq q^{h} \;, && \forall\;1\leq j\leq\ell-k\,,
\end{align*}
where we used $u\leq u^{1/2}$ for $u\leq 1$, $r_{i,m_i}-r_{k,m_k}\geq 0$ and
$r_{i,m_i}-r_{k+1,m_{k+1}}\geq 0$ for $k\leq i-1$, and
the fact that $r_{j,k}\geq h$ for all $j+k<\ell+2$.

Using previous inequalities in (\ref{eq:ctmp}) we deduce
$$
q^{-i+1}C_q(i,\GT{r},m^i)=O(q^h)+\begin{cases}
q^{r_{i,m_i}}\left\{
\frac{\qan{r_{i+1,m_{i+1}}-r_{i,m_i}}}{\qan{r_{i,m_{i+1}}-r_{i,m_i}+1}}
\right\}\rule{-4pt}{18pt}^{1/2} &\mr{if}\;0\leq i\leq\ell\,,\\
q^{r_{i,m_i}} &\mr{if}\;i=\ell+1\,.
\end{cases}
$$
Now we apply the second inequality in (\ref{eq:ineq}) to prove that
$$
q^{r_{i,m_i}}\left|1-\big\{\qan{r_{i,m_{i+1}}-r_{i,m_i}+1}\big\}^{-\frac{1}{2}}\right|
\leq q^{r_{i,m_i}}q^{2(r_{i,m_{i+1}}-r_{i,m_i})}\leq q^{r_{i,m_{i+1}}}\leq q^h \;.
$$
Recalling that $r_{1,\ell+1}=0$,
we conclude that
$|q^{-i+1}C_q(i,\GT{r},m^i)-\tilde{C}_q(i,\GT{r},m^i)|\lesssim q^h$,
where
\begin{equation}\label{eq:Ctilde}
\tilde{C}_q(i,\GT{r},m^i):=
\begin{cases}
\sqrt{1-q^{2r_{2,\ell}}}
    &\mr{if}\;\;i=1\,,\\
q^{r_{i,\ell+2-i}}\sqrt{\qan{r_{i+1,\ell+1-i}-r_{i,\ell+2-i}}}\quad
    &\mr{if}\;\;i=2,\ldots,\ell\,,\\
q^{r_{\ell+1,1}}
    &\mr{if}\;\;i=\ell+1\,.
\end{cases}
\end{equation}
Summarizing, we have
$$
\pi(z_{\ell+2-i}^*)\ket{n,h;\GT{r}}=\tilde{C}_q(i,\GT{r},m^i)\ket{n,h-1;\GT{r}^{\smash[t]{(m^i)}}}+\ldots\;,
$$
where `$\ldots$' are shift operators with weights bounded by $q^h$.

We are now in the position to prove Theorem (\ref{thm:main}).
Recall that we have the left regular representation $\pi:\A(S^{2\ell+1}_q)\to\B(\HH)$,
given by equation (\ref{eq:Lrep}), the representation
$\rho:\A(S^{2\ell+1}_q)\to\B(\hat{\HH})$ of Section \ref{sec:9.2},
and the partial isometry $W:\HH\to\hat{\HH}$ of Section \ref{sec:9.3}.

By definition of $\rho$ and $W$ one finds that
$$
W^*\rho(z_{\ell+2-i}^*)W\ket{n,h;\GT{r}}
=\tilde{C}_q(i,\GT{r},m^i)\ket{n,h-1;\GT{r}^{\smash[t]{(m^i)}}}\;,
$$
where the coefficients $\tilde{C}_q(i,\GT{r},m^i)$ are exactly the ones
defined in (\ref{eq:Ctilde}).
Then, for all $\,i=1,\ldots,\ell+1\,$, the operator $\;\pi(z_i)-W^*\rho(z_i)W\,$
is a finite sum of shift operators on $\HH$ with weights bounded by $q^h$.
This concludes the proof of the Theorem.


\backmatter

\cleardoublepage
\phantomsection
\addcontentsline{toc}{chapter}{Bibliography}

\newcommand{\etalchar}[1]{$^{#1}$}
\providecommand{\bysame}{\leavevmode\hbox to3em{\hrulefill}\thinspace}


\end{document}